\documentclass{articlenotescls}
\usepackage{gencatmacros}

\setboolean{cdbaselineon}{true}
\setboolean{showcdpunctuation}{true}

\newcommand{\cVCat}{\mathbf{\mathcal{V}\Cat}}
\newcommand{\Graph}{\mathbf{Graph}}
\newcommand{\cVGraph}{\mathbf{\mathcal{V}\Graph}}
\newcommand{\nCat}[1]{\mathbf{#1\Cat}}

\newcommand{\cZ}{\mathcal{Z}}

\newcommand{\ucC}{\underline{\mathcal{C}}}

\newcommand{\ucV}{\underline{\mathcal{V}}}

\newcommand{\uCat}{\underline{\mathbf{Cat}}}

\newcommand{\MonCat}{\mathbf{MonCat}}
\newcommand{\BrMonCat}{\mathbf{BrMonCat}}
\newcommand{\SymMonCat}{\mathbf{SymMonCat}}

\newcommand{\Mon}{\mathbf{Mon}}
\newcommand{\cMon}{\mathbf{cMon}}
\newcommand{\CoMon}{\mathbf{CoMon}}
\newcommand{\cCoMon}{\mathbf{cCoMon}}

\newcommand{\vDblCat}{\mathbf{vDblCat}}
\newcommand{\vcoDblCat}{\mathbf{vcoDblCat}}

\newcommand{\Mod}{\mathbf{Mod}}
\newcommand{\Mat}{\mathbf{Mat}}
\newcommand{\cMod}{\mathbf{cMod}}
\newcommand{\CoMod}{\mathbf{CoMod}}
\newcommand{\cCoMod}{\mathbf{cCoMod}}

\newcommand{{\Inv}}{\mathbf{Inv}}
\newcommand{\nInv}[1]{\mathbf{#1Inv}}
\newcommand{\nMon}[1]{\mathbf{#1Mon}}

\newcommand{\BiCat}{\mathbf{Bicat}}
\newcommand{\Gray}{\mathbf{Gray}}

\newcommand{\bG}{\mathbb{G}}
\newcommand{\bV}{\mathbb{V}}

\newcommand{\st}{\mathrm{st}}

\newcommand{\ps}{\mathrm{ps}}

\newcommand{\vequiv}{\mathrel{\rotatebox[origin=c]{-90}{$\equiv$}}}

\begin{document}
	\pagestyle{empty}
	\pagenumbering{Roman}
	
	\title{Relative involutions and inverses in monoidal dualities and slice constructions}
	\author{Varinderjit Mann}
	\date{\today}
	\maketitle
	
	\begin{abstract}
		Our goal in this paper is to work towards developing a framework for characterizing dualities of weak categorical structures. It is often the case that duality involutions are characterized by relations which are more general than $f^2=1$. Moreover, it is also true that many dualities cannot be expressed using involutions alone. We develop a framework of relative involutions and relative inverses to express these dualities at least within categorical and bicategorical structures. In particular, we utilize this framework to express the dualities of virtual double categories and some important constructions on them. In the second part, we explore the problem of slicing and consider how dualities of monoidal structures interact with slice constructions. We end by providing a construction for both enriched slice model categories and enriched coslice model categories. 
	\end{abstract}
	
	\tableofcontents
	\clearpage 
	
	\pagestyle{scrheadings}
	\pagenumbering{arabic}
	
	\section{Introduction} \label{sec:intro}
	\subsection{Motivation}

Our primary motivation to write this paper comes from two questions, which are broadly related under the microcosm principle. The first one concerns the dualities of categorical structures. Consider the $\op_1$ duality on $\nCat{1}$, expressed as a $2$-functor $\op_1 \colon \op_2 (\nCat{1}) \rightarrow \nCat{1}$, which satisfy the relation $\op_1 \circ \op_2 (\op_1) =1$. Clearly, this is not an involution in the traditional sense of satisfying $\op_1 \circ \op_1 = 1$, since this composition is simply not defined. Rather, it is an involution relative to a higher involution $\op_2$. This motivates the definition of a relative involution.

Now consider the duality $\op_1 \colon \op_2(\MonCat_{c}) \rightarrow \MonCat_{l}$. Here, the domain consists of colax monoidal functors, whereas the codomain consists of lax monoidal functors. This is again a duality, but it is no longer directly expressible even using relative involution as before. Instead, this motivates a further generalization to relative inverses. 

As with many algebraic structures, specifying the complete coherence data for these definitions becomes increasingly difficult in higher dimensions. Nevertheless, it is tractable and useful enough to define these for involutive bicategories and biadjunction of bicategories. In particular, this framework is sufficiently general to express the dualities of monoidal structures. Furthermore, it is also natural to consider module structures in the presence of monoidal structures. Thus, one of the main goals of this paper is to extend this framework onto virtual double categories and the corresponding module construction.  

The second question concerns slicing and its interaction with these dualities. In particular, we wish to describe a general construction of enriched slice and coslice categories, especially under the consideration of model structures. Although there are some accounts of enriched slice constructions, it is useful to have a systematic treatment that encompasses both slice and coslice constructions. Furthermore, the slice and coslice category construction is related by the $\op_1$ duality. This suggests that the enriched coslice case should be a dual of the enriched slice case in some appropriate sense. Establishing such a construction is the second main goal of this paper.

\subsection{Main results and applications}

The results can be divided into three major parts. The first part establishes a desired framework of relative involutions and relative inverses. The second part applies this framework to characterize dualities of constructions on virtual double categories. The third part discuss the slicing of module structures and establish a general construction of enriched slice and coslice model categories.

Recall that one can define a $1$-category of relative involutions, called $\nInv{0}(\cC,T,\alpha)$, inside an arbitrary lax involutive category $(\cC,T,\alpha)$ \cite{yauinvcat2020} (or see \cref{sub:1diminv}). This is an example of the microcosm principle. Motivated by this fact, we show a similar construction for adjunctions by defining relative inverses. The main result in this part is to generalize this framework one dimension higher for weak $2$-involutive bicategories and weak $2$-inverses of bicategories. This can be summarized by the following proposition.

\begin{proposition*}
	The microcosm principle holds for weak $2$-involutive bicategories and weak $2$-inverses of bicategories.
\end{proposition*}

In particular, this provides suitable definitions of relative weak $1$-inverses, relative weak $1$-involutions, and their $1$-morphisms. Note that expressing the dualities of strict $n$-categories only require relative involution or relative inverse (see \cref{prop:opdualities}). This is in contrast to the monoidal structures, which also require the $1$-morphisms. The main result here is that the $2$-functors $\Mod \colon \vDblCat \rightarrow \vDblCat$ and $\CoMod \colon \vcoDblCat \rightarrow \vcoDblCat$ are compatible with the relevant dualities.

\begin{theorem*} The following extensions hold:
	\begin{enumerate}
		\item $\Mod$ extends to a $1$-morphism of involution, respecting the horizontal duality of virtual double categories. 
		
		\item $\Mod$ and $\CoMod$ together extend to a $1$-morphism of inverses, respecting the vertical duality of virtual double categories.
	\end{enumerate}
\end{theorem*}

In the final part, we first briefly review the slice constructions for module categories. Recall that standard definitions of an enriched model category assumes the existence of a two-variable closure adjunction relating tensoring, powering, and enrichement \cite{hovey2007modelcat} (or see \cref{sub:slicecatmod}). Hence, we focus primarily on the behavior of this two-variable adjunction under slicing. We show that it can be sliced in two dual ways, one by slicing the tensoring and the other by slicing the powering. The main result is that these slicing also respect the model structures.

\begin{theorem*}
	Let $\ucC$ be a $\ucV$-enriched model categories. Let $M\in \cV$ and $Z \in \cC$ be suitable monoid and module object, respectively. Then, both $Z/\ucC$ and $\ucC/Z$ are $\ucV/M$-enriched model categories.
\end{theorem*}

\subsection{Overview}

\Cref{sec:prelim} discusses the dualities of strict $n$-categories as a motivating example of nested algebraic structures. To formulate this more precisely, it recall the microcosm principle and define the sliced microcosm principle. A relevant example of these principles is then provided in terms of pseudomonoid and pseudomodule objects.

\Cref{sec:ginvs} begins by briefly reviewing the theory of involutions and relative involutions. Correspondingly, it defines the similar theory of inverses and relative inverses. Then, the analogues of these structures are described in one dimension higher. 

\Cref{sec:dualvdcs} starts with reformulating the dualities of monoidal categories within the framework of relative inverses and involutions. This framework is then used to describe the corresponding dualities of internal monoid structures. Subsequently,  the dualities are generalized to the level of virtual double categories, encoding the action on the module data.

\Cref{sec:slicecatmod} applies the sliced microcosm principle to monoidal and module structures. First, it reviews enrichement through tensoring and states the corresponding slice construction. The results are then dualized for the case of cotensoring. Finally, these constructions are extended to the setting of enriched model categories. 
	\section{Preliminaries} \label{sec:prelim}
	\subsection{\tpdf{$\nCat{n}$, monoidal $n$-structures, and categorical dualities}}

Recall that the category $\nCat{1}=\Cat$ is canonically enriched over itself using the functor categories $\uCat(\cC,\cD)$. It also has a Kelly tensor product given by the cartesian product $\cC \times \cD$ of two categories, which forms a monoidal structure. Therefore, it is natural to ask whether this enrichment and the monoidal structure are compatible. Indeed, we can enrich the monoidal structure to a monoidal $2$-structure\footnote{The term $2$-monoidal category can be confused with $2$-tuply monoidal $1$-category which is far different from the $1$-tuply monoidal $2$-category as described in our case\cite{aguiarmon2010,yau2024bimodcats}.} in a canonical way. There exists a $\uCat$ enriched adjunction
\[\uCat(\cC \times \cD, \cE) \cong \uCat(\cC,\uCat(\cD,\cE))\cdpunct{,}\]
by using a general result which holds for any closed symmetric monoidal category category $\ucV$ \cite[Prop.~3.1.11]{yau2024bimodcats}. In this enriched adjunction, both the product and the internal hom functors are $2$-functors. Furthermore, the left unitor, the right unitor, the associator, and the symmetry are now all $2$-natural transformations. This defines, naturally, a symmetric cartesian closed monoidal $2$-structure. 

\begin{lemma}
	The tuple $(\Cat, \times, \{*\}, a_{\Cat},\lambda_{\Cat}, \rho_{\Cat})$ defines a symmetric cartesian closed monoidal structure on $\Cat$. Furthermore, this monoidal structure can be $\Cat$-enriched by canonically extending the monoidal data $\{\times , a_{\Cat} ,\lambda_{\Cat}, \rho_{\Cat}\}$ to a strict $2$-functor and strict $2$-natural transformations. All the coherence relations remain unchanged. This makes $\Cat$ a monoidal $2$-category.
\end{lemma}

As a consequence, it is natural question to ask whether the same result extends to $\nCat{n}$ or $\nCat{\omega}$.  The following general lemma answers this in affirmative. 

\begin{lemma}
	Let $\cV$ be a locally small, complete, and cocomplete cartesian closed monoidal cateogry. Then, the category $\cVCat$ (of small $\cV$-categories) is also locally small, complete, cocomplete, and cartesian closed. Furthermore, the cartesian product of $\cVCat$ is precisely the Kelly tensor product, as given by the cartesian product of sets of objects and the corresponding products of hom-objects.
\end{lemma} 

\begin{proof}
	Given any symmetric closed monoidal category $\cV$, the (co)completeness of $\cV$ implies the (co)cocompleteness of $\cVCat$ by results of \cite{wolff1974,betti1983,kellyvcatpres2001}. And, the Kelly tensor product of $\cVCat$ forms a symmetric closed monoidal structure by \cite[Sec.~2.3]{kellyect}. We are only left to verify that the kelly tensor product is indeed the cartesian product. This last part essentially follows from the monadicity of $\cVCat$ over $\cVGraph$. 
\end{proof}

In particular, this lemma provides the following desired corollary.

\begin{corollary}
	The tuple $(\nCat{n}, \times, \{*\}, a_{\nCat{n}},\lambda_{\nCat{n}}, \rho_{\nCat{n}})$ canonically defines a symmetric cartesian closed monoidal $(n+1)$-structure on $\nCat{n}$.
\end{corollary}

Finally, this result also extends to $\nCat{\omega}$, which can be defined as a form of colimit of the enrichment process. In particular, the category of all (small) $\nCat{\omega}$-enriched category is just the category $\nCat{\omega}$ itself and is cartesian closed. The interested reader maybe look up the details in \cite{leinster2004,streetncat1987}.

Using suitable universes, it is also clear that each monoidal $(n+1)$-category $\nCat{n}$ is equivalently a \enquote{monoidal object} inside the monoidal $(n+2)$-category $\nCat{(n+1)}$. This provides our first example of a \emph{nested algebraic structure}.

Another primary example is given by the duality involution of strict $n$-categories. Consider the $1$-functor $\op_1\colon \nCat{1} \rightarrow \nCat{1}$ which sends each category $\cC$ to its opposite category $\cC^{\op}$. This is a duality involution $1$-functor on $\nCat{1}$, that is, $\op_1^2=1$. However, as $\nCat{1}$ is naturally a strict $2$-category, we may wish to upgrade this duality involution to a strict $2$-functor. This is not possible. What we instead have is a strict $2$-functor $\op_1 \colon \op_2(\nCat{1}) \rightarrow \nCat{1}$, where $\op_2$ is a strict $1$-functor $\op_2 \colon \nCat{2}_{\leq 1} \rightarrow \nCat{2}_{\leq 1}$ and $(-)_{\leq i}$ denotes truncation of categorical structures. Moreover, the duality of $\op_1$ is now expressed by the equation $\op_1 \circ \op_2(\op_1)=1$.

Thus, $\op_1$ duality is fully defined only in the presence of a higher dimensional $\op_2$ duality. More generally, let $\op_{n,k}$ denote the $k$-cell duality on $n$-categories for $1 \leq k \leq n$. All the $1$-functors $(\op_{n,k})_{\leq 1} \colon n\Cat_{\leq 1} \rightarrow n\Cat_{\leq 1}$ are definable readily. See \cite{ara2015autequivinfcat,markakis2025hominfcat} for example. This leads to the following natural proposition.

\begin{proposition} \label{prop:opdualities}
	For each $1 \leq k \leq n$, the $1$-functor $(\op_{n,k})_{\leq 1} \colon n\Cat_{\leq 1} \rightarrow n\Cat_{\leq 1}$ extends to a corresponding $(n+1)$-functor $\op_{n,k} \colon \op_{n+1,k+1}(n\Cat) \rightarrow n\Cat$.
\end{proposition}

Note that both $\nCat{n}$ and $\op_{n+1,k+1}(\nCat{n})$ are enriched over $\nCat{n}$ by definition as they are $(n+1)$-categories. This means that our primary goal is to extend the given $1$-functors to corresponding $\nCat{n}$-enriched functor. Consequently, we first state the following more general lemma inspired by an example in \cite[Ex.~4.6]{shulman2018contenrich}.

\begin{lemma}
	Let $(\cV, \otimes, I)$ be a closed symmetric monoidal category equipped with a strict symmetric monoidal functor $o\colon \cV \rightarrow \cV$ satisfying $o^2 = 1$. Then, there exists a canonical $\cV$-enriched functor $o \colon o_*(\ucV) \rightarrow \ucV$, which extends the $1$-functor $o \colon \cV \rightarrow \cV$.
\end{lemma}

\begin{proof}
	$\cV$ is enriched, tensored, and cotensored over itself because it is a closed symmetric monoidal category. More precisely, $\cV$ extends to a symmetric monoidal $\cV$-enriched category \cite[Thm.~3.3.2]{yau2024bimodcats}. Thus, we can do a change of base of enrichement to produce a $\cV$-enriched category $o_*(\ucV)$. Note that the underlying category of $o_*(\ucV)$ is still $\cV$ itself because $o$ is a strict monoidal functor.

	The key idea is to note that $o_*(\ucV)$ is also tensored over $\ucV$. We define the tensoring operation using the functor $\triangleleft \colon \cV \times \cV \xrightarrow{o \times 1}  \cV \times \cV \xrightarrow{\otimes} \cV$. Our goal is now to show that there is a $\cV$-isomorphism
	\[o_*(\ucV)(A \triangleleft X, Y) \cong \ucV(A,o_*(\ucV)(X,Y))\]
	natural in all of $A,X,Y$. This is a straightforward exercise using the Yoneda lemma. Thus, we have established that $\triangleleft$ is indeed a valid tensoring operation.
	
	Finally, we note the following commuting diagram
	\[\begin{tikzcd}[ampersand replacement=\&,cramped]
		\& {\cV \times \cV} \&\&\& {\cV \times \cV} \&\& {X \otimes o(Y)} \\
		{\cV \times \cV} \\
		\& \cV \&\&\& \cV \&\& {o(X \triangleleft Y)}
		\arrow["{1 \times o}", from=1-2, to=1-5]
		\arrow["{o \times 1}"', from=1-2, to=2-1]
		\arrow["\triangleleft"', from=1-2, to=3-2]
		\arrow["{\otimes }", from=1-5, to=3-5]
		\arrow[equals, nfold, from=1-7, to=3-7]
		\arrow["\otimes"', from=2-1, to=3-2]
		\arrow["o"', from=3-2, to=3-5]
	\end{tikzcd}\cdpunct[8pt]{,}\]
	which provides a canonical identity left tensorial strength for the unenriched functor $o$. This implies the existence of a canonical $\cV$-enriched functor $o \colon o_*(\ucV) \rightarrow \ucV$ extending the $1$-functor $o \colon \cV \rightarrow \cV$. For more details, see \cite[Def.~3.1.1 and Prop~3.2.4]{segrt2012strengthsmorita}. Explicitly, the action on hom objects is given by the adjoint to composite 
	\[\begin{tikzcd}[ampersand replacement=\&,cramped]
		\&\& {o(o_*(\ucV)(X,Y) \triangleleft X)} \&\& \\
		{o_*(\ucV)(X,Y) \otimes oX} \&\&\&\& oY
		\arrow["{o(^{\triangleleft}\epsilon^Y_X)}", from=1-3, to=2-5]
		\arrow[equals, nfold, from=2-1, to=1-3]
		\arrow["{o(^{\otimes}\epsilon^Y_X)}"', from=2-1, to=2-5]
	\end{tikzcd}\cdpunct[8pt]{.}\]
\end{proof}

Now, we can use this lemma in our context to obtain the desired result.

\begin{proof}[Proof of \cref{prop:opdualities}]
	This readily follows once we note that $\op_{n+1,k+1}(\nCat{n}) = {(\op_{n,k})}_*(\nCat{n})$ and that all the $1$-functors $(\op_{n,k})_{\leq 1}$ (as defined in \cite[Sec.~5]{markakis2025hominfcat}) are all strictly symmetric and strictly involutive in the desired sense. 
\end{proof}

We can also provide a more concrete description of these $(n+1)$-functors on the hom $n$-categories. Notice that the diagram 
$\begin{tikzcd}[ampersand replacement=\&,cramped]
	{[0]} \&\& {[1]}
	\arrow[shift left=2, from=1-1, to=1-3]
	\arrow[shift right=2, from=1-1, to=1-3]
	\arrow[from=1-3, to=1-1]
\end{tikzcd}$ defines a cocategory object in $\Cat$. In particular, this allows one to compute the functor categories $\cD^{\cC}$ as $(\cD^{\cC})_0 \cong \Cat([0] \times \cC, \cD)$ and $(\cD^{\cC})_1 \cong \Cat([1] \times \cC, \cC) $. All of the source, target, identity, and composition maps for $\cD^{\cC}$ are derived from the cocategory object structure of  
$\begin{tikzcd}[ampersand replacement=\&,cramped]
	{[0]} \&\& {[1]}
	\arrow[shift left=2, from=1-1, to=1-3]
	\arrow[shift right=2, from=1-1, to=1-3]
	\arrow[from=1-3, to=1-1]
\end{tikzcd}$. 

More generally, the category of all $n$-categories $n\Cat$ has a $n$-cocategory object which can be generated from the diagram $G_0 \rightrightarrows G_1 \rightrightarrows G_2 \cdots \rightrightarrows G_n$, where each $G_i$ is the $i$-globe. This essentially follows from \cite[Prop.~3.14 and Prop.~3.15]{ara2011thetacatcocat} by considering the Globe category embedding $\bG_{n} \rightarrow \Theta_n$.

For each $1\leq k\leq n$, the globe $G_i$ has canonical isomorphisms $\op_{n,k}(G_i) \cong G_i$  induced by the corresponding duality on the $n$-Globe category $\op_{n,k} \colon \bG_n \rightarrow \bG_n$, which swaps the source map and target map that goes into the $k$th level. Indeed, the corresponding $(\op_{n,k})_{\leq 1}$ dualities on $n$-globular presheaves or $n$-categories are also precisely induced from the ones on the Globe category. See \cite[Sec.~5]{markakis2025hominfcat} for more details. Furthermore, like any arbitrary functor, the globe category dualities are compatible with the diagonal map $\bG_n \rightarrow \bG_n \times \bG_n$ and the swap map $\bG_n \times \bG_n \xrightarrow{s} \bG_n \times \bG_n$. Thus, demonstrating the strict symmetric monoidal nature of the induced dualities.

Using this, we can verify that the action of $\op_{n,k}$ on the $i$-transfor  $G_i \times \cC \rightarrow \cD$ is given by
\[\begin{tikzcd}[ampersand replacement=\&,cramped]
	{G_i \times C} \&\&\& {\op_{n,k}G_i \times \op_{n,k}C} \&\& {G_i \times \op_{n,k}C} \\
	D \&\&\& {\op_{n,k}D} \&\& {\op_{n,k}D}
	\arrow[""{name=0, anchor=center, inner sep=0}, from=1-1, to=2-1]
	\arrow[""{name=1, anchor=center, inner sep=0}, from=1-4, to=2-4]
	\arrow[""{name=2, anchor=center, inner sep=0}, from=1-6, to=2-6]
	\arrow["{(\op_{n,k})_{\leq 1}}", between={0.2}{0.8}, maps to, from=0, to=1]
	\arrow["\cong", between={0.2}{0.8}, maps to, from=1, to=2]
\end{tikzcd}\cdpunct[8pt]{.}\]
This ensures that $\op_{n,k}$ changes the direction of only the $k$-transfor. For instance, $\op_{1,1}$ acting on the $1$-functor $[1] \times \cC \rightarrow \cD$ changes only the direction of $[1]$, which is the direction of natural transformations. 

Ultimately, each lower dimensional duality is fully definable only when the higher dimensional duality is defined, at least in a $1$-truncated form. These nested structures can be better formulated and understood using the microcosm principle

\subsection{Mircocosm principle}

The initial statement of the microcosm principle goes back to Baez and Dolan, abstracting a general guiding principle from the observation that monoids are definable in any monoidal category. 

\begin{definition}[Microcosm principle, Baez and Dolan]
	The microcosm principle suggests that, in order to define a generic $(n-1)$-coherent algebraic structure inside an $n$-category, the category itself should be equipped with a categorified $n$-coherent version of the same algebraic structure.
\end{definition} 

A precise formulation of the microcosm principle is given in \cite[Def.~57]{baezhda31998} using the opetopic foundations of $n$-categories. The key idea is to define internal algebraic structures by using maps out of the terminal categorified algebraic structure. However, this result is not immediately applicable to our usual algebraic foundations. Nor does it address whether slicing with respect to an internal algebraic object is always possible. Indeed, consider a monoid object $M$ in a monoidal category $\cV$. The slice category $\cV/M$ admits a canonical monoidal category structure with the sliced tensor product $X \otimes Y \rightarrow M \otimes M \rightarrow M$. This example suggests another general principle governing these slice constructions. Accordingly, we propose the following variant of the microcosm principle.

\begin{definition}[Sliced microcosm principle]
	Let $A$ be a $(n-1)$-coherent algebraic structure inside a $n$-category $\cC$ equipped with the categorified $n$-coherent algebraic structure. The slice $n$-category $\cC/A$ admits a canonical sliced $n$-coherent algebraic structure.
\end{definition}

We emphasize that both of these are merely guiding principles rather than theorems. For instance, applying the microcosm principle on the Lawvere theory of groups doesn't yield a meaningful notion of a group object inside a \enquote{groupoidal category}. The key obstruction for groups is the necessity of the copy operation $\Delta \colon G \rightarrow  G \times G$ and the deletion operation $G \rightarrow 1$ for defining the inverse axioms. 

More generally, the microcosm principle is incompatible with algebraic theories that demands some form of copying or deletion of data. In contrast, it is well suited for theories, such as that of monoids, which require no copying or deletion. We may, in general, expect the microcosm principle to be compatible with any \enquote{linear} multisorted algebraic theory. Furthermore, one should also expect the existence of sliced algebraic structures in such a case.

\subsection{Pseudomonoids and pseudomodules}

Monoidal categories, module categories, and their dualities are important in the latter part of this paper. Therefore, it is useful to introduce more abstract definitions of pseudomonoid objects and pseudomodule objects, both for conceptual clarity and to provide another useful example of the microcosm principle.

Pseudomonoids and pseudomodules are $1$-level higher weak categorification of the idea of a monoid object and a module over a monoid object. These objects can be defined in a general monoidal bicategory $\cV$. We shall view $\cV$ as a one-object tricategory \cite[Def.~4.1]{gurski2013tricatcoh}. Furthermore, all diagrams are read inside the one object cubical tricategory $\st \cV$ (see \cite[Thm.~8.4]{gurski2013tricatcoh} and \cref{appsection:mondiag} for explanation).

\begin{definition}
	A pseudomonoid in a monoidal bicategory $\cV$ is given by the data $(M,\mu, e, a, \lambda, \rho)$ such that
	\begin{enumerate}
		\item $M \in \cV$ is an object, $\mu$ is a $1$-cell $\mu \colon M \otimes M \rightarrow M$, and $e$ is a $1$-cell $e \colon I \rightarrow M$,
		\item $a,\lambda,\rho$ are invertible $2$-cells in $\cV$ which are defined using the following diagrams:
		\[\begin{tikzcd}
			{M \otimes (M \otimes M)} && {M \otimes M} && {M \otimes I} && \\
			\\
			{(M \otimes M) \otimes M} &&&& {M \otimes M} && M \\
			\\
			{M \otimes M} && M && {I \otimes M}
			\arrow["{1 \otimes \mu}", from=1-1, to=1-3]
			\arrow["\mu", from=1-3, to=5-3]
			\arrow["{1 \otimes e}"', from=1-5, to=3-5]
			\arrow[""{name=0, anchor=center, inner sep=0}, "{\rho_{\cV}}", from=1-5, to=3-7]
			\arrow["{a_{\cV}}", from=3-1, to=1-1]
			\arrow["{\mu \otimes 1}"', from=3-1, to=5-1]
			\arrow["\mu", from=3-5, to=3-7]
			\arrow["a"', between={0.2}{0.8}, Rightarrow, from=5-1, to=1-3]
			\arrow["\mu"', from=5-1, to=5-3]
			\arrow["{e \otimes 1}", from=5-5, to=3-5]
			\arrow[""{name=1, anchor=center, inner sep=0}, "{\lambda_{\cV}}"', from=5-5, to=3-7]
			\arrow["\rho", between={0.2}{0.8}, Rightarrow, from=3-5, to=0]
			\arrow["\lambda"', between={0.2}{0.8}, Rightarrow, from=3-5, to=1]
		\end{tikzcd}\cdpunct[8pt]{,}\]
		\item the associator $a$ and the unitors $\lambda,\rho$ are subject to the pentagon axiom (\cref{app:pentagon}) and the triangle axiom (\cref{app:triangle}).
	\end{enumerate}
\end{definition}

\begin{definition}
	A (lax) $1$-morphism $(f,\chi,\iota) \colon (M,\mu,e,a,\lambda,\rho) \rightarrow (M', \mu',e',a',\lambda',\rho')$ is given by
	\begin{enumerate}
		\item a $1$-cell $f : M \rightarrow M'$,
		\item two $2$-cell $\chi$ and $\iota$ as depicted by the following diagram:
		\[\begin{tikzcd}[ampersand replacement=\&]
			{M \otimes M} \&\& {M' \otimes M'} \&\& I \&\& \\
			\\
			M \&\& {M'} \&\& M \&\& {M'}
			\arrow["{f \otimes f}", from=1-1, to=1-3]
			\arrow["\mu"', from=1-1, to=3-1]
			\arrow["\chi"', between={0.2}{0.8}, Rightarrow, nfold, from=1-3, to=3-1]
			\arrow["{\mu'}", from=1-3, to=3-3]
			\arrow["e"', from=1-5, to=3-5]
			\arrow[""{name=0, anchor=center, inner sep=0}, "{e'}", from=1-5, to=3-7]
			\arrow["f"', from=3-1, to=3-3]
			\arrow["f"', from=3-5, to=3-7]
			\arrow["\iota"', between={0.2}{0.8}, Rightarrow, nfold, from=0, to=3-5]
		\end{tikzcd}\cdpunct[8pt]{,}\]
		\item such that the two $2$-cells satisfy the morphism hexagon (\cref{app:morphhexagon}), morphism left square (\cref{app:morphlsq}), and morphism right square (\cref{app:morphrsq}).
	\end{enumerate}
\end{definition}

\begin{definition}
	A $2$-morphism $r \colon (f,\chi_{f},\iota_{f}) \Rightarrow (g,\chi_{g},\iota_{g})$ is given by
	\begin{enumerate}
		\item a $2$-cell $r \colon f \Rightarrow g$,
		\item such that the following two identities hold:
		\[\begin{tikzcd}[ampersand replacement=\&]
			{M \otimes M} \&\& {M' \otimes M'} \&\& {M \otimes M} \&\& {M' \otimes M'} \\
			\&\&\& \equiv \\
			M \&\& {M'} \&\& M \&\& {M'}
			\arrow[""{name=0, anchor=center, inner sep=0}, "{g \otimes g}"{description}, curve={height=12pt}, from=1-1, to=1-3]
			\arrow[""{name=1, anchor=center, inner sep=0}, "{f \otimes f}", curve={height=-30pt}, from=1-1, to=1-3]
			\arrow["\mu"', from=1-1, to=3-1]
			\arrow["{\chi_{g}}"{description}, between={0.2}{0.8}, Rightarrow, nfold, from=1-3, to=3-1]
			\arrow["{\mu'}", from=1-3, to=3-3]
			\arrow["{f \otimes f}", from=1-5, to=1-7]
			\arrow["\mu"', from=1-5, to=3-5]
			\arrow["{\chi_{f}}"{description}, between={0.2}{0.8}, Rightarrow, nfold, from=1-7, to=3-5]
			\arrow["{\mu'}", from=1-7, to=3-7]
			\arrow["g"', from=3-1, to=3-3]
			\arrow[""{name=2, anchor=center, inner sep=0}, "f"{description}, curve={height=-12pt}, from=3-5, to=3-7]
			\arrow[""{name=3, anchor=center, inner sep=0}, "g"{description}, curve={height=30pt}, from=3-5, to=3-7]
			\arrow["{r \otimes r}", between={0.2}{0.8}, Rightarrow, nfold, from=1, to=0]
			\arrow["r", between={0.2}{0.8}, Rightarrow, nfold, from=2, to=3]
		\end{tikzcd}\]
		\[\begin{tikzcd}[ampersand replacement=\&]
			I \&\&\&\& I \&\& \\
			\&\&\& \equiv \\
			M \&\& {M'} \&\& M \&\& {M'}
			\arrow["e"', from=1-1, to=3-1]
			\arrow[""{name=0, anchor=center, inner sep=0}, "{e'}", from=1-1, to=3-3]
			\arrow["e"', from=1-5, to=3-5]
			\arrow[""{name=1, anchor=center, inner sep=0}, "{e'}", from=1-5, to=3-7]
			\arrow[""{name=2, anchor=center, inner sep=0}, "f"{description}, curve={height=-6pt}, from=3-1, to=3-3]
			\arrow[""{name=3, anchor=center, inner sep=0}, "g"', curve={height=30pt}, from=3-1, to=3-3]
			\arrow["g"', from=3-5, to=3-7]
			\arrow["{\iota_{f}}", between={0.2}{0.8}, Rightarrow, nfold, from=0, to=3-1]
			\arrow["{\iota_{g}}", between={0.2}{0.8}, Rightarrow, nfold, from=1, to=3-5]
			\arrow["r"', between={0.2}{0.8}, Rightarrow, nfold, from=2, to=3]
		\end{tikzcd}\cdpunct[8pt]{.}\]
	\end{enumerate}
\end{definition}

\begin{lemma}
	The data of pseudomonoids and their morphisms in $\cV$ assemble into a bicategory of pseudomonoids called $\nMon{1}(\cV)$. 
\end{lemma}

\begin{proof}
	This lemma primarily follows from \cite[Sec.~3]{gurski2014itericon} by considering one object version of iconic weak enrichment over a monoidal bicategory. 
\end{proof}

\begin{example}
	As expected, if $\cV$ is the monoidal bicategory $\nCat{1}$, then the pseudomonoid objects and their morphisms are usual monoidal categories and their lax monoidal functors.
\end{example}

Now, we shall describe braided pseudomonoid objects and their morphisms. For this, we require a braided monoidal bicategory. The details are ultimately irrelevant for us but can be found carefully constructed in \cite{gurski2011monbicat}. In particular, $\cV$ is a braided monoidal bicategory if it is a monoidal bicategory equipped with a family of adjoint equivalences $b_{x,y} : X \otimes Y \leftrightarrow Y \otimes X : b_{x,y}^{-1}$ that behaves as the braiding.

\begin{definition}
	A braided pseudomonoid in a braided monoidal bicategory $\cV$ is a pseudomonoid object $M$ together with the additional data of
	\begin{enumerate}
		\item a braiding $2$-cell $b$
		\[\begin{tikzcd}[ampersand replacement=\&]
			{M \otimes M} \&\& {M \otimes M} \\
			\\
			\& M
			\arrow["{b_{M,M}}", from=1-1, to=1-3]
			\arrow[""{name=0, anchor=center, inner sep=0}, "\mu"', from=1-1, to=3-2]
			\arrow["\mu", from=1-3, to=3-2]
			\arrow["b"', between={0.2}{0.8}, Rightarrow, nfold, from=0, to=1-3]
		\end{tikzcd}\cdpunct[8pt]{,}\]
		\item such that the two hexagon axioms (\cref{app:hexagon}) hold.
	\end{enumerate}
\end{definition}

\begin{definition}\leavevmode
	\begin{enumerate}
		\item A $1$-morphism $(f,\chi,i) : (M,b) \rightarrow (M',b')$ of braided pseudomonoid is a $1$-morphism of pseudomonoid together with the following addtional axiom:
		\[\begin{tikzcd}[ampersand replacement=\&,column sep=tiny]
			\& {M \otimes M} \&\&\&\& {M \otimes M} \& \\
			{M \otimes M} \&\& {M' \otimes M'} \&\& {M \otimes M} \&\& {M' \otimes M'} \\
			\&\&\& \equiv \\
			\& {M'\otimes M'} \\
			M \&\& {M'} \&\& M \&\& {M'}
			\arrow["{f \otimes f}", from=1-2, to=2-3]
			\arrow["{f \otimes f}", from=1-6, to=2-7]
			\arrow["\mu", from=1-6, to=5-5]
			\arrow["{b_{M,M}}", from=2-1, to=1-2]
			\arrow["{f \otimes f}"{description}, from=2-1, to=4-2]
			\arrow["\mu"', from=2-1, to=5-1]
			\arrow["{\mu'}", from=2-3, to=5-3]
			\arrow["{b_{M,M}}", from=2-5, to=1-6]
			\arrow[""{name=0, anchor=center, inner sep=0}, "\mu"', from=2-5, to=5-5]
			\arrow["\chi", between={0.2}{0.8}, Rightarrow, nfold, from=2-7, to=5-5]
			\arrow["{\mu'}", from=2-7, to=5-7]
			\arrow["{b_{M',M'}}"{description}, from=4-2, to=2-3]
			\arrow["\chi"', between={0.2}{0.8}, Rightarrow, nfold, from=4-2, to=5-1]
			\arrow[""{name=1, anchor=center, inner sep=0}, "{\mu'}"{description}, from=4-2, to=5-3]
			\arrow["f"', from=5-1, to=5-3]
			\arrow["f"', from=5-5, to=5-7]
			\arrow["b"{description}, between={0.2}{0.8}, Rightarrow, nfold, from=1-6, to=0]
			\arrow["{b'}"{description}, between={0.2}{0.8}, Rightarrow, nfold, from=2-3, to=1]
		\end{tikzcd}\cdpunct[8pt]{.}\]
		
		\item A $2$-morphism of braided pseudomonoid is the same data as that of a pseudomonoid.
	\end{enumerate}
\end{definition}

\begin{lemma}
	There is a bicategory of braided pseudomonoids $\nMon{1b}(\cV)$.
\end{lemma}

\begin{proof}
	Using that the $2$-morphism of braided pseudomonoids are the same as that of pseudomonoid, we define the compositional data and the coherence constraint of the bicategory $\nMon{1b}(\cV)$ to be the same as that of $\nMon{1}(\cV)$. Thus, we only need to verify that the composite of $1$-morphism of braided pseudomonoid also satisfy the additional braiding axiom. Suppressing notations and labeling, this follows directly from the following pasting equalities:
	\[\begin{tikzcd}[ampersand replacement=\&,column sep=small]
		\& MM \&\& PP \&\&\& MM \&\&\& \\
		MM \&\& PP \&\&\& MM \&\& NN \&\& PP \\
		\& NN \\
		M \& N \& P \&\&\& M \& N \&\& P \\
		\& \vequiv \&\&\&\&\&\& \vequiv \\
		MM \&\&\&\&\& MM \\
		\& NN \&\& PP \&\&\&\& NN \&\& PP \\
		MM \& NN \& PP \&\& \equiv \& MM \& NN \\
		\\
		M \& N \& P \&\&\& M \& N \&\& P
		\arrow[from=1-2, to=1-4]
		\arrow[from=1-4, to=4-3]
		\arrow[from=1-7, to=2-8]
		\arrow[from=1-7, to=2-10]
		\arrow[from=1-7, to=4-6]
		\arrow[from=2-1, to=1-2]
		\arrow[from=2-1, to=2-3]
		\arrow[from=2-1, to=3-2]
		\arrow[from=2-1, to=4-1]
		\arrow[from=2-3, to=1-4]
		\arrow[from=2-3, to=4-3]
		\arrow[from=2-6, to=1-7]
		\arrow[from=2-6, to=4-6]
		\arrow[from=2-8, to=2-10]
		\arrow[from=2-8, to=4-7]
		\arrow[from=2-10, to=4-9]
		\arrow[from=3-2, to=2-3]
		\arrow[from=3-2, to=4-2]
		\arrow[from=4-1, to=4-2]
		\arrow[from=4-2, to=4-3]
		\arrow[from=4-6, to=4-7]
		\arrow[from=4-7, to=4-9]
		\arrow[from=6-1, to=7-2]
		\arrow[from=6-1, to=7-4]
		\arrow[from=6-6, to=7-8]
		\arrow[from=6-6, to=7-10]
		\arrow[from=7-2, to=7-4]
		\arrow[from=7-4, to=10-3]
		\arrow[from=7-8, to=7-10]
		\arrow[from=7-8, to=10-7]
		\arrow[from=7-10, to=10-9]
		\arrow[from=8-1, to=6-1]
		\arrow[from=8-1, to=8-2]
		\arrow[from=8-1, to=10-1]
		\arrow[from=8-2, to=7-2]
		\arrow[from=8-2, to=8-3]
		\arrow[from=8-2, to=10-2]
		\arrow[from=8-3, to=7-4]
		\arrow[from=8-3, to=10-3]
		\arrow[from=8-6, to=6-6]
		\arrow[from=8-6, to=8-7]
		\arrow[from=8-6, to=10-6]
		\arrow[from=8-7, to=7-8]
		\arrow[from=8-7, to=10-7]
		\arrow[from=10-1, to=10-2]
		\arrow[from=10-2, to=10-3]
		\arrow[from=10-6, to=10-7]
		\arrow[from=10-7, to=10-9]
	\end{tikzcd}\cdpunct[8pt]{.}\]
\end{proof}

Finally, we can also describe symmetric pseudomonoid object inside a sylleptic (or symmetric) monoidal bicategory. This bicategory is properly defined in \cite{gurski2013symmonbicat} and is braided monoidal bicategory together with a $2$-cell $\upsilon_{X,Y} : b_{Y,X}b_{X,Y} \Rightarrow 1$ that essentially identifies $b_{X,Y}$ with $b^{-1}_{Y,X}$. 

\begin{definition}\leavevmode
	\begin{enumerate}
		\item A symmetric pseudomonoid is a braided pseudmonoid that satisfy the following additional axiom:
		\[\begin{tikzcd}[ampersand replacement=\&,column sep=tiny]
			{M \otimes M} \&\&\&\& {M \otimes M} \&\& {M \otimes M} \&\& {M \otimes M} \\
			\&\& {M \otimes M} \&\&\& \equiv \\
			\\
			\&\& M \&\&\&\&\& M
			\arrow[""{name=0, anchor=center, inner sep=0}, "{\varepsilon_{M \otimes M}}", from=1-1, to=1-5]
			\arrow["{b_{M,M}}"{description}, from=1-1, to=2-3]
			\arrow[""{name=1, anchor=center, inner sep=0}, "\mu"', from=1-1, to=4-3]
			\arrow[""{name=2, anchor=center, inner sep=0}, "\mu", from=1-5, to=4-3]
			\arrow["{\varepsilon_{M \otimes M}}", from=1-7, to=1-9]
			\arrow[""{name=3, anchor=center, inner sep=0}, "\mu"', from=1-7, to=4-8]
			\arrow[""{name=4, anchor=center, inner sep=0}, "\mu", from=1-9, to=4-8]
			\arrow["{b_{M,M}}"{description}, from=2-3, to=1-5]
			\arrow[""{name=5, anchor=center, inner sep=0}, from=2-3, to=4-3]
			\arrow["b", between={0.2}{0.8}, Rightarrow, nfold, from=1, to=5]
			\arrow["1"', between={0.2}{0.8}, Rightarrow, nfold, from=3, to=4]
			\arrow["{\upsilon_{M,M}}", between={0.2}{0.8}, Rightarrow, nfold, from=2-3, to=0]
			\arrow["b", between={0.2}{0.8}, Rightarrow, nfold, from=5, to=2]
		\end{tikzcd}\cdpunct[8pt]{.}\]
		\item The $1$-morphisms and $2$-morphisms of symmetric pseudomonoid object has no additional data.
	\end{enumerate}
\end{definition}

\begin{lemma}
	There is a bicategory of symmetric pseudomonoids $\nMon{1s}(\cV)$, which is formed as a full sub-bicategory of braided pseudomonoids.
\end{lemma}

For the case of $\cV=\nCat{1}$, we know that $\nCat{1}$ is a symmetric monoidal $2$-category. Thus, as expected, $\nMon{1b}(\nCat{1}) = \BrMonCat_l$ and $\nMon{1s}(\nCat{1})= \SymMonCat_l$. Similarly, it is possible to consider colax $1$-morphisms in the constructions of $\nMon{1}$ and its variants. In this case, we obtain the corresponding $2$-category of monoidal categories with colax monoidal functors.

With all the relevant definitions, we can now define pseudomodules over pseudomonoid objects.

\begin{definition}
	A left pseudomodule object over a pseudomonoid object $(M,\mu,e, a_{M},\lambda_{M},\rho_{M})$ in a monoidal bicategory $\cV$ is given by the collection $(\cZ,\triangleleft,a_{\cZ},\lambda_{\cZ})$ such that:
	\begin{enumerate}
		\item $\cZ \in \cV$ is an object, $\triangleleft$ defines a left action $1$-cell $M \otimes \cZ \xrightarrow{\triangleleft} \cZ$ in $\cV$,
		\item the associator $a_{\cZ}$ and left unitor $\lambda_{\cZ}$ are invertible $2$-cells in $\cV$ which are defined using the following diagrams:
		\[\begin{tikzcd}[ampersand replacement=\&]
			{M \otimes (M \otimes \cZ)} \&\& {M \otimes \cZ} \&\& {M \otimes \cZ} \&\& \cZ \\
			\\
			{(M \otimes M) \otimes \cZ} \\
			\\
			{M \otimes \cZ} \&\& \cZ \&\& {I \otimes \cZ}
			\arrow["{1 \otimes \triangleleft}", from=1-1, to=1-3]
			\arrow["\triangleleft", from=1-3, to=5-3]
			\arrow["\triangleleft", from=1-5, to=1-7]
			\arrow["{a_{\cV}}", from=3-1, to=1-1]
			\arrow["{\mu \otimes 1}"', from=3-1, to=5-1]
			\arrow["{a_{\cZ}}"', between={0.2}{0.8}, Rightarrow, nfold, from=5-1, to=1-3]
			\arrow["\triangleleft"', from=5-1, to=5-3]
			\arrow["{e \otimes 1}", from=5-5, to=1-5]
			\arrow[""{name=0, anchor=center, inner sep=0}, "{\lambda_{\cV}}"', from=5-5, to=1-7]
			\arrow["{\lambda_{\cZ}}"', between={0.2}{0.8}, Rightarrow, nfold, from=1-5, to=0]
		\end{tikzcd}\cdpunct[8pt]{,}\]
		\item the pentagon axiom (\cref{app:pentagon}) and the triangle axiom (\cref{app:triangle}) hold with replacing the corresponding copy of $M$ with $\cZ$. 
	\end{enumerate}
\end{definition}

\begin{definition}
	A right pseudomodule object over a pseudomonoid object $(M,\mu,e, a_{M},\lambda_{M},\rho_{M})$ in a monoidal bicatgory $\cV$ is given by the collection $(\cZ,\triangleright,a_{\cZ},\rho_{\cZ})$ such that:
	\begin{enumerate}
		\item $\cZ \in \cV$ is an object, $\triangleright$ defines a right action $1$-cell $\cZ \otimes M \xrightarrow{\triangleright} \cZ$ in $\cV$,
		\item the associator $a_{\cZ}$ and right unitor $\rho_{\cZ}$ are invertible $2$-morphisms in $\cV$ which are defined using the following diagrams:
		\[\begin{tikzcd}[ampersand replacement=\&]
			{ (\cZ \otimes M) \otimes M} \&\& {\cZ \otimes M} \&\& {\cZ \otimes M} \&\& \cZ \\
			\\
			{\cZ \otimes (M \otimes M)} \\
			\\
			{\cZ \otimes M} \&\& \cZ \&\& {\cZ \otimes I}
			\arrow["{\triangleright \otimes 1}", from=1-1, to=1-3]
			\arrow["{a_{\cV}}"', from=1-1, to=3-1]
			\arrow["{a_{\cZ}}", between={0.2}{0.8}, Rightarrow, nfold, from=1-3, to=5-1]
			\arrow["\triangleright", from=1-3, to=5-3]
			\arrow["\triangleright", from=1-5, to=1-7]
			\arrow["{1 \otimes \mu}"', from=3-1, to=5-1]
			\arrow["\triangleright"', from=5-1, to=5-3]
			\arrow["{1 \otimes e}", from=5-5, to=1-5]
			\arrow[""{name=0, anchor=center, inner sep=0}, "{\rho_{\cV}}"', from=5-5, to=1-7]
			\arrow["{\rho_{\cZ}}"', between={0.2}{0.8}, Rightarrow, nfold, from=1-5, to=0]
		\end{tikzcd}\cdpunct[8pt]{,}\]
		\item the pentagon axiom (\cref{app:pentagon}) and the triangle axiom (\cref{app:triangle}) hold with replacing the corresponding copy of $M$ with $\cZ$. 
	\end{enumerate}
\end{definition}

\begin{definition}
	A bipseudomodule $\cZ \in \cV$ between two pseudomonoids $(M,\mu_{M},e_{M}, a_{M},\lambda_{M},\rho_{M})$ and $(N,\mu_{N},e_{N}, a_{N},\lambda_{N},\rho_{N})$ is given by 
	\begin{enumerate}
		\item a left pseudomodule $(\cZ,\triangleleft,a_{(\cZ,M)},\lambda_{\cZ})$ over $M$,
		\item a right pseudomodule $(\cZ,\triangleright,a_{(\cZ,N)},\rho_{\cZ})$ over $N$,
		\item an invertible balancing $2$-cell $l$ in $\cV$ which is defined using the following diagram:
		\[\begin{tikzcd}[ampersand replacement=\&,cramped]
			{ (M \otimes \cZ) \otimes N} \&\&\&\& {M \otimes (\cZ \otimes N)} \\
			\\
			{\cZ \otimes N} \&\&\&\& {M \otimes \cZ} \\
			\\
			\&\& \cZ
			\arrow["{a_{\cV}}", from=1-1, to=1-5]
			\arrow["{\triangleleft_{\cZ}\otimes 1}"', from=1-1, to=3-1]
			\arrow["{1 \otimes \triangleright_{\cZ}}", from=1-5, to=3-5]
			\arrow["l", between={0.2}{0.8}, Rightarrow, nfold, from=3-1, to=3-5]
			\arrow["{\triangleright_{\cZ}}"', from=3-1, to=5-3]
			\arrow["{\triangleleft_{\cZ}}", from=3-5, to=5-3]
		\end{tikzcd}\cdpunct[8pt]{,}\]
		\item two instances of the pentagon axiom (\cref{app:pentagon}) hold, one for $(M,\cZ,N,N)$ and the other for $(M,M,\cZ,N)$.  
	\end{enumerate}
\end{definition}

All the axioms for left, right, and bi-pseudomodule are motivated by viewing monoidal categories and module categories as parts of a two objects bicategory $K$. Furthermore, in the case of $\cV=\nCat{1}$, these provide the usual definitions.

	\section{Generalized involutions and inverses} \label{sec:ginvs}
	We aim to characterize involutions and inverses in a more general setting using the microcosm principle. In a certain sense, the data of a \textit{single} mutually invertible pair of arrows or a \textit{single} involution lies within the scope of this principle. This allows one to define relative inverses and involutions which are quite essential in describing dualities of monoidal structures. The first part primarily focus on structures internal to bicategories, whereas the latter part shall explore structures internal to tricategories.

\subsection{\tpdf{$1$-dimensional variants}} \label{sub:1diminv}

We begin by briefly reviewing $0$-dimensional involutions for motivation. As this theory is well known, we omit most of the details. The interested reader should refer to \cite{yauinvcat2020}.

\begin{definition} \leavevmode
	\begin{enumerate}
		\item  A strict involutive object in $\cC$ is given by a morphism $f\colon x \rightarrow x$ such that $f \circ f = 1$.
		\item A morphism of involutions $(x,f) \rightarrow (y,g)$ is an arrow $r\colon x \rightarrow y$ such that $gr=rf$ holds.
	\end{enumerate}
\end{definition}

Involutive objects and their morphisms define a category of involutions $\nInv{0}(\cC)$. This construction extends to a $2$-functor $\nInv{0} \colon \nCat{1} \rightarrow \nCat{1}$, which comes equipped with two canonical natural $2$-transformation, $ U \colon \nInv{0} \Rightarrow 1_{\nCat{1}}$ given by the forgetful action, and its section $i \colon 1_{\nCat{1}} \Rightarrow \nInv{0}$ given by inclusion as identity involution.

Furthermore, the category $\nInv{0}(\nInv{0}(\cC))$ has as objects a pair of involutions $f,g \colon  x \rightarrow x$ such that $fg=gf$. This implies that we have a natural $2$-transformation $\nInv{0}^2 \Rightarrow \nInv{0}$ given by $(f,g) \mapsto (fg)$. In particular, $\nInv{0}$ is a strict $2$-Monad on $\nCat{1}$.

\begin{example}[$\op_{n,k}$ duality of $\nCat{n}$]
	The $1$-functors $\op_{n,k} \colon \nCat{n}_{\leq 1} \rightarrow \nCat{n}_{\leq 1}$ satisfy the desired axiom $\op_{n,k}^2=1_{\nCat{n}}$. This is a strict $1$-involutive object.
\end{example}

The phrase strict involutive object already suggests the general definition of a weak ($1$-dimensional) involutive object. In particular, we can define weak (or pseudo-) involutive objects inside an arbitrary bicategory $K$. All the corresponding data and axioms shall be written using pasting diagrams, which are to be interpreted via the evaluation functor $e \colon \st K \rightarrow K$. See \cref{appsection:strictification} for more details.

\begin{definition}[pseudo involutive object]
	Let $K$ be a bicategory. We define a weak/pseudo involutive object by the data of
	\begin{enumerate}
		\item a $1$-cell $f \colon x \rightarrow x \in K$,
		\item an invertible $2$-cell
		\[\begin{tikzcd}[ampersand replacement=\&]
			\& x \& \\
			x \&\& x
			\arrow["f", from=1-2, to=2-3]
			\arrow["f", from=2-1, to=1-2]
			\arrow[""{name=0, anchor=center, inner sep=0}, "{\varepsilon_x}"', from=2-1, to=2-3]
			\arrow["\alpha", between={0}{0.8}, Rightarrow, nfold, from=1-2, to=0]
		\end{tikzcd}\cdpunct[8pt]{,}\]
		\item such that the $2$-cell $\alpha$ satisfy an additional axiom as depicted by the following diagram:
		\[\begin{tikzcd}[ampersand replacement=\&]
			x \&\& x \&\&\& x \&\& x \& \\
			\&\&\&\& \equiv \\
			\& x \&\& x \&\&\& x \&\& x
			\arrow[""{name=0, anchor=center, inner sep=0}, "{\varepsilon_x}", from=1-1, to=1-3]
			\arrow["f"', from=1-1, to=3-2]
			\arrow["f", from=1-3, to=3-4]
			\arrow["{\varepsilon_x}", from=1-6, to=1-8]
			\arrow["f"', from=1-6, to=3-7]
			\arrow[""{name=1, anchor=center, inner sep=0}, "f"{description}, from=1-6, to=3-9]
			\arrow["f", from=1-8, to=3-9]
			\arrow["f", from=3-2, to=1-3]
			\arrow[""{name=2, anchor=center, inner sep=0}, "{\varepsilon_x}"', from=3-2, to=3-4]
			\arrow["{\varepsilon_x}"', from=3-7, to=3-9]
			\arrow["{\alpha^{-1}}"', between={0.2}{0.8}, Rightarrow, nfold, from=0, to=3-2]
			\arrow["\alpha", between={0.2}{0.8}, Rightarrow, nfold, from=1-3, to=2]
			\arrow[between={0.2}{0.8}, Rightarrow, nfold, from=1, to=3-7]
			\arrow[between={0.2}{0.7}, Rightarrow, nfold, from=1-8, to=1]
		\end{tikzcd}\cdpunct[8pt]{.}\]
	\end{enumerate}
\end{definition}

\begin{example}
	Every strict involutive object (where $f^2=1$) in a strict $2$-category is trivially a weak involutive object.
\end{example}

One might wonder why do we require the additional identities in the specification of a weak involutive object. The key point is that the structure is incoherent without the additional axiom, analogous to how specifying a monoidal category without its pentagon or triangle identity is meaningless. This extra identity ensures that all reasonable ways of going from $f^n \Rightarrow f^m$ are equivalent. In fact, there are four equivalent way of expressing this additional axiom. The proof is relatively simple and uses the pasting of the $2$-cells $\alpha,\alpha^{-1}$, or the unitors appropriately.

\begin{lemma}
	The data of an involution in a bicategory is equivalently expressed as a $1$-cell $f\colon x \rightarrow x$ and an invertible $2$-cell $\alpha \colon f^2 \Rightarrow 1$ satisfying any one of the following axioms:
	\begin{enumerate}
		\item 
		\[\begin{tikzcd}[ampersand replacement=\&]
			x \&\& x \&\&\& x \&\& x \& \\
			\&\&\&\& \equiv \\
			\& x \&\& x \&\&\& x \&\& x
			\arrow[""{name=0, anchor=center, inner sep=0}, "{\varepsilon_x}", from=1-1, to=1-3]
			\arrow["f"', from=1-1, to=3-2]
			\arrow["f", from=1-3, to=3-4]
			\arrow["{\varepsilon_x}", from=1-6, to=1-8]
			\arrow["f"', from=1-6, to=3-7]
			\arrow[""{name=1, anchor=center, inner sep=0}, "f"{description}, from=1-6, to=3-9]
			\arrow["f", from=1-8, to=3-9]
			\arrow["f", from=3-2, to=1-3]
			\arrow[""{name=2, anchor=center, inner sep=0}, "{\varepsilon_x}"', from=3-2, to=3-4]
			\arrow["{\varepsilon_x}"', from=3-7, to=3-9]
			\arrow["{\alpha^{-1}}"', between={0.2}{0.8}, Rightarrow, nfold, from=0, to=3-2]
			\arrow["\alpha", between={0.2}{0.8}, Rightarrow, nfold, from=1-3, to=2]
			\arrow[between={0.2}{0.8}, Rightarrow, nfold, from=1, to=3-7]
			\arrow[between={0.2}{0.7}, Rightarrow, nfold, from=1-8, to=1]
		\end{tikzcd}\cdpunct[8pt]{.}\]
		\item 
		\[\begin{tikzcd}[ampersand replacement=\&]
			\&\&\& x \&\&\&\& x \&\&\& \\
			\& x \&\&\&\& \equiv \&\&\&\& x \\
			x \&\&\&\& x \&\& x \&\&\&\& x
			\arrow["f", from=1-4, to=3-5]
			\arrow["f", from=1-8, to=2-10]
			\arrow["f", from=2-2, to=1-4]
			\arrow[""{name=0, anchor=center, inner sep=0}, "{\varepsilon_x}"{description}, from=2-2, to=3-5]
			\arrow["f", from=2-10, to=3-11]
			\arrow["f", from=3-1, to=2-2]
			\arrow[""{name=1, anchor=center, inner sep=0}, "f"', from=3-1, to=3-5]
			\arrow["f", from=3-7, to=1-8]
			\arrow[""{name=2, anchor=center, inner sep=0}, "{\varepsilon_x}"{description}, from=3-7, to=2-10]
			\arrow[""{name=3, anchor=center, inner sep=0}, "f"', from=3-7, to=3-11]
			\arrow["\alpha", between={0}{0.8}, Rightarrow, nfold, from=1-4, to=0]
			\arrow["\alpha", between={0}{0.8}, Rightarrow, nfold, from=1-8, to=2]
			\arrow[between={0.2}{0.8}, Rightarrow, nfold, from=2-2, to=1]
			\arrow[between={0.2}{0.8}, Rightarrow, nfold, from=2-10, to=3]
		\end{tikzcd}\cdpunct[8pt]{.}\]
		\item 
		\[\begin{tikzcd}[ampersand replacement=\&]
			\& x \&\& x \&\&\& x \&\& x \\
			\&\&\&\& \equiv \\
			x \&\& x \&\&\& x \&\& x
			\arrow[""{name=0, anchor=center, inner sep=0}, "{\varepsilon_x}", from=1-2, to=1-4]
			\arrow["f", from=1-2, to=3-3]
			\arrow["{\varepsilon_x}", from=1-7, to=1-9]
			\arrow["f", from=3-1, to=1-2]
			\arrow[""{name=1, anchor=center, inner sep=0}, "{\varepsilon_x}"', from=3-1, to=3-3]
			\arrow["f"', from=3-3, to=1-4]
			\arrow["f", from=3-6, to=1-7]
			\arrow[""{name=2, anchor=center, inner sep=0}, "f"{description}, from=3-6, to=1-9]
			\arrow["{\varepsilon_x}"', from=3-6, to=3-8]
			\arrow["f"', from=3-8, to=1-9]
			\arrow["\alpha", between={0.2}{0.8}, Rightarrow, nfold, from=1-2, to=1]
			\arrow["{\alpha^{-1}}", between={0.2}{0.8}, Rightarrow, nfold, from=0, to=3-3]
			\arrow[between={0.2}{0.6}, Rightarrow, nfold, from=1-7, to=2]
			\arrow[between={0.4}{0.8}, Rightarrow, nfold, from=2, to=3-8]
		\end{tikzcd}\cdpunct[8pt]{.}\]
		\item 
		\[\begin{tikzcd}[ampersand replacement=\&]
			x \&\&\&\&\&\&\&\&\&\& x \\
			\& x \&\&\& x \&\& x \&\&\& x \\
			\&\&\&\&\& \equiv \\
			\&\&\& x \&\&\&\& x
			\arrow["f"', from=1-1, to=2-2]
			\arrow[""{name=0, anchor=center, inner sep=0}, "f", from=1-1, to=2-5]
			\arrow[""{name=1, anchor=center, inner sep=0}, "{\varepsilon_x}"{description}, from=2-2, to=2-5]
			\arrow["f"', from=2-2, to=4-4]
			\arrow[""{name=2, anchor=center, inner sep=0}, "f", from=2-7, to=1-11]
			\arrow[""{name=3, anchor=center, inner sep=0}, "{\varepsilon_x}"{description}, from=2-7, to=2-10]
			\arrow["f"', from=2-7, to=4-8]
			\arrow["f"', from=2-10, to=1-11]
			\arrow["f"', from=4-4, to=2-5]
			\arrow["f"', from=4-8, to=2-10]
			\arrow[between={0.2}{0.8}, Rightarrow, nfold, from=0, to=2-2]
			\arrow["{\alpha^{-1}}"', between={0.2}{0.8}, Rightarrow, nfold, from=1, to=4-4]
			\arrow[between={0.2}{0.8}, Rightarrow, nfold, from=2, to=2-10]
			\arrow["{\alpha^{-1}}"', between={0.2}{0.8}, Rightarrow, nfold, from=3, to=4-8]
		\end{tikzcd}\cdpunct[8pt]{.}\]
	\end{enumerate} 
\end{lemma}

The $2$nd and the $4$th form of the additional axiom are particularly interesting as they allow one to define a lax/colax variant for involutive object in a (lax/colax) bicategory. Thus, more generally we can work with \emph{lax involutive objects}. It is now instructive to define $1$-morphisms and $2$-morphisms of these involutions explicitly. 

\begin{definition}
	Let $(x,f,\alpha)$ and $(y,g,\beta)$ be two (lax) involutive objects in a bicategory $K$. A (lax) $1$-morphism of involution $(r,u) \colon (x,f,\alpha) \rightarrow (y,g,\beta)$ is
	\begin{enumerate}
		\item a $1$-cell $r \colon x \rightarrow y \in K$ 
		\item a $2$-cell 
		\[\begin{tikzcd}[ampersand replacement=\&]
			x \&\& y \\
			\\
			x \&\& y
			\arrow["r", from=1-1, to=1-3]
			\arrow["f"', from=1-1, to=3-1]
			\arrow["u", between={0.2}{0.8}, Rightarrow, nfold, from=1-3, to=3-1]
			\arrow["g", from=1-3, to=3-3]
			\arrow["r"', from=3-1, to=3-3]
		\end{tikzcd}\cdpunct[8pt]{,}\]
		\item such that the following axiom holds:
		\[\begin{tikzcd}[ampersand replacement=\&]
			x \&\& y \&\&\&\&\&\&\& \\
			\&\&\&\&\& x \&\& y \\
			\& x \&\& y \& \equiv \&\&\&\&\& y \\
			\&\&\&\&\& x \&\& y \\
			x \&\& y
			\arrow["r", from=1-1, to=1-3]
			\arrow["f", from=1-1, to=3-2]
			\arrow[""{name=0, anchor=center, inner sep=0}, "{\varepsilon_x}"', from=1-1, to=5-1]
			\arrow["u", between={0.2}{0.8}, Rightarrow, nfold, from=1-3, to=3-2]
			\arrow["g", from=1-3, to=3-4]
			\arrow["r", from=2-6, to=2-8]
			\arrow["{\varepsilon_x}"', from=2-6, to=4-6]
			\arrow[""{name=1, anchor=center, inner sep=0}, "r"{description}, from=2-6, to=4-8]
			\arrow["g", from=2-8, to=3-10]
			\arrow[""{name=2, anchor=center, inner sep=0}, "{\varepsilon_y}"', from=2-8, to=4-8]
			\arrow["r"', from=3-2, to=3-4]
			\arrow["f", from=3-2, to=5-1]
			\arrow["u", between={0.2}{0.8}, Rightarrow, nfold, from=3-4, to=5-1]
			\arrow["g", from=3-4, to=5-3]
			\arrow["g", from=3-10, to=4-8]
			\arrow["r"', from=4-6, to=4-8]
			\arrow["r"', from=5-1, to=5-3]
			\arrow[between={0.2}{0.8}, Rightarrow, nfold, from=1, to=4-6]
			\arrow[between={0.2}{0.8}, Rightarrow, nfold, from=2-8, to=1]
			\arrow["\alpha"', between={0.2}{0.8}, Rightarrow, nfold, from=3-2, to=0]
			\arrow["\beta"', between={0.2}{0.8}, Rightarrow, nfold, from=3-10, to=2]
		\end{tikzcd}\cdpunct[8pt]{.}\]
	\end{enumerate}
\end{definition}

\begin{definition}
	Let $(r,u)$ and $(t,v)$ be two $1$-morphisms of involution between $(x,f,\alpha)$ and $(y,g,\beta)$. A $2$-morphism $\gamma \colon (r,u) \Rightarrow (t,v)$ is
	
	\begin{enumerate}
		\item a $2$-cell $r \Rightarrow t \in K$,
		
		\item which satisfy the following axiom
		\[\begin{tikzcd}[ampersand replacement=\&]
			x \&\& y \&\& x \&\& y \\
			\&\&\& \equiv \\
			x \&\& y \&\& x \&\& y
			\arrow["r"', from=1-1, to=1-3]
			\arrow["f"', from=1-1, to=3-1]
			\arrow["u", between={0.2}{0.8}, Rightarrow, nfold, from=1-3, to=3-1]
			\arrow["g", from=1-3, to=3-3]
			\arrow[""{name=0, anchor=center, inner sep=0}, "t"', from=1-5, to=1-7]
			\arrow[""{name=1, anchor=center, inner sep=0}, "r", curve={height=-30pt}, from=1-5, to=1-7]
			\arrow["f"', from=1-5, to=3-5]
			\arrow["v", between={0.2}{0.8}, Rightarrow, nfold, from=1-7, to=3-5]
			\arrow["g", from=1-7, to=3-7]
			\arrow[""{name=2, anchor=center, inner sep=0}, "r", from=3-1, to=3-3]
			\arrow[""{name=3, anchor=center, inner sep=0}, "t"', curve={height=30pt}, from=3-1, to=3-3]
			\arrow["t", from=3-5, to=3-7]
			\arrow["\gamma", between={0.2}{0.8}, Rightarrow, nfold, from=1, to=0]
			\arrow["\gamma", between={0.2}{0.8}, Rightarrow, nfold, from=2, to=3]
		\end{tikzcd}\cdpunct[8pt]{.}\]
	\end{enumerate}
\end{definition}

As usual, all of this data assembles into a \emph{bicategory of involutions} $\nInv{1}(K)$ associated to any bicategory $K$. This follows by simply defining all the relevant compositional data strictly in $\st K$ and then descending those diagrams onto $K$ using the evaluation functor $e \colon \st K \rightarrow K$ appropriately. We do note that the identity $1$-morphisms in $\nInv{1}(K)$ are given by the $1$-cells $1_x \colon x\rightarrow x$, which are not the empty strings.

\begin{lemma}
	There is a bicategory of involutions $\nInv{1}(K)$ corresponding to each bicategory $K$.
\end{lemma}

\begin{remark}
	With a lengthy computation and using strifications appropriately, one could show that $\nInv{1}$ is a (strict) trifunctor  $\nInv{1} \colon \BiCat \rightarrow \BiCat$ on the tricategory of $\BiCat$. One could also show that there is a (strict) tritransformation $U \colon \nInv{1} \Rightarrow 1_{\BiCat}$ and a (strict) tritransformation $i \colon 1_{\BiCat} \Rightarrow \nInv{1}$. However, we shall not go into the details here.
\end{remark}

We have now provided enough discussion about the theory of involutions and want to motivate a more general definition of an involutive object. Consider again the $2$-functor $\op_1 \colon \op_2(\nCat{1}) \rightarrow \nCat{1}$. Here, the lower dimensional involution has its domain altered by a higher dimensional involution. This suggests the definition of a relative involutive object. A similar account of these can also be found in \cite{yauinvcat2020}. Hence, we shall be brief in describing the known details and provide a proof only for the missing slice construction.

\begin{definition}
	Let $(\cC,T,\alpha)$ be a pseudo (or lax) involutive object in $\nCat{1}$. We define an involutive object relative to $T$ as a morphism $f \colon Tx \rightarrow x \in \cC$ such that the following diagram commutes:
	\[\begin{tikzcd}[ampersand replacement=\&]
		{T^2x} \&\& Tx \\
		\\
		x \&\& x
		\arrow["Tf", from=1-1, to=1-3]
		\arrow["{\alpha_x}"', from=1-1, to=3-1]
		\arrow["f", from=1-3, to=3-3]
		\arrow["1"', from=3-1, to=3-3]
	\end{tikzcd}\cdpunct[8pt]{.}\]
\end{definition}

\begin{remark}
	To relate with the previous definition, we define the adjunct morphism by $f^{-1} \colon x \xrightarrow{\alpha^{-1}_x} T^2x \xrightarrow{Tf} Tx$. From this, it follows that both $f^{-1} \circ f = 1 $ and $1 = f \circ f^{-1}$ are satisfied. 
	Furthermore, if one of the equation is satisfied then the other one is automatically satisfied due to the additional axiom of $T$. However, the invertibility of $\alpha$ is essential for this, as $f^{-1}$ is not defined without $\alpha^{-1}$.
\end{remark}

\begin{lemma}
	Let $(C,T,\alpha)$ be a lax involutive object in $\nCat{1}$ and $(x,f \colon Tx \rightarrow x)$ be an involutive object relative to $T$. Then $C/x$ admits a canonical lax involutive structure using:
	\[\begin{tikzcd}[ampersand replacement=\&]
		{C/x} \&\& {C/x} \& {\left(\alpha_{/}\right)_{c} :} \& {T^2c} \&\& c \\
		c \&\& Tc \&\& {T^2x} \\
		x \&\& Tx \&\& Tx \\
		\&\& x \&\&\&\& x
		\arrow["{T/f}", from=1-1, to=1-3]
		\arrow["{\alpha_c}", from=1-5, to=1-7]
		\arrow[from=1-5, to=2-5]
		\arrow[from=1-7, to=4-7]
		\arrow[""{name=0, anchor=center, inner sep=0}, from=2-1, to=3-1]
		\arrow[""{name=1, anchor=center, inner sep=0}, from=2-3, to=3-3]
		\arrow["Tf"', from=2-5, to=3-5]
		\arrow["{\alpha_x}", from=2-5, to=4-7]
		\arrow["f", from=3-3, to=4-3]
		\arrow["f"', from=3-5, to=4-7]
		\arrow[between={0.2}{0.8}, maps to, from=0, to=1]
	\end{tikzcd}\cdpunct[8pt]{.}\]
\end{lemma}

\begin{proof}
	The following identity of diagrams implies the result immediately:
	\[\begin{tikzcd}[ampersand replacement=\&]
		{T^3c} \&\&\& Tc \&\& {T^3c} \&\&\& Tc \\
		{T^3x} \&\&\&\& \equiv \& {T^3x} \&\&\& Tx \\
		{T^2x} \&\&\&\&\& {T^2x} \\
		\&\&\& Tx \\
		Tx \&\&\& x \&\& Tx \&\&\& x
		\arrow["{T\alpha_c}", from=1-1, to=1-4]
		\arrow[from=1-1, to=2-1]
		\arrow[from=1-4, to=4-4]
		\arrow["{\alpha_{Tc}}", from=1-6, to=1-9]
		\arrow[from=1-6, to=2-6]
		\arrow[from=1-9, to=2-9]
		\arrow["{T^2f}"', from=2-1, to=3-1]
		\arrow["{T\alpha_x}"', from=2-1, to=4-4]
		\arrow["{\alpha_{Tx}}"', from=2-6, to=2-9]
		\arrow["{T^2f}"', from=2-6, to=3-6]
		\arrow["f", from=2-9, to=5-9]
		\arrow["Tf"', from=3-1, to=5-1]
		\arrow["{\alpha_x}"', from=3-1, to=5-4]
		\arrow["Tf"', from=3-6, to=5-6]
		\arrow["{\alpha_x}", from=3-6, to=5-9]
		\arrow["f", from=4-4, to=5-4]
		\arrow["f"', from=5-1, to=5-4]
		\arrow["f"', from=5-6, to=5-9]
	\end{tikzcd}\cdpunct[8pt]{.}\]
\end{proof}

As before, there is also a $1$-category $\nInv{0}(C,T,\alpha)$, whose objects are relative involutive objects and whose morphisms are maps $r \colon (x,f) \rightarrow (g,y)$ such that $g \circ Tr = r \circ f$. We can also describe this modified $\nInv{0}$ as a $2$-functor that extends the original functor.

\begin{lemma} \phantom{}
	\begin{enumerate}
		\item There exists a $2$-functor $\nInv{0} \colon \nInv{1}(\nCat{1}) \rightarrow \nCat{1}$ extending the construction of involutive objects along the identity inclusion tritransformation component $i_{\nCat{1}} \colon \nCat{1} \hookrightarrow \nInv{1}(\nCat{1})$.
		\item This further extends to the $2$-functor $\nInv{0} \colon\nInv{1}(\nCat{1}) \rightarrow \nInv{1}(\nCat{1})$ along the forgetful tritransformation component $U_{\nCat{1}} \colon \nInv{1}(\nCat{1}) \rightarrow \nCat{1}$.
		\item The natural $2$-transformation $U$ itself extends to $U \colon \nInv{0} \Rightarrow 1_{\nInv{1}(\nCat{1})}$. However, it doesn't admit any section as before. In particular, this $\nInv{0}$ is not a Monad.
	\end{enumerate}
\end{lemma}

\begin{proof}
	The proof of $(1),(2)$ and $(3)$ mainly follows from a slight modification of the arguments found in \cite[Props.~2.6.2, 2.6.3, and 2.6.5]{yauinvcat2020}.
\end{proof}

We now briefly describe the analogous definitions for inverse structures. Our goal is to provide a better perspective on these definitions that aligns with the microcosm principle and coherence principle.

\begin{definition}[pseudo inverse objects]
	Let $K$ be a bicategory. We define a tuple $(f\colon x \rightarrow y,g\colon y \rightarrow x,\eta,\epsilon)$ to be an inverse object whenever it is an internal adjoint equivalence. That is, the following two triangle identities hold:
	\[\begin{tikzcd}[ampersand replacement=\&]
		x \&\& x \&\&\& x \&\& x \& \\
		\&\&\&\& \equiv \\
		\& y \&\& y \&\&\& y \&\& y \\
		\& x \&\& x \&\&\& x \&\& x \\
		\&\&\&\& \equiv \\
		y \&\& y \&\&\& y \&\& y
		\arrow[""{name=0, anchor=center, inner sep=0}, "{\varepsilon_{x}}", from=1-1, to=1-3]
		\arrow["f"', from=1-1, to=3-2]
		\arrow["f", from=1-3, to=3-4]
		\arrow["{\varepsilon_{x}}", from=1-6, to=1-8]
		\arrow["f"', from=1-6, to=3-7]
		\arrow["f"{description}, from=1-6, to=3-9]
		\arrow["f", from=1-8, to=3-9]
		\arrow["g"{description}, from=3-2, to=1-3]
		\arrow[""{name=1, anchor=center, inner sep=0}, "{\varepsilon_{y}}"', from=3-2, to=3-4]
		\arrow["{\varepsilon_{y}}"', from=3-7, to=3-9]
		\arrow[""{name=2, anchor=center, inner sep=0}, "{\varepsilon_x}", from=4-2, to=4-4]
		\arrow["f"{description}, from=4-2, to=6-3]
		\arrow["{\varepsilon_x}", from=4-7, to=4-9]
		\arrow["g", from=6-1, to=4-2]
		\arrow[""{name=3, anchor=center, inner sep=0}, "{\varepsilon_y}"', from=6-1, to=6-3]
		\arrow["g"', from=6-3, to=4-4]
		\arrow["g", from=6-6, to=4-7]
		\arrow["g"{description}, from=6-6, to=4-9]
		\arrow["{\varepsilon_y}"', from=6-6, to=6-8]
		\arrow["g"', from=6-8, to=4-9]
		\arrow["\eta"{description}, between={0.2}{0.8}, Rightarrow, nfold, from=0, to=3-2]
		\arrow["\epsilon"{description}, between={0.2}{0.8}, Rightarrow, nfold, from=1-3, to=1]
		\arrow["\epsilon"{description}, between={0.2}{0.8}, Rightarrow, nfold, from=4-2, to=3]
		\arrow["\eta"{description}, between={0.2}{0.8}, Rightarrow, nfold, from=2, to=6-3]
	\end{tikzcd}\cdpunct[8pt]{.}\]
\end{definition}

\begin{example}[pseudo inverse objects in $\nCat{1}$]
	A pair of functors $F \colon \cC \rightarrow \cD$ and $G \colon \cD \rightarrow \cC$ specify a weak inverse object whenever they are part of an adjoint equivalence of categories. 
\end{example}

We strongly believe that this shift in perspective of viewing an adjoint equivalence as a coherent weak inverse is important and not merely pedantic. Indeed, assume that data of an adjoint equivalence of categories were genuinely different from the data of an equivalence of categories. In this hypothetical situation, working with equivalence of categories would have been no better than working with monoidal categories without the pentagon or triangle axioms.

With this idea, we can now work with general adjunctions instead of adjoint equivalences. One could view the data $\eta \colon 1 \Rightarrow GF$ as a laxly categorified version of $f$ being right inverse of $g$ (or $gf=1$), and the similar goes for $\epsilon \colon FG \Rightarrow 1$ as laxly categorified left inverse $f$. Using this viewpoint and \textit{convention}, the data of an adjunction $F \dashv G$ simultaneously expresses that $F$ is a two-sided lax inverse of $G$ and that $G$ is a two-sided colax inverse of $F$.

There is also an analogous notion of $1$-morphisms and $2$-morphisms defining the bicategory of adjunctions in $K$. To avoid notation overload, we will not label the data.

\begin{definition}
	A (lax) $1$-morphism $(x,y,f,g,\eta,\epsilon) \rightarrow (x',y',f',g',\eta',\epsilon')$ consists of
	\begin{enumerate}
		\item two $1$-cells $x \rightarrow x'$ and $y \rightarrow y'$,
		\item two $2$-cells\footnote{These $2$-cells are not mates. One may wish to closely compare this definition with \cite{chengmates2014}. However, the two immediate identities ensures that the $2$nd cell is the inverse of the mate of the first one. In particular, it is a form of Beck-Chevalley condition. Moreover, the colax variant of the $1$-morphisms will yield the dual Beck-Chevalley condition.} as depicted by the following diagram:
		\[\begin{tikzcd}[ampersand replacement=\&]
			x \&\& {x'} \& y \&\& {y'} \\
			\\
			y \&\& {y'} \& x \&\& {x'}
			\arrow[from=1-1, to=1-3]
			\arrow["f"', from=1-1, to=3-1]
			\arrow[between={0.2}{0.8}, Rightarrow, nfold, from=1-3, to=3-1]
			\arrow["{f'}", from=1-3, to=3-3]
			\arrow[from=1-4, to=1-6]
			\arrow["g"', from=1-4, to=3-4]
			\arrow[between={0.2}{0.8}, Rightarrow, nfold, from=1-6, to=3-4]
			\arrow["{g'}", from=1-6, to=3-6]
			\arrow[from=3-1, to=3-3]
			\arrow[from=3-4, to=3-6]
		\end{tikzcd}\cdpunct[8pt]{,}\]
		\item such that the following two identities holds:
		\[\begin{tikzcd}[ampersand replacement=\&]
			y \&\& {y'} \&\&\& y \&\& {y'} \& \\
			\\
			\& x \&\& {x'} \& \equiv \&\&\&\& {x'} \\
			\\
			y \&\& {y'} \&\&\& y \&\& {y'}
			\arrow[from=1-1, to=1-3]
			\arrow["g"', from=1-1, to=3-2]
			\arrow[""{name=0, anchor=center, inner sep=0}, "{\varepsilon_y}"', from=1-1, to=5-1]
			\arrow[between={0.2}{0.8}, Rightarrow, nfold, from=1-3, to=3-2]
			\arrow["{g'}", from=1-3, to=3-4]
			\arrow[from=1-6, to=1-8]
			\arrow["{\varepsilon_y}"', from=1-6, to=5-6]
			\arrow[""{name=1, anchor=center, inner sep=0}, from=1-6, to=5-8]
			\arrow["{g'}", from=1-8, to=3-9]
			\arrow[""{name=2, anchor=center, inner sep=0}, "{\varepsilon_{y'}}"', from=1-8, to=5-8]
			\arrow[from=3-2, to=3-4]
			\arrow["f"', from=3-2, to=5-1]
			\arrow[between={0.2}{0.8}, Rightarrow, nfold, from=3-4, to=5-1]
			\arrow["{f'}", from=3-4, to=5-3]
			\arrow["{f'}", from=3-9, to=5-8]
			\arrow[from=5-1, to=5-3]
			\arrow[from=5-6, to=5-8]
			\arrow[between={0.2}{0.8}, Rightarrow, nfold, from=1, to=5-6]
			\arrow[between={0.2}{0.8}, Rightarrow, nfold, from=1-8, to=1]
			\arrow["\epsilon"', between={0.2}{0.8}, Rightarrow, nfold, from=3-2, to=0]
			\arrow["{\epsilon'}"', between={0.2}{0.8}, Rightarrow, nfold, from=3-9, to=2]
		\end{tikzcd}\]
		\[\begin{tikzcd}[ampersand replacement=\&]
			\& x \&\& {x'} \&\&\& x \&\& {x'} \\
			\\
			y \&\& {y'} \&\& \equiv \& y \\
			\\
			\& x \&\& {x'} \&\&\& x \&\& {x'}
			\arrow[from=1-2, to=1-4]
			\arrow["f"', from=1-2, to=3-1]
			\arrow[between={0.2}{0.8}, Rightarrow, nfold, from=1-4, to=3-1]
			\arrow["{f'}", from=1-4, to=3-3]
			\arrow[""{name=0, anchor=center, inner sep=0}, "{\varepsilon_{x'}}", from=1-4, to=5-4]
			\arrow[from=1-7, to=1-9]
			\arrow["f"', from=1-7, to=3-6]
			\arrow[""{name=1, anchor=center, inner sep=0}, "{\varepsilon_{x}}", from=1-7, to=5-7]
			\arrow[""{name=2, anchor=center, inner sep=0}, from=1-7, to=5-9]
			\arrow["{\varepsilon_{x'}}", from=1-9, to=5-9]
			\arrow[from=3-1, to=3-3]
			\arrow["g"', from=3-1, to=5-2]
			\arrow[between={0.2}{0.8}, Rightarrow, nfold, from=3-3, to=5-2]
			\arrow["{g'}", from=3-3, to=5-4]
			\arrow["g"', from=3-6, to=5-7]
			\arrow[from=5-2, to=5-4]
			\arrow[from=5-7, to=5-9]
			\arrow["{\eta'}", between={0.2}{0.8}, Rightarrow, nfold, from=0, to=3-3]
			\arrow["\eta"', between={0.2}{0.8}, Rightarrow, nfold, from=1, to=3-6]
			\arrow[between={0.2}{0.8}, Rightarrow, nfold, from=2, to=5-7]
			\arrow[between={0.2}{0.8}, Rightarrow, nfold, from=1-9, to=2]
		\end{tikzcd}\cdpunct[8pt]{.}\]
	\end{enumerate}
\end{definition}

\begin{definition}
	A $2$-morphism consists of
	\begin{enumerate}
		\item two $2$-cells given as
		\[\begin{tikzcd}[ampersand replacement=\&]
			x \&\& {x'} \& y \&\& {y'}
			\arrow[""{name=0, anchor=center, inner sep=0}, curve={height=-12pt}, from=1-1, to=1-3]
			\arrow[""{name=1, anchor=center, inner sep=0}, curve={height=12pt}, from=1-1, to=1-3]
			\arrow[""{name=2, anchor=center, inner sep=0}, curve={height=-12pt}, from=1-4, to=1-6]
			\arrow[""{name=3, anchor=center, inner sep=0}, curve={height=12pt}, from=1-4, to=1-6]
			\arrow[between={0.2}{0.8}, Rightarrow, nfold, from=0, to=1]
			\arrow[between={0.2}{0.8}, Rightarrow, nfold, from=2, to=3]
		\end{tikzcd}\cdpunct[2pt]{,}\]
		\item such that the following two identities hold:
		\[\begin{tikzcd}[ampersand replacement=\&]
			x \&\& {x'} \&\& x \&\& {x'} \\
			\&\&\& \equiv \\
			y \&\& {y'} \&\& y \&\& {y'} \\
			\\
			y \&\& {y'} \&\& y \&\& {y'} \\
			\&\&\& \equiv \\
			x \&\& {x'} \&\& x \&\& {x'}
			\arrow[""{name=0, anchor=center, inner sep=0}, from=1-1, to=1-3]
			\arrow[""{name=1, anchor=center, inner sep=0}, curve={height=-18pt}, from=1-1, to=1-3]
			\arrow["f"', from=1-1, to=3-1]
			\arrow[between={0.2}{0.8}, Rightarrow, nfold, from=1-3, to=3-1]
			\arrow["{f'}", from=1-3, to=3-3]
			\arrow[from=1-5, to=1-7]
			\arrow["f"', from=1-5, to=3-5]
			\arrow[between={0.2}{0.8}, Rightarrow, nfold, from=1-7, to=3-5]
			\arrow["{f'}", from=1-7, to=3-7]
			\arrow[from=3-1, to=3-3]
			\arrow[""{name=2, anchor=center, inner sep=0}, from=3-5, to=3-7]
			\arrow[""{name=3, anchor=center, inner sep=0}, curve={height=18pt}, from=3-5, to=3-7]
			\arrow[""{name=4, anchor=center, inner sep=0}, from=5-1, to=5-3]
			\arrow[""{name=5, anchor=center, inner sep=0}, curve={height=-18pt}, from=5-1, to=5-3]
			\arrow["g"', from=5-1, to=7-1]
			\arrow[between={0.2}{0.8}, Rightarrow, nfold, from=5-3, to=7-1]
			\arrow["{g'}", from=5-3, to=7-3]
			\arrow[from=5-5, to=5-7]
			\arrow["g"', from=5-5, to=7-5]
			\arrow[between={0.2}{0.8}, Rightarrow, nfold, from=5-7, to=7-5]
			\arrow["{g'}", from=5-7, to=7-7]
			\arrow[from=7-1, to=7-3]
			\arrow[""{name=6, anchor=center, inner sep=0}, from=7-5, to=7-7]
			\arrow[""{name=7, anchor=center, inner sep=0}, curve={height=18pt}, from=7-5, to=7-7]
			\arrow[between={0.2}{0.8}, Rightarrow, nfold, from=1, to=0]
			\arrow[between={0.2}{0.8}, Rightarrow, nfold, from=2, to=3]
			\arrow[between={0.2}{0.8}, Rightarrow, nfold, from=5, to=4]
			\arrow[between={0.2}{0.8}, Rightarrow, nfold, from=6, to=7]
		\end{tikzcd}\cdpunct[8pt]{.}\]
	\end{enumerate}
\end{definition}

One objection the reader may have is that we never discussed the composites of adjoints in the context of our proposed \emph{bicategory of adjunctions}. To resolve this discrepancy, we shall illustrate the idea one dimension lower with a straightforward lemma.

\begin{lemma}
	There exists a double category of inverses in $\cC$ defined by the following structure:
	\begin{description}
		\item[Objects - ] are precisely the objects of $\cC$ itself,
		\item[Horizontal $1$-cells - ] are precisely the morphisms of $\cC$,
		\item[Vertical $1$-cells - ] are precisely the invertible arrows in $\cC$,
		\item[Square $2$-cells  - ] are commuting squares.
	\end{description}
\end{lemma}

From this double category, one can extract the usual $1$-category of inverses, whose objects are invertible pairs of morphisms $(f,f^{-1})$, and whose morphisms are commuting square.

This suggests that our bicategory of adjunction in $K$ should itself be a part of some higher double categorical structure, such as a Verity double bicategory \cite{morton2009doublebicat} or even an intercategory \cite{grandis2017intercatframework}. However, developing this structure is outside the scope of this paper, and we shall not pursue it further here.

Analogous to the involution case, there is also a \emph{$1$-category of relative inverse objects}. Furthermore, adjunctions admit slicing using these inverse objects.

\begin{definition}[inverse objects]
	Let $(F\colon \cC \rightarrow \cD,G\colon \cD \rightarrow \cC,\eta,\epsilon)$ be an adjunction. We define a relative inverse object to be the pair $(f \colon Fc \rightarrow d , g \colon Gd \rightarrow c)$, such that the morphisms satisfy the following identities:
	\[\begin{tikzcd}[ampersand replacement=\&]
		{G(Fc)} \&\& {G(d)} \& {F(Gd)} \&\& {F(c)} \\
		c \&\& c \& d \&\& d
		\arrow["{G(f)}", from=1-1, to=1-3]
		\arrow["g", from=1-3, to=2-3]
		\arrow["{F(g)}", from=1-4, to=1-6]
		\arrow["{\epsilon_d}"', from=1-4, to=2-4]
		\arrow["f", from=1-6, to=2-6]
		\arrow["{\eta_c}", from=2-1, to=1-1]
		\arrow["{1_c}"', from=2-1, to=2-3]
		\arrow["{1_d}"', from=2-4, to=2-6]
	\end{tikzcd}\cdpunct[8pt]{.}\]
	
\end{definition}

\begin{lemma}
	Let $F \dashv G$ be an adjunction (not necessarily equivalence) between $\cC$ and $\cD$. And let $(f \colon Fc \rightarrow d, g \colon Gd \rightarrow c)$ be the data of an inverse object relative to this adjunction. Then, there exists a canonical sliced adjunction:
	\[\begin{tikzcd}[ampersand replacement=\&]
		\& {c_1} \&\& {F(c_1)} \& \\
		\& c \&\& {F(c)} \& d \\
		\& {\cC/c} \&\& {\cD/d} \\
		\& {G(d_1)} \&\& {d_1} \\
		c \& {G(d)} \&\& d
		\arrow[""{name=0, anchor=center, inner sep=0}, from=1-2, to=2-2]
		\arrow[""{name=1, anchor=center, inner sep=0}, from=1-4, to=2-4]
		\arrow["f", from=2-4, to=2-5]
		\arrow[""{name=2, anchor=center, inner sep=0}, "{F_{/f}}", shift left=2, from=3-2, to=3-4]
		\arrow[""{name=3, anchor=center, inner sep=0}, "{G_{/g}}", shift left=2, from=3-4, to=3-2]
		\arrow[""{name=4, anchor=center, inner sep=0}, from=4-2, to=5-2]
		\arrow[""{name=5, anchor=center, inner sep=0}, from=4-4, to=5-4]
		\arrow["g", from=5-2, to=5-1]
		\arrow[between={0.2}{0.8}, maps to, from=0, to=1]
		\arrow["\dashv"{anchor=center, rotate=-90}, draw=none, from=2, to=3]
		\arrow[between={0.2}{0.8}, maps to, from=5, to=4]
	\end{tikzcd}\cdpunct[8pt]{.}\]
	Here, the action on morphisms is the one induced by $F$ and $G$.
\end{lemma}

\begin{proof}
	Let us denote the new adjoint pair using the tuple $(F_{/f}, G_{/g}, \eta_{/}, \epsilon_{/})$. Define $\eta_{/}$ and $\epsilon_{/}$ using the following two diagrams:
	\[\begin{tikzcd}[ampersand replacement=\&]
		{{\left(\eta_{/}\right)}_{\alpha} :} \& {c_1} \&\& {GF(c_1)} \& {{\left(\epsilon_{/}\right)}_{\beta} :} \& {FG(d_1)} \&\& {d_1} \\
		\&\&\& {GF(c)} \&\& {FG(d)} \\
		\& c \&\& Gd \&\& Fc \&\& d \\
		\&\& c \&\&\&\& d
		\arrow["{\eta_{c_1}}", from=1-2, to=1-4]
		\arrow["\alpha"', from=1-2, to=3-2]
		\arrow["{GF(\alpha)}", from=1-4, to=2-4]
		\arrow["{\epsilon_{d_1}}", from=1-6, to=1-8]
		\arrow["{FG(\beta)}"', from=1-6, to=2-6]
		\arrow["\beta", from=1-8, to=3-8]
		\arrow["Gf", from=2-4, to=3-4]
		\arrow["{F(g)}"', from=2-6, to=3-6]
		\arrow["{\epsilon_{d}}", from=2-6, to=3-8]
		\arrow["{\eta_c}", from=3-2, to=2-4]
		\arrow["{1_{c}}"', from=3-2, to=4-3]
		\arrow["g", from=3-4, to=4-3]
		\arrow["f"', from=3-6, to=4-7]
		\arrow["{1_{d}}", from=3-8, to=4-7]
	\end{tikzcd}\cdpunct[8pt]{.}\]
	The relevant triangle identities are then immediate by pasting these two diagrams appropriately. For instance, one of the triangle identity is given by the following commuting diagram:
	\[\begin{tikzcd}[ampersand replacement=\&]
		{F(c_1)} \&\& {FGF(c_1)} \&\& {F(c_1)} \\
		\&\& {FGF(c)} \\
		{F(c)} \&\& FGd \&\& {F(c)} \\
		\& {F(c)} \&\&\& d \\
		\&\& d
		\arrow["{F\eta_{c_1}}", from=1-1, to=1-3]
		\arrow["{F(\alpha)}"', from=1-1, to=3-1]
		\arrow["{\epsilon_{Fc_1}}", from=1-3, to=1-5]
		\arrow["{FGF(\alpha)}", from=1-3, to=2-3]
		\arrow["{F(\alpha)}", from=1-5, to=3-5]
		\arrow["{FG(f)}", from=2-3, to=3-3]
		\arrow["{\epsilon_{Fc}}", from=2-3, to=3-5]
		\arrow["{F\eta_c}", from=3-1, to=2-3]
		\arrow["{1_{F(c)}}"', from=3-1, to=4-2]
		\arrow["{F(g)}", from=3-3, to=4-2]
		\arrow["{\epsilon_{d}}", from=3-3, to=4-5]
		\arrow["f", from=3-5, to=4-5]
		\arrow["f"', from=4-2, to=5-3]
		\arrow["{1_{d}}", from=4-5, to=5-3]
	\end{tikzcd}\cdpunct[8pt]{.}\]
\end{proof}

At this stage, we have provided ample motivation and theory to view dualities and inverses through the lens of the microcosm principle. However, the original statement holds between two arbitrary dimensions $n$ and $n+1$. Therefore, a similar theory of dualities and inverses should also hold in one dimension higher. In particular, we want to define relative pseudoinvolutions and pseudoinverses.

\subsection{\tpdf{$2$-dimensional variants}} \label{sub:2diminv}

To obtain the definition of $2$-dimensional variant of involutive objects, it is best to utilize a categorification of the lax compatible form of $1$-dimensional involutive objects. A somewhat similar definition can be found in \cite[Def.~2.1]{shulman2018contenrich} for the special case of $T = \BiCat$.

\begin{definition}
	A weak $2$-involution in a tricategory $T$ consists of
	\begin{enumerate}
		\item a $1$-cell $H \colon K \rightarrow K$,
		\item an adjoint equivalence $\alpha \colon H \otimes H \iff 1_K \colon \alpha^{.} $,
		\item a $3$-cell $\zeta$ as depicted by the following diagram:
		\[\begin{tikzcd}[ampersand replacement=\&,cramped]
			K \&\& K \&\& K \&\& K \\
			\&\&\& {\Rrightarrow^{\zeta}} \\
			K \&\& K \&\& K \&\& K \\
			\& {(H \otimes H) \otimes H} \&\&\&\& {H \otimes (H \otimes H) } \\
			\\
			\& {1_K \otimes H} \&\&\&\& {H \otimes 1_K} \\
			\&\&\& H
			\arrow["H", from=1-1, to=1-3]
			\arrow[""{name=0, anchor=center, inner sep=0}, "{1_K}"{description}, from=1-1, to=3-3]
			\arrow["H", from=1-3, to=3-3]
			\arrow["H", from=1-5, to=1-7]
			\arrow["H", from=1-7, to=3-7]
			\arrow["H", from=3-1, to=1-1]
			\arrow[""{name=1, anchor=center, inner sep=0}, "H"', from=3-1, to=3-3]
			\arrow["H", from=3-5, to=1-5]
			\arrow[""{name=2, anchor=center, inner sep=0}, "{1_K}"{description}, from=3-5, to=1-7]
			\arrow[""{name=3, anchor=center, inner sep=0}, "H"{description}, from=3-5, to=3-7]
			\arrow["a", from=4-2, to=4-6]
			\arrow["{\alpha \otimes 1}"', from=4-2, to=6-2]
			\arrow["{1 \otimes \alpha}", from=4-6, to=6-6]
			\arrow["\zeta", between={0.2}{0.8}, Rightarrow, nfold, from=6-2, to=4-6]
			\arrow["l"', from=6-2, to=7-4]
			\arrow["r", from=6-6, to=7-4]
			\arrow["l"', between={0.2}{0.8}, Rightarrow, nfold, from=1-1, to=1]
			\arrow["\alpha", between={0.2}{0.8}, Rightarrow, nfold, from=1-3, to=0]
			\arrow["\alpha", between={0.2}{0.8}, Rightarrow, nfold, from=1-5, to=2]
			\arrow["r", between={0.2}{0.7}, Rightarrow, nfold, from=1-7, to=3]
		\end{tikzcd}\cdpunct[8pt]{,}\]
		\item the $3$-cell $\zeta$ satisfy the weak $2$-involution axiom (\cref{app:weak2invdef}).
		
	\end{enumerate}
\end{definition}

Although the morphisms of weak $2$-involutions are not directly relevant for our purposes, we record their definitions in \cref{appsub:weak2inv}. We now define the relative pseudo involutive objects as suggested by the microcosm principle.

\begin{definition}
	Let $(K,H,\alpha,\zeta)$ be a weak $2$-involution in the tricategory $\BiCat$. We define a relative lax involutive object by
	\begin{enumerate}
		\item a $1$-cell $f \colon Hx \rightarrow x \in B$,
		\item a $2$-cell depicted by the following diagram:
		\[\begin{tikzcd}[ampersand replacement=\&,cramped]
			{H^2x} \&\& Hx \\
			\\
			x \&\& x
			\arrow["Hf", from=1-1, to=1-3]
			\arrow["{\alpha_x}"', from=1-1, to=3-1]
			\arrow["\phi", between={0.2}{0.8}, Rightarrow, nfold, from=1-3, to=3-1]
			\arrow["f", from=1-3, to=3-3]
			\arrow["{\varepsilon_x}"', from=3-1, to=3-3]
		\end{tikzcd}\cdpunct[8pt]{,}\]
		\item such that the $2$-cell satisfy the following identity:
		\[\begin{tikzcd}[ampersand replacement=\&,column sep=scriptsize]
			{H^3x} \&\& {H^2x} \&\& Hx \&\& {H^3x} \&\& {H^2x} \&\& \\
			\&\&\&\&\& \equiv \\
			Hx \&\& x \&\& x \&\& Hx \&\& Hx \&\& x
			\arrow["{H^2f}", from=1-1, to=1-3]
			\arrow[""{name=0, anchor=center, inner sep=0}, "{\alpha_{Hx}}", from=1-1, to=3-1]
			\arrow[""{name=1, anchor=center, inner sep=0}, "{H\alpha_x}"', curve={height=30pt}, from=1-1, to=3-1]
			\arrow["Hf", from=1-3, to=1-5]
			\arrow[between={0.2}{0.8}, Rightarrow, nfold, from=1-3, to=3-1]
			\arrow["{\alpha_x}"', from=1-3, to=3-3]
			\arrow["\phi"', between={0.2}{0.8}, Rightarrow, nfold, from=1-5, to=3-3]
			\arrow["f", from=1-5, to=3-5]
			\arrow["{H^2f}", from=1-7, to=1-9]
			\arrow["{H\alpha_x}"', from=1-7, to=3-7]
			\arrow["{H\phi}", between={0.2}{0.8}, Rightarrow, nfold, from=1-9, to=3-7]
			\arrow["Hf", from=1-9, to=3-9]
			\arrow["f"', from=3-1, to=3-3]
			\arrow[""{name=2, anchor=center, inner sep=0}, "f"', curve={height=30pt}, from=3-1, to=3-5]
			\arrow["{\varepsilon_x}"', from=3-3, to=3-5]
			\arrow["{\varepsilon_{Hx}}"', from=3-7, to=3-9]
			\arrow[""{name=3, anchor=center, inner sep=0}, "f"', curve={height=30pt}, from=3-7, to=3-11]
			\arrow["f"', from=3-9, to=3-11]
			\arrow["{\zeta_x}"', between={0.2}{0.8}, Rightarrow, nfold, from=0, to=1]
			\arrow[between={0.2}{0.8}, Rightarrow, nfold, from=3-3, to=2]
			\arrow[between={0.2}{0.8}, Rightarrow, nfold, from=3-9, to=3]
		\end{tikzcd}\cdpunct[8pt]{.}\]
	\end{enumerate}
\end{definition}

\begin{example}
	The functor $\op_2 \colon \nCat{2}_{\leq 2} \rightarrow \nCat{2}_{\leq 2}$ defines a strict $2$-involutive object. Thus, $\op_1 \colon \op_2(\nCat{1}) \rightarrow \nCat{1}$ is a relative strict $1$-involutive object.
\end{example}

One can similarly define the $1$-morphisms and the $2$-morphisms of these relative involutions.

\begin{definition}		\phantom{}
	\begin{enumerate} 
		\item A $1$-morphism $(r,u)$ is given by $r \colon x \rightarrow y$, and a corresponding $2$-cell $u$ that satisfy the following identity:
		\[\begin{tikzcd}[ampersand replacement=\&,cramped,column sep=small]
			\& {H^2x} \&\& {H^2y} \&\&\&\& {H^2x} \&\& {H^2y} \& \\
			\\
			x \&\& Hx \&\& Hy \& \equiv \& x \&\& y \&\& Hy \\
			\\
			\& x \&\& y \&\&\&\& x \&\& y
			\arrow["{H^2r}", from=1-2, to=1-4]
			\arrow["{\alpha_x}"', from=1-2, to=3-1]
			\arrow["Hf", from=1-2, to=3-3]
			\arrow["Hu", between={0.2}{0.8}, Rightarrow, nfold, from=1-4, to=3-3]
			\arrow["Hg", from=1-4, to=3-5]
			\arrow["{H^2r}", from=1-8, to=1-10]
			\arrow["{\alpha_x}"', from=1-8, to=3-7]
			\arrow[between={0.2}{0.8}, Rightarrow, nfold, from=1-10, to=3-7]
			\arrow["{\alpha_y}"', from=1-10, to=3-9]
			\arrow["Hg", from=1-10, to=3-11]
			\arrow["{\varepsilon_x}"', from=3-1, to=5-2]
			\arrow["{\phi_x}"', between={0.2}{0.8}, Rightarrow, nfold, from=3-3, to=3-1]
			\arrow["Hr"', from=3-3, to=3-5]
			\arrow["f", from=3-3, to=5-2]
			\arrow["u", between={0.2}{0.8}, Rightarrow, nfold, from=3-5, to=5-2]
			\arrow["g", from=3-5, to=5-4]
			\arrow["r", from=3-7, to=3-9]
			\arrow["{\varepsilon_x}"', from=3-7, to=5-8]
			\arrow[between={0.2}{0.8}, Rightarrow, nfold, from=3-9, to=5-8]
			\arrow["{\varepsilon_y}"', from=3-9, to=5-10]
			\arrow["{\phi_y}"', between={0.2}{0.8}, Rightarrow, nfold, from=3-11, to=3-9]
			\arrow["g", from=3-11, to=5-10]
			\arrow["r"', from=5-2, to=5-4]
			\arrow["r"', from=5-8, to=5-10]
		\end{tikzcd}\cdpunct[8pt]{.}\]
		\item A $2$-morphism is a $2$-cell $\gamma \colon r \Rightarrow t$ that satisfy the following identity:
		\[\begin{tikzcd}[ampersand replacement=\&]
			Hx \&\& Hy \&\& Hx \&\& Hy \\
			\&\&\& \equiv \\
			x \&\& y \&\& x \&\& y
			\arrow["Hr", from=1-1, to=1-3]
			\arrow["f"', from=1-1, to=3-1]
			\arrow["u", between={0.2}{0.8}, Rightarrow, nfold, from=1-3, to=3-1]
			\arrow["g", from=1-3, to=3-3]
			\arrow[""{name=0, anchor=center, inner sep=0}, "Ht"', from=1-5, to=1-7]
			\arrow[""{name=1, anchor=center, inner sep=0}, "Hr", curve={height=-30pt}, from=1-5, to=1-7]
			\arrow["f"', from=1-5, to=3-5]
			\arrow["v", between={0.2}{0.8}, Rightarrow, nfold, from=1-7, to=3-5]
			\arrow["g", from=1-7, to=3-7]
			\arrow[""{name=2, anchor=center, inner sep=0}, "r", from=3-1, to=3-3]
			\arrow[""{name=3, anchor=center, inner sep=0}, "t"', curve={height=30pt}, from=3-1, to=3-3]
			\arrow["t", from=3-5, to=3-7]
			\arrow["{H\gamma}", between={0.2}{0.8}, Rightarrow, nfold, from=1, to=0]
			\arrow["\gamma", between={0.2}{0.8}, Rightarrow, nfold, from=2, to=3]
		\end{tikzcd}\cdpunct[8pt]{.}\]
	\end{enumerate}
\end{definition} 

Therefore, we obtain the bicategory $\nInv{1}(K,H,\alpha,\zeta)$. The axioms can be verified similar as before by using the strictification to encode the desired compositional data.

\begin{proposition}
	The microcosm principle holds for weak $2$-involutive bicategories.
\end{proposition}

Finally, in parallel, we can also write these higher dimensional definitions in the context of inverses. 

\begin{definition}[{\cite[Def.~2.1]{gurski2012biequivintricat}}]
	A weak $2$-adjunction in a tricategory $T$ is
	\begin{enumerate}
		\item two $1$-cell $F \colon K_1 \rightarrow K_2$ and $G \colon K_2 \rightarrow K_1$,
		
		\item two $2$-cells $\epsilon \colon FG \Rightarrow 1_{K_2}$ and $\eta \colon 1_{K_1} \Rightarrow GF$,
		
		\item two $3$-cells $\phi$ and $\psi$ as depicted by the following diagrams:
		\[\begin{tikzcd}[ampersand replacement=\&,column sep=small]
			\&\&\& F \&\&\&\&\& \\
			\\
			{K_1} \&\& {K_1} \&\&\& {K_1} \\
			\&\&\&\& {\Rrightarrow^{\phi}} \\
			\& {K_2} \&\& {K_2} \&\&\&\&\& {K_2} \\
			\\
			F \\
			\&\& {F \otimes(G \otimes F)} \&\&\&\& {(F \otimes G) \otimes F} \\
			\&\& {F \otimes 1_{K_1} } \&\&\&\& {1_{K_2} \otimes F} \\
			\&\& F \&\&\&\& F
			\arrow[between={0.2}{0.8}, Rightarrow, nfold, from=1-4, to=3-3]
			\arrow[from=1-4, to=5-4]
			\arrow[no head, from=3-1, to=1-4]
			\arrow[""{name=0, anchor=center, inner sep=0}, "{1_{K_1}}", from=3-1, to=3-3]
			\arrow["F"', from=3-1, to=5-2]
			\arrow[no head, from=3-1, to=7-1]
			\arrow["F", from=3-3, to=5-4]
			\arrow[""{name=1, anchor=center, inner sep=0}, "F"{description}, curve={height=-30pt}, from=3-6, to=5-9]
			\arrow[""{name=2, anchor=center, inner sep=0}, "F"{description}, curve={height=30pt}, from=3-6, to=5-9]
			\arrow["G"{description}, from=5-2, to=3-3]
			\arrow[""{name=3, anchor=center, inner sep=0}, "{1_{K_2}}"', from=5-2, to=5-4]
			\arrow[between={0.2}{0.8}, Rightarrow, nfold, from=5-2, to=7-1]
			\arrow[from=7-1, to=5-4]
			\arrow[""{name=4, anchor=center, inner sep=0}, "{a^{.}}", from=8-3, to=8-7]
			\arrow["{\epsilon \otimes 1}", from=8-7, to=9-7]
			\arrow["{1 \otimes \eta}", from=9-3, to=8-3]
			\arrow["l", from=9-7, to=10-7]
			\arrow["{r^{.}}", from=10-3, to=9-3]
			\arrow[""{name=5, anchor=center, inner sep=0}, "{1_F}"', from=10-3, to=10-7]
			\arrow["\eta"{description}, between={0.2}{0.8}, Rightarrow, nfold, from=0, to=5-2]
			\arrow["\epsilon"{description}, between={0.2}{0.8}, Rightarrow, nfold, from=3-3, to=3]
			\arrow["{1_F}"{description}, between={0.2}{0.8}, Rightarrow, nfold, from=1, to=2]
			\arrow["\phi", between={0.2}{0.8}, Rightarrow, nfold, from=4, to=5]
		\end{tikzcd}\]
		\[\begin{tikzcd}[ampersand replacement=\&,column sep=small]
			G \&\&\&\&\&\&\&\& \\
			\\
			\& {K_1} \&\& {K_1} \&\&\&\&\& {K_1} \\
			\&\&\&\& {\Lleftarrow^{\psi}} \\
			{K_2} \&\& {K_2} \&\&\& {K_2} \\
			\\
			\&\&\& G \\
			\&\& {(G \otimes F ) \otimes G} \&\&\&\& {G \otimes (F \otimes G)} \\
			\&\& {1_{K_1} \otimes G} \&\&\&\& {G \otimes 1_{K_2}} \\
			\&\& G \&\&\&\& G
			\arrow[between={0.2}{0.8}, Rightarrow, nfold, from=1-1, to=3-2]
			\arrow[from=1-1, to=3-4]
			\arrow[""{name=0, anchor=center, inner sep=0}, "{1_{K_1}}", from=3-2, to=3-4]
			\arrow["F"{description}, from=3-2, to=5-3]
			\arrow[no head, from=5-1, to=1-1]
			\arrow["G", from=5-1, to=3-2]
			\arrow[""{name=1, anchor=center, inner sep=0}, "{1_{K_2}}"', from=5-1, to=5-3]
			\arrow[no head, from=5-1, to=7-4]
			\arrow["G"', from=5-3, to=3-4]
			\arrow[between={0.2}{0.8}, Rightarrow, nfold, from=5-3, to=7-4]
			\arrow[""{name=2, anchor=center, inner sep=0}, "G"{description}, curve={height=-30pt}, from=5-6, to=3-9]
			\arrow[""{name=3, anchor=center, inner sep=0}, "G"{description}, curve={height=30pt}, from=5-6, to=3-9]
			\arrow[from=7-4, to=3-4]
			\arrow[""{name=4, anchor=center, inner sep=0}, "a", from=8-3, to=8-7]
			\arrow["{1 \otimes \epsilon}", from=8-7, to=9-7]
			\arrow["{\eta \otimes 1}", from=9-3, to=8-3]
			\arrow["r", from=9-7, to=10-7]
			\arrow["{l^{.}}", from=10-3, to=9-3]
			\arrow[""{name=5, anchor=center, inner sep=0}, "{1_G}"', from=10-3, to=10-7]
			\arrow["\eta"{description}, between={0.2}{0.8}, Rightarrow, nfold, from=0, to=5-3]
			\arrow["\epsilon"{description}, between={0.2}{0.8}, Rightarrow, nfold, from=3-2, to=1]
			\arrow["{1_G}"{description}, between={0.2}{0.8}, Rightarrow, nfold, from=2, to=3]
			\arrow["\psi"', between={0.2}{0.8}, Rightarrow, nfold, from=5, to=4]
		\end{tikzcd}\cdpunct[8pt]{,}\]
		\item such that two pasting diagrams as given in \cite[Def.~2.1]{gurski2012biequivintricat} or in \cref{app:weak2invs} hold. 
	\end{enumerate}
\end{definition}

As before, the morphisms of weak $2$-adjunctions can be found in \cref{appsub:weak2invs}.

\begin{definition}
	Let $(F \colon K_1 \rightarrow K_2, G: K_2 \rightarrow K_1, \epsilon,\eta,\phi,\psi)$ be the data of a weak 2-adjunction in the tricategory $\BiCat$. We define a relative (lax) weak $1$-inverse object by
	\begin{enumerate}
		\item two $1$-cells $f \colon Fk_1 \rightarrow k_2 \in K_2$ and $g \colon Gk_2 \rightarrow k_1 \in K_1$,
		
		\item two $2$-cells
		\[\begin{tikzcd}[ampersand replacement=\&]
			{FGk_2} \&\& {Fk_1} \& {GFk_1} \&\& {Gk_2} \\
			\\
			{k_2} \&\& {k_2} \& {k_1} \&\& {k_1}
			\arrow["{F(g)}", from=1-1, to=1-3]
			\arrow["{\epsilon_{k_2}}"', from=1-1, to=3-1]
			\arrow["\alpha", between={0.2}{0.8}, Rightarrow, nfold, from=1-3, to=3-1]
			\arrow["f", from=1-3, to=3-3]
			\arrow[""{name=0, anchor=center, inner sep=0}, "{G(f)}", from=1-4, to=1-6]
			\arrow["g", from=1-6, to=3-6]
			\arrow["{\varepsilon_{k_2}}"', from=3-1, to=3-3]
			\arrow["{\eta_{k_1}}", from=3-4, to=1-4]
			\arrow[""{name=1, anchor=center, inner sep=0}, "{\varepsilon_{k_1}}"', from=3-4, to=3-6]
			\arrow["\beta"', between={0.2}{0.8}, Rightarrow, nfold, from=1, to=0]
		\end{tikzcd}\cdpunct[8pt]{,}\]
		
		\item such that the $2$-cells satisfy the following identities:
		\[\begin{tikzcd}[ampersand replacement=\&,cramped]
			{Fk_1} \&\&\&\&\&\& {k_2} \\
			\\
			\&\& {FGFk_1} \&\& {FGk_2} \\
			\\
			{Fk_1} \&\&\&\& {F k_1} \&\& {k_2} \\
			\&\&\& \vequiv \\
			{Fk_1} \&\&\&\&\&\& {k_2} \\
			\\
			{Fk_1} \&\&\&\& {Fk_1} \&\& {k_2}
			\arrow[""{name=0, anchor=center, inner sep=0}, "f", from=1-1, to=1-7]
			\arrow["{\varepsilon_{k_2}}"', from=1-7, to=5-7]
			\arrow["{\epsilon_{Fk_1}}"{description}, from=3-3, to=1-1]
			\arrow[""{name=1, anchor=center, inner sep=0}, "{FG(f)}"{description}, from=3-3, to=3-5]
			\arrow["{\epsilon_{k_2}}"{description}, from=3-5, to=1-7]
			\arrow["Fg", from=3-5, to=5-5]
			\arrow[""{name=2, anchor=center, inner sep=0}, "{\varepsilon_{Fk_1}}", from=5-1, to=1-1]
			\arrow["{F\eta_{k_1}}"{description}, from=5-1, to=3-3]
			\arrow[""{name=3, anchor=center, inner sep=0}, "{\varepsilon_{Fk_1}}"', from=5-1, to=5-5]
			\arrow["\alpha"{description}, between={0.2}{0.8}, Rightarrow, nfold, from=5-5, to=1-7]
			\arrow["f"', from=5-5, to=5-7]
			\arrow["f", from=7-1, to=7-7]
			\arrow["{\varepsilon_{Fk_1}}"{description}, from=7-1, to=9-5]
			\arrow["f"{description}, from=7-1, to=9-7]
			\arrow["{\varepsilon_{k_2}}", from=7-7, to=9-7]
			\arrow["{\varepsilon_{Fk_1}}", from=9-1, to=7-1]
			\arrow["{\varepsilon_{Fk_1}}"', from=9-1, to=9-5]
			\arrow["f"', from=9-5, to=9-7]
			\arrow["\phi"{description}, between={0.2}{0.8}, Rightarrow, nfold, from=3-3, to=2]
			\arrow[between={0.2}{0.8}, Rightarrow, nfold, from=1, to=0]
			\arrow["{F\beta}"{description}, between={0.2}{0.8}, Rightarrow, nfold, from=3, to=1]
		\end{tikzcd}\]
		\[\begin{tikzcd}[ampersand replacement=\&,cramped]
			\& {Gk_2} \&\&\&\&\&\& {k_1} \\
			\\
			\& {GFGk_2} \&\&\&\& {GFk_1} \\
			\\
			\& {Gk_2} \&\&\&\& {Gk_2} \&\& {k_1} \\
			\&\&\&\& \vequiv \\
			\& {Gk_2} \&\&\&\&\&\& {k_1} \\
			{GFGk_2} \\
			\& {Gk_2} \&\&\& {Gk_2} \&\&\& {k_1}
			\arrow[""{name=0, anchor=center, inner sep=0}, "g", from=1-2, to=1-8]
			\arrow["{\eta_{Gk_2}}"', from=1-2, to=3-2]
			\arrow["{\eta_{k_1}}"{description}, from=1-8, to=3-6]
			\arrow[""{name=1, anchor=center, inner sep=0}, "{\varepsilon_{k_1}}", from=1-8, to=5-8]
			\arrow[""{name=2, anchor=center, inner sep=0}, "{GF(g)}"{description}, from=3-2, to=3-6]
			\arrow["{G\epsilon_{k_2}}"', from=3-2, to=5-2]
			\arrow["{G\alpha}"{description}, between={0.2}{0.8}, Rightarrow, nfold, from=3-6, to=5-2]
			\arrow[""{name=3, anchor=center, inner sep=0}, "Gf"{description}, from=3-6, to=5-6]
			\arrow["{\varepsilon_{Gk_2}}"', from=5-2, to=5-6]
			\arrow["g"', from=5-6, to=5-8]
			\arrow["g", from=7-2, to=7-8]
			\arrow["{\eta_{Gk_2}}"', from=7-2, to=8-1]
			\arrow[""{name=4, anchor=center, inner sep=0}, "{\varepsilon_{Gk_2}}"{description}, from=7-2, to=9-2]
			\arrow["{\varepsilon_{Gk_2}}"{description}, from=7-2, to=9-5]
			\arrow["g"{description}, from=7-2, to=9-8]
			\arrow["{\varepsilon_{k_1}}"{description}, from=7-8, to=9-8]
			\arrow["{G\epsilon_{k_2}}"', from=8-1, to=9-2]
			\arrow["{\varepsilon_{Gk_2}}"', from=9-2, to=9-5]
			\arrow["g"', from=9-5, to=9-8]
			\arrow[between={0.2}{0.8}, Rightarrow, nfold, from=0, to=2]
			\arrow["\beta"{description}, between={0.2}{0.8}, Rightarrow, nfold, from=1, to=3]
			\arrow["\psi"', between={0.2}{0.8}, Rightarrow, nfold, from=4, to=8-1]
		\end{tikzcd}\cdpunct[8pt]{.}\]
	\end{enumerate}
\end{definition}

\begin{definition}
	Let $(f\colon Fk_1 \rightarrow k_2, g \colon Gk_2 \rightarrow k_1, \alpha,\beta)$ and $(f'\colon Fk'_1 \rightarrow k'_2, g' \colon Gk'_2 \rightarrow k'_1, \alpha',\beta')$ be two relative inverse objects. The data of a (lax) $1$-morphism is given by
	\begin{enumerate}
		\item two $1$-cells $k_1 \rightarrow k'_1$ and $k_2 \rightarrow k'_2$,
		
		\item two $2$-cells
		\[\begin{tikzcd}[ampersand replacement=\&]
			{Fk_1} \&\& {Fk'_1} \&\& {Gk_2} \&\& {Gk'_2} \\
			\\
			{k_2} \&\& {k'_2} \&\& {k_1} \&\& {k'_1}
			\arrow[from=1-1, to=1-3]
			\arrow["f"', from=1-1, to=3-1]
			\arrow[between={0.2}{0.8}, Rightarrow, nfold, from=1-3, to=3-1]
			\arrow["{f'}", from=1-3, to=3-3]
			\arrow[from=1-5, to=1-7]
			\arrow["g"', from=1-5, to=3-5]
			\arrow[between={0.2}{0.8}, Rightarrow, nfold, from=1-7, to=3-5]
			\arrow["g", from=1-7, to=3-7]
			\arrow[from=3-1, to=3-3]
			\arrow[from=3-5, to=3-7]
		\end{tikzcd}\cdpunct[8pt]{,}\]
		
		\item such that the following two identities hold:
		\[\begin{tikzcd}[ampersand replacement=\&,column sep=scriptsize]
			\& {k_1} \&\& {k'_1} \&\&\& {k_1} \&\& {k_1'} \\
			\\
			{GFk_1} \&\& {GFk'_1} \&\&\& {GFk_1} \\
			\&\&\&\& \equiv \\
			{Gk_2} \&\& {Gk'_2} \&\&\& {Gk_2} \\
			\\
			\& {k_1} \&\& {k'_1} \&\&\& {k_1} \&\& {k'_1}
			\arrow[from=1-2, to=1-4]
			\arrow["{\eta_{k_1}}"', from=1-2, to=3-1]
			\arrow[between={0.2}{0.8}, Rightarrow, nfold, from=1-4, to=3-1]
			\arrow["{\eta_{k'_1}}"{description}, from=1-4, to=3-3]
			\arrow[""{name=0, anchor=center, inner sep=0}, "{\varepsilon_{k'_1}}", from=1-4, to=7-4]
			\arrow[from=1-7, to=1-9]
			\arrow["{\eta_{k_1}}"', from=1-7, to=3-6]
			\arrow[""{name=1, anchor=center, inner sep=0}, "{\varepsilon_{k_1}}", from=1-7, to=7-7]
			\arrow[from=1-7, to=7-9]
			\arrow["{\varepsilon_{k'_1}}", from=1-9, to=7-9]
			\arrow[from=3-1, to=3-3]
			\arrow["Gf"', from=3-1, to=5-1]
			\arrow[between={0.2}{0.8}, Rightarrow, nfold, from=3-3, to=5-1]
			\arrow[""{name=2, anchor=center, inner sep=0}, "{Gf'}"{description}, from=3-3, to=5-3]
			\arrow[""{name=3, anchor=center, inner sep=0}, "Gf"', from=3-6, to=5-6]
			\arrow[from=5-1, to=5-3]
			\arrow["g"', from=5-1, to=7-2]
			\arrow[between={0.2}{0.8}, Rightarrow, nfold, from=5-3, to=7-2]
			\arrow["{g'}"{description}, from=5-3, to=7-4]
			\arrow["g"', from=5-6, to=7-7]
			\arrow[from=7-2, to=7-4]
			\arrow[from=7-7, to=7-9]
			\arrow["{\beta'}"{description}, between={0.2}{0.8}, Rightarrow, nfold, from=0, to=2]
			\arrow["\beta"{description}, between={0.2}{0.8}, Rightarrow, nfold, from=1, to=3]
		\end{tikzcd}\]
		
		\[\begin{tikzcd}[ampersand replacement=\&,column sep=tiny]
			\& {FGk_2} \&\& {FGk'_2} \&\&\&\& {FGk_2} \&\& {FGk'_2} \& \\
			\\
			{k_2} \&\& {Fk_1} \&\& {Fk'_1} \& \equiv \& {k_2} \&\& {k'_2} \&\& {Fk'_1} \\
			\\
			\& {k_2} \&\& {k'_2} \&\&\&\& {k_2} \&\& {k'_2}
			\arrow[from=1-2, to=1-4]
			\arrow["{\epsilon_{k_2}}"', from=1-2, to=3-1]
			\arrow["Fg"{description}, from=1-2, to=3-3]
			\arrow[between={0.2}{0.8}, Rightarrow, nfold, from=1-4, to=3-3]
			\arrow["{Fg'}", from=1-4, to=3-5]
			\arrow[from=1-8, to=1-10]
			\arrow["{\epsilon_{k_1}}"', from=1-8, to=3-7]
			\arrow[between={0.2}{0.8}, Rightarrow, nfold, from=1-10, to=3-7]
			\arrow["{\epsilon_{k_1}}"{description}, from=1-10, to=3-9]
			\arrow["{Fg'}", from=1-10, to=3-11]
			\arrow["{\varepsilon_{k_2}}"', from=3-1, to=5-2]
			\arrow["\alpha"{description}, between={0.2}{0.8}, Rightarrow, nfold, from=3-3, to=3-1]
			\arrow[from=3-3, to=3-5]
			\arrow["f"{description}, from=3-3, to=5-2]
			\arrow[between={0.2}{0.8}, Rightarrow, nfold, from=3-5, to=5-2]
			\arrow["{f'}", from=3-5, to=5-4]
			\arrow[from=3-7, to=3-9]
			\arrow["{\varepsilon_{k_2}}"', from=3-7, to=5-8]
			\arrow[between={0.2}{0.8}, Rightarrow, nfold, from=3-9, to=5-8]
			\arrow["{\varepsilon_{k'_2}}"{description}, from=3-9, to=5-10]
			\arrow["{\alpha'}"{description}, between={0.2}{0.8}, Rightarrow, nfold, from=3-11, to=3-9]
			\arrow["{f'}", from=3-11, to=5-10]
			\arrow[from=5-2, to=5-4]
			\arrow[from=5-8, to=5-10]
		\end{tikzcd}\cdpunct[8pt]{.}\]
		
	\end{enumerate}
\end{definition}

\begin{definition}
	The data of a $2$-morphism between two $1$-morphisms is given
	\begin{enumerate}
		\item two $2$-cells
		\[\begin{tikzcd}[ampersand replacement=\&]
			{k_1} \&\& {k'_1} \&\& {k_2} \&\& {k'_2}
			\arrow[""{name=0, anchor=center, inner sep=0}, curve={height=-12pt}, from=1-1, to=1-3]
			\arrow[""{name=1, anchor=center, inner sep=0}, curve={height=12pt}, from=1-1, to=1-3]
			\arrow[""{name=2, anchor=center, inner sep=0}, curve={height=-12pt}, from=1-5, to=1-7]
			\arrow[""{name=3, anchor=center, inner sep=0}, curve={height=12pt}, from=1-5, to=1-7]
			\arrow[between={0.2}{0.8}, Rightarrow, nfold, from=0, to=1]
			\arrow[between={0.2}{0.8}, Rightarrow, nfold, from=2, to=3]
		\end{tikzcd}\cdpunct[2pt]{,}\]
		\item such that the following two identities hold:
		\[\begin{tikzcd}[ampersand replacement=\&]
			{Fk_1} \&\& {Fk'_1} \&\& {Fk_1} \&\& {Fk'_1} \\
			\&\&\& \equiv \\
			{k_2} \&\& {k'_2} \&\& {k_2} \&\& {k'_2} \\
			\\
			{Gk_2} \&\& {Gk'_2} \&\& {Gk_2} \&\& {Gk'_2} \\
			\&\&\& \equiv \\
			{k_1} \&\& {k'_1} \&\& {k_1} \&\& {k'_1}
			\arrow[""{name=0, anchor=center, inner sep=0}, curve={height=-24pt}, from=1-1, to=1-3]
			\arrow[""{name=1, anchor=center, inner sep=0}, from=1-1, to=1-3]
			\arrow["f"', from=1-1, to=3-1]
			\arrow[between={0.2}{0.8}, Rightarrow, nfold, from=1-3, to=3-1]
			\arrow["{f'}", from=1-3, to=3-3]
			\arrow[from=1-5, to=1-7]
			\arrow["f"', from=1-5, to=3-5]
			\arrow[between={0.2}{0.8}, Rightarrow, nfold, from=1-7, to=3-5]
			\arrow["{f'}", from=1-7, to=3-7]
			\arrow[from=3-1, to=3-3]
			\arrow[""{name=2, anchor=center, inner sep=0}, from=3-5, to=3-7]
			\arrow[""{name=3, anchor=center, inner sep=0}, curve={height=24pt}, from=3-5, to=3-7]
			\arrow[""{name=4, anchor=center, inner sep=0}, from=5-1, to=5-3]
			\arrow[""{name=5, anchor=center, inner sep=0}, curve={height=-24pt}, from=5-1, to=5-3]
			\arrow["g"', from=5-1, to=7-1]
			\arrow[between={0.2}{0.8}, Rightarrow, nfold, from=5-3, to=7-1]
			\arrow["{g'}", from=5-3, to=7-3]
			\arrow[from=5-5, to=5-7]
			\arrow["g"', from=5-5, to=7-5]
			\arrow[between={0.2}{0.8}, Rightarrow, nfold, from=5-7, to=7-5]
			\arrow["{g'}", from=5-7, to=7-7]
			\arrow[from=7-1, to=7-3]
			\arrow[""{name=6, anchor=center, inner sep=0}, from=7-5, to=7-7]
			\arrow[""{name=7, anchor=center, inner sep=0}, curve={height=24pt}, from=7-5, to=7-7]
			\arrow[between={0.2}{0.8}, Rightarrow, nfold, from=0, to=1]
			\arrow[between={0.2}{0.8}, Rightarrow, nfold, from=2, to=3]
			\arrow[between={0.2}{0.8}, Rightarrow, nfold, from=5, to=4]
			\arrow[between={0.2}{0.8}, Rightarrow, nfold, from=6, to=7]
		\end{tikzcd}\cdpunct[8pt]{.}\]
	\end{enumerate}
\end{definition}

\begin{proposition}
	The microcosm principle holds for weak $2$-adjunction of bicategories.
\end{proposition}

We have now developed all the desired definitions to fully express the dualities of monoidal structures, albiet at a greater level of generality than necessary. For all our purposes, we shall only require $T = \Gray$. 

As a closing remark, we note that so far we have only employed the microcosm principle in doubly nested structures. This is in accordance with how the original statement is written and demonstrated. However, there are instances where we can achieve triply/higher nested structures.

\begin{example}
	The dual vector space $1$-functor $(-)^{*} \colon \op_{1}\left(\mathbf{Vect}_{k}^{f.d.}\right) \rightarrow \mathbf{Vect}_{k}^{f.d.}$ defines an involutive object relative to the strict involutive object $\op_1 \colon \op_2(\nCat{1}) \rightarrow \nCat{1}$. The structure cells are simply $\eta_V \colon V \rightarrow V^{**} \in \mathbf{Vect}_{k}^{f.d}$ as the unit and its opposite cell $\op_1(\eta_V)$ as the counit. One could view $(-)^{*}$ as a contravariant autoequivalence on $\mathbf{Vect}_k^{f.d.}$.
\end{example}

In this example, we have a top ($2$-truncated) $3$-dimensional involution $\op_2 \colon \nCat{2}_{\leq 2} \rightarrow \nCat{2}_{\leq 2}$, inside which we define a ($1$-dimensional promoted to) $2$-dimensional involution $\op_1$, and finally inside which we define a $1$-dimensional involution $(-)^*$. In general, it is not possible to obtain such triply (or higher) nested structures just by simply considering the microcosm principle. 
	\section{Dualities of monoidal structures} \label{sec:dualvdcs}
	We wish to describe various \enquote{transport} constructions for monoids and modules across different duals of a monoidal category $\cV$. To do this, we first briefly describe these relevant dualities on $\MonCat_{x}$, where $x$ is used as a placeholder for lax, colax, or pseudo. Next, we discuss how morphisms of involutions and inverses provide a natural framework to encode this dualities. For encoding the dualities of modules, virtual double categories shall play an essential role. Moreover, these dualities themselves generalize appropriately to virtual double categories.

\subsection{\tpdf{Dualities of $\MonCat$}} \label{sub:moncatdual}

The $\op_1$ duality (and its generalized inverse $\op_1^t$) can be extended to the $2$-category of monoidal categories and provides the following relation between the lax and colax variants.

\begin{lemma}
	There exists two functors \[\op_1 \colon \op_2(\MonCat_{c}) \rightarrow \MonCat_{l}\] and \[\op_1^t \colon  \op_2(\MonCat_{l}) \rightarrow \MonCat_{c}\] extending the involution $\op_1 \colon \op_2(\nCat{1}) \rightarrow \nCat{1}$ along the forgetful $2$-functors $U_{x} \colon\MonCat_{x} \rightarrow \Cat$. These functors are strict inverses of each other relative to the $\op_2$ involution. That is, they satisfy the constraints as defined by the following diagrams:
	\[\begin{tikzcd}[ampersand replacement=\&]
		{\MonCat_c} \&\& {\MonCat_l} \& \\
		\& {\op_2(\MonCat_l)} \&\& {\op_2(\MonCat_c)} \\
		{\MonCat_c} \&\& {\MonCat_l}
		\arrow["{\op_2(\op_1)}", from=1-1, to=2-2]
		\arrow["1"', from=1-1, to=3-1]
		\arrow["{\op_2(\op^t_1)}", from=1-3, to=2-4]
		\arrow["1"', from=1-3, to=3-3]
		\arrow["{\op^t_1}", from=2-2, to=3-1]
		\arrow["{\op_1}", from=2-4, to=3-3]
	\end{tikzcd}\cdpunct[4pt]{.}\]
	Furthermore, $U_x$ can be extended to a (strict) morphism of inverses and the $2$-functors $\op_1$ and $\op_1^t$ are isomorphic when restricted onto the variant $\MonCat_{p}$.
\end{lemma}

\begin{proof}
	For each diagram defining some structural monoidal data or coherence constraint, apply the functor $\op_1 \colon \op_2(\nCat{1}) \rightarrow \nCat{1}$ to obtain the required dual diagram. More precisely, $\op_1$ is a monoidal $\nCat{1}$-enriched functor of monoidal $\nCat{1}$-enriched categories. Thus, one could also adapt \cite[Thm.~4.5.1]{yauinvcat2020} in our situation by replacing the strict structures with pseudo structures.
\end{proof}

The same duality holds for $\BrMonCat_x$ and $\SymMonCat_x$. For instance, the restriction of $\op_1$ functor yields the commuting diagram
\[\begin{tikzcd}[ampersand replacement=\&]
	{\op_2(\SymMonCat_c)} \&\& {\op_2(\BrMonCat_c)} \&\& {\op_2(\MonCat_c)} \\
	\\
	{\SymMonCat_l} \&\& {\BrMonCat_l} \&\& {\MonCat_l}
	\arrow["{\op_2(U)}", from=1-1, to=1-3]
	\arrow["{\op_1}"', from=1-1, to=3-1]
	\arrow["{\op_2(U)}", from=1-3, to=1-5]
	\arrow["{\op_1}"', from=1-3, to=3-3]
	\arrow["{\op_1}", from=1-5, to=3-5]
	\arrow["U"', from=3-1, to=3-3]
	\arrow["U"', from=3-3, to=3-5]
\end{tikzcd}\cdpunct[8pt]{.}\]
Here, $\op_1$ maps a braided monoidal category as $(\cV, b_{x,y}) \mapsto (\op_1(\cV),b^{-1}_{x,y})$. And, each square in the diagram represents a part of the corresponding morphism of inverses.

\begin{remark}
	We can combine both $\op_1$ and $\op_1^t$ as a single involution on the double category of monoidal category (\cite[Sec.~4.2.6]{grandis2019dblcat}), when considered relative to the transpose involution (\cite[Sec.~3.2.2]{grandis2019dblcat}) of the \emph{$1$-category} of (strict) double categories. In particular, there exists a double functor $\op_1 \colon t(\mathsf{\mathbb{M}onCat}) \rightarrow \mathsf{\mathbb{M}onCat}$ that satisfies $t(\op_1) \circ \op_1 = 1$ or equivalently $\op_1 \circ t(\op_1) = 1$. In particular, we have the following equivalent diagram of morphisms of involutions:
	\[\begin{tikzcd}[ampersand replacement=\&]
		{t(\mathsf{Sym\mathbb{M}onCat})} \&\& {t(\mathsf{Br\mathbb{M}onCat})} \&\& {t(\mathsf{\mathbb{M}onCat})} \\
		\\
		{\mathsf{Sym\mathbb{M}onCat}} \&\& {\mathsf{Br\mathbb{M}onCat}} \&\& {\mathsf{\mathbb{M}onCat}}
		\arrow["{t(U)}", from=1-1, to=1-3]
		\arrow["{\op_1}"', from=1-1, to=3-1]
		\arrow["{t(U)}", from=1-3, to=1-5]
		\arrow["{\op_1}"', from=1-3, to=3-3]
		\arrow["{\op_1}", from=1-5, to=3-5]
		\arrow["U"', from=3-1, to=3-3]
		\arrow["U"', from=3-3, to=3-5]
	\end{tikzcd}\cdpunct[8pt]{.}\]
\end{remark}

The other important duality is that of reversal of a monoidal category. In particular, for each monoidal category $(\cV, \otimes ,I)$, there exists a corresponding \emph{reversed monoidal category} $(\cV, \otimes_r,I)$, which is denoted by $r(\cV)$. The monoidal product of $r(\cV)$ is defined using the composite $\cV \times \cV \xrightarrow{s} \cV \times \cV \xrightarrow{\otimes} \cV$, and the monoidal unit remains the same. The associator is given by $\left(a_{r(\cV)}\right)_{x,y,z}= \left(a_{\cV}\right)^{-1}_{z,y,x}$, and the unitors by $\lambda_{r(\cV)} = \rho_{\cV}$, $\rho_{r(\cV)} = \lambda_{\cV}$.

The reversal construction $r$ also extends to an involution. That is, $r \colon \MonCat_{x} \rightarrow \MonCat_{x}$ is a $2$-functor which satisfy $r^2=1$. Furthermore, this duality restricts to the $2$-categories $\BrMonCat_x$ and $\SymMonCat_x$ appropriately, as before. That is, there are commuting squares

\[\begin{tikzcd}[ampersand replacement=\&]
	{\SymMonCat_x} \&\& {\BrMonCat_x} \&\& {\MonCat_x} \\
	\\
	{\SymMonCat_x} \&\& {\BrMonCat_x} \&\& {\MonCat_x}
	\arrow["U", from=1-1, to=1-3]
	\arrow["r"', from=1-1, to=3-1]
	\arrow["U", from=1-3, to=1-5]
	\arrow["r"', from=1-3, to=3-3]
	\arrow["r", from=1-5, to=3-5]
	\arrow["U"', from=3-1, to=3-3]
	\arrow["U"', from=3-3, to=3-5]
\end{tikzcd}\cdpunct[8pt]{,}\]
representing the corresponding morphisms of involutions. Here, $r$ sends a braiding $b_{x,y}$ to the reversed braiding $b_{y,x}$.

\begin{remark}
	Analogous double categorical extension to $\mathsf{\mathbb{M}onCat}$ (and variants) can also be formulated for the reversal duality as relative to the identity involution of strict double categories. 
\end{remark}

\subsection{Transport of monoids in virtual double categories} \label{sub:monmoddualvdc}

It is well known that the category of monoids in a monoidal category can be identified using the isomorphisms $\Mon(\cV) \cong \MonCat_l(*,\cV)$. Similarly, commutative monoids in a braided monoidal category can be identified using $\cMon(\cV)\cong \BrMonCat_l(*,\cV)$. And if $\cV$ is symmetric, we have further isomorphisms $\BrMonCat_l(*,\cV) \cong \SymMonCat_l(*,\cV)$. Our goal is to understand various ways to transport monoids and comonoids using these identifications. As mentioned before, we also wish to describe these constructions within the framework of generalized involutions/inverses, with the main motivation being the intuition for what to expect in higher dimensions.  

We start with defining the $\op$ duality functors between monoids and comonoids for a given $\cV$. 

\begin{definition}
	For each monoidal category $\cV$, define a $1$-functor $\op_{\cV} \colon \op_1(\CoMon(\cV)) \rightarrow \Mon(\op_1(\cV))$ using the commutative diagram
	\[\begin{tikzcd}[ampersand replacement=\&]
		{\op_1(\CoMon(\cV))} \&\& {\Mon(\op_1(\cV))} \\
		{\op_1(\MonCat_c(*,\cV))} \&\& {\MonCat_l(*,\op_1(\cV))}
		\arrow["{\op_{\cV}}", from=1-1, to=1-3]
		\arrow["\cong", no head, from=1-1, to=2-1]
		\arrow["\cong"', no head, from=1-3, to=2-3]
		\arrow["{\op_1 (\cong)}"', from=2-1, to=2-3]
	\end{tikzcd}\cdpunct[8pt]{.}\]
	The dual version $\op^t_{\cV} \colon \op_1(\Mon(\cV)) \rightarrow \CoMon(\op_1(\cV))$ is defined similarly.
\end{definition}

We already observe that there is no obvious way to compose $\op_{\cV}$ and $\op^t_{\cV}$ to produce the identity map as all inverses do. At the same time, it is clear that this is a form of duality. The primary issue is that we are not incorporating the functoriality in $\cV$. In particular, both $\CoMon(-)$ and $\Mon(-)$ are strict 2-functors as they are defined by a corepresentation. Therefore, we ought to extend both $\op_{\cV}$ and $\op^{t}_{\cV}$ to strict $2$-natural transformations.

\begin{lemma}
	The functors $\op_{\cV}$ and $\op^t_{\cV}$ form the components of strict $2$-natural transformations
	\[\begin{tikzcd}[ampersand replacement=\&,column sep=scriptsize]
		{\op_2(\MonCat_c)} \&\&\& {\op_2(\Cat)} \& {\op_2(\MonCat_l)} \&\&\& {\op_2(\Cat)} \\
		\\
		\\
		{\MonCat_l} \&\&\& \Cat \& {\MonCat_c} \&\&\& \Cat
		\arrow["{\op_2(\CoMon(-))}", from=1-1, to=1-4]
		\arrow["{\op_1}"', from=1-1, to=4-1]
		\arrow["\op"', between={0.2}{0.8}, Rightarrow, nfold, from=1-4, to=4-1]
		\arrow["{\op_1}", from=1-4, to=4-4]
		\arrow["{\op_2(\Mon(-))}", from=1-5, to=1-8]
		\arrow["{\op^t_1}"', from=1-5, to=4-5]
		\arrow["{\op^t}"', between={0.2}{0.8}, Rightarrow, nfold, from=1-8, to=4-5]
		\arrow["{\op_1}", from=1-8, to=4-8]
		\arrow["{\Mon(-)}"', from=4-1, to=4-4]
		\arrow["{\CoMon(-)}"', from=4-5, to=4-8]
	\end{tikzcd}\cdpunct[8pt]{.}\]
	Furthermore, these two transformations form a morphism of inverses relative to the strict $2$-involution $\op_2 \colon \nCat{2}_{\leq 2} \rightarrow \nCat{2}_{\leq 2}$.
\end{lemma}

\begin{proof}
	The $2$-naturality and strictness of the two collection of maps $\op_{\cV}$ and $\op_{\cV}^{t}$ is a direct consequence of their definition in terms of the corepresentation by the terminal monoidal category. 
	
	The two axioms for morphism of inverses are similarly straightforward. For instance, one of these axioms can be verified as indicated below:
	\[\begin{tikzcd}[ampersand replacement=\&,cramped,column sep=small]
		\& {\MonCat_c} \&\& \Cat \&\& {\MonCat_c} \\
		\\
		{\op_2^2(\MonCat_c)} \&\& {\op_2^2(\Cat)} \\
		\&\&\&\& \equiv \\
		{\op_2(\MonCat_l)} \&\& {\op_2(\Cat)} \\
		\\
		\& {\MonCat_c} \&\& \Cat \&\& \Cat \\
		{\op_1(\Mon(\op_1(\cV)))} \&\&\&\&\& {\CoMon(\op_1^2(\cV))} \\
		\\
		{\CoMon(\cV)} \&\&\&\&\& {\CoMon(\cV)}
		\arrow["{\CoMon(-)}", from=1-2, to=1-4]
		\arrow[no head, from=1-2, to=3-1]
		\arrow[between={0.2}{0.8}, equals, nfold, from=1-4, to=3-1]
		\arrow[no head, from=1-4, to=3-3]
		\arrow[""{name=0, anchor=center, inner sep=0}, "{1_{\Cat}}", from=1-4, to=7-4]
		\arrow["{\CoMon(-)}"{description}, from=1-6, to=7-6]
		\arrow["{\op_2^2(\CoMon(-))}", from=3-1, to=3-3]
		\arrow["{\op_2(\op_1)}"', from=3-1, to=5-1]
		\arrow["{\op_2(\op)}"', between={0.2}{0.8}, Rightarrow, nfold, from=3-3, to=5-1]
		\arrow[""{name=1, anchor=center, inner sep=0}, "{\op_2(\op_1)}"{description}, from=3-3, to=5-3]
		\arrow["{\op_2(\Mon(-))}", from=5-1, to=5-3]
		\arrow["{\op^t_1}"', from=5-1, to=7-2]
		\arrow["{\op^t}"', between={0.2}{0.8}, Rightarrow, nfold, from=5-3, to=7-2]
		\arrow["{\op_1}"{description}, from=5-3, to=7-4]
		\arrow["{\CoMon(-)}"', from=7-2, to=7-4]
		\arrow["{\op^t_{\op_1{\cV}}}", from=8-1, to=8-6]
		\arrow["{\CoMon(1_{\cV})}", from=8-6, to=10-6]
		\arrow["{\op_{\cV}}", from=10-1, to=8-1]
		\arrow["{1_{\CoMon(\cV)}}"', from=10-1, to=10-6]
		\arrow[between={0.3}{0.8}, equals, nfold, from=1, to=0]
	\end{tikzcd}\cdpunct[8pt]{.}\]
\end{proof}

The phenomenon of transport across $r$ and $i$ is more interesting as considering the braiding or symmetry provides additional transport structure unlike $\op_1$. First, consider that $\cV$ is just an arbitrary monoidal category. In this case, we only have the expected reversal duality.

\begin{lemma}
	There exists a strict 2-natural transformation
	\[\begin{tikzcd}[ampersand replacement=\&]
		{\MonCat_l} \&\& \Cat \\
		\\
		{\MonCat_l} \&\& \Cat
		\arrow["{\Mon(-)}", from=1-1, to=1-3]
		\arrow["r"', from=1-1, to=3-1]
		\arrow["r"', between={0.2}{0.8}, Rightarrow, from=1-3, to=3-1]
		\arrow["1", from=1-3, to=3-3]
		\arrow["{\Mon(-)}"', from=3-1, to=3-3]
	\end{tikzcd}\cdpunct[8pt]{,}\]
	which defines a morphism of involutions relative to the identity functor $1 \colon \nCat{2}_{\leq 2} \rightarrow \nCat{2}_{\leq 2}$.
\end{lemma}

Now take $\cV$ to be a braided monoidal category. The key distinction from $\op_1$ duality is that $\BrMonCat$ admits two distinct non-trivial $2$-cells 
\[\begin{tikzcd}[ampersand replacement=\&,cramped]
	{\BrMonCat_l} \&\& {\BrMonCat_l} \& {\BrMonCat_l} \&\& {\BrMonCat_l} \\
	{b_{\cV} :(r(\cV),b_{y,x})} \&\& {(\cV,b_{x,y})} \& {b^{-1}_{\cV} :(r(\cV),b_{y,x})} \&\& {(\cV,b_{x,y})}
	\arrow[""{name=0, anchor=center, inner sep=0}, "r"', shift right=3, curve={height=6pt}, from=1-1, to=1-3]
	\arrow[""{name=1, anchor=center, inner sep=0}, "1", shift left=3, curve={height=-6pt}, from=1-1, to=1-3]
	\arrow[""{name=2, anchor=center, inner sep=0}, "1", shift left=3, curve={height=-6pt}, from=1-4, to=1-6]
	\arrow[""{name=3, anchor=center, inner sep=0}, "r"', shift right=3, curve={height=6pt}, from=1-4, to=1-6]
	\arrow["{(1,b_{x,y},1)}"', from=2-1, to=2-3]
	\arrow["{(1,b^{-1}_{y,x},1)}"', from=2-4, to=2-6]
	\arrow["b"', between={0.2}{0.8}, Rightarrow, nfold, from=0, to=1]
	\arrow["{b^{-1}}"', between={0.2}{0.8}, Rightarrow, nfold, from=3, to=2]
\end{tikzcd}\]
between the identity functor and the reversal functor. And, these $2$-cell $b$ and $b^{-1}$ are part of a morphism of inverses, albeit degenerately, 
\[\begin{tikzcd}[ampersand replacement=\&,cramped]
	\& {(\cV,b_{x,y})} \&\& {\BrMonCat_l} \\
	\\
	{(r(\cV),b_{y,x})} \&\&\& {\BrMonCat_l} \\
	\\
	\& {(\cV,b_{x,y})} \&\& {\BrMonCat_{l}}
	\arrow["{r(b_{\cV})\colon(1,b_{y,x},1)}"', from=1-2, to=3-1]
	\arrow["{1_{\cV} \colon (1,1,1)}", from=1-2, to=5-2]
	\arrow[""{name=0, anchor=center, inner sep=0}, "r", shift left=3, curve={height=-6pt}, from=1-4, to=3-4]
	\arrow[""{name=1, anchor=center, inner sep=0}, "1"', shift right=3, curve={height=6pt}, from=1-4, to=3-4]
	\arrow["{b^{-1}_{\cV}\colon (1,b^{-1}_{y,x},1)}"', from=3-1, to=5-2]
	\arrow[""{name=2, anchor=center, inner sep=0}, "1"', shift right=3, curve={height=6pt}, from=3-4, to=5-4]
	\arrow[""{name=3, anchor=center, inner sep=0}, "r", shift left=3, curve={height=-6pt}, from=3-4, to=5-4]
	\arrow["b", between={0.2}{0.8}, Rightarrow, nfold, from=0, to=1]
	\arrow["{b^{-1}}", between={0.2}{0.8}, Rightarrow, nfold, from=3, to=2]
\end{tikzcd}\cdpunct[4pt]{.}\]

Expectedly, if $\cV$ is a symmetric monoidal category then both the $2$-cells $b$ and $b^{-1}$ are identical and become part of a morphism of involutions. In general, we can compose this morphism of inverses $(b,b^{-1})$ to define other variants of the reversal duality $2$-cell, such as:
\[\begin{tikzcd}[ampersand replacement=\&]
	{\BrMonCat_l} \&\& {\MonCat_l} \&\& \Cat \\
	\\
	{\BrMonCat_l} \&\& {\MonCat_l} \&\& \Cat \\
	{\Mon(\cV)} \&\& {\Mon(r(\cV))} \&\& {\Mon(\cV)}
	\arrow["U", from=1-1, to=1-3]
	\arrow[""{name=0, anchor=center, inner sep=0}, "r", from=1-1, to=3-1]
	\arrow[""{name=1, anchor=center, inner sep=0}, "1"', shift right=3, curve={height=18pt}, from=1-1, to=3-1]
	\arrow["{\Mon(-)}", from=1-3, to=1-5]
	\arrow["r"', from=1-3, to=3-3]
	\arrow["r"', between={0.2}{0.8}, Rightarrow, nfold, from=1-5, to=3-3]
	\arrow["1"', from=1-5, to=3-5]
	\arrow["U"', from=3-1, to=3-3]
	\arrow["{\Mon(-)}"', from=3-3, to=3-5]
	\arrow["{r_{\cV}}", from=4-1, to=4-3]
	\arrow["{\Mon(b_{\cV})}", from=4-3, to=4-5]
	\arrow["b"', between={0.2}{0.8}, Rightarrow, nfold, from=0, to=1]
\end{tikzcd}\cdpunct[8pt]{.}\]
To illustrate the distinction between these variants of reversal, we provide the following example.

\begin{example}
	Let $\cV$ be the symmetric monoidal category $\Set$ with the symmetry given by the swap function $s_{x,y}: x \times y \rightarrow y \times x$. In this case, the component $b_{\Set}$ is given by the braided monoidal functor $(1_{\Set},s_{x,y},1) : (r(\Set),s_{y,x}) \rightarrow (\Set,s_{x,y})$. The morphism $r_{\Set}$ merely reinterprets a given monoid $M \in \Mon(\Set)$ as a monoid in $\Mon(r(\Set))$ in the sense depicted by
	\[\begin{tikzcd}[ampersand replacement=\&]
		{M \times M} \&\& {M \times_{r} M} \\
		\& {=} \\
		M \&\& M
		\arrow["\mu"', from=1-1, to=3-1]
		\arrow["\mu", from=1-3, to=3-3]
	\end{tikzcd}\cdpunct[8pt]{.}\]
	Thus, we should view the morphism $r_{\Set}$ as a comprehension duality. In contrast, the morphism $\Mon(b_{\Set})$ then transfers this monoid back to $\Mon(\Set)$ by precomposing with the symmetry as in
	\[\begin{tikzcd}[ampersand replacement=\&]
		\&\& {M \times M} \\
		{M \times_{r} M} \\
		\&\& {M \times_{r} M} \\
		M \\
		\&\& M
		\arrow["{s_{M,M}}", from=1-3, to=3-3]
		\arrow[""{name=0, anchor=center, inner sep=0}, "\mu"', from=2-1, to=4-1]
		\arrow["\mu", from=3-3, to=5-3]
		\arrow[between={0.2}{0.8}, maps to, from=0, to=3-3]
	\end{tikzcd}\cdpunct[8pt]{.}\]
	The key takeaway is that one could replace $\Mon$ with $\mathbf{PsMon}$ and $\Set$ with $\mathbf{\Cat}$ to draw similar observation about the duality $r : \MonCat \rightarrow \MonCat$ itself. We discuss this later in \cref{sub:highermonduals}.
\end{example}

A similar theory of dualities can be developed for arbitrary bimodules in a monoidal category. In particular, we could define a $2$-functor $\Mod(-)$ that takes the input as a monoidal category and outputs the bicategory of bimodules over arbitrary monoids (as long as we restrict our attention to monoidally cocomplete categories). However, this is rather unwieldy and introduces unnecessary assumptions. A more natural framework to encode these dualities is that of virtual double categories. See \cite[Defs.~2.1, 3.1, and 7.6]{shulman2010vequip} for a complete account and also see \cite[Ch.~5]{leinster2004} for our conventions on horizontal $1$-cells. In particular, we shall generally read horiztonal $1$-cells $M \colon A \rightarrow B$ as  $(B,A)$-bimodule.

Note that in the original definition, one can vertically reverse the direction of the square multicell to obtain the dual notion of a \enquote{co}-virtual double category. As the corresponding definitions are apparent in this dual case, we shall not write them explicitly. 

The only non-trivial result which we requires is that $\vcoDblCat$ and $\vDblCat$ exhibits a vertical duality involution in the following sense.

\begin{lemma}
	There exist two $2$-functors \[\op_{V}\colon \op_2(\vcoDblCat) \rightarrow \vDblCat\] and \[\op^t_{V} \colon \op_2(\vDblCat) \rightarrow \vcoDblCat,\]
	which form a morphism of inverses relative to the $\op_2$ involution. 
\end{lemma}

\begin{proof}
	we shall mainly define the action of $\op_V$ on each co-virtual double categories $\bV$. The rest is technical but straightforward. The data of $\op_{V}(\bV)$ is given by
	\begin{description}
		\item[Objects] 
		\[\begin{tikzcd}[ampersand replacement=\&,cramped]
			\bV \&\&\&\& {\op_V(\bV)} \\
			X \&\&\&\& X
		\end{tikzcd}\]
		\item[V $1$-cell] 
		\[\begin{tikzcd}[ampersand replacement=\&,cramped]
			\bV \&\&\&\& {\op_V(\bV)} \\
			X \&\&\&\& X \\
			\\
			Y \&\&\&\& Y
			\arrow[from=2-5, to=4-5]
			\arrow[from=4-1, to=2-1]
		\end{tikzcd}\]
		\item[H $1$-cell] 
		\[\begin{tikzcd}[ampersand replacement=\&,cramped,column sep=small]
			\& \bV \&\&\&\&\&\& {\op_V(\bV)} \& \\
			X \&\& Y \&\&\&\& X \&\& Y
			\arrow["\shortmid"{marking}, from=2-1, to=2-3]
			\arrow["\shortmid"{marking}, from=2-7, to=2-9]
		\end{tikzcd}\]
		\item[Sq multicell] 
		\[\begin{tikzcd}[ampersand replacement=\&,cramped,column sep=small]
			\&\& \bV \&\&\&\&\&\& {\op_V(\bV)} \&\& \\
			{X_0} \& \cdots \&\& \cdots \& {X_n} \&\& {X_0} \& \cdots \&\& \cdots \& {X_n} \\
			\\
			{Y_0} \&\&\&\& {Y_1} \&\& {Y_0} \&\&\&\& {Y_1}
			\arrow["\shortmid"{marking}, from=2-1, to=2-2]
			\arrow[""{name=0, anchor=center, inner sep=0}, "\shortmid"{marking}, from=2-2, to=2-4]
			\arrow["\shortmid"{marking}, from=2-4, to=2-5]
			\arrow["\shortmid"{marking}, from=2-7, to=2-8]
			\arrow[from=2-7, to=4-7]
			\arrow[""{name=1, anchor=center, inner sep=0}, "\shortmid"{marking}, from=2-8, to=2-10]
			\arrow["\shortmid"{marking}, from=2-10, to=2-11]
			\arrow[from=2-11, to=4-11]
			\arrow[from=4-1, to=2-1]
			\arrow[""{name=2, anchor=center, inner sep=0}, "\shortmid"{marking}, from=4-1, to=4-5]
			\arrow[from=4-5, to=2-5]
			\arrow[""{name=3, anchor=center, inner sep=0}, "\shortmid"{marking}, from=4-7, to=4-11]
			\arrow[between={0.2}{0.8}, Rightarrow, nfold, from=1, to=3]
			\arrow[between={0.2}{0.8}, Rightarrow, nfold, from=2, to=0]
		\end{tikzcd}\]
	\end{description}
	All of the data on the right is formal data similar to how the category $\cC^{\op}$ is defined when given a category $\cC$. Therefore, we simply induce the compositions and units of $\op_V(\bV)$ from $\bV$ itself in an analogous manner. 
	
	More technically, if the underlying fc-graph of a virtual double category $\bV$ is given by $T(C_0) \xleftarrow{s} C_1 \xrightarrow{t} C_0$, we can obtain a dual co-fc-graph given by $T(C_0) \xleftarrow{t'} C_1 \xrightarrow{s'} C_0$ with $t'=s$ and $s'=t$. 
	
\end{proof}

We briefly recall that there exists a chain of embeddings $\MonCat_{l} \hookrightarrow (\BiCat_{2})_{l} \hookrightarrow \vDblCat$. The first one is well known and simply identifies one object bicategories with monoidal categories \cite{lack2008icons}. The second one, we can describe as follows:
\[\begin{tikzcd}[ampersand replacement=\&]
	\& B \&\&\&\&\& {\mathbb{B}} \&\& \\
	\&\&\&\& {x_0} \& {x_1} \& {x_2} \& \cdots \& {x_n} \\
	{x_0} \&\& {x_n} \\
	\&\&\&\& {x_0} \&\&\&\& {x_n}
	\arrow[maps to, from=1-2, to=1-7]
	\arrow["{f_1}"{inner sep=.8ex}, "\shortmid"{marking}, from=2-5, to=2-6]
	\arrow["{\varepsilon_{x_0}}"', from=2-5, to=4-5]
	\arrow["{f_2}"{inner sep=.8ex}, "\shortmid"{marking}, from=2-6, to=2-7]
	\arrow["\shortmid"{marking}, from=2-7, to=2-8]
	\arrow["\shortmid"{marking}, from=2-8, to=2-9]
	\arrow["{\varepsilon_{x_n}}", from=2-9, to=4-9]
	\arrow[""{name=0, anchor=center, inner sep=0}, "{e(f_n ,f_{n-1}, \dots ,f_1)}", curve={height=-18pt}, from=3-1, to=3-3]
	\arrow[""{name=1, anchor=center, inner sep=0}, "{g_1}"', curve={height=18pt}, from=3-1, to=3-3]
	\arrow[""{name=2, anchor=center, inner sep=0}, "{g_1}"'{inner sep=.8ex}, "\shortmid"{marking}, from=4-5, to=4-9]
	\arrow["{\alpha \colon (f_n\cdots )f_2)f_1 \Rightarrow g_1}"{description}, between={0.2}{0.8}, Rightarrow, nfold, from=2-7, to=2]
	\arrow["\alpha", between={0.2}{0.8}, Rightarrow, nfold, from=0, to=1]
\end{tikzcd}\cdpunct[8pt]{.}\]

An analogous embedding $\MonCat_{c} \hookrightarrow (\BiCat_{2})_{c} \hookrightarrow \vcoDblCat$ holds for the colax case. The first one is again clear. For the second one, we have the similar following description:
\[\begin{tikzcd}[ampersand replacement=\&]
	\& B \&\&\&\&\& {\mathbb{B}} \&\& \\
	\&\&\&\& {x_0} \& {x_1} \& {x_2} \& \cdots \& {x_n} \\
	{x_0} \&\& {x_n} \\
	\&\&\&\& {x_0} \&\&\&\& {x_n}
	\arrow[maps to, from=1-2, to=1-7]
	\arrow["{f_1}"{inner sep=.8ex}, "\shortmid"{marking}, from=2-5, to=2-6]
	\arrow["{f_2}"{inner sep=.8ex}, "\shortmid"{marking}, from=2-6, to=2-7]
	\arrow["\shortmid"{marking}, from=2-7, to=2-8]
	\arrow["\shortmid"{marking}, from=2-8, to=2-9]
	\arrow[""{name=0, anchor=center, inner sep=0}, "{e(f_n ,f_{n-1}, \dots ,f_1)}", curve={height=-18pt}, from=3-1, to=3-3]
	\arrow[""{name=1, anchor=center, inner sep=0}, "{g_1}"', curve={height=18pt}, from=3-1, to=3-3]
	\arrow["{\varepsilon_{x_0}}", from=4-5, to=2-5]
	\arrow[""{name=2, anchor=center, inner sep=0}, "{g_1}"'{inner sep=.8ex}, "\shortmid"{marking}, from=4-5, to=4-9]
	\arrow["{\varepsilon_{x_n}}"', from=4-9, to=2-9]
	\arrow["\alpha"', between={0.2}{0.8}, Rightarrow, nfold, from=1, to=0]
	\arrow["{\alpha \colon g_1 \Rightarrow (f_n\cdots )f_2)f_1}"{description}, between={0.2}{0.8}, Rightarrow, nfold, from=2, to=2-7]
\end{tikzcd}\cdpunct[8pt]{.}\]

The important fact about these two embeddings is that they canonically extend to morphisms of involutions and morphisms of inverses. 

\begin{lemma}\leavevmode
	\begin{enumerate} 
		\item There exist morphisms of inverses whose one direction is given by
		\[\begin{tikzcd}[ampersand replacement=\&]
			{\op_2(\MonCat_{c})} \&\& {\op_2(\BiCat_{2})_{c}} \&\& {\op_2(\vcoDblCat)} \\
			\\
			{\MonCat_{l}} \&\& {(\BiCat_2)_l} \&\& \vDblCat
			\arrow[hook, from=1-1, to=1-3]
			\arrow["{\op_{1}}"', from=1-1, to=3-1]
			\arrow[hook, from=1-3, to=1-5]
			\arrow["{\op_{2}}", from=1-3, to=3-3]
			\arrow["\cong"{description}, between={0.2}{0.8}, Rightarrow, nfold, from=1-5, to=3-3]
			\arrow["{\op_{V}}", from=1-5, to=3-5]
			\arrow[hook, from=3-1, to=3-3]
			\arrow[hook, from=3-3, to=3-5]
		\end{tikzcd}\cdpunct[8pt]{,}\]
		and whose dual generalized inverse is given similarly.
		
		\item There exist morphisms of involution as given by
		\[\begin{tikzcd}[ampersand replacement=\&]
			{\MonCat_{l}} \&\& {(\BiCat_2)_l} \&\& \vDblCat \\
			\\
			{\MonCat_{l}} \&\& {(\BiCat_2)_l} \&\& \vDblCat
			\arrow[hook, from=1-1, to=1-3]
			\arrow["r"', from=1-1, to=3-1]
			\arrow[hook, from=1-3, to=1-5]
			\arrow["{\op_1}"', from=1-3, to=3-3]
			\arrow["\cong"{description}, between={0.2}{0.8}, Rightarrow, nfold, from=1-5, to=3-3]
			\arrow["{\op_{H}}", from=1-5, to=3-5]
			\arrow[hook, from=3-1, to=3-3]
			\arrow[hook, from=3-3, to=3-5]
		\end{tikzcd}\cdpunct[8pt]{.}\]
	\end{enumerate}
\end{lemma}

Note that the $2$-cell appearing in the involution case primarily arises from the fact that $((f_n \circ f_{n-1}) \circ \cdots ) \circ  f_1 \neq f_n \circ (f_{n-1} \circ (\cdots \circ f_1))$ in an arbitrary bicategory. To illustrate this issue further, consider a left-biased lax monoidal category $\cV$, that is, the associator need not an invertible $2$-cell. Applying $r$ to this lax monoidal category doesn't yield a left-biased lax monoidal category, as there is no $a^{-1}$ available. At best, we obtain a canonical right-biased lax monoidal category. All of this suggests that a more thorough approach would require generalized inverses $r \colon \MonCat_l^{Lb} \rightarrow \MonCat_l^{Rb}$ and $r^t\colon \MonCat_l^{Rb}  \rightarrow \MonCat_{l}^{Lb}$. However, we are content with our approach here. 

We have now established that monoidal categories are a special case of a virtual double categories, as well as co-virtual double categories. This suggests that the dualities of monoidal categories should also extend to the virtual double categorical case. We begin by characterization the $\op_1$ duality.

\begin{theorem}
	The $2$-functors $\Mod \colon \vDblCat \rightarrow \vDblCat$ and $\CoMod \colon \vcoDblCat \rightarrow \vcoDblCat$ extend to a $1$-morphism of relative inverses by the strict $2$-natural transformations
	\[\begin{tikzcd}[ampersand replacement=\&]
		{\op_2(\vcoDblCat)} \&\&\& {\op_2(\vcoDblCat)} \\
		\\
		\vDblCat \&\&\& \vDblCat \\
		{\op_2(\vDblCat)} \&\&\& {\op_2(\vDblCat)} \\
		\\
		\vcoDblCat \&\&\& \vcoDblCat
		\arrow["{\op_2(\CoMod(-))}", from=1-1, to=1-4]
		\arrow["{\op_V}"', from=1-1, to=3-1]
		\arrow["\op"', between={0.2}{0.8}, Rightarrow, nfold, from=1-4, to=3-1]
		\arrow["{\op_V}", from=1-4, to=3-4]
		\arrow["{\Mod(-)}"', from=3-1, to=3-4]
		\arrow["{\op_2(\Mod(-))}", from=4-1, to=4-4]
		\arrow["{\op_{V}^{t}}"', from=4-1, to=6-1]
		\arrow["{\op^{t}}"', between={0.2}{0.8}, Rightarrow, nfold, from=4-4, to=6-1]
		\arrow["{\op_{V}^{t}}", from=4-4, to=6-4]
		\arrow["{\CoMod(-)}"', from=6-1, to=6-4]
	\end{tikzcd}\cdpunct[8pt]{.}\]
\end{theorem}

\begin{proof}
	As usual, we need to define the components $\op_{\bV} : \op_{V}(\CoMod(\bV)) \rightarrow \Mod(\op_{V}(\bV))$.
	\begin{description}
		\item[Objects] 
		\[\begin{tikzcd}[ampersand replacement=\&]
			\& {C \coloneqq (\bullet,C,\mu,e)} \&\&\&\& {C^{\op} \coloneqq (\bullet,C, \mu^{\op}, e^{\op})} \& \\
			\bullet \& \bullet \& \bullet \&\& \bullet \& \bullet \& \bullet \\
			\bullet \&\& \bullet \&\& \bullet \&\& \bullet \\
			\& \bullet \&\&\&\& \bullet \\
			\bullet \&\& \bullet \&\& \bullet \&\& \bullet
			\arrow[between={0.2}{0.8}, maps to, from=1-2, to=1-6]
			\arrow["C"{inner sep=.8ex}, "\shortmid"{marking}, from=2-1, to=2-2]
			\arrow["C"{inner sep=.8ex}, "\shortmid"{marking}, from=2-2, to=2-3]
			\arrow["C"{inner sep=.8ex}, "\shortmid"{marking}, from=2-5, to=2-6]
			\arrow[no head, from=2-5, to=3-5]
			\arrow["C"{inner sep=.8ex}, "\shortmid"{marking}, from=2-6, to=2-7]
			\arrow[no head, from=2-7, to=3-7]
			\arrow[no head, from=3-1, to=2-1]
			\arrow[""{name=0, anchor=center, inner sep=0}, "C"'{inner sep=.8ex}, "\shortmid"{marking}, from=3-1, to=3-3]
			\arrow[no head, from=3-3, to=2-3]
			\arrow[""{name=1, anchor=center, inner sep=0}, "C"'{inner sep=.8ex}, "\shortmid"{marking}, from=3-5, to=3-7]
			\arrow[no head, from=4-6, to=5-5]
			\arrow[no head, from=4-6, to=5-7]
			\arrow[no head, from=5-1, to=4-2]
			\arrow[""{name=2, anchor=center, inner sep=0}, "C"'{inner sep=.8ex}, "\shortmid"{marking}, from=5-1, to=5-3]
			\arrow[no head, from=5-3, to=4-2]
			\arrow[""{name=3, anchor=center, inner sep=0}, "C"'{inner sep=.8ex}, "\shortmid"{marking}, from=5-5, to=5-7]
			\arrow["{\mu^{\op}}", between={0}{0.8}, Rightarrow, nfold, from=2-6, to=1]
			\arrow["\mu", between={0.2}{1}, Rightarrow, nfold, from=0, to=2-2]
			\arrow["{e^{\op}}", between={0.2}{0.8}, Rightarrow, nfold, from=4-6, to=3]
			\arrow["e", between={0.2}{1}, Rightarrow, nfold, from=2, to=4-2]
		\end{tikzcd}\]
		\item[V $1$-cell] 
		\[\begin{tikzcd}[ampersand replacement=\&]
			\bullet \&\& \bullet \&\&\& \bullet \&\& \bullet \\
			\\
			\star \&\& \star \&\&\& \star \&\& \star
			\arrow[""{name=0, anchor=center, inner sep=0}, "C"{inner sep=.8ex}, "\shortmid"{marking}, from=1-1, to=1-3]
			\arrow[""{name=1, anchor=center, inner sep=0}, "{C^{\op}}"{inner sep=.8ex}, "\shortmid"{marking}, from=1-6, to=1-8]
			\arrow[""{name=2, anchor=center, inner sep=0}, "{f_0^{\op}}"', from=1-6, to=3-6]
			\arrow["{f_0^{\op}}", from=1-8, to=3-8]
			\arrow["{f_{0}}", from=3-1, to=1-1]
			\arrow[""{name=3, anchor=center, inner sep=0}, "D"'{inner sep=.8ex}, "\shortmid"{marking}, from=3-1, to=3-3]
			\arrow[""{name=4, anchor=center, inner sep=0}, "{f_{0}}"', from=3-3, to=1-3]
			\arrow[""{name=5, anchor=center, inner sep=0}, "{D^{\op}}"'{inner sep=.8ex}, "\shortmid"{marking}, from=3-6, to=3-8]
			\arrow["{f^{\op}}", between={0.2}{0.8}, Rightarrow, nfold, from=1, to=5]
			\arrow["f"', between={0.2}{0.8}, Rightarrow, nfold, from=3, to=0]
			\arrow[between={0.2}{0.8}, maps to, from=4, to=2]
		\end{tikzcd}\]
		\item[H $1$-cell]
		\[\begin{tikzcd}[ampersand replacement=\&,column sep=small]
			{(C,\bullet)} \&\&\&\& {(D,\star)} \&\& {(C^{\op},\bullet)} \&\&\&\& {(D^{\op},\star)} \\
			\bullet \&\&\&\& \star \&\& \bullet \&\&\&\& \star \\
			\bullet \&\& \bullet \&\& \star \&\& \bullet \&\& \bullet \&\& \star \\
			\bullet \&\&\&\& \star \&\& \bullet \&\&\&\& \star \\
			\bullet \&\& \star \&\& \star \&\& \bullet \&\& \star \&\& \star \\
			\bullet \&\&\&\& \star \&\& \bullet \&\&\&\& \star
			\arrow["{X \coloneqq(X,\lambda,\rho)}"{inner sep=.8ex}, "\shortmid"{marking}, from=1-1, to=1-5]
			\arrow[between={0.2}{0.8}, maps to, from=1-5, to=1-7]
			\arrow["{X^{\op}\coloneqq(X,\lambda^{\op},\rho^{\op})}"{inner sep=.8ex}, "\shortmid"{marking}, from=1-7, to=1-11]
			\arrow["X"{inner sep=.8ex}, "\shortmid"{marking}, from=2-1, to=2-5]
			\arrow["X"{inner sep=.8ex}, "\shortmid"{marking}, from=2-7, to=2-11]
			\arrow["C"{inner sep=.8ex}, "\shortmid"{marking}, from=3-1, to=3-3]
			\arrow["X"{inner sep=.8ex}, "\shortmid"{marking}, from=3-3, to=3-5]
			\arrow["{C^{\op}}"{inner sep=.8ex}, "\shortmid"{marking}, from=3-7, to=3-9]
			\arrow[no head, from=3-7, to=4-7]
			\arrow["{X^{\op}}"{inner sep=.8ex}, "\shortmid"{marking}, from=3-9, to=3-11]
			\arrow[no head, from=3-11, to=4-11]
			\arrow[no head, from=4-1, to=3-1]
			\arrow[""{name=0, anchor=center, inner sep=0}, "X"'{inner sep=.8ex}, "\shortmid"{marking}, from=4-1, to=4-5]
			\arrow[no head, from=4-5, to=3-5]
			\arrow[""{name=1, anchor=center, inner sep=0}, "{X^{\op}}"'{inner sep=.8ex}, "\shortmid"{marking}, from=4-7, to=4-11]
			\arrow["X"{inner sep=.8ex}, "\shortmid"{marking}, from=5-1, to=5-3]
			\arrow[no head, from=5-1, to=6-1]
			\arrow["D"{inner sep=.8ex}, "\shortmid"{marking}, from=5-3, to=5-5]
			\arrow[no head, from=5-5, to=6-5]
			\arrow["{X^{\op}}"{inner sep=.8ex}, "\shortmid"{marking}, from=5-7, to=5-9]
			\arrow[no head, from=5-7, to=6-7]
			\arrow["{D^{\op}}"{inner sep=.8ex}, "\shortmid"{marking}, from=5-9, to=5-11]
			\arrow[no head, from=5-11, to=6-11]
			\arrow[""{name=2, anchor=center, inner sep=0}, "X"'{inner sep=.8ex}, "\shortmid"{marking}, from=6-1, to=6-5]
			\arrow[""{name=3, anchor=center, inner sep=0}, "{X^{\op}}"'{inner sep=.8ex}, "\shortmid"{marking}, from=6-7, to=6-11]
			\arrow["{\rho^{\op}}"', between={0}{0.8}, Rightarrow, nfold, from=3-9, to=1]
			\arrow["\rho", between={0.2}{1}, Rightarrow, nfold, from=0, to=3-3]
			\arrow["{\lambda^{\op}}"', between={0}{0.8}, Rightarrow, nfold, from=5-9, to=3]
			\arrow["\lambda", between={0.2}{1}, Rightarrow, nfold, from=2, to=5-3]
		\end{tikzcd}\]
		\item[Sq multicell]
		\[\begin{tikzcd}[ampersand replacement=\&]
			{\bullet_1} \&\& \cdots \&\& {\bullet_n} \&\& {\bullet_1} \&\& \cdots \&\& {\bullet_n} \\
			\\
			\bullet \&\&\&\& \star \&\& \bullet \&\&\&\& \star
			\arrow["{X_1}"{inner sep=.8ex}, "\shortmid"{marking}, from=1-1, to=1-3]
			\arrow["{X_n}"{inner sep=.8ex}, "\shortmid"{marking}, from=1-3, to=1-5]
			\arrow["{X_1^{\op}}"{inner sep=.8ex}, "\shortmid"{marking}, from=1-7, to=1-9]
			\arrow[""{name=0, anchor=center, inner sep=0}, "{f^{\op}}"', from=1-7, to=3-7]
			\arrow["{X_n^{\op}}"{inner sep=.8ex}, "\shortmid"{marking}, from=1-9, to=1-11]
			\arrow["{g^{\op}}", from=1-11, to=3-11]
			\arrow["f", from=3-1, to=1-1]
			\arrow[""{name=1, anchor=center, inner sep=0}, "Y"'{inner sep=.8ex}, "\shortmid"{marking}, from=3-1, to=3-5]
			\arrow[""{name=2, anchor=center, inner sep=0}, "g"', from=3-5, to=1-5]
			\arrow[""{name=3, anchor=center, inner sep=0}, "{Y^{\op}}"'{inner sep=.8ex}, "\shortmid"{marking}, from=3-7, to=3-11]
			\arrow["{\alpha^{\op}}", between={0}{0.8}, Rightarrow, nfold, from=1-9, to=3]
			\arrow["\alpha"', between={0.2}{1}, Rightarrow, nfold, from=1, to=1-3]
			\arrow[between={0.3}{0.7}, maps to, from=2, to=0]
		\end{tikzcd}\]
	\end{description}
	The other family of components can be defined similarly.
\end{proof}

The case for reversal duality is exactly the same as above.

\begin{theorem}
	The $2$-functor $\Mod \colon \vDblCat \rightarrow \vDblCat$ extends to a $1$-morphism of involutions by the strict $2$-natural transformation
	\[\begin{tikzcd}[ampersand replacement=\&]
		\vDblCat \&\& \vDblCat \\
		\\
		\vDblCat \&\& \vDblCat
		\arrow["{\Mod(-)}", from=1-1, to=1-3]
		\arrow["{\op_H}"', from=1-1, to=3-1]
		\arrow["r", between={0.2}{0.8}, Rightarrow, nfold, from=1-3, to=3-1]
		\arrow["{\op_H}", from=1-3, to=3-3]
		\arrow["{\Mod(-)}"', from=3-1, to=3-3]
	\end{tikzcd}\cdpunct[8pt]{.}\]
\end{theorem}

\begin{proof} The components $r_{\bV} \colon \op_H(\Mod(\bV)) \rightarrow \Mod(\op_{H}(\bV))$ are defined by 
	\begin{description}
		\item[Objects] 
		\[\begin{tikzcd}[ampersand replacement=\&]
			\& {M \coloneqq (\bullet,M,\mu,e)} \&\&\&\& {M^r \coloneqq (\bullet,M^r,\mu,e)} \& \\
			\bullet \& \bullet \& \bullet \&\& \bullet \& \bullet \& \bullet \\
			\bullet \&\& \bullet \&\& \bullet \&\& \bullet \\
			\& \bullet \&\&\&\& \bullet \\
			\bullet \&\& \bullet \&\& \bullet \&\& \bullet
			\arrow[between={0.2}{0.8}, maps to, from=1-2, to=1-6]
			\arrow["M"'{inner sep=.8ex}, "\shortmid"{marking}, from=2-2, to=2-1]
			\arrow["M"'{inner sep=.8ex}, "\shortmid"{marking}, from=2-3, to=2-2]
			\arrow["{M^r}"{inner sep=.8ex}, "\shortmid"{marking}, from=2-5, to=2-6]
			\arrow[no head, from=2-5, to=3-5]
			\arrow["{M^r}"{inner sep=.8ex}, "\shortmid"{marking}, from=2-6, to=2-7]
			\arrow[no head, from=2-7, to=3-7]
			\arrow[no head, from=3-1, to=2-1]
			\arrow[no head, from=3-3, to=2-3]
			\arrow[""{name=0, anchor=center, inner sep=0}, "M"{inner sep=.8ex}, "\shortmid"{marking}, from=3-3, to=3-1]
			\arrow[""{name=1, anchor=center, inner sep=0}, "{M^r}"'{inner sep=.8ex}, "\shortmid"{marking}, from=3-5, to=3-7]
			\arrow[no head, from=4-6, to=5-5]
			\arrow[no head, from=4-6, to=5-7]
			\arrow[no head, from=5-1, to=4-2]
			\arrow[no head, from=5-3, to=4-2]
			\arrow[""{name=2, anchor=center, inner sep=0}, "M"{inner sep=.8ex}, "\shortmid"{marking}, from=5-3, to=5-1]
			\arrow[""{name=3, anchor=center, inner sep=0}, "{M^r}"'{inner sep=.8ex}, "\shortmid"{marking}, from=5-5, to=5-7]
			\arrow["\mu"', between={0}{0.8}, Rightarrow, nfold, from=2-2, to=0]
			\arrow["\mu", between={0}{0.8}, Rightarrow, nfold, from=2-6, to=1]
			\arrow["e"', between={0}{0.8}, Rightarrow, nfold, from=4-2, to=2]
			\arrow["e", between={0.2}{0.8}, Rightarrow, nfold, from=4-6, to=3]
		\end{tikzcd}\]
		\item[V $1$-cell] 
		\[\begin{tikzcd}[ampersand replacement=\&]
			\bullet \&\& \bullet \&\&\& \bullet \&\& \bullet \\
			\\
			\star \&\& \star \&\&\& \star \&\& \star
			\arrow["{f_{0}}"', from=1-1, to=3-1]
			\arrow[""{name=0, anchor=center, inner sep=0}, "M"'{inner sep=.8ex}, "\shortmid"{marking}, from=1-3, to=1-1]
			\arrow[""{name=1, anchor=center, inner sep=0}, "{f_{0}}", from=1-3, to=3-3]
			\arrow[""{name=2, anchor=center, inner sep=0}, "{M^r}"{inner sep=.8ex}, "\shortmid"{marking}, from=1-6, to=1-8]
			\arrow[""{name=3, anchor=center, inner sep=0}, "{f_0}"', from=1-6, to=3-6]
			\arrow["{f_0}", from=1-8, to=3-8]
			\arrow[""{name=4, anchor=center, inner sep=0}, "N"{inner sep=.8ex}, "\shortmid"{marking}, from=3-3, to=3-1]
			\arrow[""{name=5, anchor=center, inner sep=0}, "{N^r}"'{inner sep=.8ex}, "\shortmid"{marking}, from=3-6, to=3-8]
			\arrow["f", between={0.2}{0.8}, Rightarrow, nfold, from=0, to=4]
			\arrow[between={0.2}{0.8}, maps to, from=1, to=3]
			\arrow["f", between={0.2}{0.8}, Rightarrow, nfold, from=2, to=5]
		\end{tikzcd}\]
		\item[H $1$-cell]
		\[\begin{tikzcd}[ampersand replacement=\&,column sep=small]
			{(N,\star)} \&\&\&\& {(M,\bullet)} \&\& {(N^r,\star)} \&\&\&\& {(M^r,\bullet)} \\
			\star \&\&\&\& \bullet \&\& \star \&\&\&\& \bullet \\
			\star \&\& \bullet \&\& \bullet \&\& \star \&\& \bullet \&\& \bullet \\
			\star \&\&\&\& \bullet \&\& \star \&\&\&\& \bullet \\
			\star \&\& \star \&\& \bullet \&\& \star \&\& \star \&\& \bullet \\
			\star \&\&\&\& \bullet \&\& \star \&\&\&\& \bullet
			\arrow["{X \coloneqq(X,\lambda,\rho)}"'{inner sep=.8ex}, "\shortmid"{marking}, from=1-5, to=1-1]
			\arrow[between={0.2}{0.8}, maps to, from=1-5, to=1-7]
			\arrow["{X^{r}\coloneqq(X^r,\rho,\lambda)}"{inner sep=.8ex}, "\shortmid"{marking}, from=1-7, to=1-11]
			\arrow["X"'{inner sep=.8ex}, "\shortmid"{marking}, from=2-5, to=2-1]
			\arrow["{X^r}"{inner sep=.8ex}, "\shortmid"{marking}, from=2-7, to=2-11]
			\arrow["X"'{inner sep=.8ex}, "\shortmid"{marking}, from=3-3, to=3-1]
			\arrow["M"'{inner sep=.8ex}, "\shortmid"{marking}, from=3-5, to=3-3]
			\arrow["{X^r}"{inner sep=.8ex}, "\shortmid"{marking}, from=3-7, to=3-9]
			\arrow[no head, from=3-7, to=4-7]
			\arrow["{M^r}"{inner sep=.8ex}, "\shortmid"{marking}, from=3-9, to=3-11]
			\arrow[no head, from=3-11, to=4-11]
			\arrow[no head, from=4-1, to=3-1]
			\arrow[no head, from=4-5, to=3-5]
			\arrow[""{name=0, anchor=center, inner sep=0}, "X"{inner sep=.8ex}, "\shortmid"{marking}, from=4-5, to=4-1]
			\arrow[""{name=1, anchor=center, inner sep=0}, "{X^r}"'{inner sep=.8ex}, "\shortmid"{marking}, from=4-7, to=4-11]
			\arrow[no head, from=5-1, to=6-1]
			\arrow["N"'{inner sep=.8ex}, "\shortmid"{marking}, from=5-3, to=5-1]
			\arrow["X"'{inner sep=.8ex}, "\shortmid"{marking}, from=5-5, to=5-3]
			\arrow[no head, from=5-5, to=6-5]
			\arrow["{N^r}"{inner sep=.8ex}, "\shortmid"{marking}, from=5-7, to=5-9]
			\arrow[no head, from=5-7, to=6-7]
			\arrow["{X^r}"{inner sep=.8ex}, "\shortmid"{marking}, from=5-9, to=5-11]
			\arrow[no head, from=5-11, to=6-11]
			\arrow[""{name=2, anchor=center, inner sep=0}, "X"{inner sep=.8ex}, "\shortmid"{marking}, from=6-5, to=6-1]
			\arrow[""{name=3, anchor=center, inner sep=0}, "{X^r}"'{inner sep=.8ex}, "\shortmid"{marking}, from=6-7, to=6-11]
			\arrow["\rho"', between={0}{0.8}, Rightarrow, nfold, from=3-3, to=0]
			\arrow["\rho"', between={0}{0.8}, Rightarrow, nfold, from=3-9, to=1]
			\arrow["\lambda"', between={0}{0.8}, Rightarrow, nfold, from=5-3, to=2]
			\arrow["\lambda"', between={0}{0.8}, Rightarrow, nfold, from=5-9, to=3]
		\end{tikzcd}\]
		\item[Sq multicell]
		\[\begin{tikzcd}[ampersand replacement=\&]
			{\bullet_n} \&\& \cdots \&\& {\bullet_1} \&\& {\bullet_n} \&\& \cdots \&\& {\bullet_1} \\
			\\
			\star \&\&\&\& \bullet \&\& \star \&\&\&\& \bullet
			\arrow["g"', from=1-1, to=3-1]
			\arrow["{X_n}"'{inner sep=.8ex}, "\shortmid"{marking}, from=1-3, to=1-1]
			\arrow["{X_1}"'{inner sep=.8ex}, "\shortmid"{marking}, from=1-5, to=1-3]
			\arrow[""{name=0, anchor=center, inner sep=0}, "f", from=1-5, to=3-5]
			\arrow["{X_n^r}"{inner sep=.8ex}, "\shortmid"{marking}, from=1-7, to=1-9]
			\arrow[""{name=1, anchor=center, inner sep=0}, "g"', from=1-7, to=3-7]
			\arrow["{X_1^r}"{inner sep=.8ex}, "\shortmid"{marking}, from=1-9, to=1-11]
			\arrow["f", from=1-11, to=3-11]
			\arrow[""{name=2, anchor=center, inner sep=0}, "Y"{inner sep=.8ex}, "\shortmid"{marking}, from=3-5, to=3-1]
			\arrow[""{name=3, anchor=center, inner sep=0}, "{Y^r}"'{inner sep=.8ex}, "\shortmid"{marking}, from=3-7, to=3-11]
			\arrow["\alpha", between={0.1}{0.8}, Rightarrow, nfold, from=1-3, to=2]
			\arrow[between={0.3}{0.7}, maps to, from=0, to=1]
			\arrow["\alpha", between={0.1}{0.8}, Rightarrow, nfold, from=1-9, to=3]
		\end{tikzcd}\]
	\end{description}
\end{proof}

We shall now specialize back to the case of monoidal categories. In particular, we shall restrict $\Mod(-)$ and $\CoMod(-)$ to monoidal categories.  In what follows, we write our diagrams in strict monoidal categories $\st \cV$ and interpret the $2$-cells using the canonical evaluation map.

\begin{proposition}
	The $(\op,\op^t)$ morphism of inverses on monoidal categories extends to virtual double categories through the following strict $2$-transformation cells:
	\[\begin{tikzcd}[ampersand replacement=\&]
		{\op_2(\MonCat_c)} \&\&\& {\op_2(\vcoDblCat)} \\
		\\
		{\MonCat_l} \&\&\& \vDblCat \\
		{\op_2(\MonCat_l)} \&\&\& {\op_2(\vDblCat)} \\
		\\
		{\MonCat_c} \&\&\& \vcoDblCat
		\arrow["{\op_2(\CoMod(-))}", from=1-1, to=1-4]
		\arrow["{\op_1}"', from=1-1, to=3-1]
		\arrow["\op", between={0.2}{0.8}, Rightarrow, nfold, from=1-4, to=3-1]
		\arrow["{\op_V}", from=1-4, to=3-4]
		\arrow["{\Mod(-)}"', from=3-1, to=3-4]
		\arrow["{\op_2(\Mod(-))}", from=4-1, to=4-4]
		\arrow["{\op_1^{t}}"', from=4-1, to=6-1]
		\arrow["{\op^{t}}", between={0.2}{0.8}, Rightarrow, nfold, from=4-4, to=6-1]
		\arrow["{\op_{V}^{t}}", from=4-4, to=6-4]
		\arrow["{\CoMod(-)}"', from=6-1, to=6-4]
	\end{tikzcd}\cdpunct[8pt]{.}\]
\end{proposition}

\begin{proof}
	We primarily need to compute the components $\op_{\cV}$ and $\op^t_{\cV}$. For the component $\op_{\cV}$, this means that we need to obtain a virtual double functor $\op_{\cV}\colon \op_{V}(\CoMod(\cV)) \rightarrow \Mod(\op_{1}(\cV))$. This is given by the following data:
	\begin{description}
		\item[Objects] 
		\[\begin{tikzcd}[ampersand replacement=\&]
			{C \coloneqq (C,\mu,e)} \&\& {C^{\op} \coloneqq (C, \mu^{\op}, e^{\op})}
			\arrow[maps to, from=1-1, to=1-3]
		\end{tikzcd}\]
		\item[V $1$-cell] 
		\[\begin{tikzcd}[ampersand replacement=\&]
			{C'} \&\&\& {(C')^{\op}} \\
			\\
			C \&\&\& {C^{\op}}
			\arrow[""{name=0, anchor=center, inner sep=0}, from=1-4, to=3-4]
			\arrow[""{name=1, anchor=center, inner sep=0}, from=3-1, to=1-1]
			\arrow[between={0.2}{0.8}, maps to, from=1, to=0]
		\end{tikzcd}\]
		\item[H $1$-cell] 
		\[\begin{tikzcd}[ampersand replacement=\&]
			{C_1} \&\&\& {C_2} \&\&\& {C_1^{\op}} \&\&\& {C_2^{\op}}
			\arrow[""{name=0, anchor=center, inner sep=0}, "{X \coloneqq(X,\lambda,\rho)}"{inner sep=.8ex}, "\shortmid"{marking}, from=1-1, to=1-4]
			\arrow[""{name=1, anchor=center, inner sep=0}, "{X^{\op} \coloneqq (X ,\lambda^{\op},\rho^{\op})}"{inner sep=.8ex}, "\shortmid"{marking}, from=1-7, to=1-10]
			\arrow[between={0.4}{0.6}, maps to, from=0, to=1]
		\end{tikzcd}\]
		\item[Sq multicell] 
		\[\begin{tikzcd}[ampersand replacement=\&,cramped,column sep=scriptsize]
			{C_1'} \& \cdots \& \cdots \& {C_2'} \&\& {(C_1')^{\op}} \& \cdots \& \cdots \& {(C_2')^{\op}} \\
			\\
			{C_1} \&\&\& {C_2} \&\& {C_1^{\op}} \&\&\& {C_2^{\op}}
			\arrow["Y"{inner sep=.8ex}, "\shortmid"{marking}, from=1-1, to=1-2]
			\arrow[""{name=0, anchor=center, inner sep=0}, "\shortmid"{marking}, from=1-2, to=1-3]
			\arrow["Z"{inner sep=.8ex}, "\shortmid"{marking}, from=1-3, to=1-4]
			\arrow["{Y^{\op}}"{inner sep=.8ex}, "\shortmid"{marking}, from=1-6, to=1-7]
			\arrow[""{name=1, anchor=center, inner sep=0}, from=1-6, to=3-6]
			\arrow[""{name=2, anchor=center, inner sep=0}, "\shortmid"{marking}, from=1-7, to=1-8]
			\arrow["{Z^{\op}}"{inner sep=.8ex}, "\shortmid"{marking}, from=1-8, to=1-9]
			\arrow[from=1-9, to=3-9]
			\arrow[from=3-1, to=1-1]
			\arrow[""{name=3, anchor=center, inner sep=0}, "X"'{inner sep=.8ex}, "\shortmid"{marking}, from=3-1, to=3-4]
			\arrow[""{name=4, anchor=center, inner sep=0}, from=3-4, to=1-4]
			\arrow[""{name=5, anchor=center, inner sep=0}, "{X^{\op}}"'{inner sep=.8ex}, "\shortmid"{marking}, from=3-6, to=3-9]
			\arrow["{f^{\op} \colon Z\otimes \cdots \otimes Y\rightarrow X}"{description}, between={0.2}{0.8}, Rightarrow, nfold, from=2, to=5]
			\arrow["{f\colon X \rightarrow Z \otimes \cdots \otimes Y}"{description}, between={0.2}{0.8}, Rightarrow, nfold, from=3, to=0]
			\arrow[between={0.2}{0.8}, maps to, from=4, to=1]
		\end{tikzcd}\]
	\end{description}
	The virtual double functor $\op^{t}_{\cV} \colon \op_{V}^{t}(\Mod(\cV)) \rightarrow \CoMod(\op_1(\cV))$ can be obtained similarly.
\end{proof}

Finally, we shall similarly extend the reversal duality. 

\begin{proposition}
	The $r$ morphism of involution on monoidal categories extends to the virtual double categories case through the following strict $2$-transformation cells:
	\[\begin{tikzcd}[ampersand replacement=\&]
		{\MonCat_l} \&\& \vDblCat \\
		\\
		{\MonCat_l} \&\& \vDblCat
		\arrow["{\Mod(-)}", from=1-1, to=1-3]
		\arrow["r"', from=1-1, to=3-1]
		\arrow["r", between={0.2}{0.8}, Rightarrow, nfold, from=1-3, to=3-1]
		\arrow["{\op_H}", from=1-3, to=3-3]
		\arrow["{\Mod(-)}"', from=3-1, to=3-3]
	\end{tikzcd}\cdpunct[8pt]{.}\]
\end{proposition}

\begin{proof}
	As before, we specify the virtual double functor $r_{\cV}\colon \op_{H}(\Mod(\cV)) \rightarrow \Mod(r(\cV))$.
	\begin{description}
		\item[Objects] 
		\[\begin{tikzcd}[ampersand replacement=\&]
			{M \coloneqq (M,\mu,e)} \&\& {M^r \coloneqq (M, \mu, e)}
			\arrow[maps to, from=1-1, to=1-3]
		\end{tikzcd}\]
		\item[V $1$-cell] 
		\[\begin{tikzcd}[ampersand replacement=\&]
			M \&\&\& {M^r} \\
			\\
			{M'} \&\&\& {(M')^r}
			\arrow[""{name=0, anchor=center, inner sep=0}, from=1-1, to=3-1]
			\arrow[""{name=1, anchor=center, inner sep=0}, from=1-4, to=3-4]
			\arrow[between={0.2}{0.8}, maps to, from=0, to=1]
		\end{tikzcd}\]
		\item[H $1$-cell] 
		\[\begin{tikzcd}[ampersand replacement=\&]
			{M_2} \&\&\& {M_1} \&\&\& {(M_2)^r} \&\&\& {(M_1)^r}
			\arrow[""{name=0, anchor=center, inner sep=0}, "{X \coloneqq(X,\lambda,\rho)}"'{inner sep=.8ex}, "\shortmid"{marking}, from=1-4, to=1-1]
			\arrow[""{name=1, anchor=center, inner sep=0}, "{X^{r} \coloneqq (X ,\rho,\lambda)}"{inner sep=.8ex}, "\shortmid"{marking}, from=1-7, to=1-10]
			\arrow[between={0.4}{0.6}, maps to, from=0, to=1]
		\end{tikzcd}\]
		\item[Sq multicell] 
		\[\begin{tikzcd}[ampersand replacement=\&,cramped]
			{M_2'} \& \cdots \& \cdots \& {M_1'} \&\& {(M_2')^r} \& \cdots \& \cdots \& {(M_1')^r} \\
			\\
			{M_2} \&\&\& {M_1} \&\& {M_2^r} \&\&\& {M_1^r}
			\arrow[from=1-1, to=3-1]
			\arrow["Z"'{inner sep=.8ex}, "\shortmid"{marking}, from=1-2, to=1-1]
			\arrow[""{name=0, anchor=center, inner sep=0}, "\shortmid"{marking}, from=1-3, to=1-2]
			\arrow["Y"'{inner sep=.8ex}, "\shortmid"{marking}, from=1-4, to=1-3]
			\arrow[""{name=1, anchor=center, inner sep=0}, from=1-4, to=3-4]
			\arrow["{Z^r}"{inner sep=.8ex}, "\shortmid"{marking}, from=1-6, to=1-7]
			\arrow[""{name=2, anchor=center, inner sep=0}, from=1-6, to=3-6]
			\arrow[""{name=3, anchor=center, inner sep=0}, "\shortmid"{marking}, from=1-7, to=1-8]
			\arrow["{Y^r}"{inner sep=.8ex}, "\shortmid"{marking}, from=1-8, to=1-9]
			\arrow[from=1-9, to=3-9]
			\arrow[""{name=4, anchor=center, inner sep=0}, "X"{inner sep=.8ex}, "\shortmid"{marking}, from=3-4, to=3-1]
			\arrow[""{name=5, anchor=center, inner sep=0}, "{X^r}"'{inner sep=.8ex}, "\shortmid"{marking}, from=3-6, to=3-9]
			\arrow["{f\colon Z \otimes \cdots \otimes Y \rightarrow X}"{description}, between={0.2}{0.8}, Rightarrow, nfold, from=0, to=4]
			\arrow[between={0.2}{0.8}, maps to, from=1, to=2]
			\arrow["{f \colon Y^r \otimes_{r} \cdots \otimes_{r} Z^r \rightarrow X^r}"{description}, between={0.2}{0.8}, Rightarrow, nfold, from=3, to=5]
		\end{tikzcd}\]
	\end{description}
\end{proof}

There are many other variants to consider, such as the pseudo case or the commutative case of each of these. However, the respective results are mostly analogous.

\begin{remark}
	In the case of a braided or symmetric monoidal category $\cV$, we can again compose with the canonical braiding functor $b_{\cV}$ to obtain the following inverse form of reversal component at $\cV$:
	\[\begin{tikzcd}[ampersand replacement=\&]
		{\op_{H}(\Mod(\cV))} \&\& {\Mod(r(\cV))} \&\& {\Mod(\cV)}
		\arrow["{r_{\cV}}", from=1-1, to=1-3]
		\arrow["{\Mod(b_{\cV})}", from=1-3, to=1-5]
	\end{tikzcd}\cdpunct[2pt]{.}\]
\end{remark}

Another well known construction on virtual double categories is the $2$-functor $\Mat \colon \MonCat_{l} \rightarrow \vDblCat$. Although this functor is not our central focus, it is worth discussing the dualities of this functor. We state the following lemma on reversal duality of $\Mat$.

\begin{lemma}
	The $2$-functor $\Mat$ extends to a $1$-morphism of involutions by the strict $2$-natural transformation
	\[\begin{tikzcd}[ampersand replacement=\&,cramped]
		{\MonCat_l} \&\& \vDblCat \\
		\\
		{\MonCat_l} \&\& \vDblCat
		\arrow["{\Mat(-)}", from=1-1, to=1-3]
		\arrow["r"', from=1-1, to=3-1]
		\arrow["r", between={0.2}{0.8}, Rightarrow, nfold, from=1-3, to=3-1]
		\arrow["{\op_H}", from=1-3, to=3-3]
		\arrow["{\Mat(-)}"', from=3-1, to=3-3]
	\end{tikzcd}\cdpunct[8pt]{.}\]
\end{lemma}
\begin{proof}
	Let $(-)^t$ here denote the transpose of a matrix. The components $r_{\cV} \colon \op_H(\Mat(\cV)) \rightarrow \Mat(r(\cV))$ are defined by
	\begin{description}
		\item[Objects] 
		\[\begin{tikzcd}[ampersand replacement=\&,cramped]
			X \&\& X
			\arrow[maps to, from=1-1, to=1-3]
		\end{tikzcd}\]
		\item[V $1$-cell] 
		\[\begin{tikzcd}[ampersand replacement=\&,cramped]
			X \&\&\& X \\
			\\
			Y \&\&\& Y
			\arrow[""{name=0, anchor=center, inner sep=0}, from=1-1, to=3-1]
			\arrow[""{name=1, anchor=center, inner sep=0}, from=1-4, to=3-4]
			\arrow[between={0.2}{0.8}, maps to, from=0, to=1]
		\end{tikzcd}\]
		\item[H $1$-cell] 
		\[\begin{tikzcd}[ampersand replacement=\&,cramped]
			Y \&\&\& X \&\&\& Y \&\&\& X
			\arrow[""{name=0, anchor=center, inner sep=0}, "M"'{inner sep=.8ex}, "\shortmid"{marking}, from=1-4, to=1-1]
			\arrow[""{name=1, anchor=center, inner sep=0}, "{M^t}"{inner sep=.8ex}, "\shortmid"{marking}, from=1-7, to=1-10]
			\arrow[between={0.4}{0.6}, maps to, from=0, to=1]
		\end{tikzcd}\]
		\item[Sq multicell] 
		\[\begin{tikzcd}[ampersand replacement=\&,cramped]
			{X_n} \& \cdots \& \cdots \& {X_0} \&\& {X_n} \& \cdots \& \cdots \& {X_0} \\
			\\
			{Y_1} \&\&\& {Y_0} \&\& {Y_1} \&\&\& {M_1^r}
			\arrow[from=1-1, to=3-1]
			\arrow["{M_n}"'{inner sep=.8ex}, "\shortmid"{marking}, from=1-2, to=1-1]
			\arrow[""{name=0, anchor=center, inner sep=0}, "\shortmid"{marking}, from=1-3, to=1-2]
			\arrow["{M_1}"'{inner sep=.8ex}, "\shortmid"{marking}, from=1-4, to=1-3]
			\arrow[""{name=1, anchor=center, inner sep=0}, from=1-4, to=3-4]
			\arrow["{M_n^t}"{inner sep=.8ex}, "\shortmid"{marking}, from=1-6, to=1-7]
			\arrow[""{name=2, anchor=center, inner sep=0}, from=1-6, to=3-6]
			\arrow[""{name=3, anchor=center, inner sep=0}, "\shortmid"{marking}, from=1-7, to=1-8]
			\arrow["{M_1^t}"{inner sep=.8ex}, "\shortmid"{marking}, from=1-8, to=1-9]
			\arrow[from=1-9, to=3-9]
			\arrow[""{name=4, anchor=center, inner sep=0}, "N"{inner sep=.8ex}, "\shortmid"{marking}, from=3-4, to=3-1]
			\arrow[""{name=5, anchor=center, inner sep=0}, "{N^t}"'{inner sep=.8ex}, "\shortmid"{marking}, from=3-6, to=3-9]
			\arrow[between={0.2}{0.8}, Rightarrow, nfold, from=0, to=4]
			\arrow[between={0.2}{0.8}, maps to, from=1, to=2]
			\arrow[between={0.2}{0.8}, Rightarrow, nfold, from=3, to=5]
		\end{tikzcd}\]
	\end{description}
\end{proof}

\begin{remark}
	More generally, \cite{fujii2025famenr} propose a generalized matrices functor $\Mat \colon \vDblCat \rightarrow \vDblCat$, which arises primarily from the discrete opfibration $\Set_* \rightarrow \Set$. In particular, they define it as a composite functor 
	\[\begin{tikzcd}[ampersand replacement=\&,cramped]
		\vDblCat \&\&\& {\vDblCat/(\Set_*)_{\mathrm{hc}}} \\
		\\
		\vDblCat \&\&\& {\vDblCat/(\Set)_{\mathrm{hc}}}
		\arrow["{(\Set_*)_{\mathrm{hc}} \times (-)}", from=1-1, to=1-4]
		\arrow["\Mat"', from=1-1, to=3-1]
		\arrow["{\prod_{P_{\mathrm{hc}}}}", from=1-4, to=3-4]
		\arrow["U", from=3-4, to=3-1]
	\end{tikzcd}\cdpunct[8pt]{.}\]
	One could show that each of the individual $2$-functor appearing in this decomposition is compatible with the $\op_{H}$ duality. In particular, this would generalize the previous lemma to $\vDblCat$. 
\end{remark}

Lastly, in contrast to the horizontal duality, we note that the $2$-functor $\Mat$ doesn't respect the vertical duality of virtual double categories in the expected way. Although one could define a suitable $2$-functor $\mathbf{CoMat}\colon \vcoDblCat \rightarrow \vcoDblCat$, applying $\CoMod$ to this construction would not yield the desired variance for morphisms between cocategories. This is more readily observed by considering the category $\cV^{\op}\Cat$ of $\op_1(\cV)$-enriched categories. In particular, the action of functors on the hom-objects and on the set of objects are of opposite variance. If one is content with this disparity, we believe that one could dualize the $\Mat$ construction.

\subsection{Higher Dualities} \label{sub:highermonduals}

We briefly describe the analogous transport phenomenon in the case of $\nCat{1}$-enriched symmetric monoidal category $\nCat{1}$. 

\begin{definition}
	Let $(\cV,\otimes, I)$ be any monoidal category. We define a left $\cV$-module category $\cC$ as a left pseudomodule object in $\nCat{1}$. The definitions for right module or bimodule categories are analogous.
\end{definition}

To describe the $\op_1$ duality of monoidal categories, we need to have the variants $\nMon{1}_c$ and $\nMon{1}_l$ corresponding the construction $\nMon{1}$. This is similar to how one dimension lower, we required separate $\Mon$ and $\CoMon$ constructions. In particular, $\nMon{1}_l$ is simply the construction we have already provided involving lax $1$-morphisms. Whereas, $\nMon{1}_c$ is the colax $1$-morphism variant. 
 
\begin{lemma}
	The $\op_1$ and $\op_1^t$ duality inverses of monoidal categories are equivalently expressed as the two maps:
	\[\begin{tikzcd}[ampersand replacement=\&,column sep=scriptsize]
		{\op_1^{t} \colon} \& {\op_2(\nMon{1}_{l}(\nCat{1}))} \&\& {\nMon{1}_{c}(\op_2(\nCat{1}))} \&\&\& {\nMon{1}_{c}(\nCat{1})} \\
		{\op_1 \colon} \& {\op_{2}(\nMon{1}_{c}(\nCat{1}))} \&\& {\nMon{1}_{l}(\op_2(\nCat{1}))} \&\&\& {\nMon{1}_{l}(\nCat{1})}
		\arrow[from=1-2, to=1-4]
		\arrow["{\nMon{1}_{c}(\op_1)}", from=1-4, to=1-7]
		\arrow[from=2-2, to=2-4]
		\arrow["{\nMon{1}_{l}(\op_1)}", from=2-4, to=2-7]
	\end{tikzcd}\]
\end{lemma}
 
\begin{proof}
	The identifications $\nMon{1}_{l}(\op_2(\nCat{1})) \cong \op_2(\MonCat_{c})$ and $\nMon{1}_{c}(\op_2(\nCat{1})) \cong \op_2(\MonCat_{l})$ follows directly by the definitions. One should view these identifications as merely comprehension identifications. Whereas, $\op_1 : \op_2(\nCat{1}) \rightarrow \nCat{1}$ is a strict monoidal $2$-functor. Thus, we have transport morphisms $\nMon{1}_c(\op_1)$ and $\nMon{1}_l(\op_1)$ given by applying the $\op_1$ functor directly on data of the objects and the morphisms. This is exactly the description we provided at the beginning of this section. Now, we are stating the same description in a more formalized way.
\end{proof}

In other words, we have fully expressed the $\op_1$ and $\op_1^t$ duality of monoidal categories in terms of a hypothetical higher duality $\op_2\colon \op_3(\mathbf{MonBicat_{x}}) \rightarrow \mathbf{MonBicat_{y}}$, where $x$ and $y$ denote some suitable variants. More precisely, assuming we had all the relevant definitions, the duality can be described by the following diagram
\[\begin{tikzcd}[ampersand replacement=\&]
	{*} \&\& {\op_3(\mathbf{MonBicat}_x)} \&\& {\op_3(\BiCat_x)} \\
	\\
	{*} \&\& {\mathbf{MonBicat}_y} \&\& {\BiCat_y}
	\arrow["\Cat", from=1-1, to=1-3]
	\arrow[from=1-1, to=3-1]
	\arrow["{\nMon{1}_l}", from=1-3, to=1-5]
	\arrow["{\op_1}"{description}, between={0.2}{0.8}, Rightarrow, nfold, from=1-3, to=3-1]
	\arrow["{\op_2}"{description}, from=1-3, to=3-3]
	\arrow["{\op_2}"{description}, between={0.2}{0.8}, Rightarrow, nfold, from=1-5, to=3-3]
	\arrow["{\op_2}", from=1-5, to=3-5]
	\arrow["\Cat"', from=3-1, to=3-3]
	\arrow["{\nMon{1}_l}"', from=3-3, to=3-5]
\end{tikzcd}\cdpunct[8pt]{.}\]

It is now also straightforward to guess how pseudomodules in $\nCat{1}$ ought to transport under the dualities.

\begin{lemma}
	Let $\cV_1$ and $\cV_2$ be any two monoidal categories and let $\mathbf{Bimod}(\cV_1,\cV_2)$ be $2$-category of modules between $\cV_1$ and $\cV_2$. There exists the following involutive isomorphism induced by the action of $\op_1 \colon \op_2(\nCat{1}) \rightarrow \nCat{1}$:
	\[\begin{tikzcd}[ampersand replacement=\&,cramped]
		{\op_1 \colon} \& {\op_2(\mathbf{Bimod}(\cV_1,\cV_2))} \&\& {\mathbf{Bimod}(\op_1(\cV_1),\op_1(\cV_2))}
		\arrow[from=1-2, to=1-4]
	\end{tikzcd}\cdpunct[0pt]{.}\]
\end{lemma}

Simply note that all the defining axioms are comprised of natural isomorphisms which can be inverted. Finally, the duality for reversal is expressible similarly.

\begin{lemma} \leavevmode
	\begin{enumerate}
		\item The reversal duality $r : \MonCat_{x} \rightarrow \MonCat_{x}$ is given by the map
		\[\begin{tikzcd}[ampersand replacement=\&,cramped]
			{r \colon } \& {\nMon{1}(\nCat{1})} \&\& {\nMon1(r(\nCat{1}))} \&\&\& {\nMon{1}(\nCat{1})}
			\arrow[from=1-2, to=1-4]
			\arrow["{\nMon{1}((1,s_{\nCat{1}},1))}", from=1-4, to=1-7]
		\end{tikzcd}\cdpunct[0pt]{,}\]
		where $(1,s_{\nCat{1}},1)$ is the canonical identity symmetric monoidal strict $2$-functor whose compositor constraint is the symmetry of $\nCat{1}$. 
		\item There exists the following involutive isomorphism on the bimodule $2$-categories induced by the action of $r$:
		\[\begin{tikzcd}[ampersand replacement=\&,cramped]
			{r\colon} \& {\mathbf{Bimod}(\cV_1,\cV_2)} \&\& {\mathbf{Bimod}(r(\cV_2),r(\cV_1))}
			\arrow[from=1-2, to=1-4]
		\end{tikzcd}\cdpunct[0pt]{.}\]
	\end{enumerate}
\end{lemma}

At this point, we could also obtain similar results relating monoidal and comonoidal categories using the $\op_1$ involution of bicategories, that is, $\op_1(\mathbf{1CoMon}_l(\nCat{1})) \cong \nMon{1}_l(\op_1(\nCat{1}))$. However, these results are not essential for us, so we do not investigate them here.

	\section{Sliced modules and enriched categories} \label{sec:slicecatmod}
	The goal of this section is to describe how the microcosm principle allows one to define general enriched slice categories, especially under the consideration of model structures. Some special cases of the following results have already appeared in the literature for slice-over categories. We provide a general characterization encompassing both slice-over and slice-under categories.

\subsection{Slice monoidal categories and module} \label{sub:slicecaten}

The definition of a module category is again within the scope of the microcosm principle. Consequentially, we start with the definition of a module object inside a module category.

\begin{definition}
	Let $(\cC,\triangleleft_{\cC},a_{\cC},\lambda_{\cC})$ be a left $\cV$-module category over some monoidal category $\cV$. Let $(M,\mu,e) \in \cV$ be a monoid object. We say that an object $Z \in \cC$ is a left module object over $M$ if it has following data:
	\begin{enumerate}
		\item a left action morphism $\triangleleft_{Z} \colon M \triangleleft_{\cC} Z \rightarrow Z$,
		\item the associativity and unitality axioms hold
		\[\begin{tikzcd}[ampersand replacement=\&]
			{M \triangleleft_{\cC} (M \triangleleft_{\cC} Z)} \&\& {M \triangleleft_{\cC} Z} \& {I \triangleleft_{\cC} Z} \&\& {M \triangleleft_{\cC} Z} \\
			\\
			{(M \otimes M) \triangleleft_{\cC} Z} \&\&\&\&\& Z \\
			\\
			{M \triangleleft_{\cC} Z} \&\& Z
			\arrow["{1 \triangleleft_{\cC} (\triangleleft_{Z})}", from=1-1, to=1-3]
			\arrow["{\triangleleft_{Z}}", from=1-3, to=5-3]
			\arrow["{e \triangleleft_{\cC} Z}", from=1-4, to=1-6]
			\arrow["{\lambda_{\cC}}"', from=1-4, to=3-6]
			\arrow["{\triangleleft_{Z}}", from=1-6, to=3-6]
			\arrow["{a_{\cC}}", from=3-1, to=1-1]
			\arrow["{\mu \triangleleft 1}"', from=3-1, to=5-1]
			\arrow["{\triangleleft_{Z}}"', from=5-1, to=5-3]
		\end{tikzcd}\cdpunct[8pt]{.}\]
	\end{enumerate}
	
\end{definition}

\begin{definition}
	Let $(\cC,\triangleright_{\cC},a_{\cC},\rho_{\cC})$ be a right $\cV$-module category over some monoidal category $\cV$. Let $(M,\mu,e) \in \cV$ be a monoid object. We say that an object $Z \in \cC$ is a right module object over $M$ if it has following data:
	\begin{enumerate}
		\item a right action morphism $\triangleright_{Z} \colon Z \triangleright_{\cC} M \rightarrow Z$,
		\item the associativity and unitality axioms hold
		\[\begin{tikzcd}[ampersand replacement=\&,cramped]
			{(Z \triangleright_{\cC} M) \triangleright_{\cC} M} \&\& {Z \triangleleft_{\cC} M} \& {Z \triangleright_{\cC} I} \&\& {Z \triangleright_{\cC} M} \\
			\\
			{Z \triangleright_{\cC}(M \otimes M) } \&\&\&\&\& Z \\
			\\
			{Z \triangleright_{\cC} M} \&\& Z
			\arrow["{(\triangleright_{Z}) \triangleright_{\cC} 1}", from=1-1, to=1-3]
			\arrow["{a_{\cC}}"', from=1-1, to=3-1]
			\arrow["{\triangleright_{Z}}", from=1-3, to=5-3]
			\arrow["{Z \triangleright_{\cC} e}", from=1-4, to=1-6]
			\arrow["{\rho_{\cC}}"', from=1-4, to=3-6]
			\arrow["{\triangleright_{Z}}", from=1-6, to=3-6]
			\arrow["{1 \triangleright_{\cC} \mu}"', from=3-1, to=5-1]
			\arrow["{\triangleright_{Z}}"', from=5-1, to=5-3]
		\end{tikzcd}\cdpunct[8pt]{.}\]
	\end{enumerate}
\end{definition}

A similar definition can be stated for bimodule object in a bimodule category. One only needs to add a balancing axiom for the bimodule object relative to the balancing $2$-cell of bimodule category. 

As expected, these definitions also provides the slice constructions. Let $(\cV,\otimes,I,a_{\cV},\lambda_{\cV},\rho_{\cV})$ be a monoidal category in what follows.

\begin{lemma} \label{lem:slicingmonandmod}
	Let $(\cC,\triangleleft_{\cC},a_{\cC},\lambda_{\cC})$ be a left $\cV$-module. Let $(M, \mu, e)$ be a monoid object in $\cV$, and let $(Z, \triangleleft_{Z})$ be a left $M$-module. 
	\begin{enumerate}
		\item The slice category $\cV/M$ admits a monoidal structure denoted by $(\cV/M, \otimes_{/M}, I_{/M})$. If $\cV$ is braided or symmetric and $M$ is commutative then $\cV/M$ is also braided or symmetric. 
		\item There is a canonical left $\cV/M$-module structure on $\cC/Z$ denoted by $(\cC/Z,\triangleleft_{/})$.
	\end{enumerate}
\end{lemma}

\begin{proof}
	For the first statement, the sliced tensoring operation is defined by 
	\[f \otimes_{/M} g \coloneqq A \otimes B \xrightarrow{f \otimes g} M \otimes M \xrightarrow{\mu} M\] 
	and the unit by \[I_{/M} \coloneqq I \xrightarrow{\eta} M.\] 
	The naturality of the associator and the unitors of $\cV$ together with the monoid axioms of $M$ allows us to extend these given associator and unitors to those of $\cV/M$ directly.
	
	For the second statement, the sliced left action operation is defined exactly the same. Let $\alpha \colon X \rightarrow M \in \cV/M$ and $f \colon A \rightarrow Z \in \cC/Z$. The sliced left action is defined by \[\alpha \triangleleft_{/} f \colon X \triangleleft A \xrightarrow{\alpha \triangleleft f} M \triangleleft Z \xrightarrow{\triangleleft_{Z}} Z.\] 
	Thereafter, the extension of coherence $2$-cells follow again by the naturality and module axioms.	
\end{proof}

This lemma can be dualized in two ways: one by using the $\op_1$ duality and the other by $r$ duality. The reversal duality, however, provides nothing new.

For the $\op_1$ duality, we consider a monoid object $N$ in the monoidal category $\cV^{\op}$ and a module object $Z$ in $\cC^{\op}$. Applying the previous lemma in this situation and transporting structure back using $\op_1$ duality, we obtain the following.

\begin{lemma}
	Let $(\cC,\triangleleft_{\cC},a_{\cC},\lambda_{\cC})$ be a left $\cV$-module. Let $(N, \mu, \eta)$ be a comonoid object and $(Z,\lambda)$ be a left $N$-comodule. 
	\begin{enumerate}
		\item The coslice category $N/\cV$ admits a monoidal structure denoted by $(N/\cV, \otimes_{N/}, I_{N/})$. If $\cV$ is braided or symmetric and $N$ is commutative then $N/\cV$ is also braided or symmetric. 
		\item There is a canonical left $N/\cV$-module structure on $Z/\cV$ denoted by $(Z/\cC,\triangleleft_{Z/\cC})$.
	\end{enumerate}
\end{lemma}

We now investigate the presence of closures on left $\cV$-module structure. That is, we investigate enrichment through tensoring or powering. Consequentially, we shall now refer to the action functor $\triangleleft \colon \cV \times \cC \rightarrow \cC$ as \emph{tensoring}.

\begin{definition}
	A left $\cV$-module category $(\cC,\triangleleft)$ is called left closed if for every $c \in \cC$, the functor $- \triangleleft c \colon \cV \rightarrow \cV$ has a right adjoint $\ucC(c,-)\colon \cC \rightarrow \cV$. Thus, we have a functor $\ucC(-,-) \colon \cC^{\op} \times \cC \rightarrow \cV$, which is called the \emph{internal hom} of $\cC$.
\end{definition}

\begin{definition}
	A left $\cV$-module category $(\cC,\triangleleft)$ is called right closed if for every $x \in \cV$, the functor $x \triangleleft - \colon \cC \rightarrow \cC$ has a right adjoint $- \pitchfork x \colon \cC \rightarrow \cC$. Thus, we have a functor $- \pitchfork - \colon \cC \times \cV^{\op} \rightarrow \cC$, which is called \emph{powering}. This is usually also denoted by $b^x$ or $\{x,b\}$.
\end{definition}

\begin{definition}
	A left $\cV$-module category $(\cC,\triangleleft)$ is closed if it is both left closed and right closed.
\end{definition}

\begin{lemma}
	A closed left $\cV$-module category $(\cC,\triangleleft)$ is equivalently specified by the data of a two variable left adjunction 
	\[\begin{tikzcd}[ampersand replacement=\&,cramped]
		{(- \triangleleft - )^{\op}\colon} \& {\cV^{\op} \times \cC^{\op}} \& {\cC^{\op}} \\
		{\ucC(-,-) \colon} \& {\cC^{\op} \times \cC} \& \cV \\
		{- \pitchfork - \colon} \& {\cC \times \cV^{\op}} \& \cC
		\arrow[from=1-2, to=1-3]
		\arrow[from=2-2, to=2-3]
		\arrow[from=3-2, to=3-3]
	\end{tikzcd}\cdpunct[2pt]{.}\]
\end{lemma}
\begin{proof}
	See \cite[Def.~2.1 and Sec.~2.2]{chengmates2014} for the definition and the respective proof. In simpler terms, this is equivalent to the natural isomorphisms $\cC(x \triangleleft a,b) \cong \cV(x,\ucC(a,b)) \cong \cC(a,b \pitchfork x)$.
\end{proof}

\begin{example}
	If the canonical left $\cV$-module $(\cV, \otimes)$ is left closed then it is a left closed monoidal category in the usual sense. The same is true for the right closure of this left module structure.
\end{example}

\begin{lemma}
	There is a canonical enrichment of $\cC$ over $\cV$ using the internal homs $\ucC(-,-)$, which we call $\ucC$. Furthermore, the underlying category of $\ucC$ is isomorphic to $\cC$ itself.
\end{lemma}

\begin{proof}
	One can mimic the proof provided in \cite[Vol.~III, Sec.~3.1]{yau2024bimodcats} by defining $\mathrm{eval}^{c}_d \colon \ucC(c,d) \triangleleft c \rightarrow d$ through the counit and $\mathrm{coeval}^{d}_x \colon x \rightarrow \ucC(d,x \triangleleft d)$ through the unit. Thus, compositions are given by the adjoint to the composites
	\[\begin{tikzcd}[ampersand replacement=\&,cramped]
		{(\cC(d,e) \otimes \cC(c,d)) \triangleleft c} \&\&\& e \\
		{\cC(d,e) \triangleleft(\cC(c,d) \triangleleft c)} \&\&\& {\cC(d,e) \triangleleft d}
		\arrow[from=1-1, to=1-4]
		\arrow["a", from=1-1, to=2-1]
		\arrow["{1 \triangleleft \mathrm{eval}_d^{c}}"', from=2-1, to=2-4]
		\arrow["{\mathrm{eval}_e^{d}}"', from=2-4, to=1-4]
	\end{tikzcd}\cdpunct[8pt]{,}\]
	and unitors are given by the adjoint to the morphisms $I \triangleleft x \xrightarrow{\lambda_{x}} x$. Thereafter, the coherence follows from the $2$-cells of the left $\cV$-module structure. 
\end{proof}

Let $i\colon I \rightarrow \ucC(c,c)$ be the unitors and let $m\colon \ucC(d,e) \otimes \ucC(c,d) \rightarrow \ucC(c,e)$ be the compositors. In particular, we obtain the following two commuting diagrams
\[\begin{tikzcd}[ampersand replacement=\&,cramped]
	{(\cC(d,e) \otimes \cC(c,d)) \triangleleft c} \&\& {\ucC(c,e) \triangleleft c} \& {I \triangleleft c} \&\& {\ucC(c,c)  \triangleleft c} \\
	{\cC(d,e) \triangleleft(\cC(c,d) \triangleleft c)} \\
	{\cC(d,e) \triangleleft d} \&\& e \&\&\& c
	\arrow["{m \triangleleft 1}", from=1-1, to=1-3]
	\arrow["a", from=1-1, to=2-1]
	\arrow["{\mathrm{eval}_e^{c}}", from=1-3, to=3-3]
	\arrow["{i \triangleleft 1}", from=1-4, to=1-6]
	\arrow["\lambda"', from=1-4, to=3-6]
	\arrow["{\mathrm{eval}^c_c}", from=1-6, to=3-6]
	\arrow["{1 \triangleleft \mathrm{eval}_d^{c}}"', from=2-1, to=3-1]
	\arrow["{\mathrm{eval}_e^{d}}"', from=3-1, to=3-3]
\end{tikzcd}\cdpunct[8pt]{.}\]

The converse of this lemma is not true in general. That is, we can't obtain left module structures using left adjoints to internal homs. However, when $\cV$ is a closed symmetric monoidal category, the converse does hold, and the closed module structure lifts from $\nCat{1}$ to $\nCat{\cV}$. In particular, this defines a enriched, tensored, and cotensored category \cite[Prop.~10.1.14]{reihlcht}.

\begin{lemma} \label{lem:poweringrightmod}
	The powering functor $- \pitchfork - \colon \cC \times \cV^{\op} \rightarrow \cC$ defines a canonical right $\cV^{\op}$-module structure. 
\end{lemma}

\begin{proof}
	A left $\cV$-module structure $(\cC,\triangleleft)$ is equivalent to providing a strong monoidal functor $- \triangleleft  \colon \cV \rightarrow \nCat{1}(\cC,\cC)$. The requirement of $x \triangleleft -$ admitting right adjoint is then equivalent to $x \triangleleft -$ admitting a right dual in the endofunctor monoidal category $\nCat{1}(\cC,\cC)$. The right dual construction then defines a \enquote{partial} normal strong monoidal functor $r((\nCat{1}(\cC,\cC))^{\op}) \rightarrow \nCat{1}(\cC,\cC)$ out of the right dualizable endofunctors (see \cite[Sec.~2.10]{etingof2015tenscat}). In particular, we obtain a canonical left $r(\cV^{\op})$-module structure on $\cC$ by postcomposing with this partial right dual functor. Finally, this is equivalent to a canonical right $\cV^{\op}$-module structure on $\cC$ or a canonical right $\cV$-module structure on $\cC^{\op}$. 
\end{proof}

In concrete terms, the right $\cV^{\op}$-module structure on $\cC$ has a right unitor given by the diagram
\[\begin{tikzcd}[ampersand replacement=\&,cramped]
	{c \pitchfork I} \&\& {I \triangleleft ( c \pitchfork I)} \\
	\\
	\&\& c
	\arrow["{\lambda^{-1}_{c \pitchfork I}}", from=1-1, to=1-3]
	\arrow[from=1-1, to=3-3]
	\arrow["{\epsilon^{I}_{c}}", from=1-3, to=3-3]
\end{tikzcd}\cdpunct[2pt]{,}\]
where the inverse of this right unitor is given by a morphism $c \rightarrow (c \pitchfork I)$ that is adjoint to the left unitor $I \triangleleft c \rightarrow c$. In particular, this inverse right unitor fits in the diagram
\[\begin{tikzcd}[ampersand replacement=\&,cramped]
	{I \triangleleft c} \&\& {I \triangleleft (c \pitchfork I)} \\
	\\
	\&\& c
	\arrow[from=1-1, to=1-3]
	\arrow["\lambda"', from=1-1, to=3-3]
	\arrow["{\epsilon^I_{c}}", from=1-3, to=3-3]
\end{tikzcd}\cdpunct[2pt]{.}\]

Similarly, the associator of the right module structure on $\cC$ is given by a map $(c \pitchfork x) \pitchfork y  \rightarrow c \pitchfork (x \otimes y) \in \cC$ which fits into the diagram
\[\begin{tikzcd}[ampersand replacement=\&,cramped]
	{(x \otimes y) \triangleleft ((c \pitchfork x) \pitchfork y)} \&\& {(x \otimes y) \triangleleft (c \pitchfork (x \otimes y))} \\
	{x \triangleleft ( y \triangleleft ((c \pitchfork x)\pitchfork y))} \\
	\\
	{x \triangleleft (c \pitchfork x)} \&\& c
	\arrow[from=1-1, to=1-3]
	\arrow["a", from=1-1, to=2-1]
	\arrow["{\epsilon^{x \otimes y}_{c}}", from=1-3, to=4-3]
	\arrow["{1 \triangleleft \epsilon^y_{c \pitchfork x}}", from=2-1, to=4-1]
	\arrow["{\epsilon^x_{c}}"', from=4-1, to=4-3]
\end{tikzcd}\cdpunct[8pt]{.}\]

This lemma can be readily dualized to provide a left $\cV$-module structure by tensoring given a right $\cV^{\op}$-module structure by powering. Indeed, a right $\cV^{\op}$-module $\cC$ is equivalent to specifying a left $r(\cV)$-module $\cC^{\op}$. Now, we just apply the previous lemma to this equivalent specification to obtain the desired result. Expectedly, it is also true that an enrichment can be defined solely in terms of the powering without requiring the existence of a tensoring operation.

We can now begin discussing the slicing of closed module structures.

\begin{lemma} \label{lem:sliceovermod}
	Let $\cC$ be a left $\cV$-module, let $(M,\mu,e) \in \cV$ be a monoid object, and let $(Z,\triangleleft_Z) \in \cC$ be a left module object over the monoid $M$. And assume that $\cV$ has pullbacks.
	\begin{enumerate}
		\item If the category $\cV$ is left closed, then the category $\cV/M$ is also left closed using the sliced monoidal product. The same is true for right closure.
		\item If $\cC$ is a left closed left $\cV$-module then $\cC/Z$ is a left closed left-$\cV/M$ module. The same is true for right closure as long as $\cC$ has pullbacks.
	\end{enumerate}
\end{lemma}

\begin{proof} \phantom{}
	\begin{enumerate}
		\item We define the left closed internal hom $\ucV/M(\alpha\colon X \rightarrow M ,\beta\colon Y \rightarrow M)$ and right closed powering $\beta \pitchfork_{/} \alpha$ using the pullback diagrams
		\[\begin{tikzcd}[ampersand replacement=\&,cramped,column sep=scriptsize]
			{\ucV/M(\alpha,\beta)} \&\& {\ucV(X,Y)} \& {\beta \pitchfork_{/} \alpha} \&\&\& {Y \pitchfork X} \\
			\\
			M \& {\ucV(M,M)} \& {\ucV(X,M)} \& M \& {M \pitchfork M} \&\& {M \pitchfork X}
			\arrow[from=1-1, to=1-3]
			\arrow[from=1-1, to=3-1]
			\arrow["\ulcorner"{anchor=center, pos=0.125}, draw=none, from=1-1, to=3-2]
			\arrow["{\beta_*}", from=1-3, to=3-3]
			\arrow[from=1-4, to=1-7]
			\arrow[from=1-4, to=3-4]
			\arrow["\ulcorner"{anchor=center, pos=0.125}, draw=none, from=1-4, to=3-7]
			\arrow["{\beta \pitchfork X}", from=1-7, to=3-7]
			\arrow["\mu"', from=3-1, to=3-2]
			\arrow["{\alpha^*}"', from=3-2, to=3-3]
			\arrow["\mu"', from=3-4, to=3-5]
			\arrow["{M \pitchfork \alpha}"', from=3-5, to=3-7]
		\end{tikzcd}\]
		Thereafter, maps $\gamma \colon Z \rightarrow M$ into $\ucV/M(\alpha,\beta) \rightarrow M$ be can identified with maps $\gamma \otimes_{/M} \alpha \rightarrow \beta$ by using the left closure adjunction of $\cV$. The right closure is same.
		
		\item Let $\alpha: X \rightarrow M$, and $f,g: A,B \rightarrow Z$. The definition of internal hom and powering is similarly given pullback diagrams
		\[\begin{tikzcd}[ampersand replacement=\&,cramped]
			{\ucC/Z(f,g)} \&\& {\ucC(A,B)} \& {g \pitchfork_{/} \alpha} \&\&\& {B \pitchfork X} \\
			\\
			M \& {\ucC(Z,Z)} \& {\ucC(A,Z)} \& Z \& {Z \pitchfork M} \&\& {Z \pitchfork X}
			\arrow[from=1-1, to=1-3]
			\arrow[from=1-1, to=3-1]
			\arrow["\ulcorner"{anchor=center, pos=0.125}, draw=none, from=1-1, to=3-2]
			\arrow["{g_*}", from=1-3, to=3-3]
			\arrow[from=1-4, to=1-7]
			\arrow[from=1-4, to=3-4]
			\arrow["\ulcorner"{anchor=center, pos=0.125}, draw=none, from=1-4, to=3-7]
			\arrow["{g \pitchfork X}", from=1-7, to=3-7]
			\arrow["{\triangleleft_{Z}}"', from=3-1, to=3-2]
			\arrow["{f^*}"', from=3-2, to=3-3]
			\arrow["{\triangleleft_{Z}}"', from=3-4, to=3-5]
			\arrow["{Z \pitchfork \alpha}"', from=3-5, to=3-7]
		\end{tikzcd}\]
		The verification of adjunction is similar by using the closure adjunctions of left $\cV$-module $\cC$.
	\end{enumerate}
	
\end{proof}

A dual version of the above lemma is not possible for a comonoid $N$ as there is no transpose for the comultiplication $N \rightarrow N \otimes N$. However, we can still pursue a form of duality, if we instead start with a right $\cV^{\op}$-module $(\cC,\pitchfork)$. This right module structure is meant to encode powering. Thus, we merely ask whether this powering admits tensors using left adjoints or internal homs using mutual right adjoints. More precisely, we use the $\op_1$ duality to encode this as a right $\cV$-module $\cC^{\op}$. Then, equivalently, we ask if this right module is left closed for admitting tensors or right closed for admitting internal homs.

\begin{lemma} \label{lem:sliceundermod}
	Let $\cC^{\op}$ be a right $\cV$-module using $(- \pitchfork -)^{\op}$, let $(M,\mu,e) \in \cV$ be a monoid object, and let $(Z,\triangleright_Z) \in \cC^{\op}$ be a right module object over the monoid $M$ (equivalently, a right comodule object $(Z,\triangleright_Z) \in \cC$ over comonoid $M \in \cV^{\op}$). And assume that $\cV$ has pullbacks.
	\begin{enumerate}
		\item There is a canonical right $\cV/M$-module structure on $\cC^{\op}/Z$, equivalently on $(Z/\cC)^{\op}$. Under the $\op_1$ duality, this is a right $(\cV/M)^{\op}$-module structure on $Z/\cC$ and defines sliced powering.
		\item If $\cC^{\op}$ is a right closed right $\cV$-module then $(Z/\cC)^{\op}$ is a right closed right $\cV/M$-module. This defines a sliced internal hom $Z/\ucC(-,-) \colon (Z/\cC)^{\op} \times (Z/\cC) \rightarrow \cV/M$.
		\item If $\cC^{\op}$ is a left closed right $\cV$ module then $(Z/\cC)^{\op}$ is a left closed right $\cV/M$-module as long as $\cC$ has pushouts (equivalently $\cC^{\op}$ has pullbacks). This defines a sliced left action functor $(- \triangleleft_{/} -) \colon \cV/M \times Z/\cC \rightarrow Z/\cC$.
	\end{enumerate}
\end{lemma}

\begin{proof}
	This whole lemma is merely a direct application of the \cref{lem:sliceovermod} under appropriate dualities. Thus, we only define the functors for clarity. Let $\alpha : X \rightarrow M \in \cV/M$ and $f,g :Z \rightarrow A,B \in Z/\cC$. 
	
	The sliced powering $g^{\alpha}$ is defined by the composite
	\[\begin{tikzcd}[ampersand replacement=\&,cramped]
		Z \&\&\&\& {B \pitchfork X} \\
		\&\& {Z \pitchfork M} \\
		\&\&\&\& {Z \pitchfork X}
		\arrow["{g \pitchfork_{/} \alpha}", from=1-1, to=1-5]
		\arrow["{\triangleright_{Z}}"', from=1-1, to=2-3]
		\arrow["{Z \pitchfork \alpha}"', from=2-3, to=3-5]
		\arrow["{g \pitchfork X}"', from=3-5, to=1-5]
	\end{tikzcd}\cdpunct[2pt]{.}\]

	The internal hom $Z/\ucC(f,g)$ is defined by the pullback
	\[\begin{tikzcd}[ampersand replacement=\&,cramped]
		{Z/\ucC(f,g)} \&\& {\ucC(A,B)} \\
		\\
		M \& {\ucC(Z,Z)} \& {\ucC(Z,B)}
		\arrow[from=1-1, to=1-3]
		\arrow[from=1-1, to=3-1]
		\arrow["\ulcorner"{anchor=center, pos=0.125}, draw=none, from=1-1, to=3-3]
		\arrow["{f^*}", from=1-3, to=3-3]
		\arrow["{\triangleright_{Z}}"', from=3-1, to=3-2]
		\arrow["{g_*}"', from=3-2, to=3-3]
	\end{tikzcd}\cdpunct[4pt]{.}\]

	The sliced tensoring $\alpha \triangleleft_{Z} f$ is defined by the pushout
	\[\begin{tikzcd}[ampersand replacement=\&,cramped]
		{X \triangleleft Z} \&\& {X \triangleleft A} \\
		{M \triangleleft Z} \\
		Z \&\& {\alpha \triangleleft_{/Z} f}
		\arrow["{1 \triangleleft f}", from=1-1, to=1-3]
		\arrow["{\alpha \triangleleft 1}"', from=1-1, to=2-1]
		\arrow[from=1-3, to=3-3]
		\arrow["{\triangleright_{Z}}"', from=2-1, to=3-1]
		\arrow[from=3-1, to=3-3]
		\arrow["\ulcorner"{anchor=center, pos=0.125, rotate=180}, draw=none, from=3-3, to=1-1]
	\end{tikzcd}\cdpunct[2pt]{.}\]
\end{proof}

In the previous two lemma, we obtained the sliced enrichment structures by slicing the tensoring and by slicing the powering, respectively. Therefore, a natural question that remains is if there exists a slicing of the enrichment itself. We motivate this through the following lemma.

\begin{lemma}
	Let $\cC$ be a closed left $\cV$-module and $(M,\mu,e)$ a monoid object in $\cV$. The following data are equivalent.
	\begin{enumerate}
		\item A left module object $(Z, M \triangleleft Z \rightarrow Z)$ over the monoid object $M \in \cV$.
		
		\item A morphism of monoid objects $M \rightarrow \ucC(Z,Z)$.
		
		\item A right comodule object $(Z, Z \rightarrow Z \pitchfork M)$ over the comonoid object $M \in \cV^{\op}$.
	\end{enumerate}
\end{lemma}

\begin{proof}
	We define the monoid morphism $M \rightarrow \ucC(Z,Z)$ as adjoint to the morphism $M \triangleleft Z \rightarrow Z$. In particular, we have a commuting diagram 
	\[\begin{tikzcd}[ampersand replacement=\&,cramped]
		{M \triangleleft Z} \&\& {\ucC(Z,Z) \triangleleft Z} \\
		\\
		\&\& Z
		\arrow[from=1-1, to=1-3]
		\arrow[from=1-1, to=3-3]
		\arrow["{\mathrm{eval}^Z_Z}", from=1-3, to=3-3]
	\end{tikzcd}\cdpunct[2pt]{.}\]
	Similarly, we define the comodule map as an adjoint to the morphism $M \triangleleft Z \rightarrow Z$, which fits into the following commuting diagram
	\[\begin{tikzcd}[ampersand replacement=\&,cramped]
		{M \triangleleft Z} \&\& {M \triangleleft (Z \pitchfork M)} \\
		\\
		\&\& Z
		\arrow[from=1-1, to=1-3]
		\arrow[from=1-1, to=3-3]
		\arrow["{\epsilon^M_{Z}}", from=1-3, to=3-3]
	\end{tikzcd}\cdpunct[2pt]{.}\]
	
	The rest is simply an exercise in taking transposes of the defining axioms of the left module structure of $(Z, M \triangleleft Z \rightarrow Z)$ and using the above two diagrams to simplify.
\end{proof}

In particular, this suggests that we should be able to perform slice construction on arbitrary enriched categories $\cC$ over some base $\cV$, as long as we have a morphism of monoids $M \rightarrow \ucC(Z,Z)$. This is indeed true.

\begin{proposition}
	Consider a monoidal category $\cV$ with pullbacks. Let $\ucC$ be any enriched category and assume that there exists a morphism of monoid objects $r \colon M \rightarrow \ucC(Z,Z)$ for some $Z \in \ucC$. There exists canonical slice enriched categories $Z/\ucC$ and $\ucC/Z$ enriched over the slice monoidal category $\cV/M$. Furthermore, the underlying category of both $\ucC/Z$ and $Z/\ucC$ over $\cV/M$ (and not $\cV$) is isomorphic to the corresponding slices categories $\cC/Z$ and $Z/\cC$. 
\end{proposition}

\begin{proof}
	Let the maps $f_*$ and $f^*$ be formed by taking the unenriched representables. We only define the compositor and the unitor for the case of $\ucC/Z$ enrichment over $\cV/M$. The case for $Z/\ucC$ is identical. 
	
	The sliced internal homs are defined by the pullbacks 
	\[\begin{tikzcd}[ampersand replacement=\&,cramped]
		{\ucC/Z(f,g)} \&\& {\ucC(A,B)} \\
		\\
		M \& {\ucC(Z,Z)} \& {\ucC(A,Z)}
		\arrow[from=1-1, to=1-3]
		\arrow[from=1-1, to=3-1]
		\arrow["\ulcorner"{anchor=center, pos=0.125}, draw=none, from=1-1, to=3-2]
		\arrow["{g_*}", from=1-3, to=3-3]
		\arrow["r"', from=3-1, to=3-2]
		\arrow["{f^*}"', from=3-2, to=3-3]
	\end{tikzcd}\cdpunct[8pt]{.}\]
	
	The unitors are defined by the unique universal maps $i \colon I \rightarrow \ucC(A,A)$ as given in the diagrams
	\[\begin{tikzcd}[ampersand replacement=\&,cramped]
		I \&\&\& \\
		\& {\ucC/Z(f,f)} \&\& {\ucC(A,A)} \\
		\\
		\& M \& {\ucC(Z,Z)} \& {\ucC(A,Z)}
		\arrow["i"{description}, from=1-1, to=2-2]
		\arrow["{1_{A}}", curve={height=-12pt}, from=1-1, to=2-4]
		\arrow["e"', curve={height=12pt}, from=1-1, to=4-2]
		\arrow[from=2-2, to=2-4]
		\arrow[from=2-2, to=4-2]
		\arrow["\ulcorner"{anchor=center, pos=0.125}, draw=none, from=2-2, to=4-3]
		\arrow["{f_*}", from=2-4, to=4-4]
		\arrow["r"', from=4-2, to=4-3]
		\arrow["{f^*}"', from=4-3, to=4-4]
	\end{tikzcd}\cdpunct[8pt]{,}\]
	where the outer rectangle commutes because $r$ is a morphism of monoids.
	
	Finally, the compositors $\ucC/Z(g,h) \otimes \ucC/Z(f,g) \rightarrow \ucC/Z(f,h)$ are the unique universal maps induced by the following commuting diagrams
	\[\begin{tikzcd}[ampersand replacement=\&,cramped,column sep=small]
		{\ucC/Z(g,h) \otimes \ucC/Z(f,g)} \& {\ucC/Z(g,h) \otimes \ucC(A,B)} \& {\ucC(B,C) \otimes \ucC(A,B)} \& {\ucC(A,C)} \\
		\&\& {\ucC(B,Z) \otimes \ucC(A,B)} \\
		{M \otimes \ucC/Z(f,g)} \& {M \otimes \ucC(A,B)} \& {\ucC(Z,Z) \otimes \ucC(A,B)} \& {\ucC(A,Z)} \\
		\& {M \otimes \ucC(A,Z)} \& {\ucC(Z,Z) \otimes \ucC(A,Z)} \\
		{M \otimes M} \& {M \otimes \ucC(Z,Z)} \& {\ucC(Z,Z) \otimes \ucC(Z,Z)} \\
		\\
		M \&\&\& {\ucC(Z,Z)}
		\arrow[from=1-1, to=1-2]
		\arrow[from=1-1, to=3-1]
		\arrow[from=1-2, to=1-3]
		\arrow[from=1-2, to=3-2]
		\arrow["\circ", from=1-3, to=1-4]
		\arrow["{h_* \otimes 1}"', from=1-3, to=2-3]
		\arrow["{h_*}", from=1-4, to=3-4]
		\arrow["\circ", from=2-3, to=3-4]
		\arrow[from=3-1, to=3-2]
		\arrow[from=3-1, to=5-1]
		\arrow["{r \otimes 1}", from=3-2, to=3-3]
		\arrow["{1 \otimes g_*}"', from=3-2, to=4-2]
		\arrow["{g^* \otimes 1}", from=3-3, to=2-3]
		\arrow["{1 \otimes g_*}"', from=3-3, to=4-3]
		\arrow["{r \otimes 1}"', from=4-2, to=4-3]
		\arrow["\circ"', from=4-3, to=3-4]
		\arrow["{1 \otimes r}"', from=5-1, to=5-2]
		\arrow["\mu"', from=5-1, to=7-1]
		\arrow["{1 \otimes f^*}", from=5-2, to=4-2]
		\arrow["{r \otimes 1}"', from=5-2, to=5-3]
		\arrow["{1 \otimes f^*}", from=5-3, to=4-3]
		\arrow["\circ"', from=5-3, to=7-4]
		\arrow["r"', from=7-1, to=7-4]
		\arrow["{f^*}"', from=7-4, to=3-4]
	\end{tikzcd}\cdpunct[8pt]{.}\]
	
	For the associativity and unitality, first note that $U \colon \cV/M \rightarrow \cV$ is a faithful strict monoidal functor. Thus, $f = g$ in $\cV/M$ iff $U(f)=U(g)$ in $\cV$. Second, universal maps defined through pullbacks are uniquely characterized by the two projection maps. 
	
	Using these two facts we shall, for instance, demonstrate the left unitality. The associativity is quite similar but only technically more involved. Now, note the commuting diagram
	\[\begin{tikzcd}[ampersand replacement=\&,cramped]
		{I \otimes \ucC/Z(f,g)} \&\&\& {\ucC/Z(f,g)} \\
		\&\& {I \otimes \ucC(A,B)} \\
		\& {\ucC(B,B) \otimes \ucC/Z(f,g)} \\
		\\
		\& {\ucC/Z(g,g) \otimes \ucC(A,B)} \& {\ucC(B,B) \otimes \ucC(A,B)} \\
		{\ucC/Z(g,g) \otimes \ucC/Z(f,g)} \&\& {\ucC/Z(f,g)} \& {\ucC(A,B)}
		\arrow["\lambda", from=1-1, to=1-4]
		\arrow[from=1-1, to=2-3]
		\arrow["{1_B \otimes 1}", from=1-1, to=3-2]
		\arrow[from=1-1, to=6-1]
		\arrow[from=1-4, to=6-4]
		\arrow["{1_B \otimes 1}"', from=2-3, to=5-3]
		\arrow["\lambda", from=2-3, to=6-4]
		\arrow[from=3-2, to=5-3]
		\arrow[from=5-2, to=5-3]
		\arrow["\circ"', from=5-3, to=6-4]
		\arrow[from=6-1, to=3-2]
		\arrow[from=6-1, to=5-2]
		\arrow[from=6-1, to=6-3]
		\arrow[from=6-3, to=6-4]
	\end{tikzcd}\]
	implies that the left unitality holds under taking postcomposite with one of the projection maps. Similarly, we have the diagram
	\[\begin{tikzcd}[ampersand replacement=\&,cramped]
		{I \otimes \ucC/Z(f,g)} \&\&\&\& {\ucC/Z(g,g)  \otimes \ucC/Z(f,g)} \\
		\&\&\& {M \otimes \ucC/Z(f,g)} \\
		{\ucC/Z(f,g)} \&\& {I \otimes M} \& {M \otimes M} \& {\ucC/Z(f,g)} \\
		\\
		\\
		\&\& M
		\arrow[from=1-1, to=1-5]
		\arrow["{e \otimes 1}"', from=1-1, to=2-4]
		\arrow["\lambda"', from=1-1, to=3-1]
		\arrow[from=1-1, to=3-3]
		\arrow[from=1-5, to=2-4]
		\arrow[from=1-5, to=3-5]
		\arrow[from=2-4, to=3-4]
		\arrow[from=3-1, to=6-3]
		\arrow["{e \otimes 1}"', from=3-3, to=3-4]
		\arrow["\lambda"', from=3-3, to=6-3]
		\arrow[from=3-4, to=6-3]
		\arrow[from=3-5, to=6-3]
	\end{tikzcd}\]
	for compatibility with the other projection map. Finally, these two commuting diagrams jointly imply the desired commuting diagram 
	\[\begin{tikzcd}[ampersand replacement=\&,cramped]
		\& {\ucC/Z(g,g)  \otimes \ucC/Z(f,g)} \& \\
		\\
		{I \otimes \ucC/Z(f,g)} \&\& {\ucC/Z(f,g)}
		\arrow[from=1-2, to=3-3]
		\arrow[from=3-1, to=1-2]
		\arrow["\lambda"', from=3-1, to=3-3]
	\end{tikzcd}\]
	for left unitality. As already mentioned, the verification for associativity and the underlying category is mostly similar.

\end{proof}

Let $\ucC$ be any $\cV$ enriched category. There is a canonical monoid structure on the monoidal unit $I$ of $\cV$ given by $\mu : I \otimes I \xrightarrow{\lambda_I = \rho_I} I$ and $\eta : I \xrightarrow{1_I} I$. In particular, $\cV/I$ is a semicartesian monoidal category. That is, the induced monoidal unit $I \rightarrow I$ is a terminal object of $\cV/I$. Furthermore, for any $Z \in \cC$, the unitor map $I \rightarrow \ucC(Z,Z)$ is a morphism of monoid objects. This readily provides the following corollary.

\begin{corollary}
	There exists canonical enriched slice categories $\ucC/Z$, $Z/\ucC$ over $\cV/I$ for every $Z \in \cC$.
\end{corollary}

\begin{corollary}
	Let $\cV$ be a closed symmetric monooidal category, and let $\cC$ be enriched, tensored, and cotensored over $\cV$. Assume the relevant pullbacks and pushouts exist. Then, both the categories $\ucC/Z$ and $Z/\ucC$ are enriched, tensored, and cotensored over $\ucV/I$ for every $Z \in \ucC$. More generally, this is also true for a commutative monoid object $M$. 
\end{corollary}

\subsection{Model structures and slices} \label{sub:slicecatmod}

Monoidal model categories and enriched model categories are an important tool in abstract homotopy theory. Their defining axioms are closely linked with the two variable adjunction defining the closure of left $\cV$-modules. Hence, we can also analyze the slicing of enriched model structures.

\begin{definition}
	A monoidal model category is a \emph{closed} monoidal category  $(\cV,\otimes, I)$ together with a model structure on its underlying category $\cV$ that satisfies:
	\begin{enumerate}
		\item Unit axiom : Let $QI$ be any cofibrant replacement of $I$. The induced maps $QI \otimes X \rightarrow I \otimes X$ and $X \otimes QI \rightarrow X \otimes I$ are weak equivalences, whenever $X$ is cofibrant in $\cV$. 
		\item Pushout product axiom : Let $\alpha : X_1 \rightarrow Y_1$ and $\beta: X_2 \rightarrow Y_2$ be two arbitrary maps in  $\cV$. The pushout-product map $\alpha \bar{\otimes} \beta$ is a cofibration whenever both $\alpha$ and $\beta$ are cofibrations. Furthermore, if one of $\alpha$ or $\beta$ is also a weak equivalence then so is $\alpha \bar{\otimes} \beta$.
	\end{enumerate}
\end{definition} 

\begin{remark}
	The closedness assumption of the monoidal structure on $\cV$ is not entirely redundant. One of its important consequences is that the monoidal product $\otimes$ preserves the initial object in both of its arguments. In fact, it is sufficient to assume only this additional property for constructing the derived monoidal product using a left deformation. See \cite[Sec.~10.2]{reihlcht} or \cite[Prop.~4.3.1]{hovey2007modelcat} for more details.
\end{remark}

We now verify that this extends to a natural monoidal model structure for $\cV/M$.

\begin{proposition}
	Let $(\cV,\otimes, I)$ be a monoidal model category, $(M, \mu, \eta)$ be a monoid object inside of this monoidal category. Then, $(\cV/M, \otimes_{/M}, I_{/M})$ is also a monoidal model category equipped with the slice model structure. 
\end{proposition}

\begin{proof}
	The slice model structure is defined by the property that all of its cofibrations, fibrations, and weak equivalences are created under the forgetful functor $U : \cV/M \rightarrow \cV$. 
	\begin{enumerate}
		\item Unit Axiom : Any cofibrant replacement of $\eta : I \rightarrow M$ must also be created under the forgetful functor $U$. In particular, it must be of the form $QI \rightarrow I \rightarrow M$, where we take the functorial cofibrant replacement in $\cV$ itself. Now let $X \rightarrow M \in \cV/M$ be any cofibrant object. We obtain that the following diagram which demonstrates the required unit axiom.
		\[\begin{tikzcd}
			{QI \otimes X} && {I \otimes X} \\
			\\
			& {M \otimes X} \\
			& {M \otimes M} \\
			& M
			\arrow[from=1-1, to=1-3]
			\arrow[from=1-1, to=3-2]
			\arrow[from=1-3, to=3-2]
			\arrow[from=3-2, to=4-2]
			\arrow[from=4-2, to=5-2]
		\end{tikzcd}\]
		\item Pushout product axiom : The functor $U$ is strict monoidal and creates all colimits. Thus, the pushout product axiom of $\cV$ lifts to that of $\cV/M$ directly.  
	\end{enumerate}
\end{proof}

Having obtained the monoidal model category structure on $\cV/M$, we can now analyze the slicing of closed left $\cV$-modules under consideration of model structures. In particular, we start with a more general definition of a $\cV$-module model category instead of a enriched model category. 

\begin{definition}
	Let $(\cV,\otimes, I)$ be a monoidal model category. A left $\cV$-module model  category $\cC$ is a \emph{closed} left $\cV$-module $(\cC, \triangleleft)$ together with a model structure on $\cC$ that satisfies:
	\begin{enumerate}
		\item Tensoring unit axiom : Let $QI$ be any cofibrant replacement of $I$. The induced map $QI \triangleleft A \rightarrow I \triangleleft A$ is a weak equivalences, whenever $A$ is cofibrant in $\cC$. 
		\item Pushout product axiom : Let $\alpha : X \rightarrow Y \in \cV$ and $f : A \rightarrow B \in \cC$ be two arbitrary maps. The pushout-product map $\alpha \bar{\triangleleft} f$ is a cofibration, whenever both $\alpha$ and $f$ are cofibrations. Furthermore, if one of $\alpha$ or $f$ is a also weak equivalence  then so is $\alpha \bar{\triangleleft} f$.
	\end{enumerate}
\end{definition}

Motivated by the \cref{lem:sliceundermod}, we also have a dual definition. 

\begin{definition}
	Let $(\cV,\otimes, I)$ be a monoidal model category. A right $\cV^{\op}$-module model category $\cC$ is a \emph{closed} right $\cV$-module $(\cC^{\op}, \pitchfork^{\op})$ together with a model structure on $\cC$ that satisfies:
	\begin{enumerate}
		\item Powering unit axiom : Let $QI$ be any cofibrant replacement of $I$. The induced map $A \pitchfork I \rightarrow A \pitchfork QI$ is a weak equivalences, whenever $A$ is fibrant in $\cC$. 
		\item Pullback power axiom : Let $\alpha : X \rightarrow Y \in \cV$ and $f : A \rightarrow B \in \cC$ be two arbitrary maps. The pullback-power map $f \bar{\pitchfork} \alpha$ is a fibration, whenever $\alpha$ is a cofibration and $f$ is a fibrations. Furthermore, if one of $\alpha$ or $f$ is a also weak equivalence then so is $f \bar{\pitchfork} \alpha$.
	\end{enumerate}
\end{definition}

\begin{remark}
	Under the application of both $\op_1$ and $r$ dualities, the above definition specify precisely the tensoring unit axiom and the pushout-product axioms for the corresponding left $r(\cV)$-module structure on $\cC^{\op}$. In particular, the above definition is equivalent to specifying a left $r(\cV)$-module model category $\cC^{\op}$. The primary reason for having distinct axioms is that while $r(\cV)$ is again a monoidal model category, $\op_1(\cV)$ is not, in general.
\end{remark}

As before, there is also a correspondence between left $\cV$-module model categories and right $\cV^{\op}$-module model categories. This primarily follows from \cref{lem:poweringrightmod} together with its dual and \cite[Lemmas~4.2.2 and 4.2.7]{hovey2007modelcat}. 

\begin{theorem}
	Let $(\cV,\otimes, I)$ be a monoidal model category and let $(\cC,\triangleleft)$ be a left $\cV$-module model category. Let $(M, \mu, \eta)$ be a monoid object in $\cV$ and let $(Z,\lambda)$ be a left $M$-module in $\cC$. Then, it follows that $(\cC/Z,\triangleleft_{/Z})$ is a $\cV/M$-model category. 
\end{theorem}

\begin{proof}
	The slice model structure on $\cC/Z$ is again defined by the property that all of the cofibrations, fibrations, and weak equivalences are created under the forgetful functor $U : \cC/Z \rightarrow \cC$. 
	
	The rest of the proof is also similar to the case of a monoidal model category. Indeed, the idea of forgetful functors being \enquote{strict monoidal} is now encoded by the commuting diagram
	\[\begin{tikzcd}[ampersand replacement=\&,cramped]
		{\cV/M \times \cC/Z} \&\& {\cC/Z} \\
		\\
		{\cV \times \cC} \&\& \cC
		\arrow["{\triangleleft_{/Z}}", from=1-1, to=1-3]
		\arrow["{U \times U}"', from=1-1, to=3-1]
		\arrow["U", from=1-3, to=3-3]
		\arrow["\triangleleft"', from=3-1, to=3-3]
	\end{tikzcd}\cdpunct[8pt]{.}\]
\end{proof}

It is not possible to use the previous technique in the case of under categories $Z/\cC$. In particular, we only have the commuting diagram
\[\begin{tikzcd}[ampersand replacement=\&,cramped]
	{\cV/M \times Z/\cC} \&\& {Z/\cC} \\
	\\
	{\cV \times \cC} \&\& \cC
	\arrow["{\triangleleft_{/Z}}", from=1-1, to=1-3]
	\arrow["{U \times U}"', from=1-1, to=3-1]
	\arrow["U", from=1-3, to=3-3]
	\arrow["l"', between={0.2}{0.8}, Rightarrow, nfold, from=3-1, to=1-3]
	\arrow["\triangleleft"', from=3-1, to=3-3]
\end{tikzcd}\]
up to a $2$-cell $l\colon X \triangleleft A \rightarrow \alpha \triangleleft_{/Z} f$, encoding a \enquote{lax monoidal} functor.

We can attempt to remedy this situation by demanding that the canonical map $l$ be a natural weak equivalence. While this may potentially solve the issue, it is not a desirable assumption in practice. Instead, as before, the solution is to dualize the results with right $\cV^{\op}$-module structures.

\begin{theorem}
	Let $(\cV,\otimes, I)$ be a monoidal model category and let $(\cC,\pitchfork)$ be a right $\cV^{\op}$-module model category. Let $(M, \mu, \eta)$ be a monoid object in $\cV$ and let $(Z,\rho)$ be a right $M$-module in $\cC^{\op}$. Then, it follows that $(Z/\cC, \pitchfork_{/Z})$ is a right $(\cV/M)^{\op}$-module model category. Consequentially, this automatically implies $(Z/\cC, \triangleleft_{/Z})$ is a left $\cV/M$-module model category.
\end{theorem}

\begin{proof}
	The idea of forgetful functors being \enquote{strict monoidal} is now encoded by the following commuting diagram 
	\[\begin{tikzcd}[ampersand replacement=\&,cramped]
		{Z/\cC \times (\cV/M)^{\op}} \&\& {Z/\cC} \\
		\\
		{\cC \times \cV^{\op}} \&\& \cC
		\arrow["{\pitchfork_{/Z}}", from=1-1, to=1-3]
		\arrow["{U \times U^{\op}}"', from=1-1, to=3-1]
		\arrow["U", from=1-3, to=3-3]
		\arrow["\pitchfork"', from=3-1, to=3-3]
	\end{tikzcd}\cdpunct[8pt]{.}\]
	The functor $U \colon Z/\cC \rightarrow \cC$ creates all limit and lifts the pullback power axiom directly.
\end{proof}

Finally, we specialize to the case of a symmetric closed monoidal category for characterizing all enriched slice model categories.

\begin{theorem}
	Let $(\cV,\otimes, I)$ be a closed symmetric monoidal model category and let $\ucC$ be a $\ucV$-enriched model category with the underlying category $\cC$. Let $(M, \mu, \eta)$ be a commutative monoid object in $\cV$ and let $(Z,\lambda)$ be a left $M$-module in $\cC$ relative to tensoring. Then, it follows that both $Z/\ucC$ and $\ucC/Z$ are $\ucV/M$-enriched model categories.
\end{theorem}

	
	\newpage
	\appendix
	
	\section{Strictification functor of bicategories} \label{appsection:strictification}
	There are many sources (such as \cite{gps1995,gurski2013tricatcoh,yau2cat2021}) that discuss the strictification trifunctor $\st \colon \BiCat \rightarrow \nCat{2}_{\st,\ps}$ and its supposed right \emph{triadjoint} $i\colon \nCat{2}_{st,ps} \rightarrow \BiCat$. Here, $\nCat{2}_{\st,\ps}$ is also known as $\Gray$, and is the tricategory of strict $2$-categories, strict $2$-functors, pseudonatural transformations, and modifications. We aim to expand upon these descriptions by providing a more thorough account of this functor using concrete diagrams to fully specify the relevant trifunctorial data. Then, we shall briefly describe its usage to encode axioms of bicategorical structures using pasting diagrams. 

\subsection{Mapping of Bicategories}

Each bicategory $B \in \BiCat$ maps to a strict $2$-category called $\st B$, which is essentially determined by a pair of maps $e \colon i\circ \st B \rightarrow B$ and $s \colon B \rightarrow i \circ \st B$. For simplicity, we shall drop the $i$, unless required for better clarity. We first define the action of $e$ on the objects and $1$-cells of $\st B$. 

\begin{definition}
	The map $e$ is defined on $0$-cells and $1$-cells as follows:
	\begin{description}
		\item[$0$-cells -] for each $b \in \st B$, $e(b) = b \in B$,
		\item[$1$-cells -] for each non-empty string $(f_n,f_{n-1},\cdots,f_1)$ representing a $1$-cell in $\st B$, written right to left, $e$ is defined recursively as shown below:
		\[\begin{tikzcd}[ampersand replacement=\&,cramped,row sep=small]
			b \& {b_1} \& {b_2} \& \cdots \& {b_n} \& {\mapsto^{e}} \& b \&\& {b_n} \\
			\& {b_1} \& {b_2} \& \cdots \& {b_n} \& {\mapsto^{e}} \& {b_1} \&\& {b_n} \\
			\&\&\&\& \vdots \\
			\&\&\& {b_{n-1}} \& {b_n} \& {\mapsto^{e}} \& {b_{n-1}} \&\& {b_n}
			\arrow["{f_1}", from=1-1, to=1-2]
			\arrow["{f_2}", from=1-2, to=1-3]
			\arrow["{f_3}", from=1-3, to=1-4]
			\arrow[from=1-4, to=1-5]
			\arrow["{e(f_n,\cdots,f_2)f_1}", from=1-7, to=1-9]
			\arrow["{f_2}", from=2-2, to=2-3]
			\arrow["{f_3}", from=2-3, to=2-4]
			\arrow[from=2-4, to=2-5]
			\arrow["{e(f_n,\cdots,f_3)f_2}", from=2-7, to=2-9]
			\arrow["{f_n}", from=4-4, to=4-5]
			\arrow["{f_n}", from=4-7, to=4-9]
		\end{tikzcd}\cdpunct[2pt]{,}\]
		and for the empty string, it is defined seperately by setting $(b \xrightarrow{\varepsilon_b} b) \mapsto^{e} (1_b \colon b \rightarrow b)$.
	\end{description}
\end{definition}

Having partially defined $e$, we can now fully specify the $2$-category $\st B$ in terms of its $2$-cells, vertical composition, and horiztonal composition.

\begin{definition}
	The $2$-category $\st B$ has the following data:
	\begin{description}
		\item[Objects -] exactly the objects of $B$,
		\item[$1$-cells -] free strings of arrows $b \xrightarrow{f_1} b_1 \xrightarrow{f_2} b_2 \xrightarrow{f_2} b_3 \xrightarrow{f_3} \cdots $, where the empty string is allowed and is labeled as $b \xrightarrow{\varepsilon_{b}} b$,
		\item[Identity $1$-cells -] for each $b$, the empty string $b \xrightarrow{\varepsilon_{b}} b$ is the horizontal identity
		\item[Hcomp $1$-cells -] two strings of arrows are composed by concatenation of strings,
		\item[$2$-cells -] the $2$-cells between two strings $\vec{f}$ and $\vec{g}$ are in bijection with the $2$-cells $\alpha \colon e(\vec{f}) \Rightarrow e(\vec{g})$, 
		\[\begin{tikzcd}[ampersand replacement=\&]
			\& B \&\&\&\& {\st B} \&\&\& \\
			\\
			\&\&\&\& b \& {b_1} \& {b_2} \& \cdots \& {b_n} \\
			{b_1} \&\& {b_2} \\
			{b_1} \&\& {b_2} \\
			\&\&\&\& b \& {b'_1} \& {b'_2} \& \cdots \& {b'_{m}=b_n}
			\arrow["e"', from=1-6, to=1-2]
			\arrow["{f_1}", from=3-5, to=3-6]
			\arrow["{f_2}", from=3-6, to=3-7]
			\arrow["{f_3}", from=3-7, to=3-8]
			\arrow["{(\alpha \colon \vec{f} \Rightarrow \vec{g} \in B)}", between={0.2}{0.8}, Rightarrow, nfold, from=3-7, to=6-7]
			\arrow[from=3-8, to=3-9]
			\arrow[""{name=0, anchor=center, inner sep=0}, "{e(\vec{f})}", curve={height=-12pt}, from=4-1, to=4-3]
			\arrow[""{name=1, anchor=center, inner sep=0}, "{e(\vec{g})}"', curve={height=12pt}, from=5-1, to=5-3]
			\arrow["{g_1}"', from=6-5, to=6-6]
			\arrow["{g_2}"', from=6-6, to=6-7]
			\arrow["{g_3}"', from=6-7, to=6-8]
			\arrow[from=6-8, to=6-9]
			\arrow["\alpha", between={0.2}{0.8}, Rightarrow, nfold, from=0, to=1]
		\end{tikzcd}\cdpunct[8pt]{,}\]
		this bijection of $2$-cells also defines the action of $e$ on the $2$-cells of $\st B$,
		\item[Identity $2$-cells -] for each string $\vec{f}$, the unique $2$-cell $1_{e({\vec{f}})} \colon e(\vec{f}) \Rightarrow e(\vec{f})$ is the identity,
		\item[Vcomp $2$-cells -] the vertical composition of $2$-cells is directly lifted from the map $e$ as shown below
		\[\begin{tikzcd}[ampersand replacement=\&,column sep=small]
			\&\&\& B \&\&\&\& {\st B} \&\&\& \\
			\\
			\&\&\&\&\&\& b \& {b_1} \& {b_2} \& \cdots \& {b_n} \\
			\&\& {b_1} \&\& {b_2} \\
			\\
			{b_1} \&\& {b_2} \&\&\& b \& {b'_1} \& {b'_2} \& \cdots \& {b'_{m}=b_n} \\
			\\
			\&\& {b_1} \&\& {b_2} \\
			\&\&\&\&\&\& b \& {b''_1} \& {b''_2} \& \cdots \& {b''_{p} = b_n}
			\arrow["e"', from=1-8, to=1-4]
			\arrow["{f_1}", from=3-7, to=3-8]
			\arrow["{f_2}", from=3-8, to=3-9]
			\arrow["{f_3}", from=3-9, to=3-10]
			\arrow["{(\alpha \colon \vec{f} \Rightarrow \vec{g} \in B)}", between={0.2}{0.8}, Rightarrow, nfold, from=3-9, to=6-8]
			\arrow[""{name=0, anchor=center, inner sep=0}, from=3-10, to=3-11]
			\arrow[""{name=1, anchor=center, inner sep=0}, "{e(\vec{f})}", curve={height=-12pt}, from=4-3, to=4-5]
			\arrow[""{name=2, anchor=center, inner sep=0}, "{e(\vec{g})}"', from=6-1, to=6-3]
			\arrow["{g_1}"', from=6-6, to=6-7]
			\arrow["{g_2}"', from=6-7, to=6-8]
			\arrow["{g_3}"', from=6-8, to=6-9]
			\arrow["{(\beta \colon \vec{g} \Rightarrow \vec{h} \in B)}", between={0.2}{0.8}, Rightarrow, nfold, from=6-8, to=9-9]
			\arrow[from=6-9, to=6-10]
			\arrow[""{name=3, anchor=center, inner sep=0}, "{e(\vec{h})}"', curve={height=12pt}, from=8-3, to=8-5]
			\arrow["{h_1}"', from=9-7, to=9-8]
			\arrow["{h_2}"', from=9-8, to=9-9]
			\arrow["{h_3}"', from=9-9, to=9-10]
			\arrow[""{name=4, anchor=center, inner sep=0}, from=9-10, to=9-11]
			\arrow["{\beta\alpha}", shift left=5, curve={height=-30pt}, between={0.2}{0.8}, Rightarrow, nfold, from=0, to=4]
			\arrow["\alpha", between={0.2}{0.8}, Rightarrow, nfold, from=1, to=2]
			\arrow["{\beta\alpha}", between={0.2}{0.8}, Rightarrow, nfold, from=1, to=3]
			\arrow["\beta", between={0.2}{0.8}, Rightarrow, nfold, from=2, to=3]
		\end{tikzcd}\]
		\item[Hcomp $2$-cells -] the horiztonal composite of two $2$-cells $\alpha \colon e(\vec{f}) \Rightarrow e(\vec{g})$ and $\beta \colon e(\vec{h}) \Rightarrow e(\vec{k})$ is given as follows:
		\[\begin{tikzcd}[ampersand replacement=\&]
			b \& {b_1} \& \cdots \& {b_n=c} \& {c_1} \& \cdots \& {c_{p}=d} \\
			\\
			b \& {b'_1} \& \cdots \& {b'_{m}=b_n=c} \& {c'_1} \& \cdots \& {c'_q = c_{p}=d} \\
			\&\&\& \vequiv \\
			b \& {b_1} \& \cdots \& {b_n=c} \& {c_1} \& \cdots \& {c_p=d} \\
			\\
			b \&\&\& c \&\&\& d \\
			\\
			b \&\&\& c \&\&\& d \\
			\\
			b \& {b'_1} \& \cdots \& {b'_m=b_n=c} \& {c'_1} \& \cdots \& {c'_q=c_p=d}
			\arrow["{f_1}", from=1-1, to=1-2]
			\arrow[""{name=0, anchor=center, inner sep=0}, "{f_2}", from=1-2, to=1-3]
			\arrow[from=1-3, to=1-4]
			\arrow["{h_1}", from=1-4, to=1-5]
			\arrow[""{name=1, anchor=center, inner sep=0}, "{h_2}", from=1-5, to=1-6]
			\arrow[from=1-6, to=1-7]
			\arrow["{g_1}"', from=3-1, to=3-2]
			\arrow[""{name=2, anchor=center, inner sep=0}, "{g_2}"', from=3-2, to=3-3]
			\arrow[from=3-3, to=3-4]
			\arrow["{k_1}"', from=3-4, to=3-5]
			\arrow[""{name=3, anchor=center, inner sep=0}, "{k_2}"', from=3-5, to=3-6]
			\arrow[from=3-6, to=3-7]
			\arrow["{f_1}", from=5-1, to=5-2]
			\arrow["{f_2}", from=5-2, to=5-3]
			\arrow[from=5-3, to=5-4]
			\arrow["{h_1}", from=5-4, to=5-5]
			\arrow["{e(\vec{h},\vec{f}) \Rightarrow e(\vec{h}) e(\vec{f})}", between={0.2}{0.8}, Rightarrow, nfold, from=5-4, to=7-4]
			\arrow["{h_2}", from=5-5, to=5-6]
			\arrow[from=5-6, to=5-7]
			\arrow["{e(\vec{f})}"{description}, from=7-1, to=7-4]
			\arrow["{e(\vec{h})}"{description}, from=7-4, to=7-7]
			\arrow["{\beta \star \alpha}"{description}, between={0.2}{0.8}, Rightarrow, nfold, from=7-4, to=9-4]
			\arrow["{e(\vec{g})}"{description}, from=9-1, to=9-4]
			\arrow["{e(\vec{k})}"{description}, from=9-4, to=9-7]
			\arrow["{ e(\vec{k})e(\vec{g})\Rightarrow e(\vec{k},\vec{g})}", between={0.2}{0.8}, Rightarrow, nfold, from=9-4, to=11-4]
			\arrow["{g_1}"', from=11-1, to=11-2]
			\arrow["{g_2}"', from=11-2, to=11-3]
			\arrow[from=11-3, to=11-4]
			\arrow["{k_1}"', from=11-4, to=11-5]
			\arrow["{k_2}"', from=11-5, to=11-6]
			\arrow[from=11-6, to=11-7]
			\arrow[""{name=4, anchor=center, inner sep=0}, "{(\alpha \colon \vec{f} \Rightarrow \vec{g} \in B)}"', between={0.2}{0.8}, Rightarrow, nfold, from=0, to=2]
			\arrow[""{name=5, anchor=center, inner sep=0}, "{(\beta \colon \vec{h} \Rightarrow \vec{k} \in B)}", between={0.2}{0.8}, Rightarrow, nfold, from=1, to=3]
			\arrow["\star"{description}, draw=none, from=4, to=5]
		\end{tikzcd}\cdpunct[8pt]{,}\]
		where the $2$-cells $e(\vec{k})e(\vec{g})\Rightarrow e(\vec{k},\vec{g})$ and $e(\vec{h},\vec{f}) \Rightarrow e(\vec{h}) e(\vec{f})$ are induced uniquely from the coherence constraints, and $\beta \star \alpha$ is the horizontal composite in the bicategory $B$ itself.
	\end{description}
\end{definition}

Now, we wish to specify $e \colon i\circ \st B \rightarrow B$ and $s \colon B \rightarrow i \circ \st B$ as pseudofunctors.

\begin{definition}
	The pseudofunctor $e \colon \st B \rightarrow B$ is defined by the following data:
	\begin{description}
		\item[Objects -] for each $b\in \st B$, the pseudofunctor $e$ acts as $b \mapsto b$,
		\item[$1$-cells -] for each string $\vec{f} \in \st B$, the pseudofunctor $e$ acts as $\vec{f} \mapsto e(\vec{f})$,
		\item[$2$-cells -] each $2$-cell in $\st B$ is already defined as a $2$-cell in $B$,
		\item[Composition constraint - ] the $2$-cell $e(\vec{g}) e(\vec{f}) \Rightarrow e(\vec{g},\vec{f})$ is the unique coherence constraint $2$-cell,
		\[\begin{tikzcd}[ampersand replacement=\&]
			\& {\st B} \&\&\&\& B \& \\
			\& c \&\&\&\& c \\
			b \&\& d \& \rightsquigarrow \& b \&\& d
			\arrow["e", from=1-2, to=1-6]
			\arrow["{\vec{g}}", from=2-2, to=3-3]
			\arrow["{e(\vec{g})}", from=2-6, to=3-7]
			\arrow["{\vec{f}}", from=3-1, to=2-2]
			\arrow[""{name=0, anchor=center, inner sep=0}, "{\vec{g},\vec{f}}"', from=3-1, to=3-3]
			\arrow["{e(\vec{f})}", from=3-5, to=2-6]
			\arrow[""{name=1, anchor=center, inner sep=0}, "{e(\vec{g},\vec{f})}"', from=3-5, to=3-7]
			\arrow[between={0.2}{0.8}, equals, nfold, from=2-2, to=0]
			\arrow[between={0.2}{0.8}, Rightarrow, nfold, from=2-6, to=1]
		\end{tikzcd}\cdpunct[8pt]{,}\]
		\item[Identity constraint - ] the $2$-cell $1_b \Rightarrow e(b \xrightarrow{\varepsilon_{b}} b)$ induced by the identity $2$-cell at $1_b \colon b \rightarrow b$,
		\[\begin{tikzcd}[ampersand replacement=\&]
			\& {\st B} \&\&\&\& B \& \\
			\\
			b \&\& b \& \rightsquigarrow \& b \&\& d
			\arrow["e", from=1-2, to=1-6]
			\arrow["{\varepsilon_{b}}"', from=3-1, to=3-3]
			\arrow[""{name=0, anchor=center, inner sep=0}, "{1_{b}=e(\varepsilon_{b})}"', curve={height=12pt}, from=3-5, to=3-7]
			\arrow[""{name=1, anchor=center, inner sep=0}, "{1_{b}}", curve={height=-12pt}, from=3-5, to=3-7]
			\arrow["{1_{1_b}}", between={0.2}{0.8}, Rightarrow, nfold, from=1, to=0]
		\end{tikzcd}\cdpunct[8pt]{.}\]
	\end{description}
\end{definition}

The pseudofunctor $s \colon B \rightarrow \st B$ is supposed to be a strict section of the pseudofunctor $e \colon \st B \rightarrow B$, that is, $e \circ s = 1_{B}$. We shall specify $s$ with this view in mind.

\begin{definition}
	The pseudofunctor $s \colon B \rightarrow \st B$ is defined as follows:
	\begin{description}
		\item[Objects -] for each $b\in B$, the pseudofunctor $s$ acts as $b \mapsto b$,
		\item[$1$-cells -] for each $f \colon b \rightarrow c \in B$, the pseudofunctor $f$ acts as $f \mapsto (f)$, where $(f)$ is string of length one, for simplicity this will now be labeled as $f$ itself,
		\item[$2$-cells -] the $2$-cell $\alpha \colon f \Rightarrow g$ is mapped precisely to the $2$-cell $\alpha \colon f \Rightarrow g$ under inverse image in $e$,
		\item[Composition constraint -] the $2$-cell $(g,f) \rightarrow gf$ is induced by the identity cell $1_{gf} \colon e(g,f) \Rightarrow e(gf)$,
		\[\begin{tikzcd}[ampersand replacement=\&]
			\& B \&\&\&\& {\st B} \& \\
			\& c \&\&\&\& c \\
			b \&\& d \& \rightsquigarrow \& b \&\& d
			\arrow["f"', from=1-2, to=1-6]
			\arrow["g", from=2-2, to=3-3]
			\arrow["g", from=2-6, to=3-7]
			\arrow["f", from=3-1, to=2-2]
			\arrow[""{name=0, anchor=center, inner sep=0}, "gf"', from=3-1, to=3-3]
			\arrow["f", from=3-5, to=2-6]
			\arrow[""{name=1, anchor=center, inner sep=0}, "gf"', from=3-5, to=3-7]
			\arrow[between={0.2}{0.8}, equals, nfold, from=2-2, to=0]
			\arrow["{1_{gf}}", between={0.2}{0.8}, Rightarrow, nfold, from=2-6, to=1]
		\end{tikzcd}\cdpunct[8pt]{,}\]
		\item[Identity constraint -] the $2$-cell $\varepsilon_{b} \Rightarrow 1_b$ is induced by the unique identity cell $1_{1_b} \colon 1_b \Rightarrow 1_b$,
		\[\begin{tikzcd}[ampersand replacement=\&]
			\& B \&\&\&\& {\st B} \& \\
			\\
			b \&\& b \& \rightsquigarrow \& b \&\& b
			\arrow["f"', from=1-2, to=1-6]
			\arrow["{1_b}"', from=3-1, to=3-3]
			\arrow[""{name=0, anchor=center, inner sep=0}, "{1_b}"', curve={height=12pt}, from=3-5, to=3-7]
			\arrow[""{name=1, anchor=center, inner sep=0}, "{\varepsilon_{b}}", curve={height=-12pt}, from=3-5, to=3-7]
			\arrow["{1_{1_b}}", between={0.2}{0.8}, Rightarrow, nfold, from=1, to=0]
		\end{tikzcd}\cdpunct[8pt]{.}\]
	\end{description}
\end{definition}

It is highly important to note that the constraints cells of $s$ are non-trivial, even if both of them are represented by the identity $2$-cells under $e$. Furthermore, the equality $e \circ s = 1_{B}$ is meant as equality of pseudofunctors, which includes their coherence constraints. We will not verify the pseudo-functor axioms for both $e$ and $s$, as they are routine technical checks.

\subsection{\tpdf{Mapping of $\BiCat(A,B)$}}

Our goal now is determine the action of $\st$ on the bicategory of pseudofunctors, pseudotransformations, and modification. To this end, the following lemma shall be essential.

\begin{lemma}
	The dotted composite $\Phi$ defined in the diagram
	\[\begin{tikzcd}[ampersand replacement=\&,cramped]
		{\Gray(\st A,B)} \&\& {\BiCat(i(\st A),iB)} \& {\BiCat(i(\st A),iB) \times *} \\
		\\
		\&\& {\BiCat(A,iB)} \& {\BiCat(i(\st A),iB) \times \BiCat(A,i(\st A))}
		\arrow["i", from=1-1, to=1-3]
		\arrow["\Phi"', dashed, from=1-1, to=3-3]
		\arrow["\cong", from=1-3, to=1-4]
		\arrow["{(s_{A})^{*}}", from=1-3, to=3-3]
		\arrow["{1 \times s_{A}}", from=1-4, to=3-4]
		\arrow["\circ", from=3-4, to=3-3]
	\end{tikzcd}\]
	is a strict $2$-functor of strict $2$-categories. Here, $s_{A}$ is the constant pseudofunctor out of the terminal bicategory which picks $s_{A} \in \BiCat(A,i(\st A))$, and the functor $\circ$ is the compositor of $\BiCat$.
\end{lemma}

\begin{proof}
	We first note that $\BiCat(A,Y)$ is a strict $2$-category, whenever $Y$ is a strict $2$-category. Thus, all of the bicategories involved in defining the above composite are all strict $2$-categories. Next, we verify that $\Phi$ is a strict $2$-functor of strict $2$-categories. Note that $\Phi(\gamma \colon F \Rightarrow  G)$ is given by $\gamma \star 1_{S_{A}}$. Using the standard definition of the compositor in $\BiCat$, we can represent this using the following transformation cell in $B$:
	\[\begin{tikzcd}[ampersand replacement=\&]
		\&\& Fa \&\& Fb \&\&\&\& \\
		\&\&\&\&\&\& Fa \&\& Fb \\
		{\Phi(\gamma)} \& \equiv \& Ga \&\& Gb \& \equiv \\
		\&\&\&\&\&\& Ga \&\& Gb \\
		\&\& Ga \&\& Gb
		\arrow["Ff", from=1-3, to=1-5]
		\arrow["{\gamma_{a}}"', from=1-3, to=3-3]
		\arrow["{\gamma_{b}}", from=1-5, to=3-5]
		\arrow["Ff", from=2-7, to=2-9]
		\arrow["{\gamma_a}"', from=2-7, to=4-7]
		\arrow["{\gamma_b}", from=2-9, to=4-9]
		\arrow["{\gamma_{f}}"{description}, between={0.2}{0.8}, Rightarrow, nfold, from=3-3, to=1-5]
		\arrow["Gf"{description}, from=3-3, to=3-5]
		\arrow["{1_{Ga}}"', from=3-3, to=5-3]
		\arrow["{1_{Ga}}", from=3-5, to=5-5]
		\arrow["{\gamma_f}"{description}, between={0.2}{0.8}, Rightarrow, nfold, from=4-7, to=2-9]
		\arrow["Gf"', from=4-7, to=4-9]
		\arrow["{1_{Gf}}"{description}, between={0.2}{0.8}, Rightarrow, nfold, from=5-3, to=3-5]
		\arrow["Gf"', from=5-3, to=5-5]
	\end{tikzcd}\cdpunct[8pt]{.}\]
	Crucially, we emphasize that this equality is possible only because both $F,G \colon \st A \rightarrow B$ are strict $2$-functors of strict $2$-categories. With this description, it is now straightforward to verify both $(\gamma_2\gamma_1) \star 1_{S_{A}}=(\gamma_2 \star 1_{S_{A}})(\gamma_1 \star 1_{S_{A}})$ and $1_F \star 1_{S_{A}} = 1_{F \star S_{A}}$. Thus, the pseudofunctor $\Phi$ is indeed strict $2$-functor.
	
\end{proof}

It turns out that $\Phi$ is actually a strict $2$-isomorphism of $2$-categories. For this, we only need to show that it is a bijection on the underlying $2$-graphs. Thus, we state the following lemma.

\begin{lemma} \leavevmode
	\begin{enumerate}
		\item For every pseudofunctor $U\colon A \rightarrow Y$ with $Y \in \Gray$, there exists a unique strict functor $\tilde U \colon \st A \rightarrow Y$ such that $\tilde{U} \circ s_{A} = U$.
		\item For every pseudotransformation $\gamma\colon U \Rightarrow V \colon A \rightarrow Y$ with $Y \in \Gray$, there exists a unique pseudotransformation $\tilde \gamma$ such that $\tilde{\gamma} \star 1_{s_{A}} = \rho$.
		\item For every modification $\theta \colon \gamma \Rrightarrow \delta \colon U \Rightarrow V  \colon A \rightarrow Y$ with $Y \in \Gray$, there exists a unique modification $\tilde{\theta}$ such that $\tilde{\theta} \cdot 1_{1_{s_{A}}} = \omega$.
	\end{enumerate}
\end{lemma}

\begin{proof} 
	The only non-trivial part of the proof is defining each of the unique lifts, the rest is a routine.
	\begin{enumerate}
		\item Define $\tilde{U}$ by the following data:
		\begin{description}
			\item[Objects -] for any $a \in \st A$, $U$ maps $a \mapsto Ua \in Y$,
			\item[$1$-cells -] the non-identity $1$-cells $\vec{f}$ are mapped as
			\[\begin{tikzcd}[ampersand replacement=\&,cramped]
				a \& {a_1} \& \cdots \& {\mapsto^{\tilde{U}}} \& Ua \& {Ua_1} \& \cdots
				\arrow["{f_1}", from=1-1, to=1-2]
				\arrow["{f_2}", from=1-2, to=1-3]
				\arrow["{Uf_1}", from=1-5, to=1-6]
				\arrow["{Uf_2}", from=1-6, to=1-7]
			\end{tikzcd}\cdpunct[8pt]{,}\]
			and the identity $1$-cell is mapped to the identity $1$-cell,
			\item[$2$-cells -] the action on $2$-cells is given as
			\[\begin{tikzcd}[ampersand replacement=\&,cramped]
				a \& {a_1} \& \cdots \& b \&\&\&\&\& \\
				\\
				a \& {a'_1} \& \cdots \& b \\
				\\
				a \& {a_1} \& \cdots \& b \&\& Ua \& {Ua_1} \& \cdots \& Ub \\
				a \&\&\& b \&\& Ua \&\&\& Ub \\
				\&\&\&\& {\mapsto^{\tilde{U}}} \\
				a \&\&\& b \&\& Ua \&\&\& Ub \\
				a \& {a'_1} \& \cdots \& b \&\& Ua \& {Ua'_1} \& \cdots \& Ub
				\arrow["{f_1}", from=1-1, to=1-2]
				\arrow[""{name=0, anchor=center, inner sep=0}, "{f_2}", from=1-2, to=1-3]
				\arrow[from=1-3, to=1-4]
				\arrow["{g_1}"', from=3-1, to=3-2]
				\arrow[""{name=1, anchor=center, inner sep=0}, "{g_2}"', from=3-2, to=3-3]
				\arrow[from=3-3, to=3-4]
				\arrow["{f_1}", from=5-1, to=5-2]
				\arrow[""{name=2, anchor=center, inner sep=0}, "{f_2}", from=5-2, to=5-3]
				\arrow[from=5-3, to=5-4]
				\arrow["{Uf_1}", from=5-6, to=5-7]
				\arrow[""{name=3, anchor=center, inner sep=0}, "{Uf_2}", from=5-7, to=5-8]
				\arrow[from=5-8, to=5-9]
				\arrow[""{name=4, anchor=center, inner sep=0}, "{e(\vec{f})}"{description}, from=6-1, to=6-4]
				\arrow[""{name=5, anchor=center, inner sep=0}, "{U(e(\vec{f}))}"{description}, from=6-6, to=6-9]
				\arrow[""{name=6, anchor=center, inner sep=0}, "{e(\vec{g})}"{description}, from=8-1, to=8-4]
				\arrow[""{name=7, anchor=center, inner sep=0}, "{U(e(\vec{g}))}"{description}, from=8-6, to=8-9]
				\arrow["{g_1}"', from=9-1, to=9-2]
				\arrow[""{name=8, anchor=center, inner sep=0}, "{g_2}"', from=9-2, to=9-3]
				\arrow[from=9-3, to=9-4]
				\arrow["{Ug_1}"', from=9-6, to=9-7]
				\arrow[""{name=9, anchor=center, inner sep=0}, "{Ug_2}"', from=9-7, to=9-8]
				\arrow[from=9-8, to=9-9]
				\arrow["{(\alpha : e(\vec{f}) \Rightarrow e(\vec{g}))}", between={0.2}{0.8}, Rightarrow, nfold, from=0, to=1]
				\arrow["\equiv"{marking, allow upside down}, draw=none, from=1, to=2]
				\arrow["{\mathrm{Coh}_{S_{A}}}", between={0.2}{0.8}, Rightarrow, nfold, from=2, to=4]
				\arrow["{\mathrm{Coh}_U}", between={0.2}{0.8}, Rightarrow, nfold, from=3, to=5]
				\arrow["\alpha"{description}, between={0.2}{0.8}, Rightarrow, nfold, from=4, to=6]
				\arrow["{U(\alpha)}", between={0.2}{0.8}, Rightarrow, nfold, from=5, to=7]
				\arrow["{\mathrm{Coh}_{S_{A}}^{-1}}", between={0.2}{0.8}, Rightarrow, nfold, from=6, to=8]
				\arrow["{\mathrm{Coh}_U^{-1}}", between={0.2}{0.8}, Rightarrow, nfold, from=7, to=9]
			\end{tikzcd}\cdpunct[8pt]{,}\]
			where $\mathrm{Coh}_U$ and $\mathrm{Coh}_{S_{A}}$ represent the unique coherence constraint cells.
		\end{description}
		
		By choosing $1_{\vec{f}} \colon \vec{f} \Rightarrow e(\vec{f})$, the definition readily implies that $\tilde{U}(\mathrm{Coh}_{s_{A}}) = \mathrm{Coh}_U$. Thus, $\tilde{U} \circ s_{A} = U$ holds as desired, including the equality of coherence data.
		
		\item Define $\tilde\gamma$ by the following data:
		\begin{description}
			\item[Component $1$-cell] for each $a \in \st A$, the component $(\st \gamma)_{a}$ is defined as $\gamma_{a}$ itself,
			\item[Component $2$-cell] for each string $\vec{f} \in \st A$, the component $(\st \gamma)_{\vec{f}}$ is defined by the following composite $2$-cell in $Y$
			\[\begin{tikzcd}[ampersand replacement=\&,cramped]
				Ua \& {Ua_1} \& {Ua_2} \& \cdots \& Ub \\
				\\
				Va \& {Va_1} \& {Va_2} \& \cdots \& Vb
				\arrow["{Uf_1}", from=1-1, to=1-2]
				\arrow["{\gamma_a}"', from=1-1, to=3-1]
				\arrow["{Uf_2}", from=1-2, to=1-3]
				\arrow["{\gamma_{a_1}}"{description}, from=1-2, to=3-2]
				\arrow["{Uf_3}", from=1-3, to=1-4]
				\arrow["{\gamma_{a_2}}"{description}, from=1-3, to=3-3]
				\arrow[from=1-4, to=1-5]
				\arrow["{\gamma_b}", from=1-5, to=3-5]
				\arrow["{\gamma_{f_1}}"{description}, between={0.2}{0.8}, Rightarrow, nfold, from=3-1, to=1-2]
				\arrow["{Vf_1}"', from=3-1, to=3-2]
				\arrow["{\gamma_{f_2}}"{description}, between={0.2}{0.8}, Rightarrow, nfold, from=3-2, to=1-3]
				\arrow["{Vf_2}"', from=3-2, to=3-3]
				\arrow["{Vf_3}"', from=3-3, to=3-4]
				\arrow[from=3-4, to=3-5]
			\end{tikzcd}\cdpunct[8pt]{.}\]
		\end{description}
		
		\item Define $\tilde \theta \colon \tilde{\gamma} \Rrightarrow \tilde{\delta}$ by setting $(\tilde \theta)_{a} = \theta_a \colon \gamma_a \Rightarrow \delta_a$
	\end{enumerate}
\end{proof}

As a consequence, we obtain the desired proposition.

\begin{proposition}
		The dotted composite $\Phi$ defined in the diagram
		\[\begin{tikzcd}[ampersand replacement=\&,cramped]
			{\Gray(\st A,B)} \&\& {\BiCat(i(\st A),iB)} \& {\BiCat(i(\st A),iB) \times *} \\
			\\
			\&\& {\BiCat(A,iB)} \& {\BiCat(i(\st A),iB) \times \BiCat(A,i(\st A))}
			\arrow["i", from=1-1, to=1-3]
			\arrow["\Phi"', dashed, from=1-1, to=3-3]
			\arrow["\cong", from=1-3, to=1-4]
			\arrow["{(s_{A})^{*}}", from=1-3, to=3-3]
			\arrow["{1 \times s_{A}}", from=1-4, to=3-4]
			\arrow["\circ", from=3-4, to=3-3]
		\end{tikzcd}\]
		is a strict $2$-isomorphism of strict $2$-categories.
\end{proposition}

\begin{remark}
	One may visualize this universal property using the diagram 
	\[\begin{tikzcd}[ampersand replacement=\&]
		{i \circ \st A} \&\&\&\&\&\&\&\& iY \\
		\\
		\\
		\\
		\\
		\\
		A
		\arrow[""{name=0, anchor=center, inner sep=0}, "{\tilde{U}}", curve={height=-24pt}, from=1-1, to=1-9]
		\arrow[""{name=1, anchor=center, inner sep=0}, "{\tilde{V}}"', curve={height=24pt}, from=1-1, to=1-9]
		\arrow["{s_{A}}", from=7-1, to=1-1]
		\arrow[""{name=2, anchor=center, inner sep=0}, "U", curve={height=-24pt}, from=7-1, to=1-9]
		\arrow[""{name=3, anchor=center, inner sep=0}, "V"', curve={height=24pt}, from=7-1, to=1-9]
		\arrow[""{name=4, anchor=center, inner sep=0}, "{\tilde{\gamma}}"', curve={height=24pt}, between={0.2}{0.8}, Rightarrow, nfold, from=0, to=1]
		\arrow[""{name=5, anchor=center, inner sep=0}, "{\tilde{\delta}}", curve={height=-24pt}, between={0.2}{0.8}, Rightarrow, nfold, from=0, to=1]
		\arrow[""{name=6, anchor=center, inner sep=0}, "\gamma"', curve={height=24pt}, between={0.2}{0.8}, Rightarrow, nfold, from=2, to=3]
		\arrow[""{name=7, anchor=center, inner sep=0}, "\delta", curve={height=-24pt}, between={0.2}{0.8}, Rightarrow, nfold, from=2, to=3]
		\arrow["{\tilde{\theta}}", between={0.2}{0.8}, Rightarrow, scaling nfold=3, from=4, to=5]
		\arrow["\theta", between={0.2}{0.8}, Rightarrow, scaling nfold=3, from=6, to=7]
	\end{tikzcd}\cdpunct[8pt]{.}\]
\end{remark}

As with any ordinary adjunction, we can utilize this universal property to extend $\st$ to a trifunctor. We first extend it on $\BiCat(A,B)$.

\begin{definition}
	Define the post composition pseudofunctor $(s_{B})_{*}$, corresponding to $s_{B} \colon B \rightarrow \st B$, as the composite pseudofunctor
	\[\begin{tikzcd}[ampersand replacement=\&]
		{\BiCat(A,B)} \&\& {* \times \BiCat(A,B) } \\
		\\
		{\BiCat(A,i \circ \st B)} \&\& {\BiCat(B,i\circ\st B) \times \BiCat(A,B)}
		\arrow["\cong", from=1-1, to=1-3]
		\arrow["{(s_{B})_{*}}"', from=1-1, to=3-1]
		\arrow["{s_{B} \times 1}", from=1-3, to=3-3]
		\arrow["\circ", from=3-3, to=3-1]
	\end{tikzcd}\cdpunct[8pt]{.}\]
\end{definition}

Note that $(s_{B})_{*}$ is not a strict $2$-functor. The compositor is a modification given by unique coherence cell $(1_{s_{B}} \star \gamma_2) (1_{s_{B}} \star \gamma_1) \Rrightarrow (1_{s_{B}} \star \gamma_2\gamma_1)$. Similarly, the identity coherence modification is given by the unique coherence cell $1_{s_{B} \circ F} \Rrightarrow 1_{s_{B}} \star 1_F$.

\begin{definition}
	Define $\st_{A,B} \colon \BiCat(A,B) \rightarrow \Gray(\st A, \st B)$ as the composite pseudofunctor
	\[\begin{tikzcd}[ampersand replacement=\&]
		{\BiCat(A,B)} \&\& {\BiCat(A,i \circ \st B)} \&\& {\Gray(\st A, \st B)}
		\arrow["{(s_{B})_{*}}", from=1-1, to=1-3]
		\arrow["{\Phi^{-1}}", from=1-3, to=1-5]
	\end{tikzcd}\cdpunct[1pt]{.}\]
\end{definition}

With these definitions, one can readily check that the following diagram commutes strictly
\[\begin{tikzcd}[ampersand replacement=\&]
	{\BiCat(A,B)} \&\&\& {\BiCat(i\circ \st A, i \circ \st B)} \\
	\\
	\\
	{\BiCat(A,B)} \&\&\& {\BiCat(A,i \circ \st B)}
	\arrow["{i_{\st A,\st B} \circ \st_{A,B}}", from=1-1, to=1-4]
	\arrow["1"', from=1-1, to=4-1]
	\arrow["{(s_{A})^{*}}", from=1-4, to=4-4]
	\arrow["{(s_{B})_{*}}"', from=4-1, to=4-4]
\end{tikzcd}\cdpunct[8pt]{.}\]

Finally, the coherence data of $\st_{A,B}$ is simply the image under $\Phi^{-1}$ of the coherence constraint cells of $(s_{B})_{*}$. Thus, we obtain the following concrete description of $\st_{A,B}$.

\begin{lemma}
	The coherence data of $\st_{A,B} \colon \BiCat(A,B) \rightarrow \Gray(\st A,\st B)$ is given by
	\begin{description}
		\item[Composition constraint -] for every $\gamma_1 \colon F \Rightarrow G$ and $\gamma_2 \colon G \Rightarrow H$, $\st (\gamma_2)\st(\gamma_1) \Rrightarrow \st(\gamma_2\gamma_1)$ is given by a $2$-cell/modification in $\Gray(\st A, \st B)$, whose component at $a \in \st A$ is as depicted in:
		\[\begin{tikzcd}[ampersand replacement=\&]
			\& Fa \\
			\\
			Ga \\
			\\
			\& Ha
			\arrow["{(\gamma_1)_a}"', from=1-2, to=3-1]
			\arrow[""{name=0, anchor=center, inner sep=0}, "{(\gamma_2)_a (\gamma_1)_a}", from=1-2, to=5-2]
			\arrow["{(\gamma_2)_a}"', from=3-1, to=5-2]
			\arrow[between={0.2}{0.8}, Rightarrow, nfold, from=3-1, to=0]
		\end{tikzcd}\cdpunct[8pt]{,}\]
		\item[Identity constraint -]   for every $1_F\colon F \Rightarrow F \colon A \rightarrow B$, $1_{\st F} \Rrightarrow \st (1_F)$ is given by a $2$-cell/modification in $\Gray(\st A, \st B)$, whose component at $a \in \st A$ is as depicted in:
		\[\begin{tikzcd}[ampersand replacement=\&]
			Fa \\
			\\
			Fa
			\arrow[""{name=0, anchor=center, inner sep=0}, "{\varepsilon_{Fa}}"', curve={height=18pt}, from=1-1, to=3-1]
			\arrow[""{name=1, anchor=center, inner sep=0}, "{1_{Fa}}", curve={height=-18pt}, from=1-1, to=3-1]
			\arrow[between={0.2}{0.8}, Rightarrow, nfold, from=0, to=1]
		\end{tikzcd}\cdpunct[8pt]{.}\]
	\end{description}
\end{lemma}

\subsection{\tpdf{Coherence data of $\st$ as trifunctor}}
 
 To specify the functor $\st$ as a trifunctor, we also need to specify its identity and compositor adjoint equivalences. Note that we have not yet imposed any naturality condition on the isomorphisms $\Phi_{A,Y} \colon \Gray(\st A, Y) \rightarrow \BiCat(A, iY)$. Indeed, the optimal condition is to have $\Phi_{A,Y}$ be a strictly natural family of isomorphisms. Thus, we begin by noting that there are canonical icons 
 \[\begin{tikzcd}[ampersand replacement=\&]
 	\& {\BiCat(A,iY)} \&\& {\BiCat(A,iY)} \\
 	{\BiCat(B,iY)} \\
 	\& {\BiCat(C,iY)} \&\& {\BiCat(A,iY)}
 	\arrow["{F^*}", from=2-1, to=1-2]
 	\arrow[""{name=0, anchor=center, inner sep=0}, "{(FG)^*}"', from=3-2, to=1-2]
 	\arrow["{G^*}", from=3-2, to=2-1]
 	\arrow[""{name=1, anchor=center, inner sep=0}, "{1_{A}^*}"', curve={height=18pt}, from=3-4, to=1-4]
 	\arrow[""{name=2, anchor=center, inner sep=0}, "1", curve={height=-18pt}, from=3-4, to=1-4]
 	\arrow[between={0}{0.8}, Rightarrow, nfold, from=2-1, to=0]
 	\arrow[between={0.2}{0.8}, Rightarrow, nfold, from=2, to=1]
 \end{tikzcd}\cdpunct[8pt]{.}\]
 If we assume our family of $2$-isomorphisms $\Phi_{A,Y}$ is strictly natural, this yields similar icons on $\Gray(\st(-),Y)$ by transporting along $\Phi$. Finally, within the context of $\Gray$-enriched Yoneda lemma, this naturally motivates specific definitions of the identity and compositor adjoint equivalences (or icons). The identity adjoint equivalence is relatively straight forward. 
 
 \begin{definition}
 	The identity adjoint equivalence $\mathbf{\iota}$ is defined by the following composite of pseudo natural transformation:
 	\[\begin{tikzcd}[ampersand replacement=\&,column sep=small]
 		{*} \&\&\&\& \\
 		\\
 		\&\& {\BiCat(A,i \st A)} \&\& {\Gray(\st A, \st A)} \\
 		\\
 		{\BiCat(A,A)} \\
 		\& {\Phi^{-1}s_{A}} \&\& {\Phi^{-1}((S_A)_{*}1_A)} \& {(\Phi^{-1}(S_A)_{*})1_A} \\
 		{\iota :} \& {1_{\st A}} \&\&\& {\st_{A,A} 1_A}
 		\arrow[""{name=0, anchor=center, inner sep=0}, "{s_{A}}"{description}, from=1-1, to=3-3]
 		\arrow[""{name=1, anchor=center, inner sep=0}, "{1_{\st A}}", from=1-1, to=3-5]
 		\arrow["{1_A}"', from=1-1, to=5-1]
 		\arrow["{\Phi^{-1}}"{description}, from=3-3, to=3-5]
 		\arrow["{(s_{A})_{*}}"{description}, from=5-1, to=3-3]
 		\arrow[""{name=2, anchor=center, inner sep=0}, "{\st_{A,A}}"', from=5-1, to=3-5]
 		\arrow["{1 \star \mathrm{Coh}_{\BiCat}}", Rightarrow, nfold, from=6-2, to=6-4]
 		\arrow["{a^{.}}", equals, nfold, from=6-4, to=6-5]
 		\arrow[equals, nfold, from=6-5, to=7-5]
 		\arrow[equals, nfold, from=7-2, to=6-2]
 		\arrow[Rightarrow, nfold, from=7-2, to=7-5]
 		\arrow["1"', between={0.2}{0.8}, Rightarrow, nfold, from=1, to=3-3]
 		\arrow["{\mathrm{Coh}_{\BiCat}}"{description}, between={0.2}{0.8}, Rightarrow, nfold, from=0, to=5-1]
 		\arrow["1"', between={0.2}{0.8}, Rightarrow, nfold, from=3-3, to=2]
 	\end{tikzcd}\]
 \end{definition}

 We note that $\mathrm{Coh}_{\BiCat}$ is an invertible icon whose component at the unique $1$-cell $1_{*} \colon * \rightarrow *$ is depicted by the following diagram in $\BiCat(A, i \st A)$:
 \[\begin{tikzcd}[ampersand replacement=\&]
 	{s_{A}} \&\&\& {s_{A}} \\
 	\\
 	{s_{A}} \&\&\& {s_{A}}
 	\arrow[""{name=0, anchor=center, inner sep=0}, "{1_{s_{A}}}", Rightarrow, nfold, from=1-1, to=1-4]
 	\arrow[equals, nfold, from=1-1, to=3-1]
 	\arrow[equals, nfold, from=1-4, to=3-4]
 	\arrow[""{name=1, anchor=center, inner sep=0}, "{1_{s_{A}} \star 1_{A}}"', Rightarrow, nfold, from=3-1, to=3-4]
 	\arrow[between={0.2}{0.8}, Rightarrow, scaling nfold=3, from=1, to=0]
 \end{tikzcd}\cdpunct[8pt]{.}\]
 
 Thus, the composite invertible icon $1 \star \mathrm{Coh}_{\BiCat}$ is depicted in $\Gray(\st A, \st A)$ by
 \[\begin{tikzcd}[ampersand replacement=\&,cramped]
 	{\st(1_A)} \&\&\& {\st(1_A)} \&\& a \\
 	\\
 	{\st(1_A)} \&\&\& {\st(1_A)} \&\& a
 	\arrow[""{name=0, anchor=center, inner sep=0}, "{1_{\st 1_A}}", Rightarrow, nfold, from=1-1, to=1-4]
 	\arrow["{1_{\st 1_A}}"', equals, nfold, from=1-1, to=3-1]
 	\arrow["{1_{\st 1_A}}", equals, nfold, from=1-4, to=3-4]
 	\arrow[""{name=1, anchor=center, inner sep=0}, "{\varepsilon_{a}}", curve={height=-18pt}, from=1-6, to=3-6]
 	\arrow[""{name=2, anchor=center, inner sep=0}, "{1_{a}}"', curve={height=18pt}, from=1-6, to=3-6]
 	\arrow[""{name=3, anchor=center, inner sep=0}, "{\st(1_{1_A})}"', Rightarrow, nfold, from=3-1, to=3-4]
 	\arrow["{\mathbf{\iota}}"', between={0.2}{0.8}, Rightarrow, nfold, from=2, to=1]
 	\arrow["{\mathbf{\iota}}"', between={0.2}{0.8}, Rightarrow, scaling nfold=3, from=3, to=0]
 \end{tikzcd}\cdpunct[8pt]{.}\]
 
 Unfortunately, the compositor adjoint equivalence $\mathbf{\chi}$ is for more complicated. 
 
 \begin{definition}
 	The compositor adjoint equivalence $\mathbf{\chi}$ is defined by the following composite of pseudonatural transformations:
 	\[\begin{tikzcd}[ampersand replacement=\&,cramped]
 		{\mathbf{\chi} \colon} \& {\otimes_{\Gray} \circ (\st \times \st) } \&\& {\Phi^{-1} \circ (s_{A})^{*} \circ \otimes_{\BiCat} \circ (i \times i) \circ (\st \times \st)} \\
 		\&\&\& {\Phi^{-1} \circ \otimes_{\BiCat} \circ (1 \times (s_{A})^{*}) \circ (i \times i) \circ (\st \times \st)} \\
 		\&\&\& {\Phi^{-1} \circ \otimes_{\BiCat} \circ (1 \times (s_{B})_{*}) \circ (i \times 1) \circ (\st \times 1)} \\
 		\&\&\& {\Phi^{-1} \circ \otimes_{\BiCat}  \circ ((s_{B})^{*} \times 1)  \circ (i \times 1) \circ (\st \times 1)} \\
 		\&\&\& {\Phi^{-1} \circ \otimes_{\BiCat} \circ ((s_{C})_{*} \times 1) } \\
 		\& {\st \circ \otimes_{\BiCat}} \&\& {\Phi^{-1}  \circ (s_{C})_{*} \circ \otimes_{\BiCat}}
 		\arrow[Rightarrow, nfold, from=1-2, to=1-4]
 		\arrow[Rightarrow, nfold, from=1-2, to=6-2]
 		\arrow[Rightarrow, nfold, from=1-4, to=2-4]
 		\arrow[Rightarrow, nfold, from=2-4, to=3-4]
 		\arrow[Rightarrow, nfold, from=3-4, to=4-4]
 		\arrow[Rightarrow, nfold, from=4-4, to=5-4]
 		\arrow[Rightarrow, nfold, from=5-4, to=6-4]
 		\arrow[Rightarrow, nfold, from=6-4, to=6-2]
 	\end{tikzcd}\]
 	Here, we have suppresed the associators or the unitors of $\BiCat$ for clarity. 
 \end{definition}
 

The action of $\mathbf{\chi}$ becomes more apparent when we note that it is also an invertible icon whose components are represented by the following chain of modifications in $\Gray(\st A, \st C)$:
\[\begin{tikzcd}[ampersand replacement=\&,cramped,column sep=tiny]
	{\chi_{\alpha_{G},\alpha_{F}} :} \& {\st(\alpha_{G}) \star \st(\alpha_{F})} \&\& {G_1F_1a} \\
	{\Phi^{-1}((i\circ \st \alpha_{G} \star i \circ \st \alpha_{F}) \star 1_{s_{A}})} \& {\Phi^{-1}((i\circ (\st \alpha_{G} \star \st \alpha_{F})) \star 1_{s_{A}})} \\
	{\Phi^{-1}(i\circ \st \alpha_{G} \star (i \circ \st \alpha_{F} \star 1_{s_{A}}))} \& {\Phi^{-1}(i\circ \st \alpha_{G} \star (1_{s_{B}} \star \alpha_{F}))} \\
	{\Phi^{-1}((1_{s_{c}} \star \alpha_{G}) \star \alpha_{F})} \& {\Phi^{-1}((i\circ \st \alpha_{G} \star 1_{s_{B}}) \star \alpha_{F})} \&\& {G_2F_1a} \\
	{\Phi^{-1}(1_{s_{c}} \star (\alpha_{G} \star \alpha_{F}))} \& {\st(\alpha_{G} \star \alpha_{F})} \& {G_2F_2a}
	\arrow["{(\alpha_{G})_{F_1a}}", from=1-4, to=4-4]
	\arrow[""{name=0, anchor=center, inner sep=0}, "{G_2(\alpha_{F})_a(\alpha_{G})_{F_1a} }"', from=1-4, to=5-3]
	\arrow[Rightarrow, scaling nfold=3, no head, from=2-1, to=2-2]
	\arrow[Rightarrow, scaling nfold=3, no head, from=2-2, to=1-2]
	\arrow[Rightarrow, scaling nfold=3, from=3-1, to=2-1]
	\arrow[Rightarrow, scaling nfold=3, no head, from=3-2, to=3-1]
	\arrow[Rightarrow, scaling nfold=3, no head, from=4-1, to=4-2]
	\arrow[Rightarrow, scaling nfold=3, from=4-2, to=3-2]
	\arrow["{G_2(\alpha_{F})_a}", from=4-4, to=5-3]
	\arrow[Rightarrow, scaling nfold=3, from=5-1, to=4-1]
	\arrow[Rightarrow, scaling nfold=3, no head, from=5-2, to=5-1]
	\arrow["\chi", between={0.2}{0.8}, Rightarrow, nfold, from=0, to=4-4]
\end{tikzcd}\]

Both the tricategories $\BiCat$ and $\Gray$ can be considered $\BiCat_{2}$-enriched bicategories as deduced from \cite[Rem.~5.4]{gurski2013tricatcoh}, where $\BiCat_{2}$ is the strict $2$-category of bicategories, pseudofunctors, and Icons. Thus, the iconic structure of $\st$ implies that the trifunctor $\st$ can itself be considered as a $\BiCat_{2}$-enriched functor of $\BiCat_{2}$-enriched bicategories. In particular, any $3$-dimensional coherence data of the trifunctor $\st$ is merely induced by the coherence constraint cells. Hence, we have fully specified $\st$ as a trifunctor. 

\subsection{Encoding bicategorical pasting diagrams}

The primary usage of the strictification functor $\st$ is that it allows us to encode bicategorical pasting diagram in arbitrary bicategory $K$ using the canonical evaluation pseudofunctor $e\colon \st K \rightarrow K$. We motivate this with an example. 

Suppose we wish to specify the data of an adjunction in a bicategory $K$. Algebraically, it can be expressed using the following data:
\[\begin{tikzcd}[ampersand replacement=\&,cramped]
	\& {f \colon a \rightarrow b} \&\&\& {g \colon b \rightarrow a} \& \\
	\& {\eta  \colon 1_{a} \Rightarrow gf} \&\&\& {\epsilon  \colon fg \Rightarrow 1_{b}} \\
	\& {f1_{a}} \& {f(gf)} \& {(fg)f} \& {1_{b}f} \\
	f \&\&\&\&\& f \\
	\& {1_{a}g} \& {(gf)g} \& {g(fg)} \& {g1_{b}} \\
	g \&\&\&\&\& g
	\arrow["{1\eta}", Rightarrow, nfold, from=3-2, to=3-3]
	\arrow["{a^{-1}}", Rightarrow, nfold, from=3-3, to=3-4]
	\arrow["\epsilon1", Rightarrow, nfold, from=3-4, to=3-5]
	\arrow["l", Rightarrow, nfold, from=3-5, to=4-6]
	\arrow["{r^{-1}}", Rightarrow, nfold, from=4-1, to=3-2]
	\arrow["{1_{f}}"', Rightarrow, nfold, from=4-1, to=4-6]
	\arrow["\eta1", Rightarrow, nfold, from=5-2, to=5-3]
	\arrow["a", Rightarrow, nfold, from=5-3, to=5-4]
	\arrow["{1\epsilon}", Rightarrow, nfold, from=5-4, to=5-5]
	\arrow["r", Rightarrow, nfold, from=5-5, to=6-6]
	\arrow["{l^{-1}}", Rightarrow, nfold, from=6-1, to=5-2]
	\arrow["{1_{g}}"', Rightarrow, nfold, from=6-1, to=6-6]
\end{tikzcd}\cdpunct[8pt]{.}\]
In a strict $2$-category such as $\nCat{1}$, we could express one of the triangle identity using the following pasting diagram:
\[\begin{tikzcd}[ampersand replacement=\&,cramped]
	\& B \&\& B \&\&\& A \&\& A \\
	\&\&\&\& \equiv \\
	A \&\& A \&\&\& B \&\& B
	\arrow[""{name=0, anchor=center, inner sep=0}, "{1_{B}}", from=1-2, to=1-4]
	\arrow["G"{description}, from=1-2, to=3-3]
	\arrow["{1_{A}}", from=1-7, to=1-9]
	\arrow["F", from=3-1, to=1-2]
	\arrow[""{name=1, anchor=center, inner sep=0}, "{1_{A}}"', from=3-1, to=3-3]
	\arrow["F"', from=3-3, to=1-4]
	\arrow["F", from=3-6, to=1-7]
	\arrow["F"', from=3-6, to=1-9]
	\arrow["{1_{B}}"', from=3-6, to=3-8]
	\arrow["F"', from=3-8, to=1-9]
	\arrow["\eta"{description}, between={0.2}{0.8}, Rightarrow, nfold, from=1, to=1-2]
	\arrow["\epsilon"{description}, between={0.2}{0.8}, Rightarrow, nfold, from=3-3, to=0]
\end{tikzcd}\cdpunct[8pt]{.}\]
This motivates us to consider encoding such identities inside of $\st K$ instead of $K$ itself, where such identities are expressible using pasting diagram. To this end, we note the following relations under the evaluation map $e$ expressing both the cells $\eta$ and $\epsilon$ inside $\st K$:
\[\begin{tikzcd}[ampersand replacement=\&,cramped]
	\& {\st K} \&\&\&\& K \& \\
	\& b \\
	a \&\& a \&\& a \&\& a \\
	\& a \\
	b \&\& b \&\& b \&\& b
	\arrow["e", from=1-2, to=1-6]
	\arrow["g", from=2-2, to=3-3]
	\arrow["f", from=3-1, to=2-2]
	\arrow[""{name=0, anchor=center, inner sep=0}, "{\varepsilon_{a}}"', from=3-1, to=3-3]
	\arrow[""{name=1, anchor=center, inner sep=0}, "{1_{a}}"', curve={height=12pt}, from=3-5, to=3-7]
	\arrow[""{name=2, anchor=center, inner sep=0}, "gf", curve={height=-12pt}, from=3-5, to=3-7]
	\arrow["f", from=4-2, to=5-3]
	\arrow["g", from=5-1, to=4-2]
	\arrow[""{name=3, anchor=center, inner sep=0}, "{\varepsilon_{b}}"', from=5-1, to=5-3]
	\arrow[""{name=4, anchor=center, inner sep=0}, "{1_{b}}"', curve={height=12pt}, from=5-5, to=5-7]
	\arrow[""{name=5, anchor=center, inner sep=0}, "fg", curve={height=-12pt}, from=5-5, to=5-7]
	\arrow["\eta"', between={0.2}{0.8}, Rightarrow, nfold, from=0, to=2-2]
	\arrow["\eta"', between={0.2}{0.8}, Rightarrow, nfold, from=1, to=2]
	\arrow["\epsilon", between={0.2}{0.8}, Rightarrow, nfold, from=4-2, to=3]
	\arrow["\epsilon", between={0.2}{0.8}, Rightarrow, nfold, from=5, to=4]
\end{tikzcd}\cdpunct[8pt]{.}\]

In general, let $\gamma$ be some $2$-cell between $e(\vec{f},b)$ and $e(\vec{g},b')$, where $b$ and $b'$ signify any arbitrary choice of bracketing in the evaluation. This $2$-cell can be sent inside $\st K$ using the embedding $s_{K} \colon K \rightarrow \st K$. We do this by mapping $\gamma \in K$ to the $2$-cell $\vec{f} \xRightarrow{\mathrm{Coh}} e(\vec{f},b) \xRightarrow{\gamma} e(\vec{g},b') \xRightarrow{\mathrm{Coh}^{-1}} \vec{g}$ in $\st K$, whenever $\vec{f}$ and $\vec{g}$ contains no identity $1$-cells. If they do, we also remove them using coherence.

With this idea, we now embed one of the triangle identity of $(f,g,\eta,\epsilon)$ inside $\st k$: 
\[\begin{tikzcd}[ampersand replacement=\&]
	\&\&\& f \&\&\&\&\& \\
	a \&\& a \&\&\& a \\
	\&\&\&\& \equiv \\
	\& b \&\& b \&\&\&\&\& b \\
	f
	\arrow[between={0.2}{0.8}, Rightarrow, nfold, from=1-4, to=2-3]
	\arrow[from=1-4, to=4-4]
	\arrow[no head, from=2-1, to=1-4]
	\arrow[""{name=0, anchor=center, inner sep=0}, "{\varepsilon_{a}}"{description}, from=2-1, to=2-3]
	\arrow["f"{description}, from=2-1, to=4-2]
	\arrow[no head, from=2-1, to=5-1]
	\arrow["f"{description}, from=2-3, to=4-4]
	\arrow[""{name=1, anchor=center, inner sep=0}, "f", curve={height=-18pt}, from=2-6, to=4-9]
	\arrow[""{name=2, anchor=center, inner sep=0}, "f"', curve={height=18pt}, from=2-6, to=4-9]
	\arrow["g"{description}, from=4-2, to=2-3]
	\arrow[""{name=3, anchor=center, inner sep=0}, "{\varepsilon_{b}}"{description}, from=4-2, to=4-4]
	\arrow[between={0.2}{0.8}, Rightarrow, nfold, from=4-2, to=5-1]
	\arrow[from=5-1, to=4-4]
	\arrow["\eta", between={0.2}{0.8}, Rightarrow, nfold, from=0, to=4-2]
	\arrow["\epsilon", between={0.2}{0.8}, Rightarrow, nfold, from=2-3, to=3]
	\arrow["{1_{f}}", between={0.2}{0.8}, Rightarrow, nfold, from=1, to=2]
\end{tikzcd}\cdpunct[8pt]{.}\]
Of course, $\varepsilon_{a}$ and $\varepsilon_{b}$ are strict identities inside $\st K$. Thus, we could equivalently write
\[\begin{tikzcd}[ampersand replacement=\&]
	a \&\& a \&\&\& a \&\& a \& \\
	\&\&\&\& \equiv \\
	\& b \&\& b \&\&\& b \&\& b
	\arrow[""{name=0, anchor=center, inner sep=0}, "{\varepsilon_{a}}", from=1-1, to=1-3]
	\arrow["f"', from=1-1, to=3-2]
	\arrow["f", from=1-3, to=3-4]
	\arrow["{\varepsilon_{a}}", from=1-6, to=1-8]
	\arrow["f"', from=1-6, to=3-7]
	\arrow[""{name=1, anchor=center, inner sep=0}, "f"{description}, from=1-6, to=3-9]
	\arrow["f", from=1-8, to=3-9]
	\arrow["g"{description}, from=3-2, to=1-3]
	\arrow[""{name=2, anchor=center, inner sep=0}, "{\varepsilon_{b}}"', from=3-2, to=3-4]
	\arrow["{\varepsilon_{b}}"', from=3-7, to=3-9]
	\arrow["\eta", between={0.2}{0.8}, Rightarrow, nfold, from=0, to=3-2]
	\arrow["\epsilon", between={0.2}{0.8}, Rightarrow, nfold, from=1-3, to=2]
	\arrow[between={0.2}{0.8}, Rightarrow, nfold, from=1, to=3-7]
	\arrow[between={0.2}{0.8}, Rightarrow, nfold, from=1-8, to=1]
\end{tikzcd}\cdpunct[8pt]{}\]
to express the same triangle identity. Finally, Using the evaluation map $e \colon \st B \rightarrow B$, it is relatively straightforward to verify that this pasting diagram equality correspond to the following equality of $2$-cells inside $B$:
\[\begin{tikzcd}[ampersand replacement=\&]
	\& {f1_{a}} \& {f(gf)} \&\& {(fg)f} \& {1_{b}f} \& \\
	f \&\&\& \vequiv \&\&\& f \\
	f \&\&\&\&\&\& f
	\arrow["{1\eta}", Rightarrow, nfold, from=1-2, to=1-3]
	\arrow["{a^{-1}}", Rightarrow, nfold, from=1-3, to=1-5]
	\arrow["\epsilon1", Rightarrow, nfold, from=1-5, to=1-6]
	\arrow["l", Rightarrow, nfold, from=1-6, to=2-7]
	\arrow["{r^{-1}}", Rightarrow, nfold, from=2-1, to=1-2]
	\arrow["{1_{f}}"', Rightarrow, nfold, from=3-1, to=3-7]
\end{tikzcd}\cdpunct[8pt]{.}\]
Therefore, as a general note, whenever we write bicategorical pasting diagram inside a given bicategory $K$, we always mean to read such diagrams inside $\st K$ with the use of the evaluation functor $e \colon\st K \rightarrow K$. A more formal treatment of this idea is given using pasting schemes. See \cite[Ch.~3]{yau2cat2021} for more details. 
	\newpage
	
	\section{Coherence diagrams for pseudomonoid objects} \label{appsection:mondiag}
	Throughout this section, let $\cV$ be a \emph{monoidal bicategory}. In particular, we can choose $\cV$ to be a one-object \emph{algebraic tricategory} in the sense of \cite{gurski2006tricat, gurski2013tricatcoh}, and shall denote the relevant data of the monoidal bicategory $\cV$ as $(\cV, \otimes, I, a_{\cV}, a_{\cV}^{-1},\lambda_{\cV},\lambda_{\cV}^{-1},\rho_{\cV},\rho_{\cV}^{-1})$, omitting any labeling of the $2$-cell data. Thus, for simplicity, one may also consider $\cV$ as a $\nCat{1}$-enriched monoidal category. Under this consideration, all of the forthcoming coherence diagram are supposed to be read inside the cubical tricategory $\st\cV$. In particular, we must incorporate the fact that $\otimes$ is a pseudofunctor with non-trivial coherence constraints. Thus, for instance, $a \otimes 1$ (when seen in $\st\cV$) is supposed to be the image of $a$ combined with appropriate coherence data under the pseudofunctor $\cV \rightarrow \cV \times * \xrightarrow{1 \times M} \cV \times \cV \xrightarrow{\otimes} \cV$. Some similar or same diagrams are also available in \cite{yau2024bimodcats} for the case of $\nCat{1}$-enriched categories. More generally, some of these diagrams can be seen as one object versions of \cite[Defs.~3.2, 3.3, and 3.4]{gurski2014itericon}.

\subsection{Pentagon identity}\label{app:pentagon}

\[\begin{tikzcd}[ampersand replacement=\&,column sep=scriptsize]
	{M \otimes (M \otimes (M \otimes M))} \&\&\& {M \otimes (M \otimes M)} \\
	\& {M \otimes ((M \otimes M) \otimes M)} \& {M \otimes (M \otimes M)} \& {M \otimes M} \\
	{(M \otimes M) \otimes (M \otimes M)} \\
	\& {(M \otimes (M \otimes M)) \otimes M} \& {(M \otimes M) \otimes M} \\
	{((M \otimes M) \otimes M) \otimes M} \& {(M \otimes M) \otimes M} \& {M \otimes M} \& M
	\arrow[""{name=0, anchor=center, inner sep=0}, "{1 \otimes (1 \otimes \mu)}", from=1-1, to=1-4]
	\arrow["{1 \otimes \mu}"', from=1-4, to=2-4]
	\arrow["{1 \otimes a_{\cV}}", from=2-2, to=1-1]
	\arrow["{1 \otimes (\mu \otimes 1)}"', from=2-2, to=2-3]
	\arrow[""{name=1, anchor=center, inner sep=0}, "{1 \otimes \mu}", from=2-3, to=2-4]
	\arrow["\mu"', from=2-4, to=5-4]
	\arrow["{a_{\cV}}", from=3-1, to=1-1]
	\arrow[""{name=2, anchor=center, inner sep=0}, "{a_{\cV}}"', from=4-2, to=2-2]
	\arrow[""{name=3, anchor=center, inner sep=0}, "{(1 \otimes \mu) \otimes 1}", from=4-2, to=4-3]
	\arrow[between={0.2}{0.8}, equals, nfold, from=4-3, to=2-2]
	\arrow["{a_{\cV}}"', from=4-3, to=2-3]
	\arrow["{\mu \otimes 1}", from=4-3, to=5-3]
	\arrow["{a_{\cV}}", from=5-1, to=3-1]
	\arrow["{a_{\cV} \otimes 1}", from=5-1, to=4-2]
	\arrow["{(\mu \otimes 1) \otimes 1}"', from=5-1, to=5-2]
	\arrow["{\mu \otimes 1}"', from=5-2, to=5-3]
	\arrow[""{name=4, anchor=center, inner sep=0}, "\mu"', from=5-3, to=5-4]
	\arrow["{1 \otimes a}", between={0.2}{0.8}, Rightarrow, nfold, from=2-3, to=0]
	\arrow[between={0.2}{1}, equals, nfold, from=2, to=3-1]
	\arrow["{a \otimes 1}"', between={0.2}{0.8}, Rightarrow, nfold, from=5-2, to=3]
	\arrow["a"', between={0.2}{0.8}, Rightarrow, nfold, from=4, to=1]
\end{tikzcd}\]
\[ \vequiv \]
\[\begin{tikzcd}[ampersand replacement=\&]
	{M \otimes (M \otimes (M \otimes M))} \&\& {M \otimes (M \otimes M)} \&\& {M \otimes M} \\
	\&\& {(M \otimes M) \otimes M} \\
	{(M \otimes M) \otimes (M \otimes M)} \&\&\& {M \otimes M} \\
	\&\& {M \otimes (M \otimes M)} \\
	{((M \otimes M) \otimes M) \otimes M} \&\& {(M \otimes M) \otimes M} \& {M \otimes M} \& M
	\arrow["{1 \otimes (1 \otimes \mu)}", from=1-1, to=1-3]
	\arrow[""{name=0, anchor=center, inner sep=0}, "{1 \otimes \mu}", from=1-3, to=1-5]
	\arrow["\mu"', from=1-5, to=5-5]
	\arrow[between={0.2}{0.8}, equals, from=2-3, to=1-1]
	\arrow["{a_{\cV}}"', from=2-3, to=1-3]
	\arrow["{\mu \otimes 1}", from=2-3, to=3-4]
	\arrow["{a_{\cV}}"', from=3-1, to=1-1]
	\arrow["{(1 \otimes 1 ) \otimes \mu}", from=3-1, to=2-3]
	\arrow["{\mu \otimes (1 \otimes 1)}"', from=3-1, to=4-3]
	\arrow["\mu", from=3-4, to=5-5]
	\arrow[between={0.2}{0.8}, equals, from=4-3, to=2-3]
	\arrow["{1 \otimes \mu}"', from=4-3, to=3-4]
	\arrow[""{name=1, anchor=center, inner sep=0}, "{a_{\cV}}"', from=5-1, to=3-1]
	\arrow["{(\mu \otimes 1) \otimes 1}"', from=5-1, to=5-3]
	\arrow["{a_{\cV}}", from=5-3, to=4-3]
	\arrow["{\mu \otimes 1}"', from=5-3, to=5-4]
	\arrow["a"', between={0.2}{0.8}, Rightarrow, from=5-4, to=3-4]
	\arrow["\mu"', from=5-4, to=5-5]
	\arrow["a", between={0.2}{0.8}, Rightarrow, from=3-4, to=0]
	\arrow[between={0.2}{0.8}, equals, from=5-3, to=1]
\end{tikzcd}\]

\newpage 

\subsection{Triangle identity}\label{app:triangle}

\[\begin{tikzcd}[ampersand replacement=\&,cramped]
	{M \otimes M} \&\&\&\& M \\
	\\
	\& {M \otimes (I \otimes M)} \&\& {M \otimes (M \otimes M)} \& {} \\
	\\
	{M \otimes M} \&\&\&\& {M \otimes M} \\
	\\
	\& {(M \otimes I) \otimes M} \&\& {(M \otimes M) \otimes M}
	\arrow[""{name=0, anchor=center, inner sep=0}, "\mu", from=1-1, to=1-5]
	\arrow[""{name=1, anchor=center, inner sep=0}, "{1 \otimes \lambda_{\cV}}", from=3-2, to=1-1]
	\arrow["{1 \otimes (e \otimes 1)}"', from=3-2, to=3-4]
	\arrow["{1 \otimes \mu}"', from=3-4, to=1-1]
	\arrow[""{name=2, anchor=center, inner sep=0}, "{\varepsilon_{M \otimes M}}", from=5-1, to=1-1]
	\arrow["{\rho_{\cV}^{-1} \otimes 1}"', from=5-1, to=7-2]
	\arrow["\mu"', from=5-5, to=1-5]
	\arrow[""{name=3, anchor=center, inner sep=0}, "{a_{\cV}}"', from=7-2, to=3-2]
	\arrow["{(1 \otimes e) \otimes 1}"', from=7-2, to=7-4]
	\arrow[between={0.2}{0.8}, equals, nfold, from=7-4, to=3-2]
	\arrow["{a_{\cV}}", from=7-4, to=3-4]
	\arrow["{\mu \otimes 1}"', from=7-4, to=5-5]
	\arrow["{1 \otimes \lambda}", between={0.2}{0.8}, Rightarrow, nfold, from=3-4, to=1]
	\arrow["a"', between={0.2}{0.8}, Rightarrow, nfold, from=3-5, to=0]
	\arrow[between={0.2}{0.8}, equals, nfold, from=3, to=2]
\end{tikzcd}\]
\[\vequiv\]
\[\begin{tikzcd}[ampersand replacement=\&,cramped]
	{M \otimes M} \&\&\&\& M \\
	\\
	\\
	\\
	{M \otimes M} \&\&\&\& {M \otimes M} \\
	\\
	\& {(M \otimes I) \otimes M} \&\& {(M \otimes M) \otimes M}
	\arrow["\mu", from=1-1, to=1-5]
	\arrow["{\varepsilon_{M \otimes M}}", from=5-1, to=1-1]
	\arrow[""{name=0, anchor=center, inner sep=0}, "\mu"{description}, from=5-1, to=1-5]
	\arrow[""{name=1, anchor=center, inner sep=0}, "{\varepsilon_{M \otimes M}}"', from=5-1, to=5-5]
	\arrow["{\rho_{\cV}^{-1} \otimes 1}"', from=5-1, to=7-2]
	\arrow["\mu"', from=5-5, to=1-5]
	\arrow[""{name=2, anchor=center, inner sep=0}, "{\rho_{\cV} \otimes 1}"{description}, from=7-2, to=5-5]
	\arrow["{(1 \otimes e) \otimes 1}"', from=7-2, to=7-4]
	\arrow["{\mu \otimes 1}"', from=7-4, to=5-5]
	\arrow[between={0.2}{0.8}, equals, nfold, from=0, to=1-1]
	\arrow[between={0.2}{0.8}, equals, nfold, from=5-5, to=0]
	\arrow[between={0.2}{0.8}, equals, nfold, from=7-2, to=1]
	\arrow["{\rho \otimes 1}", between={0.2}{0.8}, Rightarrow, nfold, from=7-4, to=2]
\end{tikzcd}\]

\newpage

\subsection{Hexagon identities}\label{app:hexagon}

We assume that we have a suitable definition of a braided monoidal bicategory in sense of being a one object tricategory equipped with some coherent braiding adjoint equivalence. In particular, $\cV$ is a braided monoidal bicategory if it is a monoidal bicategory equipped with a family of adjoint equivalences $b_{x,y} : X \otimes Y \leftrightarrow Y \otimes X : b_{x,y}^{-1}$ that behaves as the braiding. The details are ultimately irrelevant for us but can be found carefully constructed in \cite{gurski2011monbicat}. Other relevant sources include \cite{yau2cat2021,pries2011monbicat}. 
\[\begin{tikzcd}[ampersand replacement=\&]
	{(M \otimes M)  \otimes M} \& {M \otimes (M \otimes M)} \&\& {(M \otimes M) \otimes M} \& {M \otimes (M \otimes M)} \\
	\\
	{M \otimes M} \& {M \otimes M} \&\& {M \otimes M} \& {M \otimes M} \\
	\\
	\&\& M
	\arrow["{a_{\cV}}", from=1-1, to=1-2]
	\arrow["{\mu \otimes 1}"', from=1-1, to=3-1]
	\arrow["{b_{M,M \otimes M}}", from=1-2, to=1-4]
	\arrow["{1 \otimes \mu}"', from=1-2, to=3-2]
	\arrow["{a_{\cV}}", from=1-4, to=1-5]
	\arrow["{\mu \otimes 1}", from=1-4, to=3-4]
	\arrow["{1 \otimes \mu}", from=1-5, to=3-5]
	\arrow["a", between={0.2}{0.8}, Rightarrow, nfold, from=3-1, to=3-2]
	\arrow["\mu"', from=3-1, to=5-3]
	\arrow[between={0.2}{0.8}, equals, nfold, from=3-2, to=1-4]
	\arrow["{b_{M,M}}", from=3-2, to=3-4]
	\arrow[""{name=0, anchor=center, inner sep=0}, "\mu"', from=3-2, to=5-3]
	\arrow["a", between={0.2}{0.8}, Rightarrow, nfold, from=3-4, to=3-5]
	\arrow["\mu", from=3-4, to=5-3]
	\arrow["\mu", from=3-5, to=5-3]
	\arrow["b", between={0.2}{0.8}, Rightarrow, nfold, from=0, to=3-4]
\end{tikzcd}\]
\[\vequiv\]
\[\begin{tikzcd}[ampersand replacement=\&]
	\& {M \otimes (M \otimes M)} \&\& {(M \otimes M) \otimes M} \& \\
	{(M \otimes M)  \otimes M} \&\&\&\& {M \otimes (M \otimes M)} \\
	\& {(M \otimes M) \otimes M} \&\& {M \otimes (M \otimes M)} \\
	\\
	{M \otimes M} \&\&\&\& {M \otimes M} \\
	\\
	\&\& M
	\arrow[""{name=0, anchor=center, inner sep=0}, "{b_{M,M \otimes M}}", from=1-2, to=1-4]
	\arrow["{a_{\cV}}", from=1-4, to=2-5]
	\arrow["{a_{\cV}}", from=2-1, to=1-2]
	\arrow["{b_{M,M} \otimes 1}"', from=2-1, to=3-2]
	\arrow[""{name=1, anchor=center, inner sep=0}, "{\mu \otimes 1}"', from=2-1, to=5-1]
	\arrow[""{name=2, anchor=center, inner sep=0}, "{1 \otimes \mu}", from=2-5, to=5-5]
	\arrow[""{name=3, anchor=center, inner sep=0}, "{a_{\cV}}", from=3-2, to=3-4]
	\arrow[""{name=4, anchor=center, inner sep=0}, "{\mu \otimes 1}", from=3-2, to=5-1]
	\arrow["{1 \otimes b_{M,M}}"', from=3-4, to=2-5]
	\arrow[""{name=5, anchor=center, inner sep=0}, "{1 \otimes \mu}"', from=3-4, to=5-5]
	\arrow["a"', between={0.2}{0.8}, Rightarrow, nfold, from=5-1, to=5-5]
	\arrow["\mu"', from=5-1, to=7-3]
	\arrow["\mu", from=5-5, to=7-3]
	\arrow["{b \otimes 1}", between={0.2}{0.8}, Rightarrow, nfold, from=1, to=4]
	\arrow[between={0.2}{0.8}, equals, nfold, from=3, to=0]
	\arrow["{1 \otimes b}", between={0.2}{0.8}, Rightarrow, nfold, from=5, to=2]
\end{tikzcd}\]

\[\begin{tikzcd}[ampersand replacement=\&]
	{M \otimes (M \otimes M)} \& {(M \otimes M) \otimes M} \&\& {M \otimes (M \otimes M)} \& {(M \otimes M) \otimes M} \\
	\\
	{M \otimes M} \& {M \otimes M} \&\& {M \otimes M} \& {M \otimes M} \\
	\\
	\&\& M
	\arrow["{a_{\cV}^{-1}}", from=1-1, to=1-2]
	\arrow["{1 \otimes \mu}"', from=1-1, to=3-1]
	\arrow["{b_{M \otimes M,M}}", from=1-2, to=1-4]
	\arrow["{\mu \otimes 1}"', from=1-2, to=3-2]
	\arrow["{a_{\cV}^{-1}}", from=1-4, to=1-5]
	\arrow["{1 \otimes \mu}", from=1-4, to=3-4]
	\arrow["{\mu \otimes 1}", from=1-5, to=3-5]
	\arrow["{a^{-1}}", between={0.2}{0.8}, Rightarrow, nfold, from=3-1, to=3-2]
	\arrow["\mu"', from=3-1, to=5-3]
	\arrow[between={0.2}{0.8}, equals, nfold, from=3-2, to=1-4]
	\arrow["{b_{M,M}}", from=3-2, to=3-4]
	\arrow[""{name=0, anchor=center, inner sep=0}, "\mu"', from=3-2, to=5-3]
	\arrow["{a^{-1}}", between={0.2}{0.8}, Rightarrow, nfold, from=3-4, to=3-5]
	\arrow["\mu", from=3-4, to=5-3]
	\arrow["\mu", from=3-5, to=5-3]
	\arrow["b", between={0.2}{0.8}, Rightarrow, nfold, from=0, to=3-4]
\end{tikzcd}\]
\[\vequiv\]
\[\begin{tikzcd}[ampersand replacement=\&]
	\& {(M \otimes M) \otimes M} \&\& {M \otimes (M \otimes M)} \& \\
	{M \otimes (M \otimes M)} \&\&\&\& {(M \otimes M) \otimes M} \\
	\& {M \otimes (M \otimes M)} \&\& {(M \otimes M) \otimes M} \\
	\\
	{M \otimes M} \&\&\&\& {M \otimes M} \\
	\\
	\&\& M
	\arrow[""{name=0, anchor=center, inner sep=0}, "{b_{M,M \otimes M}}", from=1-2, to=1-4]
	\arrow["{a_{\cV}^{-1}}", from=1-4, to=2-5]
	\arrow["{a_{\cV}^{-1}}", from=2-1, to=1-2]
	\arrow["{1 \otimes b_{M,M}}"', from=2-1, to=3-2]
	\arrow[""{name=1, anchor=center, inner sep=0}, "{1 \otimes \mu}"', from=2-1, to=5-1]
	\arrow[""{name=2, anchor=center, inner sep=0}, "{\mu \otimes 1}", from=2-5, to=5-5]
	\arrow[""{name=3, anchor=center, inner sep=0}, "{a_{\cV}^{-1}}", from=3-2, to=3-4]
	\arrow[""{name=4, anchor=center, inner sep=0}, "{1 \otimes \mu}", from=3-2, to=5-1]
	\arrow["{b_{M,M} \otimes 1}"', from=3-4, to=2-5]
	\arrow[""{name=5, anchor=center, inner sep=0}, "{\mu \otimes 1}"', from=3-4, to=5-5]
	\arrow["{a^{-1}}"', between={0.2}{0.8}, Rightarrow, nfold, from=5-1, to=5-5]
	\arrow["\mu"', from=5-1, to=7-3]
	\arrow["\mu", from=5-5, to=7-3]
	\arrow["{1 \otimes b}", between={0.2}{0.8}, Rightarrow, nfold, from=1, to=4]
	\arrow[between={0.2}{0.8}, equals, nfold, from=3, to=0]
	\arrow["{b \otimes 1}", between={0.2}{0.8}, Rightarrow, nfold, from=5, to=2]
\end{tikzcd}\]

\newpage 

\subsection{Morphism hexagon axiom}\label{app:morphhexagon}

\[\begin{tikzcd}[ampersand replacement=\&]
	\& {(M \otimes M) \otimes M} \&\& {(M' \otimes M') \otimes M'} \\
	{M \otimes (M \otimes M)} \\
	\& {M \otimes M} \&\& {M ' \otimes M'} \\
	{M \otimes M} \\
	\& M \&\& {M'}
	\arrow["{(f \otimes f) \otimes f}", from=1-2, to=1-4]
	\arrow["{a_{\cV}}"', from=1-2, to=2-1]
	\arrow["{\mu \otimes 1}"{description}, from=1-2, to=3-2]
	\arrow["{\chi \otimes 1}"{description}, between={0.2}{0.8}, Rightarrow, nfold, from=1-4, to=3-2]
	\arrow["{\mu' \otimes 1}", from=1-4, to=3-4]
	\arrow[""{name=0, anchor=center, inner sep=0}, "{1 \otimes \mu}"', from=2-1, to=4-1]
	\arrow["{f \otimes f}"', from=3-2, to=3-4]
	\arrow["\mu"{description}, from=3-2, to=5-2]
	\arrow["\chi"{description}, between={0.2}{0.8}, Rightarrow, nfold, from=3-4, to=5-2]
	\arrow["{\mu'}", from=3-4, to=5-4]
	\arrow["\mu"', from=4-1, to=5-2]
	\arrow["f"', from=5-2, to=5-4]
	\arrow["a"{description}, between={0.2}{0.8}, Rightarrow, nfold, from=3-2, to=0]
\end{tikzcd}\]
\[\vequiv\]
\[\begin{tikzcd}[ampersand replacement=\&]
	\& {(M \otimes M) \otimes M} \&\& {(M' \otimes M') \otimes M'} \\
	{M \otimes (M \otimes M)} \&\& {M' \otimes (M' \otimes M')} \\
	\&\&\& {M ' \otimes M'} \\
	{M \otimes M} \&\& {M' \otimes M'} \\
	\& M \&\& {M'}
	\arrow[""{name=0, anchor=center, inner sep=0}, "{(f \otimes f) \otimes f}", from=1-2, to=1-4]
	\arrow["{a_{\cV}}"', from=1-2, to=2-1]
	\arrow["{a_{\cV}}"{description}, from=1-4, to=2-3]
	\arrow["{\mu' \otimes 1}", from=1-4, to=3-4]
	\arrow[""{name=1, anchor=center, inner sep=0}, "{f \otimes (f \otimes f)}"{description}, from=2-1, to=2-3]
	\arrow["{1 \otimes \mu}"', from=2-1, to=4-1]
	\arrow["{1 \otimes \chi}"{description}, between={0.2}{0.8}, Rightarrow, nfold, from=2-3, to=4-1]
	\arrow[""{name=2, anchor=center, inner sep=0}, "{1 \otimes \mu'}"{description}, from=2-3, to=4-3]
	\arrow["{\mu'}", from=3-4, to=5-4]
	\arrow["{f  \otimes f}"{description}, from=4-1, to=4-3]
	\arrow["\mu"', from=4-1, to=5-2]
	\arrow["\chi"{description}, between={0.2}{0.8}, Rightarrow, nfold, from=4-3, to=5-2]
	\arrow["{\mu'}"{description}, from=4-3, to=5-4]
	\arrow["f"', from=5-2, to=5-4]
	\arrow[between={0.2}{0.8}, equals, nfold, from=0, to=1]
	\arrow["{a'}"{description}, between={0.2}{0.8}, Rightarrow, nfold, from=3-4, to=2]
\end{tikzcd}\]

\newpage

\subsection{Morphism left square} \label{app:morphlsq}
\[\begin{tikzcd}[ampersand replacement=\&]
	{I \otimes M} \&\&\& {I \otimes M'} \& \\
	\&\& {M' \otimes M} \\
	\\
	\& {M \otimes M} \&\&\& {M' \otimes M'} \\
	\\
	\\
	M \&\&\& {M'}
	\arrow["{1 \otimes f}", from=1-1, to=1-4]
	\arrow[""{name=0, anchor=center, inner sep=0}, "{e' \otimes 1}"{description}, from=1-1, to=2-3]
	\arrow["{e \otimes 1}"{description}, from=1-1, to=4-2]
	\arrow[""{name=1, anchor=center, inner sep=0}, "{\lambda_{\cV}}"', from=1-1, to=7-1]
	\arrow[between={0.2}{0.8}, equals, nfold, from=1-4, to=2-3]
	\arrow["{e' \otimes 1}", from=1-4, to=4-5]
	\arrow[""{name=2, anchor=center, inner sep=0}, "{1 \otimes f}"{description}, from=2-3, to=4-5]
	\arrow["{f \otimes 1}"{description}, from=4-2, to=2-3]
	\arrow[""{name=3, anchor=center, inner sep=0}, "{f \otimes f}"{description}, from=4-2, to=4-5]
	\arrow["\mu"{description}, from=4-2, to=7-1]
	\arrow["\chi"{description}, between={0.2}{0.8}, Rightarrow, nfold, from=4-5, to=7-1]
	\arrow["{\mu'}", from=4-5, to=7-4]
	\arrow["f"', from=7-1, to=7-4]
	\arrow["{\iota \otimes 1}"{description}, between={0.2}{0.8}, Rightarrow, nfold, from=0, to=4-2]
	\arrow[between={0.2}{0.8}, equals, nfold, from=2, to=3]
	\arrow["\lambda", between={0.2}{0.8}, Rightarrow, nfold, from=4-2, to=1]
\end{tikzcd}\]
\[\vequiv\]
\[\begin{tikzcd}[ampersand replacement=\&]
	{I \otimes M} \&\&\& {I \otimes M'} \& \\
	\\
	\\
	\&\&\&\& {M' \otimes M'} \\
	\\
	\\
	M \&\&\& {M'}
	\arrow["{1 \otimes f}", from=1-1, to=1-4]
	\arrow["{\lambda_{\cV}}"', from=1-1, to=7-1]
	\arrow["{e' \otimes 1}", from=1-4, to=4-5]
	\arrow[between={0.2}{0.8}, equals, nfold, from=1-4, to=7-1]
	\arrow[""{name=0, anchor=center, inner sep=0}, "{\lambda_{\cV}}"{description}, from=1-4, to=7-4]
	\arrow["{\mu'}", from=4-5, to=7-4]
	\arrow["f"', from=7-1, to=7-4]
	\arrow["{\lambda'}", between={0.2}{0.8}, Rightarrow, nfold, from=4-5, to=0]
\end{tikzcd}\]

\newpage

\subsection{Morphism right square} \label{app:morphrsq}

\[\begin{tikzcd}[ampersand replacement=\&]
	{M \otimes I} \&\&\& {M' \otimes I} \& \\
	\&\& {M \otimes M'} \\
	\\
	\& {M \otimes M} \&\&\& {M' \otimes M'} \\
	\\
	\\
	M \&\&\& {M'}
	\arrow["{f \otimes 1}", from=1-1, to=1-4]
	\arrow[""{name=0, anchor=center, inner sep=0}, "{1 \otimes e'}"{description}, from=1-1, to=2-3]
	\arrow["{1 \otimes e}"{description}, from=1-1, to=4-2]
	\arrow[""{name=1, anchor=center, inner sep=0}, "{\rho_{\cV}}"', from=1-1, to=7-1]
	\arrow[between={0.2}{0.8}, equals, nfold, from=1-4, to=2-3]
	\arrow["{1 \otimes e'}", from=1-4, to=4-5]
	\arrow[""{name=2, anchor=center, inner sep=0}, "{f \otimes 1}"{description}, from=2-3, to=4-5]
	\arrow["{1 \otimes f}"{description}, from=4-2, to=2-3]
	\arrow[""{name=3, anchor=center, inner sep=0}, "{f \otimes f}"{description}, from=4-2, to=4-5]
	\arrow["\mu"{description}, from=4-2, to=7-1]
	\arrow["\chi"{description}, between={0.2}{0.8}, Rightarrow, nfold, from=4-5, to=7-1]
	\arrow["{\mu'}", from=4-5, to=7-4]
	\arrow["f"', from=7-1, to=7-4]
	\arrow["{1 \otimes \iota}"{description}, between={0.2}{0.8}, Rightarrow, nfold, from=0, to=4-2]
	\arrow[between={0.2}{0.8}, equals, nfold, from=2, to=3]
	\arrow["\rho", between={0.2}{0.8}, Rightarrow, nfold, from=4-2, to=1]
\end{tikzcd}\]
\[\vequiv\]
\[\begin{tikzcd}[ampersand replacement=\&]
	{M \otimes I} \&\&\& {M' \otimes I} \& \\
	\\
	\\
	\&\&\&\& {M' \otimes M'} \\
	\\
	\\
	M \&\&\& {M'}
	\arrow["{f \otimes 1}", from=1-1, to=1-4]
	\arrow["{\rho_{\cV}}"', from=1-1, to=7-1]
	\arrow["{e' \otimes 1}", from=1-4, to=4-5]
	\arrow[between={0.2}{0.8}, equals, nfold, from=1-4, to=7-1]
	\arrow[""{name=0, anchor=center, inner sep=0}, "{\rho_{\cV}}"{description}, from=1-4, to=7-4]
	\arrow["{\mu'}", from=4-5, to=7-4]
	\arrow["f"', from=7-1, to=7-4]
	\arrow["{\rho'}", between={0.2}{0.8}, Rightarrow, nfold, from=4-5, to=0]
\end{tikzcd}\]

	\newpage
	
	\section{Higher morphisms of weak involutions and inverses}\label{appsection:highermorphdefn}
	\subsection{Morphisms of weak $2$-involutions} \label{appsub:weak2inv}

\begin{definition}
	Let $(K,H,\alpha,\zeta)$ and $(K',H',\alpha',\zeta')$ be two weak $2$-involutions. A $1$-morphism of these is given by
	\begin{enumerate}
		\item a $1$-cell $R \colon K \rightarrow K'$,
		
		\item an adjoint equivalence $u,u^{.}$ given by the following $2$-cell:
		\[\begin{tikzcd}[ampersand replacement=\&]
			K \&\& {K'} \\
			\\
			K \&\& {K'}
			\arrow["R", from=1-1, to=1-3]
			\arrow["H"', from=1-1, to=3-1]
			\arrow["u", between={0.2}{0.8}, Rightarrow, nfold, from=1-3, to=3-1]
			\arrow["{H'}", from=1-3, to=3-3]
			\arrow["R"', from=3-1, to=3-3]
		\end{tikzcd}\cdpunct[8pt]{,}\]
		\item a $3$-cell $\delta$ as depicted by the following diagram:
		\[\begin{tikzcd}[ampersand replacement=\&,column sep=small]
			\&\& {K'} \&\&\&\& {K'} \& \\
			\\
			\&\& K \& {K'} \& {\Rrightarrow^{\delta}} \& K \&\& {K'} \\
			K \\
			\& K \& {K'} \&\&\&\& {K'} \\
			\& R \\
			\&\& {(H'H')R} \&\& {1_{K'}R} \&\& R \\
			\& {H'(H'R)} \&\&\&\&\&\& {R1_K} \\
			\& {H'(RH)} \&\&\&\&\&\& {R(HH)} \\
			\&\& {(H'R)H} \&\&\&\& {(RH)H}
			\arrow["u", between={0.2}{0.8}, Rightarrow, nfold, from=1-3, to=3-3]
			\arrow["{H'}", from=1-3, to=3-4]
			\arrow["{H'}", from=1-7, to=3-8]
			\arrow[""{name=0, anchor=center, inner sep=0}, "{1_{K'}}"{description}, from=1-7, to=5-7]
			\arrow["R"{description}, from=3-3, to=3-4]
			\arrow["H"{description}, from=3-3, to=5-2]
			\arrow["u", between={0.2}{0.8}, Rightarrow, nfold, from=3-4, to=5-2]
			\arrow["{H'}", from=3-4, to=5-3]
			\arrow["R", from=3-6, to=1-7]
			\arrow["R"', from=3-6, to=5-7]
			\arrow["{H'}", from=3-8, to=5-7]
			\arrow["R", from=4-1, to=1-3]
			\arrow["H"{description}, from=4-1, to=3-3]
			\arrow[""{name=1, anchor=center, inner sep=0}, "{1_{K}}"{description}, from=4-1, to=5-2]
			\arrow[no head, from=4-1, to=6-2]
			\arrow["R"{description}, from=5-2, to=5-3]
			\arrow[between={0.2}{0.8}, Rightarrow, nfold, from=5-2, to=6-2]
			\arrow[from=6-2, to=5-3]
			\arrow["{\alpha' 1}", from=7-3, to=7-5]
			\arrow["a"', from=7-3, to=8-2]
			\arrow["l", from=7-5, to=7-7]
			\arrow["1u"', from=8-2, to=9-2]
			\arrow["r"', from=8-8, to=7-7]
			\arrow["{a^{.}}"', from=9-2, to=10-3]
			\arrow["{1\alpha}"', from=9-8, to=8-8]
			\arrow[""{name=2, anchor=center, inner sep=0}, "u1"', from=10-3, to=10-7]
			\arrow["a"', from=10-7, to=9-8]
			\arrow["\alpha"{description}, between={0.2}{0.8}, Rightarrow, nfold, from=3-3, to=1]
			\arrow["{\alpha'}"{description}, between={0.2}{0.8}, Rightarrow, nfold, from=3-8, to=0]
			\arrow["\delta"', between={0.2}{0.8}, Rightarrow, nfold, from=2, to=7-5]
		\end{tikzcd}\cdpunct[8pt]{,}\]
		\item the $3$-cell $\delta$ satisfy the weak $2$-involution $1$-morphism axiom (\cref{app:weak2invmorph}).
		
	\end{enumerate}
\end{definition}

\begin{definition}
	A $2$-morphism $(\gamma,\omega) \colon (R,u,\delta_{u}) \Rightarrow (S,v,\delta_{v})$ is given by
	\begin{enumerate}
		\item a $2$-cell $\gamma \colon R \rightarrow S$,
		\item a $3$-cell $\omega$ as depicted by the following diagram:
		\[\begin{tikzcd}[ampersand replacement=\&]
			K \&\& {K'} \&\& K \&\& {K'} \\
			\&\&\& {\Rrightarrow^{\omega}} \\
			K \&\& {K'} \&\& K \&\& {K'}
			\arrow[""{name=0, anchor=center, inner sep=0}, "S"{description}, from=1-1, to=1-3]
			\arrow[""{name=1, anchor=center, inner sep=0}, "R", curve={height=-24pt}, from=1-1, to=1-3]
			\arrow["H"', from=1-1, to=3-1]
			\arrow["v"{description}, between={0.2}{0.8}, Rightarrow, nfold, from=1-3, to=3-1]
			\arrow["{H'}", from=1-3, to=3-3]
			\arrow["R"{description}, from=1-5, to=1-7]
			\arrow["H"', from=1-5, to=3-5]
			\arrow["u"{description}, between={0.2}{0.8}, Rightarrow, nfold, from=1-7, to=3-5]
			\arrow["{H'}", from=1-7, to=3-7]
			\arrow["S"', from=3-1, to=3-3]
			\arrow[""{name=2, anchor=center, inner sep=0}, "R"{description}, from=3-5, to=3-7]
			\arrow[""{name=3, anchor=center, inner sep=0}, "S"', curve={height=24pt}, from=3-5, to=3-7]
			\arrow["\gamma", between={0.2}{0.8}, Rightarrow, nfold, from=1, to=0]
			\arrow["\gamma", between={0.2}{0.8}, Rightarrow, nfold, from=2, to=3]
		\end{tikzcd}\cdpunct[8pt]{,}\]
		\item the $3$-cell $\omega$ satisfy the weak $2$-involution $2$-morphism axiom (\cref{app:weak2inv2morph}).
	\end{enumerate}
\end{definition}

\begin{definition}
	A $3$-morphism between $\rho \colon (\gamma,\omega_{\gamma}) \Rrightarrow (\lambda, \omega_{\lambda})$ is given by
	\begin{enumerate}
		\item a $3$-cell $\upsilon \colon \gamma \Rrightarrow \lambda$,
		
		\item the following identity holds:
		\[\begin{tikzcd}[ampersand replacement=\&]
			{H'R} \&\& {H'S} \&\& {H'R} \&\& {H'S} \\
			\&\&\& \equiv \\
			RH \&\& SH \&\& RH \&\& SH
			\arrow["{1\gamma}", from=1-1, to=1-3]
			\arrow["u"', from=1-1, to=3-1]
			\arrow["{\omega_{\gamma}}"{description}, between={0.2}{0.8}, Rightarrow, nfold, from=1-3, to=3-1]
			\arrow["v", from=1-3, to=3-3]
			\arrow[""{name=0, anchor=center, inner sep=0}, "{1\gamma}", curve={height=-30pt}, from=1-5, to=1-7]
			\arrow[""{name=1, anchor=center, inner sep=0}, "{1\lambda}"{description}, from=1-5, to=1-7]
			\arrow["u"', from=1-5, to=3-5]
			\arrow["{\omega_{\lambda}}"{description}, between={0.2}{0.8}, Rightarrow, nfold, from=1-7, to=3-5]
			\arrow["v", from=1-7, to=3-7]
			\arrow[""{name=2, anchor=center, inner sep=0}, "\gamma1"{description}, from=3-1, to=3-3]
			\arrow[""{name=3, anchor=center, inner sep=0}, "\lambda1"{description}, curve={height=30pt}, from=3-1, to=3-3]
			\arrow["\lambda1"', from=3-5, to=3-7]
			\arrow["\upsilon", between={0.2}{0.8}, Rightarrow, nfold, from=0, to=1]
			\arrow["\upsilon", between={0.2}{0.8}, Rightarrow, nfold, from=2, to=3]
		\end{tikzcd}\cdpunct[8pt]{.}\]
	\end{enumerate}
\end{definition}

\clearpage

\subsection{Morphisms of weak $2$-inverses} \label{appsub:weak2invs}

\begin{definition}
	A $1$-morphism $(K_1,K_2,F,G,\epsilon,\eta,\phi,\psi) \rightarrow (K'_1,K'_2,F',G',\epsilon',\eta',\phi',\psi')$ consists of
	\begin{enumerate}
		\item two $1$-cells $R_1\colon K_1 \rightarrow K'_1$ and $R_2\colon K_2 \rightarrow K'_2$,
		
		\item two $2$-cells
		\[\begin{tikzcd}[ampersand replacement=\&]
			{K_1} \&\& {K'_1} \&\& {K_2} \&\& {K'_2} \\
			\\
			{K_2} \&\& {K'_2} \&\& {K_1} \&\& {K'_1}
			\arrow["{R_1}", from=1-1, to=1-3]
			\arrow["F"', from=1-1, to=3-1]
			\arrow["{u_{12}}", between={0.2}{0.8}, Rightarrow, nfold, from=1-3, to=3-1]
			\arrow["{F'}", from=1-3, to=3-3]
			\arrow["{R_2}", from=1-5, to=1-7]
			\arrow["G"', from=1-5, to=3-5]
			\arrow["{u_{21}}", between={0.2}{0.8}, Rightarrow, nfold, from=1-7, to=3-5]
			\arrow["{G'}", from=1-7, to=3-7]
			\arrow["{R_2}"', from=3-1, to=3-3]
			\arrow["{R_1}"', from=3-5, to=3-7]
		\end{tikzcd}\cdpunct[8pt]{,}\]
		
		\item two $3$-cells as depicted by the following:
		\[\begin{tikzcd}[ampersand replacement=\&,cramped,column sep=tiny]
			\&\& {K_1} \&\&\&\&\&\& {K_1} \&\& \\
			\\
			{K_2} \&\& {K'_1} \&\&\& \Lleftarrow \& {K_2} \\
			\\
			\& {K'_2} \&\&\& {K'_1} \&\&\&\& {K_1} \&\& {K'_1} \\
			\\
			{K_1} \\
			\&\& {(G'F')R_1} \&\& {G'(F'R_1)} \&\& {G'(R_2F)} \&\& {(G'R_2)F} \\
			\\
			\&\& {1_{K'_1} R_1} \\
			\\
			\&\& {R_1} \&\& {R_11_{K_1}} \&\& {R_1(GF)} \&\& {(R_1G)F}
			\arrow["F"', from=1-3, to=3-1]
			\arrow["{R_1}"{description}, from=1-3, to=3-3]
			\arrow["{R_1}", from=1-3, to=5-5]
			\arrow["F"', from=1-9, to=3-7]
			\arrow[""{name=0, anchor=center, inner sep=0}, "{1_{K_1}}"{description}, from=1-9, to=5-9]
			\arrow["{R_1}", from=1-9, to=5-11]
			\arrow["{R_2}"{description}, from=3-1, to=5-2]
			\arrow["G"', from=3-1, to=7-1]
			\arrow["{u_{12}}"{description}, between={0.2}{0.8}, Rightarrow, nfold, from=3-3, to=3-1]
			\arrow["{F'}"{description}, from=3-3, to=5-2]
			\arrow[""{name=1, anchor=center, inner sep=0}, "{1_{K'_1}}"{description}, from=3-3, to=5-5]
			\arrow["G"', from=3-7, to=5-9]
			\arrow["{G'}"{description}, from=5-2, to=5-5]
			\arrow["{u_{21}}"{description}, between={0.2}{0.8}, Rightarrow, nfold, from=5-2, to=7-1]
			\arrow["{R_1}"', from=5-9, to=5-11]
			\arrow["{R_1}"', from=7-1, to=5-5]
			\arrow["a", from=8-3, to=8-5]
			\arrow[""{name=2, anchor=center, inner sep=0}, "{1u_{12}}", from=8-5, to=8-7]
			\arrow["{a^{.}}", from=8-7, to=8-9]
			\arrow["{u_{21}1}", from=8-9, to=12-9]
			\arrow["{\eta'1}", from=10-3, to=8-3]
			\arrow["{l^{.}}", from=12-3, to=10-3]
			\arrow["{r^{.}}"', from=12-3, to=12-5]
			\arrow[""{name=3, anchor=center, inner sep=0}, "{1\eta}"', from=12-5, to=12-7]
			\arrow["{a^{.}}"', from=12-7, to=12-9]
			\arrow["\eta"', between={0.2}{0.8}, Rightarrow, nfold, from=0, to=3-7]
			\arrow["{\eta'}"{description}, between={0.2}{1}, Rightarrow, nfold, from=1, to=5-2]
			\arrow[between={0.2}{0.8}, Rightarrow, nfold, from=3, to=2]
		\end{tikzcd}\]
		\[\begin{tikzcd}[ampersand replacement=\&,cramped,column sep=tiny]
			{K_2} \&\&\& {K'_2} \&\&\&\& {K_2} \&\& {K'_2} \&\& \\
			\&\&\& {K_1} \&\& {K'_1} \& \Rrightarrow \\
			\&\& {K_2} \&\&\&\&\&\&\&\&\& {K'_1} \\
			\\
			\&\& {K'_2} \&\&\&\&\&\&\& {K'_2} \\
			\&\&\& {F'(G'R_2)} \&\& {F'(R_1G)} \&\& {(F'R_1)G} \&\& {(R_2F)G} \\
			\&\&\&\&\&\&\&\&\& {R_2(FG)} \\
			\&\&\&\&\&\&\&\&\& {R_21_{K_2}} \\
			\&\&\& {(F'G')R_2} \&\&\& {1_{K'_2}R_2} \&\&\& {R_2}
			\arrow["{R_2}", from=1-1, to=1-4]
			\arrow["G"{description}, from=1-1, to=2-4]
			\arrow[""{name=0, anchor=center, inner sep=0}, "{1_{K_2}}"{description}, from=1-1, to=3-3]
			\arrow["{R_2}"', from=1-1, to=5-3]
			\arrow["{u_{21}}", between={0.2}{0.8}, Rightarrow, nfold, from=1-4, to=2-4]
			\arrow["{G'}", from=1-4, to=2-6]
			\arrow["{R_2}", from=1-8, to=1-10]
			\arrow["{R_2}"', from=1-8, to=5-10]
			\arrow["{G'}", from=1-10, to=3-12]
			\arrow[""{name=1, anchor=center, inner sep=0}, "{1_{K'_2}}"{description}, from=1-10, to=5-10]
			\arrow["{R_1}"{description}, from=2-4, to=2-6]
			\arrow["F"{description}, from=2-4, to=3-3]
			\arrow["{u_{12}}", between={0.2}{0.8}, Rightarrow, nfold, from=2-6, to=3-3]
			\arrow["{F'}", from=2-6, to=5-3]
			\arrow["{R_2}"{description}, from=3-3, to=5-3]
			\arrow["{F'}", from=3-12, to=5-10]
			\arrow["{1u_{21}}", from=6-4, to=6-6]
			\arrow[""{name=2, anchor=center, inner sep=0}, "{a^{.}}", from=6-6, to=6-8]
			\arrow["{u_{12}1}", from=6-8, to=6-10]
			\arrow["a", from=6-10, to=7-10]
			\arrow["{1\epsilon}", from=7-10, to=8-10]
			\arrow["r", from=8-10, to=9-10]
			\arrow["a", from=9-4, to=6-4]
			\arrow["{\epsilon'1}"', from=9-4, to=9-7]
			\arrow["l"', from=9-7, to=9-10]
			\arrow["\epsilon"{description}, between={0.2}{0.8}, Rightarrow, nfold, from=2-4, to=0]
			\arrow["{\epsilon'}", between={0.2}{0.8}, Rightarrow, nfold, from=3-12, to=1]
			\arrow[between={0.2}{0.8}, Rightarrow, nfold, from=2, to=9-7]
		\end{tikzcd}\]
		\item such that the two weak $2$-inverse $1$-morphism axioms hold (\cref{app:weak2invsmorph}).
	\end{enumerate}
\end{definition}

\begin{definition}
	A $2$-morphism $(R_1,R_2,u_{12},u_{21},\Lleftarrow, \Rrightarrow) \Rightarrow (S_1,S_2,v_{12},v_{21},\Lleftarrow, \Rrightarrow)$ consists of
	\begin{enumerate}
		\item two $2$-cells $\gamma_1 \colon R_1 \Rightarrow S_1$ and $\gamma_2 \colon R_2 \Rightarrow S_2$,
		
		\item two $3$-cells
		\[\begin{tikzcd}[ampersand replacement=\&]
			{K_1} \&\& {K'_1} \&\& {K_1} \&\& {K'_1} \\
			\&\&\& \Rrightarrow \\
			{K_2} \&\& {K'_2} \&\& {K_2} \&\& {K'_2}
			\arrow[""{name=0, anchor=center, inner sep=0}, "{R_1}", curve={height=-30pt}, from=1-1, to=1-3]
			\arrow[""{name=1, anchor=center, inner sep=0}, "{S_1}"{description}, from=1-1, to=1-3]
			\arrow["F"', from=1-1, to=3-1]
			\arrow["{v_{12}}", between={0.2}{0.8}, Rightarrow, nfold, from=1-3, to=3-1]
			\arrow["{F'}", from=1-3, to=3-3]
			\arrow["{R_1}", from=1-5, to=1-7]
			\arrow["F"', from=1-5, to=3-5]
			\arrow["{u_{12}}", between={0.2}{0.8}, Rightarrow, nfold, from=1-7, to=3-5]
			\arrow["{F'}", from=1-7, to=3-7]
			\arrow["{S_2}"', from=3-1, to=3-3]
			\arrow[""{name=2, anchor=center, inner sep=0}, "{R_2}"{description}, from=3-5, to=3-7]
			\arrow[""{name=3, anchor=center, inner sep=0}, "{S_2}"', curve={height=30pt}, from=3-5, to=3-7]
			\arrow["{\gamma_1}", between={0.2}{0.8}, Rightarrow, nfold, from=0, to=1]
			\arrow["{\gamma_2}"', between={0.2}{0.8}, Rightarrow, nfold, from=2, to=3]
		\end{tikzcd}\]
		\[\begin{tikzcd}[ampersand replacement=\&]
			{K_2} \&\& {K'_2} \&\& {K_2} \&\& {K'_2} \\
			\&\&\& \Rrightarrow \\
			{K_1} \&\& {K'_1} \&\& {K_1} \&\& {K'_1}
			\arrow[""{name=0, anchor=center, inner sep=0}, "{S_2}"{description}, from=1-1, to=1-3]
			\arrow[""{name=1, anchor=center, inner sep=0}, "{R_2}", curve={height=-30pt}, from=1-1, to=1-3]
			\arrow["G"', from=1-1, to=3-1]
			\arrow["{v_{21}}", between={0.2}{0.8}, Rightarrow, nfold, from=1-3, to=3-1]
			\arrow["{G'}", from=1-3, to=3-3]
			\arrow["{R_2}", from=1-5, to=1-7]
			\arrow["G"', from=1-5, to=3-5]
			\arrow["{u_{21}}", between={0.2}{0.8}, Rightarrow, nfold, from=1-7, to=3-5]
			\arrow["{G'}", from=1-7, to=3-7]
			\arrow["{S_1}"', from=3-1, to=3-3]
			\arrow[""{name=2, anchor=center, inner sep=0}, "{R_1}"{description}, from=3-5, to=3-7]
			\arrow[""{name=3, anchor=center, inner sep=0}, "{S_1}"', curve={height=30pt}, from=3-5, to=3-7]
			\arrow["{\gamma_2}", between={0.2}{0.8}, Rightarrow, nfold, from=1, to=0]
			\arrow["{\gamma_1}", between={0.2}{0.8}, Rightarrow, nfold, from=2, to=3]
		\end{tikzcd}\cdpunct[8pt]{,}\]
		
		\item such that the two weak $2$-inverse $2$-morphism axioms hold (\cref{app:weak2invs2morph}).
		
	\end{enumerate}
\end{definition}

\begin{definition}
	A $3$-morphism $(\gamma_1,\gamma_2,\omega_{12},\omega_{21}) \Rrightarrow (\lambda_1,\lambda_2,\omega_{12},\omega_{21})$ consists of
	\begin{enumerate}
		\item two $3$-cells $\upsilon_1 \colon \gamma_1 \Rrightarrow \lambda_1$ and $\upsilon_2\colon \gamma_2 \Rrightarrow \lambda_2$.
		
		\item such that the following two identities holds:
		\[\begin{tikzcd}[ampersand replacement=\&]
			{G'R_2} \&\& {G'S_2} \&\& {G'R_2} \&\& {G'S_2} \\
			\&\&\& \equiv \\
			{R_1G} \&\& {S_1G} \&\& {R_1G} \&\& {S_1G}
			\arrow[""{name=0, anchor=center, inner sep=0}, "{1\lambda_2}"{description}, from=1-1, to=1-3]
			\arrow[""{name=1, anchor=center, inner sep=0}, "{1\gamma_2}", curve={height=-30pt}, from=1-1, to=1-3]
			\arrow["{u_{21}}"', from=1-1, to=3-1]
			\arrow[between={0.2}{0.8}, Rightarrow, nfold, from=1-3, to=3-1]
			\arrow["{v_{21}}", from=1-3, to=3-3]
			\arrow["{1\gamma_2}", from=1-5, to=1-7]
			\arrow["{u_{21}}"', from=1-5, to=3-5]
			\arrow[between={0.2}{0.8}, Rightarrow, nfold, from=1-7, to=3-5]
			\arrow["{v_{21}}", from=1-7, to=3-7]
			\arrow["{\lambda_11}"', from=3-1, to=3-3]
			\arrow[""{name=2, anchor=center, inner sep=0}, "{\gamma_11}"{description}, from=3-5, to=3-7]
			\arrow[""{name=3, anchor=center, inner sep=0}, "{\lambda_11}"', curve={height=30pt}, from=3-5, to=3-7]
			\arrow["{1\upsilon_2}", between={0.2}{0.8}, Rightarrow, nfold, from=1, to=0]
			\arrow["{\upsilon_11}"', between={0.2}{0.8}, Rightarrow, nfold, from=2, to=3]
		\end{tikzcd}\]
		\[\begin{tikzcd}[ampersand replacement=\&]
			{F'R_1} \&\& {F'S_1} \&\& {F'R_1} \&\& {F'S_1} \\
			\&\&\& \equiv \\
			{R_2F} \&\& {S_2F} \&\& {R_2F} \&\& {S_2F}
			\arrow[""{name=0, anchor=center, inner sep=0}, "{1\lambda_1}"{description}, from=1-1, to=1-3]
			\arrow[""{name=1, anchor=center, inner sep=0}, "{1\gamma_1}", curve={height=-30pt}, from=1-1, to=1-3]
			\arrow["{u_{12}}"', from=1-1, to=3-1]
			\arrow[between={0.2}{0.8}, Rightarrow, nfold, from=1-3, to=3-1]
			\arrow["{v_{12}}", from=1-3, to=3-3]
			\arrow["{1\gamma_1}", from=1-5, to=1-7]
			\arrow["{u_{12}}"', from=1-5, to=3-5]
			\arrow[between={0.2}{0.8}, Rightarrow, nfold, from=1-7, to=3-5]
			\arrow["{v_{12}}", from=1-7, to=3-7]
			\arrow["{\lambda_21}"', from=3-1, to=3-3]
			\arrow[""{name=2, anchor=center, inner sep=0}, "{\lambda_21}"', curve={height=30pt}, from=3-5, to=3-7]
			\arrow[""{name=3, anchor=center, inner sep=0}, "{\gamma_21}"{description}, from=3-5, to=3-7]
			\arrow["{1\upsilon_1}", between={0.2}{0.8}, Rightarrow, nfold, from=1, to=0]
			\arrow["{\upsilon_21}"', between={0.2}{0.8}, Rightarrow, nfold, from=3, to=2]
		\end{tikzcd}\cdpunct[8pt]{.}\]
	\end{enumerate}
\end{definition}
	\newpage
	
	\section{Coherence diagrams for weak objects} \label{appsection:wobjdiag}

\subsection{Weak $2$-involution} \label{app:weak2invdef}

\[\begin{tikzcd}[ampersand replacement=\&]
	{((HH)H)H} \&\& {(HH)(HH)} \&\& {H(H(HH))} \\
	\\
	\& {(H(HH))H} \&\& {H((HH)H)} \\
	\& {(H1)H} \&\& {H(1H)} \\
	{(1H)H} \& HH \&\& HH \& {H(H1)} \\
	\&\& 1
	\arrow["a", from=1-1, to=1-3]
	\arrow["a1", from=1-1, to=3-2]
	\arrow[""{name=0, anchor=center, inner sep=0}, "{(\alpha1)1}"', from=1-1, to=5-1]
	\arrow["a", from=1-3, to=1-5]
	\arrow["{1(1\alpha)}", from=1-5, to=5-5]
	\arrow["a", from=3-2, to=3-4]
	\arrow["{(1\alpha)1}", from=3-2, to=4-2]
	\arrow["1a", from=3-4, to=1-5]
	\arrow["{1(\alpha1)}"', from=3-4, to=4-4]
	\arrow["{1\zeta}", between={0.2}{0.8}, Rightarrow, nfold, from=3-4, to=5-5]
	\arrow["a", from=4-2, to=4-4]
	\arrow["r1", from=4-2, to=5-2]
	\arrow["1l"', from=4-4, to=5-4]
	\arrow["l1"', from=5-1, to=5-2]
	\arrow["{\varepsilon_{HH}}", from=5-2, to=5-4]
	\arrow["\alpha"', from=5-2, to=6-3]
	\arrow["\alpha", from=5-4, to=6-3]
	\arrow["1r", from=5-5, to=5-4]
	\arrow["\zeta1"', between={0.2}{1}, Rightarrow, nfold, from=0, to=4-2]
\end{tikzcd}\]
\[\vequiv\]
\[\begin{tikzcd}[ampersand replacement=\&]
	{((HH)H)H} \&\& {(HH)(HH)} \&\& {H(H(HH))} \\
	{(1H)H} \& {1(HH)} \&\& {(HH)1} \& {H(H1)} \\
	\& HH \& 11 \& HH \\
	\\
	\&\& 1
	\arrow["a", from=1-1, to=1-3]
	\arrow["{(\alpha1)1}"', from=1-1, to=2-1]
	\arrow["a", from=1-3, to=1-5]
	\arrow["\alpha1"', from=1-3, to=2-2]
	\arrow["{1\alpha}", from=1-3, to=2-4]
	\arrow["{1(1\alpha)}", from=1-5, to=2-5]
	\arrow["a", from=2-1, to=2-2]
	\arrow["l1"', from=2-1, to=3-2]
	\arrow["l"', from=2-2, to=3-2]
	\arrow["{1\alpha}"', from=2-2, to=3-3]
	\arrow["a", from=2-4, to=2-5]
	\arrow["\alpha1", from=2-4, to=3-3]
	\arrow["r", from=2-4, to=3-4]
	\arrow["1r", from=2-5, to=3-4]
	\arrow["\alpha"', from=3-2, to=5-3]
	\arrow["l"', curve={height=12pt}, from=3-3, to=5-3]
	\arrow["r", curve={height=-12pt}, from=3-3, to=5-3]
	\arrow["\alpha", from=3-4, to=5-3]
\end{tikzcd}\]
\clearpage

\subsection{Weak $2$-involution $1$-morphism} \label{app:weak2invmorph}

\[\begin{adjustbox}{
		angle=90,
		max totalsize={\linewidth}{0.95\textheight},
		valign=t
	}
	\begin{tikzcd}[ampersand replacement=\&,cramped,column sep=tiny]
		\& {H'((H'H')R)} \&\&\& {H'(1R)} \&\&\&\&\&\&\&\& \\
		{(H'(H'H'))R} \&\&\&\&\&\&\& {H'R} \\
		{((H'H')H')R} \& {H'(H'(H'R))} \&\& {H'((H'R)H)} \\
		\&\& {H'(H'(RH))} \&\&\& {H'(R1)} \&\& {(RH)1} \\
		{(H'H')(H'R)} \&\&\&\& {H'((RH)H)} \&\& {(H'R)1} \&\& {R(H1)} \&\&\&\& RH \\
		\& {(H'H')(RH)} \&\&\&\& {H'(R(HH))} \\
		\& {((H'H')R)H} \& {(H'(H'R))H} \&\&\&\& {(RH)(HH)} \\
		\&\&\& {(H'(RH))H} \&\& {(H'R)(HH)} \&\& {R(H(HH))} \&\& {R(1H)} \\
		\&\&\&\&\&\&\& {R((HH)H)} \\
		\&\&\&\& {((H'R)H)H} \&\&\&\&\& {(R1)H} \\
		\&\&\&\&\& {((RH)H)H} \& {(R(HH))H}
		\arrow[""{name=0, anchor=center, inner sep=0}, "{1(\alpha'1)}", from=1-2, to=1-5]
		\arrow["1a"{description}, from=1-2, to=3-2]
		\arrow["1l", from=1-5, to=2-8]
		\arrow["a", from=2-1, to=1-2]
		\arrow["u", from=2-8, to=5-13]
		\arrow["a1", from=3-1, to=2-1]
		\arrow["a"{description}, from=3-1, to=5-1]
		\arrow["{1(1u)}"{description}, from=3-2, to=4-3]
		\arrow["{1(u1)}"{description}, from=3-4, to=5-5]
		\arrow["{1a^{.}}"{description}, from=4-3, to=3-4]
		\arrow["1r"{description}, from=4-6, to=2-8]
		\arrow["{a^{.}}"{description}, from=4-6, to=5-7]
		\arrow["a"{description}, from=4-8, to=5-9]
		\arrow["r"{description}, from=4-8, to=5-13]
		\arrow["a"{description}, from=5-1, to=3-2]
		\arrow["{(11)u}"', from=5-1, to=6-2]
		\arrow["1a"{description}, from=5-5, to=6-6]
		\arrow["r"{description}, from=5-7, to=2-8]
		\arrow["u1"{description}, from=5-7, to=4-8]
		\arrow["1r"{description}, from=5-9, to=5-13]
		\arrow["a"{description}, from=6-2, to=4-3]
		\arrow["{a^{.}}"', from=6-2, to=7-2]
		\arrow["{1(1\alpha)}"{description}, from=6-6, to=4-6]
		\arrow["{a^{.}}"{description}, from=6-6, to=8-6]
		\arrow["a1"', from=7-2, to=7-3]
		\arrow["a"{description}, from=7-3, to=3-4]
		\arrow["{(1u)1}"', from=7-3, to=8-4]
		\arrow["{(11)\alpha}"{description}, from=7-7, to=4-8]
		\arrow["a"{description}, from=7-7, to=8-8]
		\arrow["a"{description}, from=8-4, to=5-5]
		\arrow["{a^{.}1}"', from=8-4, to=10-5]
		\arrow["{(11)\alpha}"{description}, from=8-6, to=5-7]
		\arrow["{u(11)}"{description}, from=8-6, to=7-7]
		\arrow["{1(1\alpha)}"{description}, from=8-8, to=5-9]
		\arrow[""{name=1, anchor=center, inner sep=0}, "1l"{description}, from=8-10, to=5-13]
		\arrow["1a"', from=9-8, to=8-8]
		\arrow["{1(\alpha1)}"{description}, from=9-8, to=8-10]
		\arrow["a"{description}, from=10-5, to=8-6]
		\arrow["{(u1)1}"', from=10-5, to=11-6]
		\arrow["r1"', from=10-10, to=5-13]
		\arrow[from=10-10, to=8-10]
		\arrow["a"{description}, from=11-6, to=7-7]
		\arrow["a1"', from=11-6, to=11-7]
		\arrow["a"{description}, from=11-7, to=9-8]
		\arrow["{(1\alpha)1}"', from=11-7, to=10-10]
		\arrow["{1\delta}", between={0.2}{0.8}, Rightarrow, nfold, from=3-4, to=0]
		\arrow["{1\zeta}", between={0.2}{0.8}, Rightarrow, nfold, from=1, to=5-9]
	\end{tikzcd}	
\end{adjustbox}
\]
\clearpage
\[\vequiv\]
\[\begin{adjustbox}{
		angle=90,
		max totalsize={\linewidth}{0.95\textheight},
		valign=t
	}
	\begin{tikzcd}[ampersand replacement=\&,cramped,column sep=tiny]
		\&\& {H'((H'H')R)} \&\&\&\&\& \\
		\& {(H'(H'H'))R} \&\&\& {H'(1R)} \\
		\&\& {(H'1)R} \&\&\&\& {H'R} \\
		{((H'H')H')R} \&\& {(1H')R} \&\&\&\&\& RH \\
		\&\&\& {1(H'R)} \\
		{(H'H')(H'R)} \&\&\&\& {1(RH)} \&\&\& {(R1)H} \\
		\& {(H'H')(RH)} \&\&\&\& {(1R)H} \\
		\& {((H'H')R)H} \&\&\&\&\&\& {(R(HH))H} \\
		\&\& {(H'(H'R))H} \&\&\&\& {((RH)H)H)} \\
		\&\&\& {(H'(RH))H} \& {((H'R)H)H}
		\arrow["{1(\alpha'1)}", from=1-3, to=2-5]
		\arrow["a", from=2-2, to=1-3]
		\arrow["{(1\alpha')1}"{description}, from=2-2, to=3-3]
		\arrow["1l", from=2-5, to=3-7]
		\arrow["a"{description}, from=3-3, to=2-5]
		\arrow["r1"{description}, from=3-3, to=3-7]
		\arrow["u", from=3-7, to=4-8]
		\arrow["a1", from=4-1, to=2-2]
		\arrow["{(\alpha'1)1}"{description}, from=4-1, to=4-3]
		\arrow["a"', from=4-1, to=6-1]
		\arrow["{\zeta'1}", between={0.2}{0.8}, Rightarrow, nfold, from=4-3, to=3-3]
		\arrow["l1"{description}, from=4-3, to=3-7]
		\arrow["a"{description}, from=4-3, to=5-4]
		\arrow["l"{description}, from=5-4, to=3-7]
		\arrow["1u"{description}, from=5-4, to=6-5]
		\arrow["{\alpha'(11)}"{description}, from=6-1, to=5-4]
		\arrow["{(11)u}"', from=6-1, to=7-2]
		\arrow["l"{description}, from=6-5, to=4-8]
		\arrow["{a^{.}}"{description}, from=6-5, to=7-6]
		\arrow["r1"', from=6-8, to=4-8]
		\arrow["{\alpha'(11)}"{description}, from=7-2, to=6-5]
		\arrow["{a^{.}}"', from=7-2, to=8-2]
		\arrow["l1"{description}, from=7-6, to=4-8]
		\arrow["{(\alpha'1)1}"{description}, from=8-2, to=7-6]
		\arrow["a1"', from=8-2, to=9-3]
		\arrow["{(1\alpha)1}"', from=8-8, to=6-8]
		\arrow["{(1u)1}"', from=9-3, to=10-4]
		\arrow[""{name=0, anchor=center, inner sep=0}, "a1"', from=9-7, to=8-8]
		\arrow["{a^{.}1}"', from=10-4, to=10-5]
		\arrow["{(u1)1}"', from=10-5, to=9-7]
		\arrow["\delta1"', between={0.2}{0.8}, Rightarrow, nfold, from=0, to=7-6]
	\end{tikzcd}
\end{adjustbox}
\]
\clearpage

\subsection{Weak $2$-involution $2$-morphism} \label{app:weak2inv2morph}


\[\begin{adjustbox}{
		angle=90,
		max totalsize={\linewidth}{0.95\textheight},
		valign=t
	}
	\begin{tikzcd}[ampersand replacement=\&,cramped,column sep=tiny]
		\& {(H'H')S} \&\& {H'(H'S)} \&\& {H'(SH)} \&\& {(H'S)H} \&\& {(SH)H} \&\& {S(HH)} \&\& S1 \&\& S \\
		\\
		\\
		{(H'H')R} \&\& {H'(H'R)} \&\& {H'(RH)} \&\& {(H'R)H} \&\& {(RH)H} \&\& {R(HH)} \&\& R1 \&\& R \\
		\\
		\&\&\&\&\&\&\& 1R
		\arrow["a", from=1-2, to=1-4]
		\arrow["1v", from=1-4, to=1-6]
		\arrow["{1\omega}"', between={0.2}{0.8}, Rightarrow, nfold, from=1-4, to=4-5]
		\arrow["{a^{.}}", from=1-6, to=1-8]
		\arrow["v1", from=1-8, to=1-10]
		\arrow["\omega1"', between={0.2}{0.8}, Rightarrow, nfold, from=1-8, to=4-9]
		\arrow["a", from=1-10, to=1-12]
		\arrow["{1\alpha}", from=1-12, to=1-14]
		\arrow["r", from=1-14, to=1-16]
		\arrow["{(11)\gamma}", from=4-1, to=1-2]
		\arrow["a"{description}, from=4-1, to=4-3]
		\arrow["{\alpha'1}"', from=4-1, to=6-8]
		\arrow["{1(1\gamma)}"{description}, from=4-3, to=1-4]
		\arrow["1u"{description}, from=4-3, to=4-5]
		\arrow["{1(\gamma1)}"{description}, from=4-5, to=1-6]
		\arrow["{a^{.}}"{description}, from=4-5, to=4-7]
		\arrow["{(1\gamma)1}"{description}, from=4-7, to=1-8]
		\arrow[""{name=0, anchor=center, inner sep=0}, "u1"{description}, from=4-7, to=4-9]
		\arrow["{(\gamma1)1}"{description}, from=4-9, to=1-10]
		\arrow["a"{description}, from=4-9, to=4-11]
		\arrow["{\gamma(11)}"{description}, from=4-11, to=1-12]
		\arrow["{1\alpha}"{description}, from=4-11, to=4-13]
		\arrow["\gamma1"{description}, from=4-13, to=1-14]
		\arrow["r"{description}, from=4-13, to=4-15]
		\arrow["\gamma"', from=4-15, to=1-16]
		\arrow["l"', from=6-8, to=4-15]
		\arrow["{\delta_u}"', between={0.2}{0.8}, Rightarrow, nfold, from=0, to=6-8]
	\end{tikzcd}
\end{adjustbox}
\]
\clearpage
\[\vequiv\]
\[\begin{adjustbox}{
		angle=90,
		max totalsize={\linewidth}{0.90\textheight},
		valign=t
	}
	\begin{tikzcd}[ampersand replacement=\&,cramped]
		\& {(H'H')S} \&\& {H'(H'S)} \&\& {H'(SH)} \&\& {(H'S)H} \&\& {(SH)H} \&\& {S(HH)} \&\& S1 \&\& S \\
		\\
		\&\&\&\&\&\&\&\& 1S \\
		{(H'H')R} \&\&\&\&\&\&\&\&\&\&\&\&\&\& R \\
		\\
		\&\&\&\&\&\&\& 1R
		\arrow["a", from=1-2, to=1-4]
		\arrow["{\alpha'1}"{description}, from=1-2, to=3-9]
		\arrow["1v", from=1-4, to=1-6]
		\arrow["{a^{.}}", from=1-6, to=1-8]
		\arrow[""{name=0, anchor=center, inner sep=0}, "v1", from=1-8, to=1-10]
		\arrow["a", from=1-10, to=1-12]
		\arrow["{1\alpha}", from=1-12, to=1-14]
		\arrow["r", from=1-14, to=1-16]
		\arrow["l"{description}, from=3-9, to=1-16]
		\arrow["{(11)\gamma}", from=4-1, to=1-2]
		\arrow["{\alpha'1}"', from=4-1, to=6-8]
		\arrow["\gamma"', from=4-15, to=1-16]
		\arrow["{1\gamma}"{description}, from=6-8, to=3-9]
		\arrow["l"', from=6-8, to=4-15]
		\arrow["{\delta_{v}}"', between={0.2}{0.8}, Rightarrow, nfold, from=0, to=3-9]
	\end{tikzcd}
\end{adjustbox}
\]

\clearpage
\subsection{Weak $2$-inverse} \label{app:weak2invs}

\[\begin{tikzcd}[ampersand replacement=\&,cramped,column sep=scriptsize]
	\&\&\&\&\&\& FG \&\&\& \\
	\&\&\&\&\& {(1F)G} \\
	\\
	\&\&\&\&\& {((FG)F)G} \\
	\&\& {(F1)G} \&\& {(F(GF))G} \&\& {1(FG)} \\
	FG \&\&\&\&\& {(FG)(FG)} \&\& 11 \&\& 1 \\
	\&\& {F(1G)} \&\& {F((GF)G)} \&\& {(FG)1} \\
	\&\&\&\&\& {F(G(FG))} \\
	\\
	\&\&\&\&\& {F(G1)} \\
	\&\&\&\&\&\& FG
	\arrow["{l^.}"{description}, from=1-7, to=5-7]
	\arrow["\epsilon", from=1-7, to=6-10]
	\arrow["l1"{description}, from=2-6, to=1-7]
	\arrow["a"{description}, from=2-6, to=5-7]
	\arrow["{(\epsilon1)1}"{description}, from=4-6, to=2-6]
	\arrow["a"{description}, from=4-6, to=6-6]
	\arrow["{(1\eta)1}"{description}, from=5-3, to=5-5]
	\arrow["a"{description}, from=5-3, to=7-3]
	\arrow[""{name=0, anchor=center, inner sep=0}, "{a^.1}"{description}, from=5-5, to=4-6]
	\arrow["a"{description}, from=5-5, to=7-5]
	\arrow["{1\epsilon}"{description}, from=5-7, to=6-8]
	\arrow[""{name=1, anchor=center, inner sep=0}, "1", curve={height=-30pt}, from=6-1, to=1-7]
	\arrow["{r^.1}"{description}, from=6-1, to=5-3]
	\arrow["{1l^.}"{description}, from=6-1, to=7-3]
	\arrow[""{name=2, anchor=center, inner sep=0}, "1"', curve={height=30pt}, from=6-1, to=11-7]
	\arrow["\epsilon1"{description}, from=6-6, to=5-7]
	\arrow["{1\epsilon}"{description}, from=6-6, to=7-7]
	\arrow["l"{description}, curve={height=-12pt}, from=6-8, to=6-10]
	\arrow["r"{description}, curve={height=12pt}, from=6-8, to=6-10]
	\arrow["{1(\eta1)}"{description}, from=7-3, to=7-5]
	\arrow[""{name=3, anchor=center, inner sep=0}, "1a"{description}, from=7-5, to=8-6]
	\arrow["\epsilon1"{description}, from=7-7, to=6-8]
	\arrow["{a^.}"{description}, from=8-6, to=6-6]
	\arrow["{1(1\epsilon)}"{description}, from=8-6, to=10-6]
	\arrow["{a^.}"{description}, from=10-6, to=7-7]
	\arrow["1r", from=10-6, to=11-7]
	\arrow["\epsilon"', from=11-7, to=6-10]
	\arrow["{r^.}"{description}, from=11-7, to=7-7]
	\arrow["\phi1", between={0.2}{0.8}, Rightarrow, nfold, from=0, to=1]
	\arrow["{1\psi}"', between={0.2}{0.8}, Rightarrow, nfold, from=2, to=3]
\end{tikzcd}\]
\[\vequiv\]
\[\begin{tikzcd}[ampersand replacement=\&,cramped]
	\&\& FG \&\& \\
	\\
	FG \&\&\&\& 1 \\
	\\
	\&\& FG
	\arrow["\epsilon", from=1-3, to=3-5]
	\arrow["1", from=3-1, to=1-3]
	\arrow["1"', from=3-1, to=5-3]
	\arrow["\epsilon"', from=5-3, to=3-5]
\end{tikzcd}\]

\clearpage

\[\begin{tikzcd}[ampersand replacement=\&,cramped,column sep=scriptsize]
	\&\&\& GF \&\&\&\&\&\& \\
	\&\&\&\& {(1G)F} \\
	\\
	\&\&\&\& {((GF)G)F} \\
	\&\&\& {1(GF)} \&\& {(G(FG))F} \&\& {(G1)F} \\
	1 \&\& 11 \&\& {(GF)(GF)} \&\&\&\&\& GF \\
	\&\&\& {(GF)1} \&\& {G((FG)F)} \&\& {G(1F)} \\
	\&\&\&\& {G(F(GF))} \\
	\\
	\&\&\&\& {G(F1)} \\
	\&\&\& GF
	\arrow["{l^.1}"{description}, from=1-4, to=2-5]
	\arrow[""{name=0, anchor=center, inner sep=0}, "1", curve={height=-30pt}, from=1-4, to=6-10]
	\arrow["{(\eta1)1}"{description}, from=2-5, to=4-5]
	\arrow[""{name=1, anchor=center, inner sep=0}, "a1"{description}, from=4-5, to=5-6]
	\arrow["l"{description}, from=5-4, to=1-4]
	\arrow["{a^.}"{description}, from=5-4, to=2-5]
	\arrow["\eta1"{description}, from=5-4, to=6-5]
	\arrow["{(1\epsilon)1}"{description}, from=5-6, to=5-8]
	\arrow["a"{description}, from=5-6, to=7-6]
	\arrow["r1"{description}, from=5-8, to=6-10]
	\arrow["a"{description}, from=5-8, to=7-8]
	\arrow["\eta", from=6-1, to=1-4]
	\arrow["{l^.}"{description}, curve={height=-12pt}, from=6-1, to=6-3]
	\arrow["{r^.}"{description}, curve={height=12pt}, from=6-1, to=6-3]
	\arrow["\eta"', from=6-1, to=11-4]
	\arrow["{1\eta}"{description}, from=6-3, to=5-4]
	\arrow["\eta1"{description}, from=6-3, to=7-4]
	\arrow["{a^.}"{description}, from=6-5, to=4-5]
	\arrow["a"{description}, from=6-5, to=8-5]
	\arrow["{1\eta}"{description}, from=7-4, to=6-5]
	\arrow["a"{description}, from=7-4, to=10-5]
	\arrow["r"{description}, from=7-4, to=11-4]
	\arrow["{1(\epsilon1)}"{description}, from=7-6, to=7-8]
	\arrow["1l"{description}, from=7-8, to=6-10]
	\arrow[""{name=2, anchor=center, inner sep=0}, "{1a^.}"{description}, from=8-5, to=7-6]
	\arrow["{1(1\eta)}"{description}, from=10-5, to=8-5]
	\arrow[""{name=3, anchor=center, inner sep=0}, "1"', curve={height=30pt}, from=11-4, to=6-10]
	\arrow["{1r^.}"{description}, from=11-4, to=10-5]
	\arrow["\psi1", between={0.2}{0.8}, Rightarrow, nfold, from=0, to=1]
	\arrow["{1\phi}", between={0.2}{0.8}, Rightarrow, nfold, from=2, to=3]
\end{tikzcd}\]
\[\vequiv\]
\[\begin{tikzcd}[ampersand replacement=\&,cramped]
	\&\& GF \&\& \\
	\\
	1 \&\&\&\& GF \\
	\\
	\&\& GF
	\arrow["1", from=1-3, to=3-5]
	\arrow["\eta", from=3-1, to=1-3]
	\arrow["\eta"', from=3-1, to=5-3]
	\arrow["1"', from=5-3, to=3-5]
\end{tikzcd}\]

\clearpage

\subsection{Weak $2$-inverse $1$-morphism} \label{app:weak2invsmorph}

\[\begin{adjustbox}{
		angle=90,
		max totalsize={\linewidth}{0.90\textheight},
		valign=t
	}
	\begin{tikzcd}[ampersand replacement=\&,cramped,column sep=scriptsize,row sep=2.25em]
		\&\&\& {F'((R_1G)F)} \& {(F'(R_1G))F} \& {((F'R_1)G)F} \&\&\&\&\& \\
		\&\& {F'(R_1(GF))} \&\& {(F'R_1)(GF)} \&\& {((R_2F)G)F} \\
		\& {F'(R_11)} \&\& {(F'R_1)1} \&\& {(R_2F)(GF)} \&\& {(R_2(FG))F} \\
		\&\&\&\& {(R_2F)1} \&\& {R_2(F(GF))} \&\& {(R_21)F} \\
		{F'R_1} \&\&\& {R_2F} \&\& {R_2(F1)} \&\& {R_2((FG)F)} \\
		\&\&\&\&\&\&\&\& {R_2(1F)} \\
		\&\&\&\&\&\&\&\&\&\& {R_2F}
		\arrow["{a^.}", from=1-4, to=1-5]
		\arrow["{a^.1}", from=1-5, to=1-6]
		\arrow["{(u_{12}1)1}", from=1-6, to=2-7]
		\arrow["{1a^.}", from=2-3, to=1-4]
		\arrow["{a^.}"{description}, from=2-3, to=2-5]
		\arrow["{a^.}"{description}, from=2-5, to=1-6]
		\arrow["{u_{12}(11)}"{description}, from=2-5, to=3-6]
		\arrow["a1", from=2-7, to=3-8]
		\arrow["{1(1\eta)}", from=3-2, to=2-3]
		\arrow["{a^.}"{description}, from=3-2, to=3-4]
		\arrow["{(11)\eta}"{description}, from=3-4, to=2-5]
		\arrow["{u_{12}1}"{description}, from=3-4, to=4-5]
		\arrow["{a^.}"{description}, from=3-6, to=2-7]
		\arrow["a"{description}, from=3-6, to=4-7]
		\arrow["{(1\epsilon)1}", from=3-8, to=4-9]
		\arrow["a"{description}, from=3-8, to=5-8]
		\arrow["{(11)\eta}"{description}, from=4-5, to=3-6]
		\arrow["a"{description}, from=4-5, to=5-6]
		\arrow["{1a^.}"{description}, from=4-7, to=5-8]
		\arrow["a"{description}, from=4-9, to=6-9]
		\arrow["r1", from=4-9, to=7-11]
		\arrow["{1r^.}", from=5-1, to=3-2]
		\arrow["{r^.}"{description}, from=5-1, to=3-4]
		\arrow["{u_{12}}"{description}, from=5-1, to=5-4]
		\arrow["{u_{12}}"', curve={height=30pt}, from=5-1, to=7-11]
		\arrow["{r^.}"{description}, from=5-4, to=4-5]
		\arrow["{1r^.}"{description}, from=5-4, to=5-6]
		\arrow[""{name=0, anchor=center, inner sep=0}, "1"{description}, from=5-4, to=7-11]
		\arrow["{1(1\eta)}"{description}, from=5-6, to=4-7]
		\arrow["{1(\epsilon1)}"{description}, from=5-8, to=6-9]
		\arrow["1l"{description}, from=6-9, to=7-11]
		\arrow["{1\phi}"{description}, between={0.2}{0.8}, Rightarrow, nfold, from=4-7, to=0]
	\end{tikzcd}
\end{adjustbox}
\]

\clearpage
\[\vequiv\]

\[\begin{adjustbox}{
		angle=90,
		max totalsize={\linewidth}{0.90\textheight},
		valign=t
	}
	\begin{tikzcd}[ampersand replacement=\&,cramped]
		\&\&\&\& {F'(R_1(GF))} \&\&\&\& {F'((R_1G)F)} \&\& \\
		\&\&\&\&\&\&\&\&\& {(F'(R_1G))F} \\
		\&\&\&\&\&\&\& {F'((G'R_2)F)} \\
		\&\&\&\&\&\&\&\& {(F'(G'R_2))F} \&\& {((F'R_1)G)F} \\
		\& {F'(R_11)} \&\&\& {F'(G'(F'R_1))} \&\& {F'(G'(R_2F))} \&\&\&\& {((R_2F)G)F} \\
		\&\&\&\&\& {(F'G')(F'R_1)} \&\&\& {((F'G')R_2)F} \&\& {(R_2(FG))F} \\
		\&\& {F'((G'F')R_1)} \&\& {((F'G')F')R_1} \&\&\& {(F'G')(R_2F)} \&\& {(1R_2)F} \& {(R_21)F} \\
		\&\&\& {(F'(G'F'))R_1} \&\& {(1F')R_1} \& {1(F'R_1)} \&\& {1(R_2F)} \\
		\& {F'(1R_1)} \&\&\&\&\&\& {F'R_1} \&\&\& {R_2F} \\
		\&\& {(F'1)R_1} \\
		{F'R_1}
		\arrow["{1a^.}", from=1-5, to=1-9]
		\arrow["{1(-)}", between={0.2}{0.8}, Rightarrow, nfold, from=1-5, to=5-5]
		\arrow["{a^.}", from=1-9, to=2-10]
		\arrow["{a^.1}", from=2-10, to=4-11]
		\arrow["{(-)1}"', between={0.2}{0.8}, Rightarrow, nfold, from=2-10, to=7-10]
		\arrow["{1(u_{21}1)}"{description}, from=3-8, to=1-9]
		\arrow["{a^.}"{description}, from=3-8, to=4-9]
		\arrow["{(1u_{21})1}"{description}, from=4-9, to=2-10]
		\arrow["{(u_{12}1)1}", from=4-11, to=5-11]
		\arrow["{1(1\eta)}", from=5-2, to=1-5]
		\arrow["{1(1u_{12})}"{description}, from=5-5, to=5-7]
		\arrow["{a^.}"{description}, from=5-5, to=6-6]
		\arrow["{1a^.}"{description}, from=5-7, to=3-8]
		\arrow["{a^.}"{description}, from=5-7, to=7-8]
		\arrow["a1", from=5-11, to=6-11]
		\arrow["{(11)u_{12}}"{description}, from=6-6, to=7-8]
		\arrow["{\epsilon'(11)}"{description}, from=6-6, to=8-7]
		\arrow["a1"{description}, from=6-9, to=4-9]
		\arrow["{(\epsilon'1)1}"{description}, from=6-9, to=7-10]
		\arrow["{(1\epsilon)1}", from=6-11, to=7-11]
		\arrow["1a"{description}, from=7-3, to=5-5]
		\arrow["{a^.}"{description}, from=7-3, to=8-4]
		\arrow["a"{description}, from=7-5, to=6-6]
		\arrow["{(\epsilon'1)1}"{description}, from=7-5, to=8-6]
		\arrow["{a^.}"{description}, from=7-8, to=6-9]
		\arrow["{\epsilon'(11)}"{description}, from=7-8, to=8-9]
		\arrow["l1"{description}, from=7-10, to=9-11]
		\arrow["r1", from=7-11, to=9-11]
		\arrow[""{name=0, anchor=center, inner sep=0}, "{a^.1}"{description}, from=8-4, to=7-5]
		\arrow["a"{description}, from=8-6, to=8-7]
		\arrow["{1u_{12}}"{description}, from=8-7, to=8-9]
		\arrow["l"{description}, from=8-7, to=9-8]
		\arrow["{a^.}"{description}, from=8-9, to=7-10]
		\arrow["l"{description}, from=8-9, to=9-11]
		\arrow["{1(\eta'1)}"{description}, from=9-2, to=7-3]
		\arrow["{a^.}"{description}, from=9-2, to=10-3]
		\arrow["{u_{12}}"{description}, from=9-8, to=9-11]
		\arrow["{(1\eta')1}"{description}, from=10-3, to=8-4]
		\arrow["{1r^.}", from=11-1, to=5-2]
		\arrow["{1l^.}"{description}, from=11-1, to=9-2]
		\arrow[""{name=1, anchor=center, inner sep=0}, "1"{description}, from=11-1, to=9-8]
		\arrow["{u_{12}}"', curve={height=30pt}, from=11-1, to=9-11]
		\arrow["{r^.1}"{description}, from=11-1, to=10-3]
		\arrow["{\phi'1}"{description}, between={0.2}{0.8}, Rightarrow, nfold, from=0, to=1]
	\end{tikzcd}
\end{adjustbox}
\]

\clearpage

\[\begin{adjustbox}{
		angle=90,
		max totalsize={\linewidth}{0.90\textheight},
		valign=t
	}
	\begin{tikzcd}[ampersand replacement=\&,cramped,row sep=large]
		{G'R_2} \&\&\&\&\&\&\&\&\&\& {R_1G} \\
		\&\&\&\&\&\&\&\&\& {G'R_2} \\
		\& {(1G')R_2} \&\&\&\&\& {(G'1)R_2} \\
		\&\& {((G'F')G')R_2} \&\&\&\&\& {G'(1R_2)} \\
		\&\&\&\& {(G'(F'G'))R_2} \&\& {G'((F'G')R_2)}
		\arrow["{u_{21}}", from=1-1, to=1-11]
		\arrow[""{name=0, anchor=center, inner sep=0}, "1"{description}, from=1-1, to=2-10]
		\arrow["{l^.1}"', from=1-1, to=3-2]
		\arrow["{u_{21}}"', from=2-10, to=1-11]
		\arrow["{(\eta'1)1}"', from=3-2, to=4-3]
		\arrow["r1"{description}, from=3-7, to=2-10]
		\arrow["a"{description}, from=3-7, to=4-8]
		\arrow[""{name=1, anchor=center, inner sep=0}, "a1"', from=4-3, to=5-5]
		\arrow["1l"', from=4-8, to=2-10]
		\arrow["{(1\epsilon')1}"{description}, from=5-5, to=3-7]
		\arrow["a"', from=5-5, to=5-7]
		\arrow["{1(\epsilon'1)}"', from=5-7, to=4-8]
		\arrow["{\psi'1}", between={0.2}{0.8}, Rightarrow, nfold, from=0, to=1]
	\end{tikzcd}
\end{adjustbox}
\]

\clearpage
\[\vequiv\]
\[\begin{adjustbox}{
		angle=90,
		max totalsize={\linewidth}{0.90\textheight},
		valign=t
	}
	\begin{tikzcd}[ampersand replacement=\&,cramped,column sep=scriptsize,row sep=2.25em]
		{G'R_2} \&\&\&\&\&\&\&\&\&\&\&\&\&\&\&\& {R_1G} \\
		\\
		\&\&\& {R_1G} \&\&\&\&\&\&\&\&\& {R_1(G1)} \&\& {(R_1G)1} \\
		\&\&\&\&\&\& {R_1(1G)} \&\& {R_1((GF)G)} \&\& {R_1(G(FG))} \\
		{(1G')R_2} \&\&\&\& {(1R_1)G} \& {(R_11)G} \&\& {(R_1(GF))G} \&\&\&\& {(R_1G)(FG)} \&\& {(G'R_2)1} \&\&\& {G'R_2} \\
		\& {1(G'R_2)} \&\& {1(R_1G)} \&\& {((G'F')R_1)G} \&\&\&\&\& {((R_1G)F)G} \&\&\&\& {G'(R_21)} \\
		{((G'F')G')R_2} \&\&\&\&\&\& {(G'(F'R_1))G} \&\&\&\&\& {(G'R_2)(FG)} \\
		\& {(G'F')(G'R_2)} \&\& {(G'F')(R_1G)} \&\&\&\&\& {(G'(R_2F))G} \&\& {((G'R_2)F)G} \&\& {G'(R_2(FG))} \\
		\&\&\&\&\&\& {G'((F'R_1)G)} \\
		\&\&\&\&\& {G'(F'(R_1G))} \&\&\& {G'((R_2F)G)} \\
		\& {(G'(F'G'))R_2} \&\&\& {G'(F'(G'R_2))} \\
		\&\&\& {G'((F'G')R_2)} \&\&\&\&\&\&\&\&\& {G'(1R_2)}
		\arrow["{u_{21}}", from=1-1, to=1-17]
		\arrow["{u_{21}}"{description}, from=1-1, to=3-4]
		\arrow["{l^.1}"', from=1-1, to=5-1]
		\arrow["{l^.}"{description}, from=1-1, to=6-2]
		\arrow[""{name=0, anchor=center, inner sep=0}, from=3-4, to=1-17]
		\arrow["{1l^.}"{description}, from=3-4, to=4-7]
		\arrow["{l^.1}"{description}, from=3-4, to=5-5]
		\arrow["{r^.1}"{description}, from=3-4, to=5-6]
		\arrow["{l^.}"{description}, from=3-4, to=6-4]
		\arrow["1r"{description}, from=3-13, to=1-17]
		\arrow["{a^.}"{description}, from=3-13, to=3-15]
		\arrow["r"{description}, from=3-15, to=1-17]
		\arrow["{1(\eta1)}"{description}, from=4-7, to=4-9]
		\arrow[""{name=1, anchor=center, inner sep=0}, "1a"{description}, from=4-9, to=4-11]
		\arrow["{1(1\epsilon)}"{description}, from=4-11, to=3-13]
		\arrow["{a^.}"{description}, from=4-11, to=5-12]
		\arrow["a"{description}, from=5-1, to=6-2]
		\arrow["{(\eta'1)1}"', from=5-1, to=7-1]
		\arrow["{(\eta'1)1}"{description}, from=5-5, to=6-6]
		\arrow["a"{description}, from=5-6, to=4-7]
		\arrow[""{name=2, anchor=center, inner sep=0}, "{(1\eta)1}"{description}, from=5-6, to=5-8]
		\arrow["a"{description}, from=5-8, to=4-9]
		\arrow["{a^.1}"{description}, from=5-8, to=6-11]
		\arrow["{(11)\epsilon}"{description}, from=5-12, to=3-15]
		\arrow["{u_{21}1}"{description}, from=5-14, to=3-15]
		\arrow["r"{description}, from=5-14, to=5-17]
		\arrow["a"{description}, from=5-14, to=6-15]
		\arrow["{u_{21}}"', from=5-17, to=1-17]
		\arrow["{1u_{21}}"{description}, from=6-2, to=6-4]
		\arrow["{\eta'(11)}"{description}, from=6-2, to=8-2]
		\arrow["{a^.}"{description}, from=6-4, to=5-5]
		\arrow["{\eta'(11)}"{description}, from=6-4, to=8-4]
		\arrow["a1"{description}, from=6-6, to=7-7]
		\arrow["a"{description}, from=6-11, to=5-12]
		\arrow["1r"{description}, from=6-15, to=5-17]
		\arrow["a"{description}, from=7-1, to=8-2]
		\arrow["a1"', from=7-1, to=11-2]
		\arrow["{(1u_{12})1}"{description}, from=7-7, to=8-9]
		\arrow["a"{description}, from=7-7, to=9-7]
		\arrow["{u_{21}(11)}"{description}, from=7-12, to=5-12]
		\arrow["{(11)\epsilon}"{description}, from=7-12, to=5-14]
		\arrow["a"{description}, from=7-12, to=8-13]
		\arrow["{(11)u_{21}}"{description}, from=8-2, to=8-4]
		\arrow["a"{description}, from=8-2, to=11-5]
		\arrow["{a^.}"{description}, from=8-4, to=6-6]
		\arrow["a"{description}, from=8-4, to=10-6]
		\arrow["{a^.1}"{description}, from=8-9, to=8-11]
		\arrow["a"{description}, from=8-9, to=10-9]
		\arrow["{(u_{21}1)1}"{description}, from=8-11, to=6-11]
		\arrow["a"{description}, from=8-11, to=7-12]
		\arrow["{1(1\epsilon)}"{description}, from=8-13, to=6-15]
		\arrow["{1(u_{12}1)}"{description}, from=9-7, to=10-9]
		\arrow["{1a^.}"{description}, from=10-6, to=9-7]
		\arrow[""{name=3, anchor=center, inner sep=0}, "1a"{description}, from=10-9, to=8-13]
		\arrow["a"', from=11-2, to=12-4]
		\arrow["{1(1u_{21})}"{description}, from=11-5, to=10-6]
		\arrow["1a"{description}, from=12-4, to=11-5]
		\arrow["{1(\epsilon'1)}"', from=12-4, to=12-13]
		\arrow["1l"', from=12-13, to=5-17]
		\arrow["{1\psi}"{description}, between={0.2}{0.8}, Rightarrow, nfold, from=0, to=1]
		\arrow["{(-)1}"{description}, between={0.2}{1}, Rightarrow, nfold, from=2, to=7-7]
		\arrow["{1(-)}"{description}, between={0.2}{0.8}, Rightarrow, nfold, from=3, to=12-13]
	\end{tikzcd}
\end{adjustbox}
\]

\subsection{Weak $2$-inverse $2$-morphism} \label{app:weak2invs2morph}

\[\begin{adjustbox}{
		angle=90,
		max totalsize={\linewidth}{0.90\textheight},
		valign=t
	}
	\begin{tikzcd}[ampersand replacement=\&,cramped,column sep=scriptsize,row sep=2.25em]
		\&\&\&\& {S_1 1} \&\&\&\& {S_1(GF)} \&\&\&\& \\
		\\
		{S_1} \&\& {1S_1} \&\& {(G'F')S_1} \&\& {G'(F'S_1)} \&\& {G'(S_2F)} \&\& {(G'S_2)F} \&\& {(S_1G)F} \\
		\\
		{R_1} \&\& {1R_1} \&\& {(G'F')R_1} \&\& {G'(F'R_1)} \&\& {G'(R_2F)} \&\& {(G'R_2)F} \&\& {(R_1G)F}
		\arrow[""{name=0, anchor=center, inner sep=0}, "{1\eta}", from=1-5, to=1-9]
		\arrow["{a^.}", from=1-9, to=3-13]
		\arrow["{r^.}", from=3-1, to=1-5]
		\arrow["{l^.}"{description}, from=3-1, to=3-3]
		\arrow["{\eta'1}"{description}, from=3-3, to=3-5]
		\arrow["a"{description}, from=3-5, to=3-7]
		\arrow["{1v_{12}}"{description}, from=3-7, to=3-9]
		\arrow["{1(-)}"{description}, between={0.2}{0.8}, Rightarrow, nfold, from=3-7, to=5-9]
		\arrow["{a^.}"{description}, from=3-9, to=3-11]
		\arrow["{v_{21}1}"{description}, from=3-11, to=3-13]
		\arrow["{(-)1}"{description}, between={0.2}{0.8}, Rightarrow, nfold, from=3-11, to=5-13]
		\arrow["{\gamma_1}", from=5-1, to=3-1]
		\arrow["{l^.}"', from=5-1, to=5-3]
		\arrow["{1\gamma_1}"{description}, from=5-3, to=3-3]
		\arrow["{\eta'1}"', from=5-3, to=5-5]
		\arrow["{(11)\gamma_1}"{description}, from=5-5, to=3-5]
		\arrow["a"', from=5-5, to=5-7]
		\arrow["{1(1\gamma_1)}"{description}, from=5-7, to=3-7]
		\arrow["{1u_{12}}"', from=5-7, to=5-9]
		\arrow["{1(\gamma_21)}"{description}, from=5-9, to=3-9]
		\arrow["{a^.}"', from=5-9, to=5-11]
		\arrow["{(1 \gamma_2)1}"{description}, from=5-11, to=3-11]
		\arrow["{u_{21}1}"', from=5-11, to=5-13]
		\arrow["{(\gamma_11)1}"', from=5-13, to=3-13]
		\arrow[between={0.2}{0.8}, Rightarrow, nfold, from=0, to=3-7]
	\end{tikzcd}
\end{adjustbox}
\]

\clearpage
\[\vequiv\]

\[\begin{adjustbox}{
		angle=90,
		max totalsize={\linewidth}{0.90\textheight},
		valign=t
	}
	\begin{tikzcd}[ampersand replacement=\&,cramped,column sep=scriptsize,row sep=2.25em]
		\&\&\&\& {S_1 1} \&\&\&\& {S_1(GF)} \&\&\&\& \\
		\\
		{S_1} \&\&\&\& {R_11} \&\&\&\& {R_1(GF)} \&\&\&\& {(S_1G)F} \\
		\\
		{R_1} \&\& {1R_1} \&\& {(G'F')R_1} \&\& {G'(F'R_1)} \&\& {G'(R_2F)} \&\& {(G'R_2)F} \&\& {(R_1G)F}
		\arrow["{1\eta}", from=1-5, to=1-9]
		\arrow["{a^.}", from=1-9, to=3-13]
		\arrow["{r^.}", from=3-1, to=1-5]
		\arrow["{\gamma_11}"{description}, from=3-5, to=1-5]
		\arrow[""{name=0, anchor=center, inner sep=0}, "{1\eta}"{description}, from=3-5, to=3-9]
		\arrow["{\gamma_1(11)}"{description}, from=3-9, to=1-9]
		\arrow["{a^.}"{description}, from=3-9, to=5-13]
		\arrow["{\gamma_1}", from=5-1, to=3-1]
		\arrow["{r^.}"{description}, from=5-1, to=3-5]
		\arrow["{l^.}"', from=5-1, to=5-3]
		\arrow["{\eta'1}"', from=5-3, to=5-5]
		\arrow["a"', from=5-5, to=5-7]
		\arrow["{1u_{12}}"', from=5-7, to=5-9]
		\arrow["{a^.}"', from=5-9, to=5-11]
		\arrow["{u_{21}1}"', from=5-11, to=5-13]
		\arrow["{(\gamma_11)1}"', from=5-13, to=3-13]
		\arrow[between={0.2}{0.8}, Rightarrow, nfold, from=0, to=5-7]
	\end{tikzcd}
\end{adjustbox}
\]

\clearpage

\[\begin{adjustbox}{
		angle=90,
		max totalsize={\linewidth}{0.90\textheight},
		valign=t
	}
	\begin{tikzcd}[ampersand replacement=\&,cramped,column sep=small,row sep=2.25em]
		{(F'G')S_2} \&\& {F'(G'S_2)} \&\& {F'(S_1G)} \&\& {(F'S_1)G} \&\& {(S_2F)G} \&\& {S_2(FG)} \&\& {S_21} \&\& {S_2} \\
		\\
		{(F'G')R_2} \&\& {F'(G'R_2)} \&\& {F'(R_1G)} \&\& {(F'R_1)G} \&\& {(R_2F)G} \&\& {R_2(FG)} \&\& {R_21} \&\& {R_2} \\
		\\
		\&\&\&\&\&\&\& {1R_2}
		\arrow["a", from=1-1, to=1-3]
		\arrow["{1v_{21}}", from=1-3, to=1-5]
		\arrow["{1(-)}"{description}, between={0.2}{0.8}, Rightarrow, nfold, from=1-3, to=3-5]
		\arrow["{a^.}", from=1-5, to=1-7]
		\arrow["{v_{12}1}", from=1-7, to=1-9]
		\arrow["{(-)1}"{description}, between={0.2}{0.8}, Rightarrow, nfold, from=1-7, to=3-9]
		\arrow["a", from=1-9, to=1-11]
		\arrow["{1\epsilon}", from=1-11, to=1-13]
		\arrow["r", from=1-13, to=1-15]
		\arrow["{(11)\gamma_2}", from=3-1, to=1-1]
		\arrow["a"{description}, from=3-1, to=3-3]
		\arrow["{\epsilon'1}"', from=3-1, to=5-8]
		\arrow["{1(1\gamma_2)}"{description}, from=3-3, to=1-3]
		\arrow["{1u_{21}}"{description}, from=3-3, to=3-5]
		\arrow["{1(\gamma_11)}"{description}, from=3-5, to=1-5]
		\arrow["{a^.}"{description}, from=3-5, to=3-7]
		\arrow["{(1\gamma_1)1}"{description}, from=3-7, to=1-7]
		\arrow[""{name=0, anchor=center, inner sep=0}, "{u_{12}1}"{description}, from=3-7, to=3-9]
		\arrow["{(\gamma_21)1}"{description}, from=3-9, to=1-9]
		\arrow["a"{description}, from=3-9, to=3-11]
		\arrow["{\gamma_2(11)}"{description}, from=3-11, to=1-11]
		\arrow["{1\epsilon}"{description}, from=3-11, to=3-13]
		\arrow["{\gamma_21}"', from=3-13, to=1-13]
		\arrow["r"{description}, from=3-13, to=3-15]
		\arrow["{\gamma_2}"', from=3-15, to=1-15]
		\arrow["l"', from=5-8, to=3-15]
		\arrow[between={0.2}{0.8}, Rightarrow, nfold, from=0, to=5-8]
	\end{tikzcd}
\end{adjustbox}
\]

\clearpage
\[\vequiv\]

\[\begin{adjustbox}{
		angle=90,
		max totalsize={\linewidth}{0.90\textheight},
		valign=t
	}
	\begin{tikzcd}[ampersand replacement=\&,cramped,column sep=small,row sep=2.25em]
		{(F'G')S_2} \&\& {F'(G'S_2)} \&\& {F'(S_1G)} \&\& {(F'S_1)G} \&\& {(S_2F)G} \&\& {S_2(FG)} \&\& {S_21} \&\& {S_2} \\
		\\
		{(F'G')R_2} \&\&\&\&\&\&\& {1S_2} \&\&\&\&\&\&\& {R_2} \\
		\\
		\&\&\&\&\&\&\& {1R_2}
		\arrow["a", from=1-1, to=1-3]
		\arrow["{\epsilon'1}"{description}, from=1-1, to=3-8]
		\arrow["{1v_{21}}", from=1-3, to=1-5]
		\arrow["{a^.}", from=1-5, to=1-7]
		\arrow[""{name=0, anchor=center, inner sep=0}, "{v_{12}1}", from=1-7, to=1-9]
		\arrow["a", from=1-9, to=1-11]
		\arrow["{1\epsilon}", from=1-11, to=1-13]
		\arrow["r", from=1-13, to=1-15]
		\arrow["{(11)\gamma_2}", from=3-1, to=1-1]
		\arrow["{\epsilon'1}"', from=3-1, to=5-8]
		\arrow["l"{description}, from=3-8, to=1-15]
		\arrow["{\gamma_2}"', from=3-15, to=1-15]
		\arrow["{1\gamma_2}"{description}, from=5-8, to=3-8]
		\arrow["l"', from=5-8, to=3-15]
		\arrow[between={0.2}{0.8}, Rightarrow, nfold, from=0, to=3-8]
	\end{tikzcd}
\end{adjustbox}
\]

	\newpage
	
	
	\clearpage 
	
	\bibliography{master} 

@Book{reihlcht,
  author    = {Riehl, Emily},
  publisher = {Cambridge University Press},
  title     = {{Categorical Homotopy Theory}},
  year      = {2014},
  isbn      = {9781107261457},
  month     = may,
  doi       = {10.1017/cbo9781107261457},
}

@Article{betti1983,
  author    = {Betti, Renato and Carboni, Aurelio and Street, Ross and Walters, Robert},
  journal   = {Journal of Pure and Applied Algebra},
  title     = {{Variation through enrichment}},
  year      = {1983},
  issn      = {0022-4049},
  month     = aug,
  number    = {2},
  pages     = {109--127},
  volume    = {29},
  doi       = {10.1016/0022-4049(83)90100-7},
  publisher = {Elsevier BV},
}

@Article{wolff1974,
  author    = {Wolff, Harvey},
  journal   = {Journal of Pure and Applied Algebra},
  title     = {{V-cat and V-graph}},
  year      = {1974},
  issn      = {0022-4049},
  month     = apr,
  number    = {2},
  pages     = {123--135},
  volume    = {4},
  doi       = {10.1016/0022-4049(74)90018-8},
  publisher = {Elsevier BV},
}

@Book{kellyect,
  author    = {G. M. Kelly},
  publisher = {Reprints in Theory and Applications of Categories},
  title     = {{Basic Concepts of Enriched Category Theory}},
  year      = {2005},
  note      = {Reprint of the 1982 Cambridge University Press edition},
  url       = {http://www.tac.mta.ca/tac/reprints/articles/10/tr10.pdf},
}

@Book{yau2cat2021,
  author    = {Johnson, Niles and Yau, Donald},
  publisher = {Oxford University Press USA - OSO},
  title     = {{2-Dimensional Categories}},
  year      = {2021},
  isbn      = {9780192645678},
}

@Book{aguiarmon2010,
  author    = {Aguiar, Marcelo and Mahajan, Swapneel},
  publisher = {American Mathematical Society and Centre de Recherches Mathématiques},
  title     = {{Monoidal functors, species, and Hopf algebras}},
  year      = {2010},
  isbn      = {978-0-8218-4776-3},
  series    = {CRM Monograph Series},
  volume    = {29},
  pagetotal = {789},
}

@Book{leinster2004,
  author     = {Leinster, Tom},
  publisher  = {Cambridge University Press},
  title      = {{Higher Operads, Higher Categories}},
  year       = {2004},
  series     = {London Mathematical Society Lecture Note Series},
  collection = {London Mathematical Society Lecture Note Series},
  doi        = {https://doi.org/10.1017/CBO9780511525896},
  place      = {Cambridge},
}

@Book{yauinvcat2020,
  author    = {Yau, Donald},
  publisher = {Springer International Publishing},
  title     = {{Involutive Category Theory}},
  year      = {2020},
  isbn      = {9783030612030},
  doi       = {10.1007/978-3-030-61203-0},
  issn      = {1617-9692},
  journal   = {Lecture Notes in Mathematics},
}

@Article{baezhda31998,
  author    = {Baez, John C. and Dolan, James},
  journal   = {Advances in Mathematics},
  title     = {{Higher-Dimensional Algebra III.n-Categories and the Algebra of Opetopes}},
  year      = {1998},
  issn      = {0001-8708},
  month     = may,
  number    = {2},
  pages     = {145--206},
  volume    = {135},
  doi       = {10.1006/aima.1997.1695},
  publisher = {Elsevier BV},
}

@Article{kellyvcatpres2001,
  author  = {Kelly, G. M. and Lack, Stephen},
  journal = {Theory and Applications of Categories},
  title   = {{V-Cat is locally presentable or locally bounded if V is so}},
  year    = {2001},
  issn    = {1201-561X},
  pages   = {555--575},
  volume  = {08},
  doi     = {10.70930/tac/to984gmg},
}

@Article{streetncat1987,
  author    = {Street, Ross},
  journal   = {Journal of Pure and Applied Algebra},
  title     = {{The algebra of oriented simplexes}},
  year      = {1987},
  issn      = {0022-4049},
  month     = dec,
  number    = {3},
  pages     = {283--335},
  volume    = {49},
  doi       = {10.1016/0022-4049(87)90137-x},
  publisher = {Elsevier BV},
}

@PhdThesis{gurski2006tricat,
  author = {Gurski, Nick},
  school = {University of Chicago},
  title  = {{An algebraic theory of tricategories}},
  year   = {2006},
  url    = {https://ncatlab.org/nlab/files/Gurski-AlgebraicTricategories.pdf},
}

@Book{gurski2013tricatcoh,
  author    = {Gurski, Nick},
  publisher = {Cambridge University Press},
  title     = {{Coherence in Three-Dimensional Category Theory}},
  year      = {2013},
  isbn      = {9781107034891},
  month     = Mar,
  doi       = {10.1017/cbo9781139542333},
}

@Article{gurski2012biequivintricat,
  author    = {Gurski, Nick},
  journal   = {Theory and Applications of Categories},
  title     = {Biequivalences in tricategories},
  year      = {2012},
  issn      = {1201-561X},
  pages     = {349--384},
  volume    = {26},
  doi       = {10.70930/tac/5bwbxsdl},
  publisher = {Mount Allison University},
}

@Article{gps1995,
  author    = {Gordon, R. and Power, A. J. and Street, Ross},
  journal   = {Memoirs of the American Mathematical Society},
  title     = {Coherence for tricategories},
  year      = {1995},
  issn      = {0065-9266},
  number    = {558},
  volume    = {117},
  doi       = {10.1090/memo/0558},
  publisher = {American Mathematical Society (AMS)},
}

@Article{chengmates2014,
  author    = {Cheng, Eugenia and Gurski, Nick and Riehl, Emily},
  journal   = {Journal of K-Theory},
  title     = {Cyclic multicategories, multivariable adjunctions and mates},
  year      = {2014},
  issn      = {1865-5394},
  month     = Jan,
  number    = {2},
  pages     = {337--396},
  volume    = {13},
  doi       = {10.1017/is013012007jkt250},
  publisher = {Cambridge University Press (CUP)},
}

@Article{gurski2011monbicat,
  author    = {Gurski, Nick},
  journal   = {Advances in Mathematics},
  title     = {Loop spaces, and coherence for monoidal and braided monoidal bicategories},
  year      = {2011},
  issn      = {0001-8708},
  month     = Mar,
  number    = {5},
  pages     = {4225--4265},
  volume    = {226},
  doi       = {10.1016/j.aim.2010.12.007},
  publisher = {Elsevier BV},
}

@Misc{pries2011monbicat,
  author    = {Schommer-Pries, Christopher J.},
  title     = {The Classification of Two-Dimensional Extended Topological Field Theories},
  year      = {2011},
  copyright = {arXiv.org perpetual, non-exclusive license},
  doi       = {10.48550/ARXIV.1112.1000},
  publisher = {arXiv},
  url       = {https://arxiv.org/abs/1112.1000},
}

@Article{gurski2013symmonbicat,
  author    = {Gurski, Nick and Osorno, Angélica M.},
  journal   = {Advances in Mathematics},
  title     = {Infinite loop spaces, and coherence for symmetric monoidal bicategories},
  year      = {2013},
  issn      = {0001-8708},
  month     = Oct,
  pages     = {1--32},
  volume    = {246},
  doi       = {10.1016/j.aim.2013.06.028},
  publisher = {Elsevier BV},
}

@Book{yau2024bimodcats,
  author    = {Johnson, Niles and Yau, Donald},
  publisher = {American Mathematical Society},
  title     = {{Bimonoidal Categories, En-Monoidal Categories, and Algebraic K-Theory}},
  year      = {2024},
  isbn      = {9781470478117},
  month     = Oct,
  doi       = {10.1090/surv/285},
  issn      = {2331-7159},
  journal   = {Mathematical Surveys and Monographs},
}

@Article{grandis2017intercatframework,
  author    = {Grandis, Marco and Paré, Robert},
  journal   = {Journal of Pure and Applied Algebra},
  title     = {Intercategories: A framework for three-dimensional category theory},
  year      = {2017},
  issn      = {0022-4049},
  month     = May,
  number    = {5},
  pages     = {999--1054},
  volume    = {221},
  doi       = {10.1016/j.jpaa.2016.08.002},
  publisher = {Elsevier BV},
}

@Book{grandis2019dblcat,
  author    = {Grandis, Marco},
  publisher = {WORLD SCIENTIFIC},
  title     = {Higher Dimensional Categories: From Double to Multiple Categories},
  year      = {2019},
  isbn      = {9789811205118},
  month     = Apr,
  doi       = {10.1142/11406},
}

@Article{shulman2010vequip,
  author    = {Cruttwell, G.S.H. and Shulman, Michael A.},
  journal   = {Theory and Applications of Categories},
  title     = {A unified framework for generalized multicategories.},
  year      = {2010},
  pages     = {580-655},
  volume    = {24},
  language  = {eng},
  publisher = {Mount Allison University, Department of Mathematics and Computer Science, Sackville},
  url       = {http://eudml.org/doc/228027},
}

@Article{gurski2014itericon,
  author  = {Cheng, Eugenia and Gurski, Nick},
  journal = {Theory and Applications of Categories},
  title   = {Iterated icons},
  year    = {2014},
  number  = {32},
  pages   = {929--977},
  volume  = {29},
  eprint  = {https://www.tac.mta.ca/tac/volumes/29/32/29-32.pdf},
}

@Book{etingof2015tenscat,
  author    = {Etingof, Pavel and Shlomo, Gelaki and Dmitri, Nikshych and Victor, Ostrik},
  publisher = {American Mathematical Society},
  title     = {Tensor categories},
  year      = {2015},
  address   = {Providence, Rhode Island},
  isbn      = {9781470434410},
  number    = {205},
  series    = {Mathematical surveys and monographs},
  pagetotal = {344},
  ppn_gvk   = {1616794232},
  url       = {https://math.mit.edu/~etingof/egnobookfinal.pdf},
}

@Book{hovey2007modelcat,
  author    = {Hovey, Mark},
  publisher = {American Mathematical Society},
  title     = {{Model Categories}},
  year      = {2007},
  isbn      = {9781470412906},
  month     = Oct,
  doi       = {10.1090/surv/063},
  issn      = {2331-7159},
  journal   = {Mathematical Surveys and Monographs},
}

@Article{ara2015autequivinfcat,
  author    = {Ara, Dimitri and Groth, Moritz and Gutiérrez, Javier J.},
  journal   = {Mathematische Zeitschrift},
  title     = {On autoequivalences of the {$(\infty,1)$}-category of {$\infty$}-operads},
  year      = {2015},
  issn      = {1432-1823},
  month     = sep,
  number    = {3-4},
  pages     = {807--848},
  volume    = {281},
  doi       = {10.1007/s00209-015-1509-5},
  publisher = {Springer Science and Business Media LLC},
}

@Article{markakis2025hominfcat,
  author    = {Benjamin, Thibaut and Markakis, Ioannis},
  journal   = {Journal of the London Mathematical Society},
  title     = {Hom {$\omega$}‐categories of a computad are free},
  year      = {2025},
  issn      = {1469-7750},
  month     = nov,
  number    = {6},
  volume    = {112},
  doi       = {10.1112/jlms.70367},
  publisher = {Wiley},
}

@PhdThesis{segrt2012strengthsmorita,
  author      = {Segrt Ratkovic, Kruna},
  school      = {{Universit{\'e} de Nice Sophia Antipolis - UFR}},
  title       = {{Morita theory in enriched context}},
  year        = {2012},
  month       = feb,
  hal_id      = {tel-00785301},
  hal_version = {v1},
  url         = {https://theses.hal.science/tel-00785301},
}

@Article{ara2011thetacatcocat,
  author    = {Ara, Dimitri},
  journal   = {Applied Categorical Structures},
  title     = {{The Groupoidal Analogue $\widetilde{\Theta}$ to Joyal’s Category $\Theta$ is a Test Category}},
  year      = {2011},
  issn      = {1572-9095},
  month     = jul,
  number    = {6},
  pages     = {603--649},
  volume    = {20},
  doi       = {10.1007/s10485-011-9255-2},
  publisher = {Springer Science and Business Media LLC},
}

@Article{shulman2018contenrich,
  author    = {Shulman, Michael},
  journal   = {Theory and Applications of Categories},
  title     = {Contravariance through enrichment},
  year      = {2018},
  issn      = {1201-561X},
  pages     = {95--130},
  volume    = {33},
  doi       = {10.70930/tac/s8y6o655},
  publisher = {Mount Allison University},
}

@Article{morton2009doublebicat,
  author  = {Morton, Jeffrey C.},
  journal = {Journal of Homotopy and Related Structures},
  title   = {{Double Bicategories and Double Cospans}},
  year    = {2009},
  number  = {1},
  pages   = {389--428},
  volume  = {4},
  url     = {https://tcms.org.ge/Journals/JHRS/xvolumes/2009/n1a19/v4n1a19.pdf},
}

@Article{lack2008icons,
  author    = {Lack, Stephen},
  journal   = {Applied Categorical Structures},
  title     = {Icons},
  year      = {2008},
  issn      = {1572-9095},
  month     = Apr,
  number    = {3},
  pages     = {289--307},
  volume    = {18},
  doi       = {10.1007/s10485-008-9136-5},
  publisher = {Springer Science and Business Media LLC},
}

@Misc{fujii2025famenr,
  author    = {Fujii, Soichiro and Lack, Stephen},
  title     = {The familial nature of enrichment over virtual double categories},
  year      = {2025},
  copyright = {arXiv.org perpetual, non-exclusive license},
  doi       = {10.48550/ARXIV.2507.05529},
  publisher = {arXiv},
}
	
\end{document}